\documentclass[11pt,leqno]{article}
\usepackage[english]{babel}
\usepackage{amsmath}
\usepackage{amssymb}
\usepackage{mathrsfs}
\usepackage[all]{xy}
\UseComputerModernTips
\usepackage{latexsym}
\usepackage{amsfonts}

\usepackage[dvips]{graphicx}

\usepackage{amsthm,enumerate}
\newcommand{\earrow}
{\mathrel{\text{\smash{\lower.477ex\hbox{$\overset{\text{\normalsize \smash{\lower.55ex\hbox{$\sim\,$}}}}{\to}$}}}}}

\numberwithin{equation}{section}

\newcommand{\A}{\mathcal{A}}

\newcommand{\C}{\mathcal{C}}

\newcommand{\D}{\mathcal{D}}

\renewcommand{\mod}{\mathrm{Mod}}

\newcommand{\op}{\mathrm{Op}}

\newcommand{\rc}{\mathbb{R}\textrm{-}\mathrm{c}}

\newcommand{\Powl}{\hat{\mathcal{P}}(\Z_\ell)}

\newcommand{\CC}{\mathbb{C}}

\newcommand{\R}{\mathbb{R}}

\newcommand{\Z}{\mathbb{Z}}

\newcommand{\N}{\mathbb{N}}

\newcommand{\OO}{\mathcal{O}}

\newcommand{\iso}{\stackrel{\sim}{\to}}
\newcommand{\dbt}{\mathcal{D}\mathit{b}^t}

\newcommand{\db}{\mathcal{D}\mathit{b}}

\newcommand{\ot}{\mathcal{O}^t}

\newcommand{\OW}{\OO^\mathrm{w}}
\newcommand{\OWX}{\OO^\mathrm{w}_X}

\newcommand{\CW}{\C^{{\infty ,\mathrm{w}}}}
\newcommand{\wtens}{\overset{\mathrm{w}}{\otimes}}

\newcommand{\rh}{\mathit{R}\mathcal{H}\mathit{om}}

\newcommand{\ho}{\mathcal{H}\mathit{om}}

\newcommand{\id}{\mathrm{id}}

\newcommand{\supp}{\mathrm{supp}}

\newcommand{\RP}{\mathbb{R}^{{\scriptscriptstyle{+}}}}

\newcommand{\imin}[1]{#1^{-1}}

\newcommand{\lind}[1]{\underset{#1}{\underrightarrow{\lim}}}

\newcommand{\lpro}[1]{\underset{#1}{\underleftarrow{\lim}}}

\newcommand{\exs}[3]{0 \to {#1} \to {#2} \to {#3} \to 0}

\newcommand{\lexs}[3]{0 \to {#1} \to {#2} \to {#3}}

\newcommand{\dt}[3]{{#1} \to {#2} \to {#3} \stackrel{+}{\to}}

\newtheorem{teo}{Theorem}[section]%[chapter]
\newtheorem{cor}[teo]{Corollary}%[chapter]
\newtheorem{prop}[teo]{Proposition}%[chapter]
\newtheorem{lem}[teo]{Lemma}%[chapter]
\newtheorem{nt}[teo]{Notation}%[chapter]
\newtheorem{es}[teo]{Example}%[chapter]
\newtheorem{oss}[teo]{Remark}%[chapter]
\newtheorem{df}[teo]{Definition}%[chapter]

\usepackage{latexsym}

\newcommand{\subsetl}{{\widehat{\mathcal{P}}_\ell}}

\newcommand{\RRpair}{{\mathcal{H}_{\tau} \times \mathbb{R}_{\ge 0}}}
\newcommand{\RRTpair}{{\mathcal{H}_{\tau,\lambda} \times \mathbb{R}_{\ge 0}}}
\newcommand{\Ssec}{{S_H(V,\,\{\epsilon\})}}
\newcommand{\compl}{\complement}
\newcommand{\comD}{{\sigma_A}}
\newcommand{\intA}{{\comD A_\chi}}
\newcommand{\semiG}{{\mathcal{G}}}

\newcommand{\ckM}{{\bigcup_k M}}
\newcommand{\ckZ}{{\bigcup_k Z}}

\author{Naofumi Honda and Luca Prelli}
\title{Generalization of multi-specializations and multi-asymptotics}
\date{}

\begin{document}

\maketitle\thispagestyle{empty}

\begin{abstract}
The aim of this paper is to give a new description of the
geometry appearing in the multi-specialization along a general
family of submanifolds of a real analytic manifold (including some important cases as clean intersection or a
simultaneously linearizable family of Lagrangian submanifolds in a cotangent bundle) and then,
to extend several properties of the multi-specialization.
The notion of multi-asymptotic expansions is also extended.
 In the local model more general cases are studied:
locally we can construct new sheaves of multi-asymptotically developable
functions closely related with asymptotics along a subvariety with
a simple singularity such as a cusp.
\end{abstract}

\tableofcontents

\addcontentsline{toc}{section}{\textbf{Introduction}}

\section*{Introduction}

We have established, in \cite{HP13} and \cite{HPY}, the theory of
multi-specialization
and multi-microlocalization along a simultaneously linearizable family $\chi$ of
closed submanifolds (see Definition \ref{df:sim-linear} and Proposition
\ref{prop:canonical-from})
on which we additionally
impose some moderate condition
for their geometrical configuration, that is, the condition H2 of \cite{HP13} for $\chi$.
A typical example  of such a configuration (with the condition H2)
of two submanifolds is
that they transversally intersect or that one submanifold completely
contains the other.
However, some important cases do not satisfy this geometrical condition such as
two submanifolds which cleanly intersect but not transversally or a
simultaneously linearizable family of Lagrangian submanifolds in a cotangent bundle,
where we encounter the latter situation
in the micro-support estimation of a multi-microlocalized object (Theorem 3.6 \cite{HPY}).

The primary purpose of our paper is to give a new description of the
geometry appearing in the multi-specialization along a general
linearizable family $\chi$ (without the condition H2), and then,
to extend several properties of the multi-specialization
functor already established in \cite{HP13} to this general $\chi$.
One should be aware that, however, the geometry
appearing in the multi-specialization for a general $\chi$ is
quite different from that for a $\chi$ with the condition H2.
For a usual specialization along a closed submanifold $M$,
a (locally defined) conic cone along the submanifold $M$
can be considered as a fundamental geometrical object of this
specialization in the sense that
both the normal cone $C_M(\cdot)$ along $M$
and a fiber formula for the specialization functor $\nu_M(\cdot)$ along $M$
can be described by using a family of these cones.
When we consider the multi-specialization along
a simultaneously linearizable family $\chi$ of submanifolds
$\{M_1,\,\dots,\,M_\ell\}$ in $X$ with the condition H2,
the situation is quite similar as that for the usual specialization. The fundamental
geometrical object in this case is a multicone that is, by definition, the
intersection of $\ell$ conic cones
along submanifolds $M_1$, $M_2$, $\dots$, $M_\ell$, respectively.
And the multi-normal cone $C_\chi(\cdot)$ along $\chi$
and a fiber formula for the
multi-specialization functor $\nu_\chi(\cdot)$ along $\chi$
can be described by using a family of these multicones
as we have shown in \cite{HP13}.

On the other hand a multicone for a general linearizable family $\chi$
can be no longer obtained by a set-theoretical operation of conic cones
along submanifolds. In fact, as we see
in Section \ref{section:Multi-cones for a general case},
it is defined by a semigroup generated by
monomials associated with local actions $\mu_k$'s $(k=1,\dots,\ell)$,
where each action $\mu_k$ is induced from the one on the conormal bundle $T_{M_k}X$
of $M_k$ with a local identification $X \simeq T_{M_k}X$ near $M_k$.
%Hence we cannot
%expect linear structure of a multicone anymore.
Furthermore, when we drop the condition H2,
we cannot either
expect the zero section $S_\chi$ of the multi-normal deformation $\widetilde{X}_\chi$
along $\chi$,
i.e., the subset $\{t=0\}$ in $\widetilde{X}_\chi$, to
be a vector bundle over the base space $M=M_1 \cap \dots \cap M_\ell$.
In fact, let us consider $\chi$ which
consists of two submanifolds
intersecting cleanly but not transversally. Then,
as Example {\ref{es:non-vector-bundle}} shows,
$S_\chi$ is not a vector bundle over $M$ even for this simple case.
Although several unexpected facts are observed in geometry of the
multi-specialization for a general $\chi$ as we have already seen,
Theorems \ref{teo:multi-normal-cone-multi-cone} and
\ref{sections of multispec}
and Equation \eqref{eq:fibers} in this paper
say that, by using our newly introduced multicones,
we can characterize the multi-normal cone
and establish a fiber formula for the multi-specialization functor in a usual way.

The shape of a multicone becomes a little bit complicated
for a general $\chi$, however, it
still has many good properties. It is, in particular,
infinitesimally stable under the actions $\mu_k$'s.
This fact suggests us that an asymptotic expansion
can be defined on this multicone since
it is, roughly speaking,
nothing but a formal Taylor expansion
along orbits generated by the actions $\mu_k$'s.
As a matter of fact, in Section \ref{section:Multi-asymptotic expansions},
the notion of multi-asymptotic expansions
along $\chi$ introduced in \cite{HP13} is successfully extended
to the one for a general simultaneously linearizable family $\chi$.
Further, our generalization of the multi-asymptotic expansion can
manipulate the case that each action $\mu_k$ is not necessarily linear.
Hence, by choosing a family of appropriate non-linear actions,
locally we can construct a new sheaf of multi-asymptotically developable
functions which is closely tied with asymptotics along a subvariety with
a simple singularity such as a cusp.

We study, in the Section \ref{section:Multi-specialization}, the functorial nature of multi-asymptotics, proving some vanishing theorems which provide a general Borel-Ritt exact sequence in this setting. Similar vanishing results hold for the multispecialization of tempered holomorphic functions.
We also study functorial operations of multi-specializations
such as a direct image and an inverse image
in our general settings.
By an operation associated with a desingularization map, for example,
we can show isomorphisms between several solution sheaves of
holomorphic functions multi-asymptotically
developable along different kinds of families of closed complex submanifolds.
Such an application to $\mathcal{D}$-modules will be
given in a forthcoming paper.

Examples of our previous constructions are given in Section \ref{section:Examples}, in particular in $\CC^2$ and $\CC^3$. We also give, in the two manifolds case, a classification for asymptotics provided by the matrix associated to the multi-normal deformation.

We end this work with an Appendix in which we introduce the notion of multi-conic sheaves. Using o-minimal geometry we also obtain a decomposition theorem for subanalytic open sets needed for study of the geometry of sections of multi-conic sheaves.

%The plan of our paper is as follows: In Section 2, after establishment of
%some basic facts on a linearizable family $\chi$ of
%closed submanifolds,
%we quickly review the construction
%of multi-normal deformation $\widetilde{X}_\chi$
%along $\chi$. Then we explain the local model
%of multi-normal deformation in detail and give some examples.
%We also study the zero-section $S_\chi := \{t = 0\} \subset \widetilde{X}_\chi$
%of $\widetilde{X}_\chi$, where
%%As already stated in this introduction,
%$S_\chi$ generally has no vector bundle structure over the base space.
%we first introduce a suitable
%class of a family $\chi$ for which the zero-section $S_\chi$ becomes a vector bundle
%over the base space, and we give concrete structure of a vector bundle for such a $\chi$.
%
%We define, in Section 3, a multicone along $\chi$
%after the introduction of
%the semigroup generated by monomials associated with the actions
%$\mu_k$'s.
%Then

\section{Multi-normal deformation}\label{section:Multi-normal deformation}

In \cite{HP13} the notion of multi-normal deformation was introduced.
Here we consider a generalization where we replace the condition H2 with a weaker one.
We refer to \cite{KS90} for the classical normal deformation with respect to one submanifold.

\subsection{A simultaneously linearizable family of manifolds}\label{subsection:A simultaneously linearizable family of manifolds}

Let $X$ be a real analytic manifold with $\dim X = n$. Hereafter
all the maninfolds appearing in this paper
are assumed to be countable at infinity.
Let $\chi = \{M_1,\dots,M_\ell\}$ be a family of closed real analytic submanifolds of $X$.
Set $M := M_1 \cap \dots \cap M_\ell$.
Throughout this note, we always assume that
all the submanifolds are connected and
that $M$ is also connected.
Recall that, for closed submanifolds $M$ and $M'$, we say that
$M$ and $M'$ {\it{intersect cleanly}} at $q \in M \cap M'$
if and only if $M \cap M'$ is a closed submanifold near $q$ and
\begin{equation}
 (T^*_{M}X)_q + (T^*_{M'}X)_q =
 (T^*_{M \cap M'}X)_q
\end{equation}
holds.

\begin{df}{\label{df:sim-linear}}
We say that $\chi$ is {\it{simultaneously linearizable}} at $q \in M$
if $M_j$ and $M_k$ $(1 \le j,k \le \ell)$ intersect cleanly at $q$
and if
there exist a vector subspace $V$ of $(T^*X)_q$ and its decomposition
$V = V_1 \oplus \dots \oplus V_m$ by vector subspaces such that every $(T^*_{M_j}X)_q$ is
a direct sum of some $V_k$'s.

\end{df}

Note that
if $\chi$ satisfies the condition H2 of \cite{HP13}, then it is simultaneously linearizable.
\begin{prop}{\label{prop:canonical-from}}
Let $\chi = \{M_1,\dots,M_\ell\}$ be a family of closed submanifolds of $X$
and $q \in M$.
Then the following conditions are equivalent.
\begin{enumerate}
 \item[i.] $\chi$ is simultaneously linearizable at $q$.
 \item[ii.]
 There exist a neighborhood $W$ of $q$,
a neighborhood $U$ of $0$ with a system of local coordinates
$(x_1, \dots, x_n)$, the isomorphism $\psi: W \to U$ and
subsets $I_j$'s of $\{1, \dots, n\}$ $(j=1,\dots,\ell)$ such
that each $\psi(M_j \cap W)$ is defined by equations $x_i = 0$ $(i \in I_j)$.
\end{enumerate}
\end{prop}
The proof of the proposition is rather long, and hence it is given in Appendix \ref{D}.
Note that the decomposition
$V = V_1 \oplus \dots \oplus V_m$
appearing in Definition \ref{df:sim-linear}
is not unique. However,
if we consider minimal decomposition such that Definition \ref{df:sim-linear} holds,
the number of $V_k$'s of such a decomposition is unique,
and this minimal number $m$ is called {\it{ rank of the family $\chi$ at $q$}}.
Furthermore, when $\chi$ is simultaneously linearizable at any point in $M$,
as $M$ is connected by assumption,
the rank of the family $\chi$ does not depend on the point
under consideration.

%\begin{df}
%We say that a simultaneously linearizable family $\chi$ is of normal type
%if the rank of the family $\chi$ is equal to the number of submanifolds in $\chi$.
%\end{df}

By Proposition \ref{prop:canonical-from}, we have locally the family $\{I_j\}_{j=1}^\ell$
of subsets in $\{1,2,\dots,n\}$ so that each $M_j$ is defined
by $x_k = 0$ for $k \in I_j$. Let us consider the family $\{\widehat{I}_j\}_{j=1}^m$
of equivalence classes of $\bigcup_{j=1}^\ell I_j$ by the equivalence relation
\begin{equation}{\label{eq:equiavlent-def}}
i_1 \sim i_2 \iff \,\,
\text{"\,$i_1 \in I_j \iff i_2 \in I_j$ for all $j=1,2,\dots,\ell$\,"}.
\end{equation}
%Then $\chi$ is of normal type if and only if the number of equivalence classes defined above is equal to $\ell$.
Note that the number $m$ of equivalence classes coincides with the
the rank of the family $\chi$.
For convenience, we set
\begin{equation}
\widehat{I}_0 = I_0 := \{1,\dots,n\} \setminus (\bigcup_{j=1}^\ell I_j).
\end{equation}

\subsection{Geometry of multi-normal deformation}\label{subsection:Geometry of multi-normal deformation}

First recall the classical construction of \cite{KS90} of the normal deformation
of $X$ along $M_1$. We denote it by $\widetilde{X}_{M_1}$ and
we denote by $t_1\in\R$ the deformation parameter.
Let $\Omega_{M_1}=\{t_1>0\}$ and let us identify $\imin s(0)$ with $T_{M_1}X$.
We have the commutative diagram
\begin{equation}\label{eq:normal-deformation}
\xymatrix{T_{M_1}X \ar[r]^{s_{M_1}} \ar[d]^{\tau_{M_1}} & \widetilde{X}_{M_1} \ar[d]^{p_{M_1}} &
\Omega_{M_1} \ar[l]_{\ \ i_{\Omega_{M_1}}} \ar[dl]^{\widetilde{p}_{M_1}} \\
M \ar[r]^{i_{M_1}} & X. & }
\end{equation}
Set $\widetilde{\Omega}_{M_1}=\{(x;\,t_1)\;;\;t_1\neq 0\}$.
Let $\chi = \{M_1,\,\dots,\,M_\ell\}$ be the simultaneously
linearizable family given in the previous subsection, and define
\[
\widetilde{M}_2:=\overline{(p_{M_1}|_{\widetilde{\Omega}_{M_1}})^{-1}M_2}.
\]
Then $\widetilde{M}_2$ is a closed smooth submanifold of $\widetilde{X}_{M_1}$
as $\chi$ is simultaneously linearizable.\\

Now we can define the normal deformation along $M_1,M_2$ as
\[
\widetilde{X}_{M_1,M_2}:=(\widetilde{X}_{M_1})^\sim_{\widetilde{M}_2}.
\]
Then we can define recursively the normal deformation along
%$M_1,\ldots,M_\ell$
$\chi$ as
\[
\widetilde{X}=\widetilde{X}_{M_1,\ldots,M_\ell}:=(\widetilde{X}_{M_1,\ldots,M_{\ell-1}})^\sim_{\widetilde{M}_\ell}.
\]
Set $S_\chi=\{t_1,\ldots,t_\ell=0\}$, $M=\bigcap_{i=1}^\ell M_i$ and $\Omega_\chi=\{t_1,\ldots,t_\ell>0\}$. Then we have the commutative diagram
\begin{equation}\label{eq:multi-normal-deformation}
\xymatrix{S_\chi \ar[r]^s \ar[d]^\tau & \widetilde{X} \ar[d]^p &
\Omega_\chi \ar[l]_{\ \ i_\Omega} \ar[dl]^{\widetilde{p}} \\
M \ar[r]^{i_M} & X. & }
\end{equation}
In what follows, we always assume:
\begin{enumerate}
\item[A.1] $\chi$ is a simultaneously linearizable family of manifolds  at any point in $M$.
%\item[A.2]  The associated multi-action $\mu_\chi: \widetilde{X} \times (\mathbb{R}^+)^\ell
%	\to \widetilde{X}$ is non-degenerate.
\end{enumerate}

We denote $\mu_\chi: \widetilde{X} \times (\RP)^\ell \to \widetilde{X}$ the
multi-action associated with $\chi$.
Note that, as $\mu_\chi(S_\chi,\,\lambda) \subset S_\chi$ holds for $\lambda \in (\RP)^\ell$,  $\mu_\chi$ induces the action on $S_\chi$, which is denoted by
the same symbol hereafter. Then, for a generic point $p \in S_\chi$,
the orbit $\mu_\chi(p,\, (\RP)^\ell) \subset S_\chi$
becomes a real analytic manifold,  whose maximal dimension is called {\it{the dimension of
an orbit of $\mu_\chi$}}.

\begin{df}{\label{def:maximal_rank}}
We say that the multi-action $\mu_\chi$ has maximal rank if
the dimension of an orbit of $\mu_\chi$ is equal to either the number of
submanifolds in $\chi$ $($usually denoted by $\ell$$)$ or
the rank of the family $\chi$ $($usually denoted by $m$$)$.
In particular, $\mu_\chi$ is said to be non-degenerate if
the dimension of the orbit of $\mu_\chi$ coincides with $\ell$
and transitive if it coincides with $m$.
\end{df}
Then we have:
\begin{lem}
The following conditions are equivalent.
\begin{enumerate}
\item The multi-action $\mu_\chi$ is non-degenerate.
\item
$\mu_\chi(x,\lambda) = \mu_\chi(x, \lambda')$ for any $x \in S_\chi$ implies $\lambda = \lambda'$.
\end{enumerate}
\end{lem}

Now we give the local description of the multi-normal deformation if $\chi$ satisfies
the condition A.1.
Let $\{I_j\}_{j=1}^\ell$ be the family of subsets in $\{1,\dots,n\}$ in
Proposition \ref{prop:canonical-from}, and let $\{\widehat{I}_k\}_{k=1}^m$
be the equivalence classes in $\bigcup_{j=1}^\ell I_j$ by the equivalence relation
(\ref{eq:equiavlent-def}), where $m$ coincides with the rank of the family $\chi$. We define the monomial $\varphi_k$ $(k=1,\dots, m)$
of the variables
$t_1$, $\dots$, $t_\ell$ by
\begin{equation}{\label{eq:def-monomial}}
	\varphi_k(t_1,\dots,t_\ell) = \prod_{\{j \in \{1,\dots,\ell\};\,\widehat{I}_k \subset I_j\}} t_j.
\end{equation}
We also define the $\ell \times m$-matrix $A_\chi = (a_{jk})$ with
$$
a_{jk} =
\left\{
\begin{array}{ll}
1 \qquad & \widehat{I}_k \subset I_j, \\
0 \qquad & \text{otherwise}.
\end{array}
\right.
$$
Then we have, for $t = (t_1, \dots, t_\ell) \in (\mathbb{R}^+)^\ell$,
$$
(\log \varphi_1,\dots, \log \varphi_m) =
(\log t_1, \dots, \log t_\ell) A_\chi.
$$
%We also assume, for simplicity, the following
%\begin{enumerate}
%\item[A.3] The rank of $A_\chi$ is $m$ (that also implies $m \leq \ell$).
%\end{enumerate}
%Note that condition A.2 implies that the rank of $A_\chi$ is $\ell$
%(that also implies $\ell \leq m$). \\

Using these $\varphi_k$'s, the projection $p: \widetilde{X} \to X$ and the
action $\mu_\chi$ on $\widetilde{X}$ are locally defined as follows:
We take
\begin{equation}{\label{eq:local-coordinates}}
(x^{(0)}, x^{(1)}, \dots, x^{(m)}) \quad\text{and}\quad
(x^{(0)}, x^{(1)}, \dots, x^{(m)};\, t_1, \dots, t_\ell)
\end{equation}
the local coordinates of $X$ and those of $\widetilde{X}$ respectively. Here
$x^{(j)}$ $(j > 0)$ denotes the set of coordinates $x_i$'s ($i \in \widehat{I}_j$)
and $x^{(0)}$ denotes the rest of coordinates.

Then the projection $p$ is given by, for
$(x;\,t) = (x^{(0)},\, x^{(1)},\, \dots, x^{(m)};\, t_1, \dots,t_\ell) \in \widetilde{X}$,
$$
p(x;\,t) = \left(x^{(0)},\, \varphi_1(t) x^{(1)},\, \dots,\, \varphi_m(t)x^{(m)}\right).
$$
The action $\mu_\chi$ on $\widetilde{X}$ is defined locally by,
for $\lambda = (\lambda_1, \dots, \lambda_\ell) \in (\RP)^\ell$,
\begin{equation}\label{eq:action_mu_chi}
\mu_\chi((x;\,t), \lambda ) =
\left(x^{(0)},\, \varphi_1(\lambda) x^{(1)},\, \dots,\, \varphi_m(\lambda)x^{(m)};\,
t_1/\lambda_1,\, \dots,\, t_\ell/\lambda_\ell\right).
\end{equation}
It is the composition of the actions $\mu_j$ associated to the $I_j$'s, $j=1,\dots,\ell$, which are defined locally by,
for $\lambda_j \in \RP$,
\begin{equation}\label{eq:action_mu_j}
\mu_j((x;\,t), \lambda_j ) =
\left(x^{(0)},\, \lambda_j^{a_{j1}} x^{(1)},\, \dots,\,
\lambda_{j}^{a_{jm}}x^{(m)};\,
t_1,\, \dots,\, t_j/\lambda_j,\, \dots,\, t_\ell\right).
\end{equation}
%where $\lambda_{ji}=\lambda_j$ if $\widehat{I_i} \subseteq I_j$ and $\lambda_{ji}=1$ otherwise.

From these local observations, the following lemma easily follows:
\begin{lem}{\label{lem:maximal_rank}}
	$\mu_\chi$ has maximal rank if and only if $A_\chi$ has
	maximal rank, that is,
$$
\operatorname{Rank}\,A_\chi = \min\, \{\ell,\,m\}.
$$
Furthermore, $\mu_\chi$ is non-degenerate if $\operatorname{Rank}\,A_\chi = \ell$,
and it is transitive if $\operatorname{Rank}\,A_\chi = m$,
\end{lem}

%Let us consider the diagram \eqref{eq:multi-normal-deformation}.
%In local coordinates let $I_1,\ldots ,I_\ell \subseteq \{1,\ldots,n\}$ such that $M_i=\{x_k=0\;;\;k\in I_i\}$. For $j \in \{0,\dots,\ell\}$ set
%\[
%\JJ_j=\{k\in\{1,\ldots,\ell\}\;;\;\widehat{I}_j \subseteq I_k\}, \ \ \ \ t_{\JJ_j}=\prod_{k\in \JJ_j}t_k,
%\]
%where $t_1,\ldots,t_\ell\in\R$ and $t_{\JJ_0}=1$.
%Then $p:\widetilde{X}\to X$ is defined by
%\[
%(x^{(0)},x^{(1)}\ldots,x^{(\ell)};\,t_1,\ldots,t_\ell) \mapsto (t_{\JJ_0}x^{(0)},t_{\JJ_1}x^{(1)},\ldots,t_{\JJ_\ell}x^{(\ell)}).
%\]
Now we introduce an important geometrical object.
\begin{df} Let $Z$ be a subset of $X$.
The multi-normal cone to $Z$ along $\chi$ is the set
$
C_\chi(Z)=\overline{\imin{\widetilde{p}}(Z)} \cap S_\chi.
$
\end{df}

The following result has been proved in \cite{HPY}

\begin{lem} \label{lem: multi-normal cone} Let $p=(p^{(0)},p^{(1)},\dots,p^{(m)};\,0,\dots,0) \in S_\chi$ and let $Z \subset X$. The following conditions are equivalent:
\begin{enumerate}
\item $p \in C_\chi(Z)$.
\item There exist two sequences
$$
\left\{c_\kappa = (c_{1,\kappa},\dots,c_{\ell, \kappa})\right\}_{\kappa \in \mathbb{N}} \subset (\RP)^\ell, \quad
\left\{\left(q^{(0)}_\kappa,q^{(1)}_\kappa,\dots,q^{(m)}_\kappa\right)\right\}_{\kappa \in \mathbb{N}} \subset Z$$
such that $\varphi_j(c_\kappa) q_\kappa^{(j)} \to p^{(j)}$ $(j=1,\dots,m)$
and $c_{j,\kappa} \to +\infty$ $(j=1,\dots,\ell)$ when $\kappa \to \infty$.
\end{enumerate}
\end{lem}

\subsection{A simultaneously linearizable family of transitive type}\label{subsection:A simultaneously linearizable family of normal type}

We first study the case where the zero section $S_\chi$
of the multi-normal deformation $\widetilde{X}$ along $\chi$
becomes a vector bundle over the base space $M := M_1 \cap \cdots \cap M_\ell$.

\begin{df}
We say that a simultaneously linearizable family $\chi$ is of transitive type
if the associated multi-action $\mu_\chi$ is transitive.
Furthermore, it is said to be of normal type if
$\mu_\chi$ is non-degenerate and transitive.
\end{df}

Let $m$ be the rank of the family $\chi$, and let $\ell$ be the number of
submanifolds in $\chi$.
Note that
$\chi$ is of normal type if and only if,
under  the local description of the multi-normal deformation
given in the previous subsection,
the $A_\chi$ is a square matrix,
i.e., $\ell = m$ and it is invertible.
Furthermore, it follows from the definition that
$\ell \ge m$ holds when $\chi$ is of transitive type.

\begin{teo} \label{th:bundle-structure-of-zero-section}
Assume $\chi$ is of transitive type. Then $S_\chi$ becomes
a vector bundle over $M$.
\end{teo}
\begin{proof}
Suppose that $X$ has
a system of local coordinates blocks
$$
(x^{(0)},x^{(1)},\dots,x^{(m)})
$$
described in the previous subsection, where $m$ is the rank of the family $\chi$.
%satisfying the conditions described in Proposition \ref{prop:canonical-from}.
Let
$$
f(x) = (f^{(0)}(x),\, f^{(1)}(x), \dots, f^{(m)}(x))
$$
be a coordinate transformation on $X$
which sends each $M_j$ to $M_j$ $(j=1,\dots,\ell)$.
Let, for $k \ge 1$,
$$
f^{(k)}(x) = \sum_{\alpha \in \mathbb{Z}_{\ge 0}^n} c^{(k)}_{\alpha}x^\alpha \qquad
(c^{(k)}_{\alpha} \in \mathbb{R}^{\#\widehat{I}_k})
$$
be the Taylor expansion of $f^{(k)}$ at the origin and set
$f^{(k)}_{\alpha}(x) := c^{(k)}_{\alpha} x^\alpha$.

\

We first show the claim that
\begin{equation}{\label{eq:division-basic}}
\text{$f^{(k)}_\alpha(x^{(0)}, \varphi_1(t)x^{(1)}, \dots,
\varphi_m(t)x^{(m)})$ is
divided by the monomial $\varphi_k(t)$.}
\end{equation}
Let $\alpha$ with $c^{(k)}_\alpha \ne 0$.
By definition, $f^{(k)}(x)$ vanishes on $M_j$ for $j$ with
$\widehat{I}_k \subset I_j$, and hence,
$f^{(k)}_{\alpha}(x)$ also vanishes on $M_j$ for such a $j$, which entails
that there exists $i \in I_j$ with $\alpha_i > 0$.
Hence $f^{(k)}_{\alpha}(x^{(0)},\varphi_1(t)x^{(1)},\dots, \varphi_m(t)^{(m)})$
has the factor $t_j$ for $j$ with $\widehat{I}_k \subset I_j$, from which
the claim follows.

\

It follows from the claim \eqref{eq:division-basic} that
the induced map $\widetilde{f}$ of $f$ on $\widetilde{X}$ is well-defined
since their definitions are
$$
\begin{aligned}
\tilde{f}^{(k)}(x) &= \left.\dfrac{1}{\varphi_k(t)} f(x^{(0)}, \varphi_1(t)x^{(1)},\dots, \varphi_m(t)^{(m)})\right|_{t = 0} \qquad (k=1,2,\dots,m),\\
\tilde{f}^{(0)}(x) &= f(x^{(0)}, 0, \dots, 0),
\end{aligned}
$$
and $t_j = t_j$.
Now assume that there exist some $k \in \{1,\dots,m\}$ and
$\alpha \in \mathbb{Z}_{\ge 0}^n$ satisfying
$$
\varphi_k(t) = \prod_{j=1}^m (\varphi_j(t))^{\alpha^{(j)}},
$$
where we set $\alpha^{(j)} = \sum_{i \in \widehat{I}_j} \alpha_i$.
Let $e_k$ be the unit vector with the $k$-th element being $1$.
The above equation implies
$$
A_\chi e_k = A_\chi {}^t(\alpha^{(1)}, \dots, \alpha^{(m)})
$$
which is equal to
$$
A_\chi ({}^t(\alpha^{(1)}, \dots, \alpha^{(m)}) - e_k) = 0.
$$
Since $\operatorname{Rank} A_\chi = m$ holds, we have
$({}^t(\alpha^{(1)}, \dots, \alpha^{(m)}) - e_k) = 0$, that is,
allowable $\alpha$'s are multi-indices
in which the only one element $\alpha_i$ ($i \in \widehat{I}_k$) is equal to 1.
Therefore the associated map of $f$ on $\widetilde{X}$ is the one of vector bundles
with the fiber coordinates $(x^{(1)},\dots, x^{(m)})$ whose transformation law is given by
the block diagonal map
$$
\dfrac{\partial f^{(1)}}{\partial x^{(1)}}(x^{(0)},0,\dots,0)
\oplus \dots \oplus
\dfrac{\partial f^{(m)}}{\partial x^{(m)}}(x^{(0)},0,\dots,0).
$$
Here $\dfrac{\partial f^{(k)}}{\partial x^{(k)}}$ denotes the Jacobian matrix
of the map $f^{(k)}$ with respect to variables $x_i$'s ($i \in \widehat{I}_k$).
\end{proof}

We continue to assume $\chi$ to be of transitive type.
We now give an explicit description of $S_\chi$. We first define $m$ subsets
$B_1$, $\dots$, $B_m$ of $\{1,\dots,\ell\}$ in the following way. First
take a point $p \in M$ and take a system of local coordinates blocks $(x^{(0)},\,
x^{(1)},\, \dots, x^{(m)})$ of $X$ given in the previous subsection, where
$m$ is the rank of the family $\chi$.
Then we set, for $1 \le k \le m$,
$$
B_k := \{j \in \{1,\dots,\ell\};\, \widehat{I}_k \subset I_j\}.
$$
Since $M$ is connected,
a different choice of a point $p$ and that of coordinates gives
the same family of subsets $B_k$'s up to permutation.
Then we have:
\begin{prop}{\label{prop:bundle-description}}
For $k =1,2,\dots,m$, set
$$
N_k
:= \left\{
\begin{array}{ll}
X	     &\text{if $B_k = \{1,\dots,\ell\}$}, \\
\underset{j \notin B_k}{\bigcap} M_j \qquad &\text{otherwise.}
\end{array}
\right.
$$
Then $S_\chi$ is a direct sum of
the following vector bundles over $M$:
$$
\dfrac{TN_k \underset{X}{\times}{M}}
{\displaystyle\sum_{j \in B_k} \left(T(M_j \cap N_k) \underset{X}{\times}{M}\right) + TM}
\qquad (k=1,\dots,m).
$$
\end{prop}
\begin{proof}
Let $(x^{(0)},\, x^{(1)},\, \dots, x^{(m)})$
be a system of local coordinates blocks given in
the previous subsection, and let $(x^{(0)};\, \xi^{(1)},\, \dots, \xi^{(m)})$ be
that of $S_\chi$. Then the local coordinates of
$TN_k \underset{X}{\times}{M}$ is given by
$$
(x^{(0)};\, \xi^{(k')})_{k' \in \Xi}
$$
where $\Xi$ consists of indices $k'$'s such that $\widehat{I}_{k'} \cap I_j = \emptyset$
holds for any $j \notin B_k$.  Note that $\widehat{I}_k \cap I_j = \emptyset$ also holds
for any $j \notin B_k$.

Let $k' \in \Xi$ with $k' \ne k$. Then,
since $\widehat{I}_k$ and $\widehat{I}_{k'}$ are different equivalence classes,
it follows from the definition of the equivalence relation that
there exists $j \in B_k$ with $\widehat{I}_{k'} \cap I_j = \emptyset$.
As a matter of fact, if such a $j$ does not exist, then for any $j \in B_k$, we have
$\widehat{I}_{k'} \subset I_j$. As $\widehat{I}_k \subset I_j$ also holds for
any $j \in B_k$,
we have $\widehat{I}_k = \widehat{I}_{k'}$, which contradicts $k \ne k'$.

Therefore we conclude that the local coordinates of
$$
\dfrac{TN_k \underset{X}{\times}{M}}
{\displaystyle\sum_{j \in B_k} \left(T(M_j \cap N_k) \underset{X}{\times}{M}\right) + TM}
$$
are given by $(x^{(0)};\, \xi^{(k)})$.
Now let $f$ be a coordinate transformation considered in the proof of
Theorem \ref{th:bundle-structure-of-zero-section}. Then the corresponding
coordinate transformation of the above bundle is clearly given by the matrix
$\dfrac{\partial f^{(k)}}{\partial x^{(k)}}(x^{(0)},\,0,\,\dots,\,0)$, and that of $S_\chi$ is given by
$$
\dfrac{\partial f^{(1)}}{\partial x^{(1)}}(x^{(0)},0,\dots,0)
\oplus \dots \oplus
\dfrac{\partial f^{(m)}}{\partial x^{(m)}}(x^{(0)},0,\dots,0).
$$
as in the proof of Theorem \ref{th:bundle-structure-of-zero-section}.
Hence the bundle associated with $B_k$ is a subbundle of $S_\chi$ over $M$.
This completes the proof.
\end{proof}
We give some examples.
\begin{es}
Let $X=\mathbb{R}^n$ with coordinates $(x^{(0)},\, x^{(1)},\,x^{(2)},\,x^{(3)})$.
Define $M_1 = \{x^{(1)} = x^{(3)}= 0\}$,
$M_2 = \{x^{(2)} = x^{(3)}= 0\}$ and $M_3 = \{x^{(1)} = x^{(2)} = x^{(3)}= 0\}$.
Then
$$
S_\chi = T_MM_1 \oplus T_MM_2 \oplus \dfrac{TX \underset{X}{\times} M}
{TM_1 \underset{X}{\times} M + TM_2 \underset{X}{\times} M}.
$$
\end{es}
\begin{es}
Let $X=\mathbb{R}^n$ with coordinates $(x^{(0)},\, x^{(1)},\,x^{(2)},\,x^{(3)})$.
Define $M_1 = \{x^{(1)} = x^{(2)}= 0\}$,
$M_2 = \{x^{(2)} = x^{(3)}= 0\}$ and $M_3 = \{x^{(1)} = x^{(3)}= 0\}$.
Then
$$
S_\chi = T_MM_1 \oplus T_MM_2 \oplus T_MM_3.
$$
\end{es}

\begin{oss}
When $\chi$ is of transitive type,
the zero section $S_\chi$ becomes a vector bundle over $M$ as we have already seen.
However, in general, the simultaneously linearizable condition is not
enough to assure the existence of a vector bundle structure on
$S_\chi$. Another important exceptional case where
 $S_\chi$ has a vector bundle structure has been
studied in \cite{HPY}.
\end{oss}

If $\chi$ is not of transitive type, generally $S_\chi$ is not a vector bundle over
$M = \bigcap_k\, M_k$ as the following example shows:
\begin{es}{\label{es:non-vector-bundle}}
Let $X = \mathbb{R}^3$ with coordinates $(x_1,\, x_2,\, x_3)$, and
let $\chi = \{M_1,\, M_2\}$ be a family of closed submanifolds in $X$ defined by
$M_1 = \{x_2 = x_3 = 0\}$ and $M_2 = \{x_1 = x_3 = 0\}$. Then
$S_\chi$ is locally isomorphic to $\mathbb{R}^3$ with coordinates
$(\xi_1,\, \xi_2,\, \xi_3)$. Let $f = (f_1,\, f_2,\, f_3):
X \to Y$ be a coordinates transformation on $X$ to
its copy $Y$ with coordinates $(y_1, y_2, y_3)$
which sends $M_1$ and $M_2$ to their copies defined by the same equations
$\{y_2 = y_3 = 0\}$ and $\{y_1 = y_3 = 0\}$ respectively. Then the associated coordinates
transformation from $S_\chi$ to its copy $S_\chi$ with coordinates $(\eta_1,\, \eta_2,\, \eta_3)$ is given by
$$
\begin{aligned}
\eta_1 &= \dfrac{\partial f_1}{\partial x_1}(0) \xi_1, \\
\eta_2 &= \dfrac{\partial f_2}{\partial x_2}(0) \xi_2, \\
\eta_3 &= \dfrac{\partial f_3}{\partial x_3}(0) \xi_3 +
\dfrac{\partial^2 f_3}{\partial x_1 \partial x_2}(0) \xi_1 \xi_2.
\end{aligned}
$$
Hence $S_\chi$ is not a vector bundle over $M = \{0\}$.
\end{es}

\subsection{The local model}\label{subsection:The local model}

Throughout the paper we will often study the local model, i.e., the case in which $X=\R^n$ with coordinates $x=(x_1,\dots,x_n)$ and the submanifolds $M_j \in \chi$ ($j=1,\dots,\ell$) are defined by $\{x_i=0,\,i \in I_j\}$ with $I_j \subseteq \{1,\dots,n\}$. In this setting we are able to study a more general situation we are going to explain below.

First of all, we may allow $I_j=I_k$ if $j \neq k$. The action associated with $\chi$ is now defined by $\mu_j(x,\lambda)=(\lambda^{a_{j1}}x_1,\dots,\lambda^{a_{jn}}x_n)$ with $a_{ji}$ non-negative rational numbers, $a_{ji} \neq 0$ if $i \in I_j$, $a_{ji}=0$ otherwise. We introduce such a situation from now.

Let $m+1$ be the number of coordinates blocks,
that is, coordinates $(x)$ of $\mathbb{R}^n$ are divided into ($m+1$)-coordinates blocks
$$
(x^{(0)},\,x^{(1)}, \dots, x^{(m)}).
$$
Under these coordinates blocks, each $M_j$ is assumed to be linearized, that is,
for each $M_j$, there exists the subset $K_j \subset \{1,2,\dots,m\}$ such that
\begin{equation}
M_j = \{x^{(k)} = 0\,\, (k \in K_j)\} \qquad (j=1,2,\dots,\ell).
\end{equation}
Let $A_\chi = (a_{ij})$
be an $\ell \times m$ matrix with $a_{ij} \in \mathbb{Q}_{\ge 0}$, and
we take the positive rational number $\comD$ so that
$\comD a_{ij} \in \mathbb{Z}$ and all the $\comD a_{ij}$
($1 \le i \le \ell$, $1 \le j \le m$) have no common divisors.

Then we can define a more general normal deformation $\widetilde{X}=\R^n \times \R^\ell$ with the map $p:\widetilde{X} \to X$ defined by
\begin{equation}
p(x;\,t)=(x^{(0)},\,\varphi_1(t)^\comD x^{(1)},\, \dots,\, \varphi_m(t)^\comD x^{(m)})
\end{equation}
with
\begin{equation}
\varphi_k(t)=\prod_{j=1}^\ell t_j^{a_{jk}} \qquad (k = 1,2,\dots,m).
\end{equation}
Comparing with the matrix $A_\chi$, when $t \in (\RP)^\ell$ we have
$$
(\log\varphi_1,\dots,\log\varphi_m) = (\log t_1,\dots,\log t_\ell)A_\chi.
$$
In this setting, the action
$\mu: \widetilde{X}\times (\RP)^\ell \to \widetilde{X}$
is defined by the same equation as the one in \eqref{eq:action_mu_chi}, that is,
\begin{equation}
\mu((x;\,t), \lambda ) =
\left(x^{(0)},\, \varphi_1(\lambda) x^{(1)},\, \dots,\, \varphi_m(\lambda)x^{(m)};\,\,
\dfrac{t_1}{\lambda_1^{1/\comD}},\, \dots,\, \dfrac{t_\ell}{\lambda_\ell^{1/\comD}}\right),
\end{equation}
and $\mu_j: \widetilde{X}\times \RP \to \widetilde{X}$
for each $j$ is obtained by
\begin{equation}
\mu_j((x;\,t), \lambda_j ) =
\left(x^{(0)},\, \lambda_{j}^{a_{j1}} x^{(1)},\, \dots,\, \lambda_{j}^{a_{jm}}x^{(m)};
\,\, t_1,\, \dots,\, \dfrac{t_j}{1/\lambda_j^\comD},\, \dots,\, t_\ell\right)
\end{equation}
for $\lambda_j \in \RP$. Note that the actions $\mu$ and $\mu_j$ on $X$
(which are denoted by the same symbols as those on $\widetilde{X}$)
are given by
\begin{equation}
\begin{aligned}
\mu(x, \lambda ) &=
\left(x^{(0)},\, \varphi_1(\lambda) x^{(1)},\, \dots,\, \varphi_m(\lambda)x^{(m)}\right), \\
\mu_j(x, \lambda_j ) &=
\left(x^{(0)},\, \lambda_{j}^{a_{j1}} x^{(1)},\, \dots,\, \lambda_{j}^{a_{jm}}x^{(m)}
\right).
\end{aligned}
\end{equation}
Then we assume:
%\begin{itemize}
%\item The action $\mu$ on $\widetilde{X}$ is non-degenerate, i.e.,
%$\operatorname{Rank} A_\chi = \ell$.
%\end{itemize}
%And we also assume:
\begin{itemize}
\item[A2.] Each $M_j$ coincides with the set of fixed points of the action $\mu_j$ on $X$.
\end{itemize}
That is, we have
\begin{equation}
x \in M_j \iff \mu_j(x,\,\lambda_j) = x\,\, \text{for any $\lambda_j \in \RP$}.
\end{equation}
In terms of $A_\chi$, this condition is equivalent to saying that
\begin{equation}
a_{jk} \ne 0 \iff x^{(k)} = 0\,\, \text{on $M_j$}.
\end{equation}

\begin{es} Let us see some examples in $\R^3$ with variables $(x_1,x_2,x_3)$.
\begin{itemize}
\item Let $M_1=\{x_1=0\}$, $M_2=\{x_2=x_3=0\}$ and consider the action on $X$ by
$$
\mu(x,t)=(x_1t_1,x_2t_2,x_3t_2).
$$
Then we can take $x^{(1)}=x_1$, $x^{(2)}=(x_2,x_3)$ as the coordinates blocks
and the associated matrix is
$$
A_\chi=\left(
\begin{array}{cc}
1 & 0 \\
0 & 1
\end{array}
\right).
$$
\item Let $M_1=\{x_1=x_2=0\}$, $M_2=\{x_2=0\}$ and consider the action on $X$ by
$$
\mu(x,t)=(x_1t_1,x_2t_1t_2,x_3).
$$
Then we can take $x^{(1)}=x_1$, $x^{(2)}=x_2$, $x^{(3)}=x_3$ as the coordinates blocks and the associated matrix is
$$
A_\chi=\left(
\begin{array}{ccc}
1 & 1 & 0\\
0 & 1 & 0
\end{array}
\right).
$$
\item Let $M_1=M_2=\{0\}$ and consider the action on $X$ by
$$
\mu(x,t)=(x_1t_1^3t_2,x_2t_1^3t_2,x_3t_1^2t_2).
$$
Then we can take $x^{(1)}=(x_1,x_2)$, $x^{(2)}=x_3$ as the coordinates blocks
and the associated matrix is
$$
A_\chi=\left(
\begin{array}{cc}
3 & 2 \\
1 & 1
\end{array}
\right).
$$
\item Let $M_1=\{x_1=x_2=0\}$, $M_2=\{x_2=x_3=0\}$ and consider the action on $X$ by
$$
\mu(x,t)=(x_1t_1,x_2t_1t_2,x_3t_2).
$$
Then we can take $x^{(1)}=x_1$, $x^{(2)}=x_2$, $x^{(3)}=x_3$ as the coordinates blocks
and the associated matrix is
$$
A_\chi=\left(
\begin{array}{ccc}
1 & 1 & 0 \\
0 & 1 & 1
\end{array}
\right).
$$
\end{itemize}
\end{es}

Finally, in a local model, we use the same terminologies as in Definition \ref{def:maximal_rank} and Lemma \ref{lem:maximal_rank}. That is,
\begin{df}
We say that $\mu_\chi$ has {\it{maximal rank}}
if $\operatorname{Rank} A_\chi= \min\,\{\ell,m\}$
holds for the associated matrix $A_\chi$ of $\mu_\chi$. It is said to be
{\it{non-degenerate}} if
$\operatorname{Rank} A_\chi = \ell$ and {\it{transitive}} if $\operatorname{Rank} A_\chi = m$.
\end{df}

\section{Multicones for a general case}\label{section:Multi-cones for a general case}

 In this section we find a family of multicones associated to a multi-normal deformation. We mainly work in the local model introduced in Subsection \ref{subsection:The local model}.
%always assume A-1 and A-2. If $\chi$ is not of normal type, generally $S_\chi$ is not a vector bundle over
%$M = \cap M_k$.
%However $S_\chi$ has a conic structure induced from the action
%$\mu: (\mathbb{R}_+)^\ell \times \widetilde{X} \to \widetilde{X}$,
%and hence, we can define some important notions on $S_\chi$ from this action.

\subsection{Fixed points of $S_\chi$}\label{subsection:fixed points}

As we have seen in the previous section, if $\chi$ is of transitive type,
$S_\chi$ is the Cartesian product of vector bundles over $M$. Hence we can
consider the subbundle consisting of points
$\xi = \xi^{(1)}\times_M \cdots \times_M \xi^{(m)}$ where some of $\xi^{(k)}$'s
belong to the zero section of the corresponding vector bundles. Here we introduce
its counterpart for a general case.

For $a = (a_1, \dots, a_\ell) \in \mathbb{R}^\ell$, we define the action
$\mu^a: S_\chi \times \RP \to S_\chi$ by
$$
\mu^a(p,\, t) := \mu(p,\, t^{a_1},\dots,t^{a_\ell}) \qquad (t \in \RP).
$$
\begin{df}
We say that a point $p \in S_\chi$ is a fixed point (with respect to the action $\mu$)
if there exists
$a = (a_1, \dots, a_\ell) \in \mathbb{R}^\ell$ for which
$\mu^a$ is not the identity action on $S_\chi$
and $p$ is a fixed point of $\mu^a$,
that is, $\mu^a(\cdot,\,t) \ne \operatorname{id}_{S_\chi}$ for some $t > 0$
and $\mu^a(p,\, t) = p$ for any $t > 0$.
\end{df}

We give some examples.

\begin{es}
Let $\chi$ be a linearizable family of closed submanifolds.
We assume $\chi$ to be of transitive type.
Then $S_\chi$ has locally a system of coordinates blocks
$
(x^{(0)};\, \xi^{(1)},\, \dots,\, \xi^{(m)})
$
as described in the previous section,
where $m$ is the rank of the family $\chi$.
In this coordinates, the set of fixed points is given by
$
\{|\xi^{(1)}||\xi^{(2)}| \cdots |\xi^{(m)}| = 0\}
$.
Hence a point is located outside the set of fixed points if and only if
all the $\xi^{(k)}$'s are different from the origin.
\end{es}
\begin{es}
Let us consider the situation observed in Example {\ref{es:non-vector-bundle}}.
In this case, the set of fixed points is given by
$$
\{\xi_1 = \xi_2 = 0 \text { or } \xi_2 = \xi_3 = 0 \text{ or } \xi_1 = \xi_3 = 0\}.
$$
Hence a point is outside the set of fixed points if and only if
at least two $\xi_i$'s are non-zero.
\end{es}

Let us consider a local model described in Subsection \ref{subsection:The local model},
where $m+1$ denotes the number of coordinates blocks and $\ell$ designates the number
of manifolds in $\chi$.
Set $L := \operatorname{Rank} A_\chi$. Note that $L \le \min\,\{\ell,m\}$ holds.
\begin{lem}{\label{lem:outside-fixed-points-coordinate}}
Let $p = (x^{(0)},\,\xi^{(1)},\,\dots,\, \xi^{(m)}) \in S_\chi$ be a point outside
the set of fixed points.
Let $1 \le j_1 < \dots < j_L \le \ell$ and assume that
the $j_i$-th rows $(1 \le i \le L)$ are linearly independent.
Then there exist $1 \le k_1 < \dots < k_L \le m$
which satisfy the following conditions.
\begin{enumerate}
\item $\xi^{(k_i)} \ne 0$ $(1 \le i \le L)$.
\item The $L \times L$ submatrix made from the $j_i$-th rows and
the $k_i$-th columns $(1 \le i \le L)$ in $A_\chi$ is invertible.
\end{enumerate}
\end{lem}
\begin{proof}
By taking $j_i$-th rows $(1 \le i \le L)$ in $A_\chi$,
we may assume that the size of $A_\chi$ is $L \times m$. By
a permutation of columns, we may further assume $\xi^{(k)} = 0$ $(k > m')$ for
some $m'$. Let $\widetilde{A}$ be the $L \times m'$ submatrix of $A_\chi$ consisting of
the first $m'$-columns. If any choice of $L$-columns in $\widetilde{A}$ gives
 linearly dependent column vectors, then the rank of $\widetilde{A}$ is less
than $L$, and hence, we can find a non-zero row vector $a = (a_1,\dots, a_L) \in
\mathbb{R}^L$ with $a \widetilde{A} = 0$. Then the associated action $\mu^a$ is
different from $\operatorname{id}_{S_\chi}$ and $\mu^a(p,t) = p$ holds
for any $t > 0$. This contradicts the fact that $p$ is not a fixed point.
Hence we can find $L$-columns which satisfy the claim of the lemma.
\end{proof}

\begin{cor}{\label{cor:outside-fixed-points-coordinate}}
Let $p = (x^{(0)},\,\xi^{(1)},\,\dots,\, \xi^{(m)}) \in S_\chi$ be a point outside
the set of fixed points. Then,
by an appropriate permutation of coordinates blocks and that of
the parameters in the action $\mu$, we have the followings:
\begin{enumerate}
\item The first $L$-directions of $p$ are different from the origin, i.e.,
$\xi^{(k)} \ne 0$ $(1 \le k \le L)$.
\item The $L \times L$ submatrix made from the first $L$-columns and
the first $L$-rows in $A_\chi$ is invertible.
\end{enumerate}
\end{cor}

%In what follows, we consider a local model, and thus, we study a little bit general situation as specified below:
\subsection{Semigroups of generators}\label{subsection:Semigroups of generators}
From now on, we consider the problem under the local model described
in Section \ref{subsection:The local model}. Let $m+1$ be the number of coordinates blocks
and $\ell$ the number of submanifolds in $\chi$.
Set $L := \operatorname{Rank} A_\chi$ $(1 \le L \le \min\{\ell,m\})$.
By an appropriate permutation of coordinates blocks and that of parameters of the action
$\mu$, we assume that
\begin{equation}{\label{eq:L_L_invertibility}}
\begin{aligned}
&\text{the first $L \times L$ submatrix in $A_\chi = (a_{jk})$ is invertible,} \\
&\quad\text{that is, }\,\, \operatorname{det}\,\, (a_{jk})_{1\le j,k \le L} \ne 0.
\end{aligned}
\end{equation}
Take a point
\begin{equation}{\label{eq:point_xi_in_S_chi}}
p= (x^{(0)};\, \xi) =
(x^{(0)};\, \xi^{(1)},\,\dots,\, \xi^{(m)}) \in
S_\chi.
\end{equation}
Then we also assume that
\begin{equation}{\label{eq:xi_0_for_zero_vector}}
\begin{aligned}
&\text{if the $k$-th column of $A_\chi$ is a zero vector,} \\
&\quad\text{then the corresponding $\xi^{(k)}$ is zero.}
\end{aligned}
\end{equation}
By the assumption \eqref{eq:L_L_invertibility},
we may assume that $p$ is normalized as
\begin{equation}
\text{$|\xi^{(k)}| = 1$ for
$1 \le k \le L$ with $\xi^{(k)} \ne 0$}.
\end{equation}
Note that we cannot assume
this normalization for indices $k > L$.

Recall that the action on $X$ is given by
$$
(x^{(0)},\, \varphi_1(\lambda)x^{(1)},\dots,
\varphi_m(\lambda)x^{(m)})
$$
for $\lambda = (\lambda_1,\, \dots,\, \lambda_\ell)$,
where $m$-functions $\varphi_1$, $\dots$, $\varphi_m$ are monomials of $\lambda$
defined by
$$
\varphi_k(\lambda) = \lambda_1^{a_{1k}} \cdots \lambda_\ell^{a_{\ell k}}
\qquad (k=1,\dots,m).
$$
Here $a_{jk}$'s are entries of $A_\chi = (a_{jk})$ with $a_{jk} \in \mathbb{Q}_{\ge 0}$.

Let $\tau$ denote $m$-variables $(\tau_1,\,\dots,\,\tau_m)$, and we
denote by $\tau'$ the first $L$-variables of $\tau$ and by $\tau''$ the rest, that is,
we have $\tau = (\tau',\, \tau'')$ with
\begin{equation}
\tau' = (\tau_1, \dots, \tau_L)\,\, \text{ and }
\tau'' = (\tau_{L+1}, \dots, \tau_m).
\end{equation}
In the same way, we denote by $\lambda'$ the first $L$-variables of $\lambda$ and
by $\lambda''$ the rest, that is,
$\lambda = (\lambda',\, \lambda'')$ with
\begin{equation}
\lambda' = (\lambda_1, \dots, \lambda_L)\,\, \text{ and }
\lambda'' = (\lambda_{L+1}, \dots, \lambda_\ell).
\end{equation}
Consider the $L$-equations of $\lambda'$
\begin{equation}
\tau' = \varphi'(\lambda',\, \lambda''),
\end{equation}
where $\varphi'$ denotes $(\varphi_1,\, \dots,\, \varphi_L)$.
Then, by the assumption, the system of the equations can be solved
for the variables
$\lambda' = (\lambda_1, \dots, \lambda_L)$. That is, we can find
\begin{equation}{\label{eq:def_varphi_inv}}
\lambda' = \varphi^{-1}(\tau',\, \lambda'') =
( \varphi_1^{-1}(\tau',\,\lambda''),\,\dots,\,
\varphi_L^{-1}(\tau',\,\lambda''))
\end{equation}
which satisfies
\begin{equation}
\tau_k = \varphi_k(\varphi^{-1}(\tau',\,\lambda''),\,\lambda'')
\qquad (k = 1,\dots, L).
\end{equation}

\

Let $\mathcal{H}_\tau$
(resp. $\mathcal{H}_{\tau,\lambda}$)
be the set of non-zero monomials of rational powers of $\tau$
(resp. $\tau$ and $\lambda$)
with coefficients $1$, that is, each element has a form
\begin{equation}{\label{eq:form_monomials}}
\begin{aligned}
&\tau_1^{\alpha_1}\cdots\tau_m^{\alpha_m}
\quad(\alpha_1,\dots,\alpha_m \in \mathbb{Q}) \\
&\big(\text{resp. }\,\,
\tau_1^{\alpha_1}\cdots\tau_m^{\alpha_m}
\lambda_1^{\beta_1}\cdots\lambda_\ell^{\beta_\ell}
\quad(\alpha_1,\dots,\alpha_m,\,\beta_1,\dots,\beta_\ell \in \mathbb{Q})\big).
\end{aligned}
\end{equation}
For $f \in \mathcal{H}_{\tau,\lambda}$ given in \eqref{eq:form_monomials}, we denote by $\nu_k(f)$
(resp. $\nu^{\lambda}_j(f)$) the exponent of the variable $\tau_k$ (resp.
$\lambda_j$) in $f$, i.e.,
\begin{equation}
\nu_k(f) = \alpha_k\,\, \text{ and }\,\, \nu^{\lambda}_j(f) = \beta_j.
\end{equation}
%for $f \in \mathcal{H}_{\tau,\lambda}$.
Then
we regard $\RRTpair$ to be
a semigroup with its multiplication defined by,
for pairs $(f,\,v)$ and $(g,\,w) \in \RRTpair$
\begin{equation}
(f,\,v) (g,\,w) = (fg,\, vw).
\end{equation}
Let $\mathcal{A}$ and $\mathcal{B}$ be semigroups in $\RRTpair$.
Then we define
the semigroup $\mathcal{A} \vee \mathcal{B}$
as the minimal semigroup in $\RRTpair$
containing $\mathcal{A}$ and $\mathcal{B}$, that is,
\begin{equation}
\mathcal{A} \vee \mathcal{B} =
\{f,\, g,\, fg;\, f \in \mathcal{A},\, g \in \mathcal{B}\}.
\end{equation}
We also define the semigroup
\begin{equation}
\mathcal{A} * \mathcal{B} := \{f g;\, f \in \mathcal{A},\, g \in \mathcal{B}\}.
\end{equation}
For example, $\{(1,\,0)\}$ is a semigroup in $\RRTpair$ which consists
of only one element $(1,\,0)$. Then
$$
\mathcal{B} * \{(1,\,0)\} = \{(f,\,0);\, (f,v) \in \mathcal{B}\}.
$$
We can also regard $\RRpair \subset \RRTpair$ to be
a semigroup and define
the corresponding operations on $\RRpair$ induced from those on $\RRTpair$.

\

Set, for $p = (x^{(0)};\, \xi) = (x^{(0)};\,\xi^{(1)},\dots,\xi^{(m)}) \in S_\chi$ given in {\eqref{eq:point_xi_in_S_chi}},
\begin{equation}
J_Z(\xi) := \{k \in \{1,2,\dots,m\};\, \xi^{(k)} = 0\},
\end{equation}
and assume
\begin{equation}{\label{eq:def_k_1_k_q}}
\begin{aligned}
&J_Z(\xi) \cap \{1,2, \dots, L\} = \{k_1,\,k_2,\,\dots,\,k_q\} \\
&\qquad(1 \le k_1 < k_2 < \dots < k_q \le L),
\end{aligned}
\end{equation}
where $q := \#\left(J_Z(\xi) \cap \{1,2,\dots,L\}\right) \ge 0$.
We define, for a subset $A$ in $\RRTpair$,
\begin{equation}
%A_{\ge} = \left\{(f,\,v) \in A;\,
%f \in \mathcal{H}_\tau,\,\,
%\nu_{k}(f) \ge 0\,\,\text{for any $k \in J_Z(\xi) \cap \{1,\dots,L\}$}
%\right\},
A_{\ge} := \big\{(f,\,v) \in A;\,
\nu_{k_i}(f) \ge 0\,\,(1 \le i \le q)
\big\}
\, \bigcap\, (\RRpair),
\end{equation}
\begin{equation}
%A_{=} = \left\{(f,\,v) \in A;\,
%f \in \mathcal{H}_\tau,\,\,
%\nu_{k}(f) = 0\,\,\text{for any $k \in J_Z(\xi) \cap \{1,\dots,L\}$}
%\right\}
A_{=} := \big\{(f,\,v) \in A;\,
\nu_{k_i}(f) = 0\,\,(1 \le i \le q)
\big\}\, \bigcap\, (\RRpair),
\end{equation}
and
\begin{equation}
A_{>} := A_{\ge} \setminus A_{=} \,\,\subset (\RRpair).
\end{equation}
Note that, if $\mathcal{B}$ is a semigroup in $\RRTpair$,
then $\mathcal{B}_{\ge}$, $\mathcal{B}_{=}$ and $\mathcal{B}_{>}$
are semigroups in $\RRpair$.
For a subset $A \subset \RRTpair$,
we introduce its fraction ${\rm{Q}}(A) \subset \RRTpair$ by
\begin{equation}
A \,\bigcup\, \left\{(f^{-1},\,v^{-1})\right\}_{\{(f,\,v) \in A;\,v \ne 0\}}.
\end{equation}

\begin{df}
Let $\mathcal{A}$ and $\mathcal{B}$ be semigroups in $\RRpair$.
We say that $\mathcal{A}$ and $\mathcal{B}$ are equivalent if and only if
there exist natural numbers $N_1$ and $N_2$ with
\begin{equation}
\langle \mathcal{A}\rangle_{N_1}\,\, \subset\,\,
\langle \mathcal{B}\rangle_{N_2}\,\, \subset\,\,
\langle \mathcal{A}\rangle.
\end{equation}
Here $\langle \mathcal{A}\rangle_N$ denotes the semigroup
\begin{equation}
\langle \mathcal{A}\rangle_{N}\,\,
=\left\{(f,\,v)^N;\, (f,\,v) \in \mathcal{A}\right\} \subset \mathcal{A}.
\end{equation}
\end{df}

\begin{oss}
If
$\langle \mathcal{G} \rangle_N = \mathcal{A}$  $(N \in \mathbb{N})$ holds
for semigroups $\mathcal{A} \subset \mathcal{G}$,
we sometimes say that the radical of $\mathcal{A}$ is $\mathcal{G}$.
\end{oss}

Now recall that $\varphi^{-1}(\tau',\lambda'')$ was defined
in \eqref{eq:def_varphi_inv}. Then we define
\begin{equation}
\psi_k(\tau,\,\lambda'') :=
\dfrac{\tau_k}{\varphi_k(\varphi^{-1}(\tau',\,\lambda''),\,\lambda'')}
\qquad (L < k \le m).
\end{equation}
\begin{lem}{\label{lem:independent_psi_lambda}}
Let $L < k \le m$.
Then $\varphi_k(\varphi^{-1}(\tau',\,\lambda''),\,\lambda'')$ depends
only on the variables $\tau'$.
In particular, $\psi_k(\tau,\,\lambda'')$
does not depend on the variables $\lambda''$.
\end{lem}
\begin{proof}
The $k$-th column of $A_\chi$ is a linear combination of the first $L$-columns
in $A_\chi$, and thus, there exist $\alpha_{kj} \in \mathbb{Q}$ satisfying
$$
\varphi_k = \underset{1 \le j \le L}{\prod} \varphi_j^{\alpha_{kj}}
\qquad (L < k \le m).
$$
Hence we have
$$
\begin{aligned}
\varphi_k(\varphi^{-1}(\tau',\lambda''), \lambda'')
&= \underset{1 \le j \le L}{\prod}
\varphi_j(\varphi^{-1}(\tau',\lambda''),\lambda'')^{\alpha_{kj}} \\
& = \underset{1 \le j \le L}{\prod}\tau_j^{\alpha_{kj}}.
\end{aligned}
$$
\end{proof}
By the lemma,
$\psi_k(\tau,\,\lambda'')$
is denoted by $\psi_k(\tau)$ $(L < k \le m)$ in what follows.
Now we introduce the finite subset
$G$ in $\RRTpair$ defined by
\begin{equation}{\label{eq:def_M}}
G :=
Q\left(
\left\{\left(\varphi^{-1}_j,\,0\right)\right\}_{1 \le j \le L}\,
\bigcup\,
\left\{\left(\psi_k,\,|\xi^{(k)}|\right)\right\}_{L < k \le m}\,
\bigcup\,
\left\{\left(\lambda_j,\,0\right)\right\}_{L < j \le \ell}\right).
\end{equation}
\begin{oss}{\label{oss:general_G}}
In the definition of $G$, we use the fact that the point $\xi$ is normalized,
that is, $|\xi^{(k)}| = 1$ for $1 \le k \le L$ with $\xi^{(k)} \ne 0$.
When $\xi$ is not normalized, we define $G$ by
\begin{equation}\label{eq:def_M_general}
G :=
Q\left(
\left\{\left(\varphi^{-1}_j,\,0\right)\right\}_{1 \le j \le L}\,
\bigcup\,
\big\{\left(\psi_k,\,n_k(\xi)\right)\big\}_{L < k \le m}\,
\bigcup\,
\left\{\left(\lambda_j,\,0\right)\right\}_{L < j \le \ell}\right).
\end{equation}
Here we set
\begin{equation}
n_k(\xi) := \psi_k
\left(|\xi^{(*)}|_{+,1},\, |\xi^{(*)}|_{+,2},\,\dots,\,
|\xi^{(*)}|_{+,m}\right),
\end{equation}
where $|\xi^{(*)}|_{+,k} = |\xi^{(k)}|$ if $k \notin J_Z(\xi) \cap \{1,\dots,L\}$
and $|\xi^{(*)}|_{+,k} = 1$ otherwise.
\end{oss}
\begin{oss}
If the associated action $\mu$ has maximal rank, the definition of $G$
becomes much simpler. In fact, when the action $\mu$ is non-degenerate, we have
\begin{equation}
G =
Q\left(
\left\{\left(\varphi^{-1}_j(\tau'),\,0\right)\right\}_{1 \le j \le \ell}\,
\bigcup\,
\big\{\left(\psi_k(\tau),\,n_k(\xi)\right)\big\}_{\ell < k \le m}
\right).
\end{equation}
Note that, in this case, $G$ is contained in $\RRpair$, that is,
the variables $\lambda$ do not appear.
On the other hand, when the action $\mu$ is transitive, we have
\begin{equation}
G =
Q\left(
\left\{\left(\varphi^{-1}_j(\tau,\,\lambda''),\,0\right)\right\}_{1 \le j \le m}\,
\bigcup\,
\left\{\left(\lambda_j,\,0\right)\right\}_{m < j \le \ell}\right).
\end{equation}
In this case, the function $\psi_k$ does not appear.
\end{oss}

We then define the semigroup in $\RRpair$ by
\begin{equation}{\label{eq:def_S}}
\semiG :=
\left[G\right]_{=} \bigvee
\left(\left[G\right]_{>} * \left\{(1,\,0)\right\}\right).
\end{equation}
Here $[B]$ denotes the semigroup in $\RRTpair$ generated
by a subset $B$, that is, if $B$ consists of
$(f_s,\,v_s)$ $(s \in \Lambda)$ for a finite subset $\Lambda$, then
\begin{equation}
[B] :=
\left\{
\underset{s \in \Lambda}{\prod} (f_s,\,v_s)^{\alpha_s} \in
\RRTpair;\,
\alpha_s \in \mathbb{Z}_{\ge 0}\,\,\, (s \in \Lambda)\right\}.
\end{equation}
Note that $[B]$ contains the unit $(1,\,1)$ because we allow
all the $\alpha_s$ to be zero.

The semigroup $\semiG$ seems to depend on the choice of
linearly independent $L$-columns in $A_\chi$.
However, by the following lemma, we see that it is essentially independent
of such a choice for our purpose.
%Let $\mathcal{H}_{\tau,\lambda}$ be the set of non-zero rational monomials of
%the variables $\tau$ and $\lambda$ with their coefficients being $1$.
%Then we can regard $\mathcal{H}\times \mathbb{R}_{\ge 0}$ as
%a sub-semigroup of $\mathcal{H}_{\tau,\lambda} \times \mathbb{R}_{\ge 0}$.
Define the finite subset $\widehat{G} \subset \RRTpair$ by
\begin{equation}
Q\left(\left(\dfrac{\tau_1}{\varphi_1(\lambda)},\,|\xi^{(1)}|\right),\,\cdots,\,
\left(\dfrac{\tau_m}{\varphi_m(\lambda)},\,|\xi^{(m)}|\right),\,
(\lambda_1,\,0),\,\cdots,\, (\lambda_\ell,\,0)\right),
\end{equation}
and set
\begin{equation}
\widehat{\semiG} := (\RRpair)\, \bigcap\, [\,\widehat{G}\,].
\end{equation}
\begin{prop}{\label{prop:equivalence_G_hat_G}}
The semigroups $\semiG$ and $\widehat{\semiG}$ are equivalent.
To be more precise, we have
\begin{equation}
\langle \semiG \rangle_N\,\, \subset\,\, \widehat{\semiG}
\,\,\subset\,\, \semiG
\end{equation}
for some $N \in \mathbb{N}$.
%where $[\,\widehat{G}\,]$ is the semigroup generated by
%$\widehat{G}$ in $\RRTpair$.
\end{prop}
In this proof, we adopt \eqref{eq:def_M_general} as the definition of $G$, that is,
the normalization of $\xi$ is not assumed.
Set
$$
g_k(\tau,\lambda) := \dfrac{\tau_k}{\varphi_k(\lambda)} \quad (k=1,\dots,m).
$$
Then the proposition is a consequence of the following 3 lemmas.

\begin{lem}
We have
$$
(\RRpair)\, \bigcap\, [\,\widehat{G}\,] =
\widehat{\semiG}
\subset \semiG.
$$
\end{lem}

\begin{proof}
Take $(f,\,v) \in (\RRpair)\, \bigcap\, [\,\widehat{G}\,]$.
Then, by definition, there exist $\alpha_k, \beta_j \in \mathbb{Z}$
such that
$$
(f(\tau),\,v) = \underset{1 \le k \le m}{\prod} (g_k,\,|\xi^{(k)}|)^{\alpha_k}\,\,
\underset{1 \le j \le \ell}{\prod} (\lambda_j,\,0)^{\beta_j}.
$$
Here $\alpha_k \ge 0$ if $k \in J_Z(\xi)$ and $\beta_j \ge 0$ ($1 \le j \le \ell$).
Note that it follows from construction of $\widehat{G}$
that $\nu_k(f) \ge 0$ for any $k \in J_Z(\xi)$.
Since $f(\tau)$ does not depend on $\lambda$, by putting
$\lambda' = \varphi^{-1}(\tau', \lambda'')$
into the above equation, we obtain
%$f(\tau)$ is a product
%of $\varphi^{-1}_k$ ($1 \le k \le \ell$) and $\psi_k$ ($\ell < k \le m$).
%That is, there exist $\gamma_k \in \mathbb{Z}$ such that
$$
f = \underset{1 \le j \le L}{\prod} \varphi^{-1}_j(\tau',\,\lambda'')^{\beta_j}\,\,
\underset{L < k \le m}{\prod} \psi_k(\tau)^{\alpha_k}\,\,
\underset{L < j \le \ell}{\prod} \lambda_j^{\beta_j}.
$$
%where $\gamma_k \ge 0$ if $1 \le k \le \ell$ or $\ell < k$ with $k \in J_Z(\xi)$.
%In fact, we have $\gamma_k = \beta_k$ for $1 \le k \le \ell$ and
%$\gamma_k = \alpha_k$ for $\ell < k \le m$.
Set
$$
(f,\,v') = \underset{1 \le j \le L}{\prod} (\varphi^{-1}_j,\,0)^{\beta_j}\,\,
\underset{L < k \le m}{\prod} (\psi_k,\, n_k(\xi))^{\alpha_k}\,\,
\underset{L < j \le \ell}{\prod} (\lambda_j,\, 0)^{\beta_j}.
$$
Note $f \in \mathcal{H}_\tau$.
Then it follows from the definition of $\semiG$ that we have
$(f,\,v') \in \semiG$ if $\nu_k(f) = 0$
for all $k \in J_Z(\xi) \cap \{1,\dots,L\}$ and
$(f,\,0) = (f,\,v') * (1,0) \in \semiG$ if $\nu_k(f) > 0$
for some $k \in J_Z(\xi) \cap \{1,\dots,L\}$.
Hence, to show $(f,\,v) \in \semiG$, it suffices to prove $v = v'$
for the former case and $v = 0$ for the latter case.

Assume $\nu_k(f) > 0$ for some $k \in J_Z(\xi) \cap \{1,\dots,L\}$.
Since only $g_k$ contains the variable $\tau_k$,
the corresponding $\alpha_k$
must be greater than $0$, and thus, $v$ becomes zero. This implies
$(f,\,v) \in [G]_{>} * [(1,0)]$.

Assume $\nu_k(f) = 0$ for any $k \in J_Z(\xi) \cap \{1,\dots,L\}$.
Since $\alpha_k = 0$ for any $k \in J_Z(\xi) \cap \{1,\dots,L\}$, if
$v = 0$, then some $\beta_j$ or
some $\alpha_k$ with $k \in J_Z(\xi) \cap \{L+1,\,\dots,\,m\}$ is non-zero. Hence,
we have $v'=0$ and we conclude
$(f,\,v) \in [G]_{=}$ when $v = 0$.

When $v \ne 0$, we have $\beta_j = 0$ for $1 \le j \le \ell$ and
$\alpha_k = 0$ for $k \in J_Z(\xi)$, which implies $v' \ne 0$.
Then $(f,\,v) \in [G]_{=}$ immediately follows from the following elementary lemma.
Hence we obtained the desired inclusion.
\end{proof}

\begin{lem}{\label{eq:lemma-v}}
Let $(h, w) \in \semiG$
with $w \ne 0$.
If $h$ can be written in the form
$\underset{1 \le k \le m}{\prod} \tau_k^{\delta_k}$
for some $\delta_k \in \mathbb{Q}$ $(1 \le k \le m)$,
then we have
$w = \underset{1 \le k \le m}{\prod} |\xi^{(k)}|^{\delta_k}$.
The same fact holds for an element
$(h, w) \in \widehat{\mathcal{G}}$
with $w \ne 0$.
\end{lem}
\begin{proof}
Since $w\ne 0$, we have $(h,w) \in [G]_=$ and it is a
product of integer powers of $(\psi_k,\,n_k(\xi))$, i.e.,
there exist $\alpha_k \in \mathbb{Z}$ such that
\begin{equation}{\label{eq:tmp1-factor-lemma}}
(h,w) = \underset{L<k\le m}{\prod} (\psi_k,\,n_k(\xi))^{\alpha_k}.
\end{equation}
We temporarily regard $|\xi^{(1)}|$, $\dots$, $|\xi^{(m)}|$
as variables and, for $1 \le i \le m$,
we define $\nu^\xi_i$ in the similar way as that for $\nu_i$, i.e.,
$\nu^\xi_i(g)$
denotes the exponent of $|\xi^{(i)}|$ in a polynomial $g$ of variables
$|\xi^{(1)}|$, $\dots$, $|\xi^{(m)}|$. Then, by the definition of $n_k(\xi)$,
we have, for $1 \le k \le m$,
$$
\nu_i(\psi_k) = \nu^\xi_i(n_k(\xi)) \qquad (i \notin J_Z(\xi) \cap
\{1,\dots,L\})
$$
and
$$
\nu^\xi_i(n_k(\xi)) = 0 \qquad (i \in J_Z(\xi) \cap \{1,\dots,L\}).
$$
Then, by \eqref{eq:tmp1-factor-lemma}, we get
$$
\nu_i(h) = \nu^\xi_i(w(\xi)) \qquad (i \notin (J_Z(\xi) \cap \{1,\dots,L\}).
$$
Since $(h,w) \in [G]_=$ implies $\nu_i(h) = 0$ for
$i \in J_Z(\xi) \cap \{1,\dots,L\}$, we finally obtain
$$
\nu_i(h) = \nu^\xi_i(w(\xi)) \qquad (1 \le i \le m),
$$
from which we get the first claim.

Since $(h,w) \in \widehat{\mathcal{G}}$ with $w \ne 0$ implies
that $h$ is a product of integer powers of $\tau_k/\varphi_k(\lambda)$'s
$(1 \le k \le m)$, i.e., no factors $(\lambda_j,\,0)$ appear, we can easily conclude the second claim.
\end{proof}

Now we will prove the converse inclusion.
\begin{lem}
We have
$$
\langle \semiG\rangle_{N} \subset \widehat{\semiG} = (\RRpair)\, \bigcap\, [\,\widehat{G}\,]
$$
for some $N \in \mathbb{N}$.
\end{lem}
\begin{proof}
It follows from the relations
$$
\varphi^{-1}_j(\varphi'(\lambda', \lambda''), \lambda'')  = \lambda_j
\qquad (1 \le j \le L)
$$
that
$\dfrac{\varphi_j^{-1}(\tau',\,\lambda'')}{\lambda_j}$ ($1 \le j \le L$)
is expressed by a product of rational powers of $g_1$, $\dots$, $g_L$,
that is, there exist $\gamma_{jk} \in \mathbb{Q}$ such that
\begin{equation}\label{eq:varphi_lambda_product}
\dfrac{\varphi_j^{-1}(\tau',\,\lambda'')}{\lambda_j}
= \underset{1 \le k \le L}{\prod} g_k^{\gamma_{jk}}
%\,\, \underset{L < k \le \ell}{\prod} \lambda_k^{\kappa_{jk}}
\qquad (1 \le j \le L),
\end{equation}
equivalently we have
\begin{equation}\label{eq:varphi_lambda_product_1}
\varphi_j^{-1}(\tau',\,\lambda'')
= \left(\underset{1 \le k \le L}{\prod} g_k^{\gamma_{jk}}
%\,\, \underset{L < k \le \ell}{\prod} \lambda_k^{\kappa_{jk}}
  \right) \lambda_j
\qquad (1 \le j \le L).
\end{equation}
By using \eqref{eq:varphi_lambda_product},
we see that $\psi_k$ is also expressed by a product of rational powers of
$g_1$, $\dots$, $g_m$, that is,
there exist $\gamma'_{jk} \in \mathbb{Q}$ such that
\begin{equation}\label{eq:varphi_psi_product_1}
\psi_j
= \underset{1 \le k \le m}{\prod} g_k^{\gamma'_{jk}}
%\underset{L < k \le m}{\prod} \lambda_k^{\kappa'_{jk}}
\qquad (L < j \le m).
\end{equation}
Let $(f,\,v) \in \semiG$. By definition,
there exist $\alpha_k$ and $\beta_j \in \mathbb{Z}$
such that
$$
(f,\,v) =
\underset{1 \le k \le L}{\prod}(\varphi^{-1}_k,\,0)^{\alpha_k}\,\,
\underset{L < k \le m}{\prod}(\psi_k,\,n_k(\xi))^{\alpha_k}\,\,
\underset{L < j \le \ell}{\prod}(\lambda_j,\,0)^{\beta_j}
$$
holds if $\nu_k(f) = 0$ for any $k \in J_Z(\xi) \cap \{1,2,\dots,L\}$, and
$$
(f,\,v) =
\underset{1 \le k \le L}{\prod}(\varphi^{-1}_k,\,0)^{\alpha_k}\,\,
\underset{L < k \le m}{\prod}(\psi_k,\,n_k(\xi))^{\alpha_k}\,\,
\underset{L < j \le \ell}{\prod}(\lambda_j,\,0)^{\beta_j}\,\,*\,\, (1,\,0)
$$
holds if $\nu_k(f) > 0$ for some $k \in J_Z(\xi) \cap \{1,2,\dots,L\}$.
Here $\beta_j \ge 0$ for any $j$ and
$\alpha_k \ge 0$ if $1 \le k \le L$ or if $k > L$ with $k \in J_Z(\xi)$.

Then, by putting \eqref{eq:varphi_lambda_product_1} and \eqref{eq:varphi_psi_product_1}
into the above $f$,
we can find
$\delta_k \in \mathbb{Q}$ satisfying
$$
f = \underset{1 \le k \le m}{\prod} g_k^{\delta_k}
\underset{1 \le k \le L}{\prod} \lambda_k^{\alpha_k}
\underset{L < j \le \ell}{\prod} \lambda_j^{\beta_j}.
$$
Since only $g_k$ contains the variable $\tau_k$, we know
$\delta_k = \nu_k(f)$ ($1 \le k \le m$). In particular,
we have $\delta_k \ge 0$ for $k \in J_Z(\xi)$.
Set
$$
(f,\,v')
 = \underset{1 \le k \le m}{\prod} (g_k,\,|\xi^{(k)}|)^{\delta_k}
\underset{1 \le k \le L}{\prod} (\lambda_k,\,0)^{\alpha_k}
\underset{L < j \le \ell}{\prod} (\lambda_j,\,0)^{\beta_j}.
$$
By  definition,
if the natural number $N$ is taken to be the common denominator of
rationals $\gamma_{jk}$'s and $\gamma'_{jk}$'s,
we find $(f^N,\,v'^N) \in (\RRpair)\, \bigcap\, [\,\widehat{G}\,]$.
Hence to show $(f^N,\,v^N) \in (\RRpair)\, \bigcap\, [\,\widehat{G}\,]$,
it is sufficient to prove $v = v'$.

Now assume $v = 0$. If $\nu_k(f) > 0$ for some $k \in J_Z(\xi) \cap \{1,2,\dots,L\}$,
clearly we get $v' = 0$ because of $\delta_k = \nu_k(f) > 0$.
Hence we assume $\nu_k(f) = 0$
for any $k \in J_Z(\xi) \cap \{1,2,\dots,L\}$.
Then $v = 0$ implies either $\alpha_k > 0$ for some $1 \le k \le L$ or
$\beta_j > 0$ for some $L < j \le \ell$ or
$\alpha_k > 0$ for some $L < k \le m$ with $k \in J_Z(\xi)$.
The last case implies $\delta_k > 0$ for the same $k$ as $\tau_k$
only appears in $g_k$.
Hence we have $v' = 0$ and we conclude that $(f^N,\,v^N)
\in (\RRpair)\, \bigcap\, [\, \widehat{G}\,]$ in this case.

Next assume $v \ne 0$. Then we have $\beta_j = 0$ for $L < j \le \ell$ and
$\alpha_k = 0$ for any $1 \le k \le L$ or any
$L < k$ with $k \in J_Z(\xi)$. Furthermore, we have $\nu_k(f) = 0$
for any $k \in J_Z(\xi) \cap \{1,2,\dots,L\}$
because $v\ne 0$ means $(f,\,v) \in [G]_=$.
Hence, in particular, we obtain $\delta_k = \nu_k(f) = 0$ for $k \in J_Z(\xi)$ and
$$
(f,\,v') = \underset{1 \le k \le m, k \notin J_Z(\xi)}{\prod} (g_k,\,|\xi^{(k)}|)^{\delta_k},
$$
which implies $v' \ne 0$.  Then $v = v'$ follows from Lemma {\ref{eq:lemma-v}}.
This completes the proof.
\end{proof}

%Set
%\begin{equation}
%G_0 :=
%\left\{\left(\varphi^{-1}_j(\tau),\,0\right)\right\}_{1 \le j \le \ell}
%\bigcup
%\left\{\left(\psi_k(\tau),\,0\right)\right\}_{\ell < k \le m,\, k \in J_Z(\xi)}.
%\end{equation}
We can also show the following finiteness property.
\begin{prop}
$\semiG$ and $\widehat{\semiG}$ are finitely generated.
\end{prop}
\begin{proof}
We only show finiteness of $\semiG$ because
the proof goes in the same way for the case of $\widehat{\semiG}$.
We denote by $\pi: \RRpair \to \mathcal{H}_\tau$ the projection which forgets $v$
of $(f,\,v) \in \RRpair$.
Clearly we have $\pi(\semiG) = \pi([G]_\ge)$, and
it follows from Gordan's lemma that $\pi([G]_\ge)$ is finitely generated.

Hence it suffices to show
$$
\pi|_{\semiG}: \semiG \to \pi([G]_\ge)
$$
is an isomorphism of semigroups.
We make the inverse of $\pi|_\semiG$.
Let $h \in \pi([G]_\ge)$.
Then $h$ has the unique form
$$
h = (\varphi^{-1}_1)^{\alpha_1} \cdots (\varphi^{-1}_L)^{\alpha_L}
\psi_{L+1}^{\alpha_{L+1}} \cdots \psi_m^{\alpha_m}
\lambda_{L+1}^{\beta_{L+1}}\cdots \lambda_\ell^{\beta_\ell},
$$
where $\beta_j \in \mathbb{Z}_{\ge 0}$ $(L < j \le \ell)$,
$\alpha_k \in \mathbb{Z}_{\ge 0}$ for either $k \in J_Z(\zeta)$ or $1 \le k \le L$,
and $\alpha_k \in \mathbb{Z}$ otherwise. As a matter of fact, the uniqueness
comes from the condition \eqref{eq:L_L_invertibility}
and the fact that $\tau_k$ ($k > L$) is contained only in $\psi_k$.

Now we define
$$
v := \left\{
\begin{aligned}
&0 \qquad\qquad\qquad
\begin{aligned}
&\text{if $\alpha_k > 0$ for some $1 \le k \le L$} \\
&\text{or if $\beta_j > 0$ for some $L < j \le \ell$} \\
&\text{or if $\nu_k(h) > 0$ for some $k \in J_Z(\xi) \cap \{1,2,\dots,L\}$}, \\
\end{aligned}
\\
&\underset{L < k \le m}{\prod}n_{k}(\xi)^{\alpha_{k}} \qquad
 \text{otherwise}.
\end{aligned}
\right.
$$
Then we see $(h,\,v) \in \semiG$. It is easily checked that
this correspondence from $\pi([G]_\ge)$ to $\semiG$
gives a morphism of semigroup and becomes the inverse of $\pi|_{\semiG}$.
\end{proof}

The following lemma is useful through the paper.
\begin{lem}{\label{lem:fundamental_S}}
We have the followings:
\begin{enumerate}
\item For $k \in J_Z(\xi)$ and $(f,\,v)$ in $\semiG$ or
$\widehat{\semiG}$, we have $\nu_k(f) \ge 0$.
\item There exists $N \in \mathbb{N}$ such that
$(\tau_k^N,\,0)$ belongs to $\semiG$ and $\widehat{\semiG}$ $(1 \le k \le m)$.
\end{enumerate}
\end{lem}
\begin{proof}
Since $\semiG$ and $\widehat{\semiG}$ are equivalent,
it suffices to show the claims for $\widehat{\semiG}$.

The claim 1.~comes immediately from the definition of $\widehat{G}$
because every pair in $\widehat{G}$ does not contain a negative power of the variable
$\tau_k$ for $k \in J_Z(\xi)$.

Let us show the claim 2. If $\varphi_k = 1$,
we have $\xi_k = 0$ by the assumption \eqref{eq:xi_0_for_zero_vector},
and thus, $(\tau_k,\, 0) \in \widehat{\semiG}$.
If $\varphi_k \ne 1$,
by noticing $(\lambda_j,\,0) \in \widehat{G}$ ($j=1,\dots,\ell$)
and $(\tau_k/\varphi_k(\lambda),\, |\xi^{(k)}|) \in \widehat{G}$,
we have $(\tau_k^N,\,0) \in \widehat{\semiG}$ for sufficiently
large $N \in \mathbb{N}$.
This completes the proof.
\end{proof}

\subsection{Families of cones}\label{subsection:Families of cones}
Let $p = (x^{(0)};\, \xi) = (x^{(0)}; \xi^{(1)},\dots, \xi^{(m)})$ be a
point in $S_\chi$.
Set
\begin{equation}
\begin{aligned}
	|x^{(*)}| &:= (|x^{(1)}|,\,|x^{(2)}|,\,\cdots,\, |x^{(L)}|,\,\cdots,\,|x^{(m)}|), \\
|x^{(*)}|' &:= (|x^{(1)}|,\,|x^{(2)}|,\,\cdots,\, |x^{(L)}|).
\end{aligned}
\end{equation}
Let $H$ be a finite subset in $\RRpair$.
We denote by $C(\xi,\, H)$ a family of subsets of $X$ in the form
\begin{equation}
\begin{aligned}
& \quad \left\{(x^{(0)},x^{(1)},\dots,x^{(m)}) \in X;
\begin{aligned}
&x^{(k)} \in W_k\quad(k=1,\dots,m), \\
&v - \epsilon < g\left(|x^{(*)}|\right) < v + \epsilon \\
&\text{for $(g,\,v) \in H$}
\end{aligned}
\right\},
\end{aligned}
\end{equation}
where $\epsilon > 0$ and each $W_k$ satisfies either 1.~or 2.~below.
\begin{enumerate}
\item
$W_k$ is an $\mathbb{R}^+$-conic open subset
in $\mathbb{R}^{n_k} \setminus \{0\}$ containing the point $\xi^{(k)}$
if $k \notin J_Z(\xi)$.
\item $W_k$ is just $\mathbb{R}^{n_k}$ if $k \in J_Z(\xi)$.
\end{enumerate}

\begin{lem}
Assume that $[H]$ and $\semiG$
defined in \eqref{eq:def_S} are equivalent.
Then $C(\xi,\, H)$ is well-defined, and it is an
open subset and bounded with respect to
the variables $x^{(1)}$, $\dots$, $x^{(m)}$.
\end{lem}
\begin{proof}
it follows from the claim 1.~of Lemma \ref{lem:fundamental_S} and
the definition of $W_k$ that,
for any $(g,\,v) \in H$,
the function $g\left(|x^{(*)}|\right)$ is well-defined and
continuous on the region
$(x^{(1)},\,\dots,x^{(m)})\, \in W_1 \times \dots \times W_m$.
Hence $C(\xi,\,H)$ is an open subset.
Boundness is a consequence of the claim 2.~of the same lemma.
\end{proof}

Then we set
\begin{equation}
C(\xi) = \bigcup_{H} C(\xi,\,H),
\end{equation}
where $H$ runs through
the family of finite subsets in $\RRpair$
such that $[H]$ and $\semiG$ are equivalent.
Let $\pi: S_\chi \to M = M_1 \cap \dots \cap M_\ell$ be the canonical projection.

\begin{df}
We say that an open subset $V \subset X$ is a germ of a multicone at $p \in S_\chi$
if there exist an open neighborhood $\Omega$ of $\pi(p)$ in $X$ and $W \in C(\xi)$ satisfying
\begin{equation}
(W \cap \Omega) \subset V.
\end{equation}
\end{df}

\begin{oss}
As the semigroup $\semiG$ is finitely generated,
we can take
generators of $\semiG$ as $H$.
Then it is easy to see that
$C(\xi,\, H)$ becomes a cofinal family of $C(\xi)$.
Furthermore, since $\semiG$ and $\widehat{\semiG}$ are equivalent,
we can obtain the same family of multicones if
we replace $\semiG$ with $\widehat{\semiG}$ in the above definition.
\end{oss}

\

Remember $q := \#(J_Z(\xi) \cap \{1,\dots,L\})$.
We now give the explicit method to obtain a finite set $F^q$ in $\RRpair$ for which
$C(\xi, F^q)$ gives a cofinal family of $C(\xi)$.

Let $F$ be a finite set in $\RRTpair$ and $k \in \{1,2,\dots,L\}$.
We will consider an operation $\mathcal{L}_k$ on $F$ which generates
a finite subset in $\RRTpair$:
The subset $\mathcal{L}_k(F) \subset \RRTpair$ consists of the following elements:
\begin{enumerate}
	\item[F1.] $(f,\,v)$ for $(f,\,v) \in F$ with $\nu_k(f) = 0$.
	\item[F2.] $(f,\,0)$ for $(f,\,v) \in F$ with $\nu_k(f) > 0$.
	\item[F3.] $(f^{a}g^{b},\, v^{a}w^{b})$ where
		$(f,v)$ and $(g,w)$ in $F$
		with $\nu_k(f) > 0$ and $\nu_k(g) < 0$, and
		the natural numbers $a$ and $b$ are taken to be
		prime to each other with
		$a |\nu_k(f)| = b |\nu_k(g)|$.
\end{enumerate}

The following lemma is an easy consequence of the above operations.
\begin{lem}
For any $(f,v) \in \mathcal{L}_{k}(F)$,
we have $\nu_k(f) \ge 0$. Furthermore, for $(f,v) \in \mathcal{L}_k(F)$ with
$v \ne 0$, we have $\nu_k(f) = 0$, i.e., $f$ does not depend on $\tau_k$.
\end{lem}

Let $j \in \{1,2,\dots,\ell\}$.
We also introduce the similar operation $\mathcal{L}^\lambda_j(F)$ with respect to
the variable $\lambda_j$. The subset $\mathcal{L}^\lambda_j(F) \subset \RRTpair$
consists of the following elements:
\begin{enumerate}
	\item[F1.] $(f,\,v)$ for $(f,\,v) \in F$ with $\nu^\lambda_j(f) = 0$.
	\item[F3.] $(f^{a}g^{b},\,
		v^{a}w^{b})$ where
		$(f,v)$ and $(g,w)$ in $F$
		with $\nu^\lambda_j(f) > 0$ and $\nu^\lambda_j(g) < 0$, and
		the natural numbers $a$ and $b$ are taken to be prime
		to each other with
		$a |\nu^{\lambda}_j(f)| = b |\nu^{\lambda}_j(g)|$.
\end{enumerate}
Note that, in this case, there is no counterpart of the operation F2 in $\mathcal{L}_k$.
Hence, unlike the operation $\mathcal{L}_k$,
the variable $\lambda_j$ is
eliminated completely by the operation $\mathcal{L}^\lambda_j$.

Recall the finite subset $G \subset \RRTpair$ defined by \eqref{eq:def_M} or
\eqref{eq:def_M_general} when $\xi$ is not normalized. Then,
we apply the operation $\mathcal{L}^\lambda_j$'s ($L < j \le \ell$) successively to $G$,
we obtain the finite subset $F^0$ whose elements are independent of the variables
$\lambda$ because an element of $G$ is independent of the variables $\lambda'$
and the operation $\mathcal{L}^\lambda_j$ eliminates the variable $\lambda_j$ ($L < j \le \ell$), that is,
\begin{equation}
F^0 := (\mathcal{L}^\lambda_\ell \circ \cdots \circ \mathcal{L}^\lambda_{L+1})(G)\,\,
\subset \RRpair.
\end{equation}
Next we apply the operations $\mathcal{L}_{k_1}$, $\dots$,
$\mathcal{L}_{k_q}$ successively to $F^0$,
we get the subset $F^q \subset \RRpair$, i.e.,
\begin{equation}
F^q := (\mathcal{L}_{k_q} \circ \cdots \circ \mathcal{L}_{k_1})(F^0)\,\,
\subset \RRpair,
\end{equation}
where $q = \# (J_Z(\xi) \cap \{1,2,\dots,L\})$ and
$J_Z(\xi) \cap \{1,2,\dots,L\} = \{k_1,\dots,k_q\}$ with
$1 \le k_1 < \dots < k_q \le L$.
\begin{oss}
For convenience, we also set
\begin{equation}
F^{0,r} := (\mathcal{L}^\lambda_r \circ \cdots \circ \mathcal{L}^\lambda_{L+1})(G)\,\,
\subset \RRTpair \qquad (L < r \le \ell),
\end{equation}
and
\begin{equation}
F^s := (\mathcal{L}_{k_s} \circ \cdots \circ \mathcal{L}_{k_1})(F^0)
\subset \RRpair \qquad (1 \le s \le q).
\end{equation}
\end{oss}
\begin{oss}
If $\mu_\chi$ is non-degenerate, as $L = \ell$ holds, we get $F^0 = G$, i.e.,
an element of $G$ is independent of the variables $\lambda$.
Hence we have
\begin{equation}
F^q := (\mathcal{L}_{k_q} \circ \cdots \circ \mathcal{L}_{k_1})(G).
\end{equation}
On the other hand, when $p = (x^{(0)};\,\xi) \in S_\chi$
is outside fixed points, by Corollary \ref{cor:outside-fixed-points-coordinate}, we can
take coordinates blocks so that the conditions 1.~and 2.~of the corollary hold.
Then, under such coordinates blocks, since $J_Z(\xi) \cap \{1,2,\dots,L\}$ becomes
empty, we have $q = 0$ which implies
$$
F^q = F^0 = (\mathcal{L}^\lambda_\ell \circ \cdots \circ \mathcal{L}^\lambda_{L+1})(G).
$$
In particular, if $\mu_\chi$ is non-degenerate and $p$ is outside fixed points,
under the above coordinates blocks, we have
$$
F^q = F^0 = G.
$$
\end{oss}

The following lemma easily follows from the previous lemma.
\begin{lem}{\label{lem:non-zero-alpha}}
For any $(f,v) \in F^q$ and $k \in J_Z(\xi)$,
we have $\nu_{k}(f) \ge 0$.
Furthermore, for $(f,v) \in F^q$ with
$v \ne 0$ and $k \in J_Z(\xi)$, we have
$\nu_{k}(f) = 0$.
\end{lem}
$F^q$ is not a family of generators of $\semiG$ in general, however,
it enjoys the following good property from a geometrical point of view.
\begin{prop}\label{prop:radical-of-S}
\begin{enumerate}
\item We have $[F^q] \subset \semiG$. We also have
${\rm{Q}}(F^s) = F^s$ for any $1 \le s \le q$ and
${\rm{Q}}(F^{0,r}) = F^{0,r}$ for $L < r \le \ell$.
\item $[F^q]$ and $\semiG$ are equivalent. In particular,
the radical of $[F^q]$ is $\semiG$.
That is, there exists $N \in \mathbb{N}$ such that $(f,\,v)^N \in [F^q]$
for any $(f,\,v) \in \semiG$.
\end{enumerate}
\end{prop}
\begin{proof}
These claims can be proved by the induction on the number of operations
$\mathcal{L}^\lambda_j$'s and $\mathcal{L}_{k_i}$'s.
Since argument used in the inductive step on the number of
the operations $\mathcal{L}^\lambda_j$'s is
the same as that for the operations $\mathcal{L}_{k_i}$'s,
for simplicity, we show the claims only for the case $L = \ell$, i.e., $F^0 = G$.

It is easy to see that, for a subset $A \subset \RRTpair$, we have
$$
({\rm{Q}}\circ \mathcal{L}_k \circ {\rm{Q}})(A) = (\mathcal{L}_k \circ {\rm{Q}})(A).
$$
In what follows, we denote by $A_{>^{\tilde{q}}}$, $\semiG^{\tilde{q}}$ and etc.~the corresponding objects when $q = \tilde{q}$.

Let us show the claim 1.~by the induction on $q$.
When $q=0$, the claim clearly holds.
Suppose that the claim is true for $q = \tilde{q}$,
and we will show the claim when $q = \tilde{q}+1$.
We have
$$
\begin{aligned}
F^{\tilde{q}+1} &= \mathcal{L}_{k_{\tilde{q}+1}}(F^{\tilde{q}})
= (\mathcal{L}_{k_{\tilde{q}+1}} \circ {\rm{Q}})(F^{\tilde{q}})
= ({\rm{Q}} \circ \mathcal{L}_{k_{\tilde{q}+1}} \circ {\rm{Q}})(F^{\tilde{q}}) \\
&= ({\rm{Q}} \circ \mathcal{L}_{k_{\tilde{q}+1}})(F^{\tilde{q}})
 = {\rm{Q}}(F^{\tilde{q}+1}),
\end{aligned}
$$
which show the first claim of 1.~for $q = \tilde{q} + 1$.
Since we have
$$
F^{\tilde{q}+1} = \mathcal{L}_{k_{\tilde{q}+1}}(F^{\tilde{q}})
\subset \mathcal{L}_{k_{\tilde{q}+1}}(\semiG^{\tilde{q}})
$$
by the induction hypothesis and
since we can easily confirm
$\mathcal{L}_{k_{\tilde{q}+1}}(\semiG^{\tilde{q}})
\subset \semiG^{\tilde{q}+1}$,
we have obtained $[F^{\tilde{q}+1}] \subset \semiG^{\tilde{q}+1}$, which
shows the second claim of 1.~to be true for $q = \tilde{q}+1$.
Hence, by the induction, we have shown the claim 1.

\

Next we will show the claim 2.~by the induction on $q \ge 0$.
Assume first $q = 0$. In this case,
$
\semiG = [G] = [F^0]
$
holds, and hence, the claim 2.~is true.

Suppose that the claim 2.~is true for $q = \tilde{q} - 1 \ge 0$.
We will show the claim 2.~for $q = \tilde{q}$.
Let $(f,\,v) \in \semiG^{\tilde{q}}$ and assume
$F^{\tilde{q}-1} = \{(h_r,\,v_r)\}_r$.

Let us first consider the case $v \ne 0$. In this case,
we have $(f,\,v) \in [G]_{=^{\tilde{q}}}$, and hence,
$(f,\,v) \in \semiG^{\tilde{q}-1}$ also holds.
It follows from induction hypothesis that
there exists $N_{\tilde{q}-1} \in \mathbb{N}$
and $\{\alpha_r\}_r$ of non-negative integers such that
$$
(f,\,v)^{N_{\tilde{q}-1}} = \prod_r (h_r,\,v_r)^{\alpha_r}.
$$
Then $(f,\,v) \in [G]_{=^{\tilde{q}}}$ implies
$$
\sum_r \alpha_r \nu_{k_{\tilde{q}}}(h_r) = 0.
$$
Now let us consider the following procedure:
Choose an index $r^*$ such that $|\alpha_{r^*} \nu_{k_{\tilde{q}}}(h_{r^*})|$
is minimum in the set of non-zero $|\alpha_{r} \nu_{k_{\tilde{q}}}(h_{r})|$'s, and then,
find an index $r'$ and a positive rational number $\beta_{r'}$ such that
$$
0 < \beta_{r'} \le \alpha_{r'},\quad
\alpha_{r^*} \nu_{k_{\tilde{q}}}(h_{r^*})
+ \beta_{r'} \nu_{k_{\tilde{q}}}(h_{r'}) = 0.
$$
Note that such an index $r'$ and $\beta_{r'}$ necessarily exist by the choice of $r^*$.
By repeated applications of the procedure described above,
we can find a finite number of $(g_r,\,w_r) \in F^{\tilde{q}}$
with $\nu_{k_{\tilde{q}}}(g_r) = 0$
and positive rational numbers $\beta_r$ such that
$$
(f,\,v)^{N_{\tilde{q}-1}} = \prod_r (g_r,\,w_r)^{\beta_r},
$$
where each $g_r$ is generated by the operations F1 or F3.

Next we consider the case $v = 0$.
By induction hypothesis,
there exist $N_{\tilde{q}-1} \in \mathbb{N}$
and a family $\{\alpha_r\}_r$ of non-negative integers such that we have either
$$
(f,\,v)^{N_{\tilde{q}-1}} = \prod_r (h_r,\,v_r)^{\alpha_r} * (1,0), \qquad
\sum_r \alpha_r \nu_{k_{\tilde{q}}}(h_r) > 0
$$
if $\nu_{k_i}(f) = 0$ for $1 \le i \le \tilde{q}-1$ and
$\nu_{k_{\tilde{q}}}(f) >0$ hold, or
$$
(f,\,v)^{N_{\tilde{q}-1}} = \prod_r (h_r,\,v_r)^{\alpha_r}, \qquad
\sum_r \alpha_r \nu_{k_{\tilde{q}}}(h_r) \ge 0
$$
otherwise.
In both the above cases, by applying the same argument, we can find
a finite number of $(g_r,\,w_r) \in F^{\tilde{q}}$ with
$\nu_{k_{\tilde{q}}}(g_r) = 0$,
$(\tilde{g}_s,\,0) \in F^{\tilde{q}}$ with $\nu_{k_{\tilde{q}}}(\tilde{g}_s) > 0$,
positive rational numbers $\beta_r$ and $\tilde{\beta}_s$ such that
$$
(f,\,v)^{N_{\tilde{q}-1}} =
\prod_r (g_r,\,w_r)^{\beta_r} \prod_s (\tilde{g}_s,\, 0)^{\tilde{\beta}_s},
$$
where each $g_r$ is generated by the operation F1 or F3 and
each $\tilde{g}_s$ is a consequence of the operation F2.

As a conclusion, for each $(f,\,v) \in \semiG^{\tilde{q}}$,
there exists $N \in \mathbb{N}$ such that $(f,\,v)^N$ belongs to $[F^{\tilde{q}}]$.
Furthermore, since $\semiG^{\tilde{q}}$ is finitely generated,
we can take such an $N$ uniformly.
Therefore the claim 2.~is true for $q = \tilde{q}$, and thus, it is true for all $q$.
This completes the proof.
\end{proof}

As an immediate corollary, we obtain:
\begin{cor}
$C(\xi,\,F^q)$ gives a cofinal family of $C(\xi)$.
\end{cor}

\begin{es}
Now let us consider the following simple example.
Let $m=3$ and $\ell=2$ where
the action $\mu$ is defined by
$$
(x^{(0)},\,\lambda_1 x^{(1)},\, \lambda_2x^{(2)},\, \lambda_1 \lambda_2 x^{(3)})
\qquad \lambda = (\lambda_1, \lambda_2) \in (\RP)^2.
$$
Hence we have
$$
\varphi_1(\lambda) = \lambda_1,\,\,
\varphi_2(\lambda) = \lambda_2,\,\,
\varphi_3(\lambda) = \lambda_1\lambda_2
$$
and
$$
\varphi^{-1}_1(\tau) = \tau_1,\,\, \varphi^{-1}_2(\tau) = \tau_2.
$$
Take $p = (0;\,0,0,\xi^{(3)}) \in S_\chi$ with $\xi^{(3)} \ne 0$.
Then the initial set $F^0$ is given by
$$
\left\{(\tau_1,\,0),\,\, (\tau_2,\,0),\,\,
(\tau_3/(\tau_1\tau_2),\, |\xi^{(3)}|),\,\,
((\tau_1\tau_2)/\tau_3,\, 1/|\xi^{(3)}|)
\right\}.
$$
Therefore $F^1$ is given by
$$
\left\{(\tau_1,\,0),\,\, (\tau_2,\,0),\,\, ((\tau_1 \tau_2)/\tau_3,\,0),\,\,
(\tau_3/\tau_2,\, 0),\,\, (1,1)\right\}.
$$
Then $F^2$ is
$$
\left\{(\tau_1,\,0),\,\, (\tau_2,\,0),\,\, ((\tau_1 \tau_2)/\tau_3,\,0),\,\,
(\tau_3,\, 0),\,\,(1,1) \right\}.
$$
Hence the most important inequality of a multicone in this case is
$$
0-\epsilon < \dfrac{|x^{(1)}||x^{(2)}|}{|x^{(3)}|} < 0+\epsilon
\,\,\iff \,\,
|x^{(1)}||x^{(2)}| < \epsilon |x^{(3)}|.
$$
As a conclusion, a cofinal family of multicones is given by
$$
\left\{
(x^{(0)},\, x^{(1)},x^{(2)},x^{(3)});\,
\begin{aligned}
&x^{(k)} \in W_3, \\
&|x^{(1)}| < \epsilon,\, |x^{(2)}| < \epsilon, \\
& |x^{(1)}||x^{(2)}| < \epsilon |x^{(3)}|
\end{aligned}
\right\},
$$
where $\epsilon > 0$ and $W_3$ is a proper convex cone containing the direction $\xi^{(3)}$.
\end{es}

We have the following characterization of the multi-normal cone by
using multicones thus defined.
\begin{teo}\label{teo:multi-normal-cone-multi-cone}
Let $p=(x^{(0)};\, \xi) \in S_\chi$ and $Z$ a subset in $X$.
Then the following conditions are equivalent.
\begin{enumerate}
\item[i)] $p \notin C_\chi(Z)$.
\item[ii)] There exist a $V \in C(\xi)$ and an open neighborhood $U$
	of $\pi(p)$ in $X$ such that $Z \cap V \cap U = \emptyset$.
\item[iii)] There exist a $V \in C(\xi,\, F^q)$ and an open neighborhood $U$
	of $\pi(p)$ in $X$ such that $Z \cap V \cap U = \emptyset$.
\end{enumerate}
\end{teo}
To prove the theorem, it suffices to show the following equivalence:
\begin{enumerate}
\item $p \in C_\chi(Z)$.
\item $Z \cap V \cap U$ is non-empty for any $V \in C(\xi)$ and
any open neighborhood $U$ of $\pi(p)$.
\item $Z \cap V \cap U$ is non-empty for any $V \in C(\xi,\, F^q)$ and
any open neighborhood $U$ of $\pi(p)$.
\end{enumerate}
We will prove these equivalences step by step in several lemmas below.
Clearly the implication from 2.~to 3 holds. Hence
we first show the implication from 1.~to 2.
\begin{lem}
If $p \in C_\chi(Z)$, then
$Z \cap V \cap U$ is non-empty for any $V \in C(\xi)$ and
any open neighborhood $U$ of $\pi(p)$.
\end{lem}
\begin{proof}
Let
$$
p_i = (\tilde{x}_i,t_i) =
(\tilde{x}^{(0)}_i, \tilde{x}^{(1)}_i, \dots, \tilde{x}^{(m)}_i, t_{1,i},
\dots, t_{\ell,i}) \in \widetilde{X}
\qquad (i=1,2,\dots)
$$
be a sequence satisfying that $\widetilde{p}(p_i) \in Z$ and
$$
\tilde{x}^{(0)}_i \to x^{(0)},\quad
\tilde{x}^{(k)}_i \to \xi^{(k)}\,\, (k=1,\dots,m),\quad
t_{j,i} \to 0\,\, (j=1,\dots,\ell),
$$
when $i \to \infty$. Set
$$
x_i = (x_i^{(0)},\,x_i^{(1)},\,\dots,\,x_i^{(m)}) := \tilde{p}(p_i) \in Z \subset X.
$$
By Proposition \ref{prop:equivalence_G_hat_G},
it is enough to prove, for any $(f,\,v) \in (\RRpair) \cap [\,\widehat{G}\,]$,
\begin{equation}{\label{eq:xxx-1}}
f(|x_i^{(1)}|,\,\dots,\, |x_i^{(m)}|) \to v
\qquad (i \to \infty).
\end{equation}
Set $g_k := \tau_k/\varphi_k(\lambda)$ $(k=1,2,\dots,m)$. Then, by  definition
of $\widehat{G}$, there exist $\alpha_k, \beta_j \in \mathbb{Z}$ such that
$$
f(\tau) = \underset{1 \le k \le m}{\prod} g_k(\tau,\lambda)^{\alpha_k}
\underset{1 \le j \le \ell}{\prod} \lambda_j^{\beta_j},
$$
where $\beta_j \ge 0$ and $\alpha_k \ge 0$ if $k \in J_Z(\xi)$.
Now, by putting $\tau_k = |x_i^{(k)}|$ and $\lambda_j = t_{j,i}^\comD$ into $f$,
we have
$$
f(|x_i^{(*)}|)
= \underset{1 \le k \le m}{\prod} g_k(|x_i^{(*)}|,\, t_i^\comD)^{\alpha_k}
\underset{1 \le j \le \ell}{\prod} t_{j,i}^{\comD\beta_j},
$$
where $t_i^\comD := (t_{1,i}^\comD,\,\dots,\, t_{\ell,i}^\comD)$.
Note that
$$
g_k(|x_i^{(*)}|,\, t_i^\comD) = \dfrac{\varphi_k(t_i)^\comD|\tilde{x}^{(k)}_{i}|}
{\varphi_k(t_i^\comD)} = |\tilde{x}^{(k)}_i| \to |\xi^{(k)}| \qquad (k=1,2,\dots,m)
$$
holds when $i \to \infty$.

If $v = 0$, then we have either $\beta_j > 0$ for some $j$ or
$\alpha_k > 0$ for some $k \in J_Z(\xi)$. Hence we get
$$
f(|x_i^{(*)}|)
= \underset{1 \le k \le m}{\prod} g_k(|x_i^{(*)}|,\, t_i^\comD)^{\alpha_k}
\underset{1 \le j \le \ell}{\prod} t_{j,i}^{\comD\beta_j} \to 0 = v.
$$

If $v \ne 0$, then we have $\beta_j = 0$ and $\alpha_k = 0$ for any $k \in J_Z(\xi)$.
Therefore we get
$$
f(|x_i^{(*)}|)
= \underset{1 \le k \le m, k \notin J_Z(\xi)}{\prod} g_k(|x_i^{(*)}|,\, t_i^\comD)^{\alpha_k}
\to
\underset{1 \le k \le m, k \notin J_Z(\xi)}{\prod} |\xi^{(k)}|^{\alpha_k}.
$$
It follows from Lemma {\ref{eq:lemma-v}} that the last term is equal to $v$.
This completes the proof.
\end{proof}

Let us show the implication from 3.~to 1.
We need to prepare several lemmas.
\begin{lem}{\label{lem:sequence_converge}}
Let $a_1, \dots, a_m \in \mathbb{R}^\times$ be non-zero real numbers and
$v_1, \dots, v_m \in \mathbb{R}_{\ge 0}$ non-negative real numbers,
and let $\{\kappa_{k,i}\}_{i=1}^\infty$
$(k=1,\dots,m)$ be $m$-sequences of
positive real numbers. Then the following two conditions
are equivalent.
\begin{enumerate}
\item There exists a sequence $\{\epsilon_i\}_{i=1}^\infty$ of positive
real numbers satisfying, when $i \to \infty$,
$$
\epsilon_i \to 0, \qquad (\epsilon_i)^{-a_k} \kappa_{k,i} \to v_k \quad (k=1,\dots,m).
$$
\item The $m$-sequences satisfy, when $i \to \infty$,
$$
\left\{
\begin{array}{ll}
\kappa_{p,i} \to 0  \quad &\left(p \in O_+ \cup N_+\right),
\\
\\
1/\kappa_{p,i} \to 0  \quad &\left(p \in N_-\right),
\\
\\
(\kappa_{p,i})^{-a_q}(\kappa_{q,i})^{a_p} \to
(v_{p,i})^{-a_q}(v_{q,i})^{a_p}
\quad &\left(
\begin{aligned}
&p \in O_+ \cup N_+ \, \text{ and}\\
&q \in O_- \cup N_-
\end{aligned}
\right),
\\
\\
(\kappa_{p,i})^{|a_q|}(\kappa_{q,i})^{-|a_p|} \to
(v_{p,i})^{|a_q|}(v_{q,i})^{-|a_p|}
\quad &\left(
\begin{aligned}
&p \in O_+ \cup N_+, q \in N_+ \text{ or }\\
&p \in O_- \cup N_-, q \in N_-
\end{aligned}
\right),
\end{array}
\right.
$$
where we set
$O_+ = \{k;\, v_k = 0,\, a_k > 0\}$,
$O_- = \{k;\, v_k = 0,\, a_k < 0\}$,
$N_+ = \{k;\, v_k > 0,\, a_k > 0\}$ and
$N_- = \{k;\, v_k > 0,\, a_k < 0\}$.
\end{enumerate}
\end{lem}
\begin{proof}
The implication from 1.~to 2.~is trivial. Let us show the converse implication.

First assume that $N_+ \cup N_-$ is empty.
By considering $(\kappa_{j,m})^{1/|a_j|}$, we may assume $a_j = \pm 1$.
Determine two sequences of positive real numbers
$$
c_m := \max_{p \in O_+}\, \kappa_{p,m}, \quad
d_m := \max_{q \in O_-}\, \kappa_{q,m}.
$$
Then, by the conditions, we have
$$
1/c_m \to \infty,\quad d_m \prec 1/c_m.
$$
Here $\alpha_m \prec \beta_m$ means $\alpha_m/\beta_m \to 0$ ($m \to \infty$).
Now we can find a sequence $\{\epsilon_m\}$ so that
$$
\epsilon_m \to 0,\quad d_m \prec 1/\epsilon_m \prec 1/c_m,
$$
which satisfies the required condition.

Now assume $N_+ \cup N_-$ to be non-empty.
If $N_+$ is non-empty, we take some $p \in N_+$ and
set $\epsilon_i := (\kappa_{p,i}/v_p)^{1/a_p}$. It is easy to see that 1.~holds
for $\{\epsilon_i\}$. The other case can be shown in the same way.
\end{proof}

Then we introduce the modified operation $\tilde{\mathcal{L}}_{k}$ of $\mathcal{L}_k$
and $\tilde{\mathcal{L}}^\lambda_j$ of $\mathcal{L}^\lambda_j$,
which generate conditions appearing in the above lemma.
For a finite subset $F$ in $\RRTpair$, we define some subsets
$F_{0,>}$, $F_{0,<}$, $F_{\times,>}$ and
$F_{\times,<}$ of $F$ by
$$
\begin{aligned}
F_{0,>} &:=
\{(f,v) \in F;\, v = 0,\, \nu_k(f) > 0\}, \\
F_{0,<} &:=
\{(f,v) \in F;\, v = 0,\, \nu_k(f) < 0\}, \\
F_{\times,>} &:=
\{(f,v) \in F;\, v > 0,\, \nu_k(f) > 0\}, \\
F_{\times,<} &:=
\{(f,v) \in F;\, v > 0,\, \nu_k(f) < 0\}.
\end{aligned}
$$
Then,
the set $\tilde{\mathcal{L}}_k(F)$ consists of the following elements:
\begin{enumerate}
	\item[$\rm{\tilde{F}1}$.] $(f,\,v)$ for $(f,\,v) \in F$ with $\nu_k(f) = 0$.
	\item[$\rm{\tilde{F}2}$.] $(f,\,0)$
		where $(f,v) \in F_{0, >} \cup F_{\times, >}$ or
	        $f = 1/g$ for some $(g,v) \in F_{\times, <}$.
	\item[$\rm{\tilde{F}3}$.] $(f^{a}g^{b},\, v^{a}w^{b}$) where we take a combination of pairs $(f,v)$ and $(g,w)$ in sets with indices $>$ and $<$ respectively,
		that is, $(f,v) \in F_{0,>} \cup F_{\times,>}$
		and $(g,w) \in F_{0,<} \cup F_{\times,<}$. And see below for
		$a, b \in \mathbb{N}$.
	\item[$\rm{\tilde{F}4}$.] $(f^{a}g^{-b},\,
		v^{a}w^{-b})$
		where a pair satisfies
	\begin{enumerate}
	\item $(f,v) \in F_{0,>} \cup F_{\times,>}$ and
		$(g,w) \in F_{\times,>}$.
	\item $(f,v) \in F_{0,<} \cup F_{\times,<}$ and
		$(g,w) \in F_{\times,<}$.
	\end{enumerate}
\end{enumerate}
Here the natural numbers $a$ and $b$ appearing in $\rm{\tilde{F}3}$ and
$\rm{\tilde{F}4}$ are taken to be prime to each other with
$
a | \nu_k(f)| = b | \nu_k(g) |.
$

The modified operation $\tilde{\mathcal{L}}^\lambda_j$ is defined in the similar way:
We first define subsets $F_{0,>}$, $F_{0,<}$, $F_{\times,>}$ and
$F_{\times,<}$ of $F$ by replacing $\nu_k(f)$ with $\nu^\lambda_j(f)$ in the above
definitions of corresponding subsets for $\tilde{\mathcal{L}}_k$.
Then, the set $\tilde{\mathcal{L}}^\lambda_j(F)$ consists of the following elements:
\begin{enumerate}
	\item[$\rm{\tilde{F}1}$.] $(f,\,v)$ for $(f,\,v) \in F$ with
		$\nu^\lambda_j(f) = 0$.
	\item[$\rm{\tilde{F}2}$.]
		$\left(\lambda_j^{a}\,f^b,\,0\right)$
		where $(f,v) \in F_{0, <} \cup F_{\times, <}$ or
	        $f = 1/g$ for some $(g,v) \in F_{\times, >}$, and
		the natural numbers $a$ and $b$
		are taken to be prime to each other
		with $|\nu^{\lambda}_j(f)| = a/b$.
	\item[$\rm{\tilde{F}3}$.] $(f^{a}g^{b},\,
		v^{a}w^{b})$
		where we take a combination
		of pairs $(f,v)$ and $(g,w)$
		in sets with indices $>$ and $<$ respectively,
		that is, $(f,v) \in F_{0,>} \cup F_{\times,>}$
		and $(g,w) \in F_{0,<} \cup F_{\times,<}$, and see below for
		$a, b \in \mathbb{N}$.
	\item[$\rm{\tilde{F}4}$.] $(f^{a}g^{-b},\, v^{a}w^{-b})$
		where a pair satisfies
	\begin{enumerate}
	\item $(f,v) \in F_{0,>} \cup F_{\times,>}$ and
		$(g,w) \in F_{\times,>}$.
	\item $(f,v) \in F_{0,<} \cup F_{\times,<}$ and
		$(g,w) \in F_{\times,<}$.
	\end{enumerate}
\end{enumerate}
Here the natural numbers $a$ and $b$ appearing in $\rm{\tilde{F}3}$ and
$\rm{\tilde{F}4}$ are taken to be prime to each other with
$
a | \nu^\lambda_j(f)| = b | \nu^\lambda_j(g) |.
$

\begin{oss}
The above $\rm{\tilde{F}1}$, $\rm{\tilde{F}3}$ and $\rm{\tilde{F}4}$
are the same as those in $\tilde{\mathcal{L}}_k$,
where we just replace $\nu_k(f)$ and $\nu_k(g)$ with
$\nu^\lambda_j(f)$ and $\nu^\lambda_j(g)$ in their definitions, respectively.
While $\rm{\tilde{F}2}$ in $\tilde{\mathcal{L}}^\lambda_j$
is different from the corresponding one in $\tilde{\mathcal{L}}_k$.
\end{oss}

Then we can easily confirm the following lemma:
\begin{lem}
\begin{enumerate}
\item Assume $F = {\rm{Q}}(F)$. Then we have
$\mathcal{L}_k(F) = \tilde{\mathcal{L}}_k(F)$.
\item Assume $F = {\rm{Q}}(F)$ and $(\lambda_j,\,0) \in F$.
Then we have $\mathcal{L}^\lambda_j(F) = \tilde{\mathcal{L}}^\lambda_j(F)$.
\end{enumerate}
\end{lem}
It follows from the lemma and facts
$$
G = {\rm{Q}}(G), \quad
Q\circ\mathcal{L}^\lambda_{j}\circ Q = \mathcal{L}^\lambda_{j}\circ Q, \quad
Q\circ\mathcal{L}_{k_i}\circ Q = \mathcal{L}_{k_i}\circ Q
$$
and
$$
Q\circ\tilde{\mathcal{L}}^\lambda_{j}\circ Q =
\tilde{\mathcal{L}}^\lambda_{j}\circ Q, \quad
Q\circ\tilde{\mathcal{L}}_{k_i}\circ Q = \tilde{\mathcal{L}}_{k_i}\circ Q
$$
that we get
$$
F^0 =
(\mathcal{L}^\lambda_{\ell} \circ \cdots \circ \mathcal{L}^\lambda_{L+1}) (G)
=
(\tilde{\mathcal{L}}^\lambda_{\ell} \circ \cdots \circ
\tilde{\mathcal{L}}^\lambda_{L+1}) (G)
$$
and
$$
F^q =
(\mathcal{L}_{k_q} \circ \cdots \circ
\mathcal{L}_{k_1}) (F^0)
=
(\tilde{\mathcal{L}}_{k_q} \circ \cdots \circ
\tilde{\mathcal{L}}_{k_1}) (F^0).
$$
As a conclusion, we may consider $F^q$ to be generated by the modified operations
$\tilde{\mathcal{L}}^\lambda_j$'s ($L < j \le m$) and
$\tilde{\mathcal{L}}_{k_i}$'s ($1 \le i \le q$).

\begin{proof}
Now we are ready to show the implication from 3. to 1.
We may assume $x^{(0)} = 0$.  We can find a sequence in $Z$
$$
x_i = (x^{(0)}_i, x^{(1)}_i, \dots, x^{(m)}_i) \in Z \qquad (i=1,2,\dots)
$$
which satisfies, when $i \to \infty$,
\begin{enumerate}
\item $x_i \to 0$.
\item $x_i^{(k)}/|x_i^{(k)}| \to \xi^{(k)}/|\xi^{(k)}|$ for $k \notin J_Z(\xi)$.
\item For any $(f,\,v) \in F^q$,
$$
f(|x_i^{(1)}|, \dots, |x_i^{(m)}|) \to v,\quad (i \to \infty).
$$
\end{enumerate}
By considering a subsequence, we may assume, for each $k \in \{1,2,\dots,m\}$,
either $x_i^{(k)} \ne 0$ for all $i$ or $x_i^{(k)} = 0$ for all $i$.
Note that latter case occurs only for $k \in J_Z(\xi)$.

Set $\tau_{k,i} = |x_i^{(k)}|$ for $1 \le k \le m$. Now we define
$\tilde{\tau}_{k,i}$ as follows:
For $k \in \{1,2,\dots,m\}$ with $\tau_{k,i} \ne 0$ ($i = 1,2,\dots$),
we set $\tilde{\tau}_{k,i} = \tau_{k,i}$.
For $k \in \{1,2,\dots,m\}$ with $\tau_{k,i}$ being identically zero ($i=1,2,\dots$),
we take, as a \{$\tilde{\tau}_{k,i}\}_{i=1}^\infty$,
a sequence of positive real numbers satisfying
$\tilde{\tau}_{k,i} \to 0$ ($i \to \infty$) and
\begin{equation}{\label{eq:xxx-xxx-1}}
f(\tilde{\tau}_{1,i}, \dots, \tilde{\tau}_{m,i}) \to v, \quad (i \to \infty)
\end{equation}
for any $(f,\,v) \in F^q$. Note that it is possible because
such a $k$ belongs to $J_Z(\xi)$ and we have $\nu_k(f) \ge 0$ for $k \in J_Z(\xi)$.
Furthermore, if $\nu_k(f) > 0$ for some $k \in J_Z(\xi)$, then the corresponding
$v$ must be zero by Lemma \ref{lem:non-zero-alpha}.
Therefore it suffices to take $\tilde{\tau}_{k,i}$
rapidly decreasing to $0$ ($i \to \infty$).

Recall that $J_Z(\xi) \cap \{1,2,\dots,L\}
= \{k_1, \dots, k_{q}\}$ ($1 \le k_1 < \dots < k_q \le L$).
Then, as $F^q$ is generated by
applying the operation $\tilde{\mathcal{L}}_{k_q}$ to $F^{q-1}$,
by the conditions \eqref{eq:xxx-xxx-1} and Lemma \ref{lem:sequence_converge},
we find a sequence $\{\epsilon_{k_{q},i}\}$ of positive real
numbers such that, when $i \to \infty$,
\begin{equation}{\label{eq:xxx-xxx-2}}
\epsilon_{k_{q},i} \to 0,\quad
f\left(\tilde{\tau}_{1,i},\, \dots,\,
\overset{\text{$k_{q}$-th}}{\left(\tilde{\tau}_{k_{q},i}/\epsilon_{k_{q},i}\right)},\, \dots,\, \tilde{\tau}_{m,i}\right) \to v,
\end{equation}
for any $(f,\,v) \in F^{q-1}$.
Hence, thanks to the conditions \eqref{eq:xxx-xxx-2},
by applying the lemma again,
we find a sequence $\{\epsilon_{k_{q-1},i}\}$ of positive real
numbers such that $\epsilon_{k_{q-1},i} \to 0$  ($i \to \infty$) and
$$
f\left(\tilde{\tau}_{1,i},\, \dots,\,
\overset{\text{$k_{q-1}$-th}}
{\left(\tilde{\tau}_{k_{q-1},i}/\epsilon_{k_{q-1},i}\right)},\, \dots,\,
\overset{\text{$k_q$-th}}
{\left(\tilde{\tau}_{k_{q},i}/\epsilon_{k_{q},i}\right)},\, \dots,\,
\tilde{\tau}_{m,i}\right) \to v
\quad (i \to \infty)
$$
for any $(f,\,v) \in F^{q-2}$. By repeated applications of the lemma in this way,
we finally find the $q$-sequences $\{\epsilon_{k,i}\}_{i=1}^\infty$
($k \in J_Z(\xi) \cap \{1,2,\dots,L\}$) of positive real numbers
which satisfies, when $i \to \infty$,
$$
\begin{array}{ll}
\epsilon_{k,i} \to 0\,\,\, &(k \in J_Z(\xi) \cap \{1,2,\dots,L\}),
\\
f(\hat{\tau}_{1,i},\, \dots,\, \hat{\tau}_{m,i}) \to v
\qquad& (f,\,v) \in F^0,
\end{array}
$$
where $\hat{\tau}_{k,i}$
denotes $\tilde{\tau}_{k,i}/\epsilon_{k,i}$
if $k \in J_Z(\xi) \cap \{1,2,\dots,L\}$
and $\hat{\tau}_{k,i} = \tilde{\tau}_{k,i}$ for other $k$.
Set
$$
\hat{\tau}_i = (\hat{\tau}_{1,i},\,\dots,\, \hat{\tau}_{m,i})\, \text{ and }\,
\hat{\tau}'_i = (\hat{\tau}_{1,i},\, \dots,\, \hat{\tau}_{L,i}).
$$
Now applying the same argument to the variables
$\lambda'' = (\lambda_{L+1},\dots,\lambda_\ell)$ repeatedly, we can find
sequences
$$
\lambda''_i = (\lambda_{L+1,i},\,\dots,\lambda_{\ell,i})
\quad (i=1,2,\dots)
$$
of positive real numbers such that, when $i \to \infty$,
$$
\lambda''_i \to 0,\qquad
f(\hat{\tau}_i,\, \lambda''_i) \to v
$$
for any $(f,\,v) \in G$.  Finally we define the sequence in $\widetilde{X}$ by
$$
\begin{aligned}
\widetilde{x}_i :=
\big(& x^{(0)}_i,\, \left(x^{(1)}_i/\hat{\tau}_{1,i}\right),\, \dots,\,
\big(x_i^{(L)}/\hat{\tau}_{L,i}\big), \\
&\big(x_i^{(L+1)}/\varphi_{L+1}(t_i)\big),\, \dots,\,
\big(x_i^{(m)}/\varphi_m(t_i)\big),\,
t_{1,i},\,\dots,\, t_{\ell,i}
\big) \in \widetilde{X},
\end{aligned}
$$
where
$$
\begin{array}{lll}
	t_{j,i} &:= \varphi^{-1}_j(\hat{\tau}'_{i}, \lambda''_i)^{1/\comD}
\qquad &(1 \le j \le L), \\
\\
t_{j,i} &:= \lambda_{j,i}^{1/\comD}
\qquad &(L < j \le \ell).
\end{array}
$$
%and
%$$
%\sigma_{k,i} := \varphi_k(t_{1,i},\,t_{2,i},\, \dots,\, t_{\ell,i})
%\qquad (L < k \le m).
%$$
Let $\hat{p} = (0;\, \hat{\xi}^{(1)}, \dots, \hat{\xi}^{(\ell)})$
be the normalized point of $p$. That is,
$$
\hat{\xi}^{(k)} :=
\left\{
\begin{array}{ll}
0
\qquad
&(k \in J_Z(\xi) \cap \{1,\dots,L\}), \\

\dfrac{\xi^{(k)}}{|\xi^{(k)}|}
\qquad
&(1 \le k \le L,\, \xi^{(k)} \ne 0), \\

\dfrac{n_k(\xi) \xi^{(k)}}{|\xi^{(k)}|}
\qquad &(k > L).
\end{array}
\right.
$$
Note that $p$ and $\hat{p}$ belong to the same orbit of $(\RP)^\ell$-actions
on $S_\chi$.
Then we find $\widetilde{p}(\widetilde{x}_i) = x_i \in Z$ and
$\widetilde{x}_i \to \hat{p}$ ($i \to \infty$), which implies
$\hat{p} \in C_\chi(Z)$.
Since $C_\chi(Z)$ is an $(\RP)^\ell$-conic subset, we have finally obtained
$p \in C_\chi(Z)$.  This completes the proof.
\end{proof}

By the same argument as that in the first part of the proof of the theorem, we see that
a multicone enjoys the following good property.
\begin{cor}\label{cor: stability under contraction}
Let $H$ be a finite subset in $\RRpair$ for which $[H]$ and $\semiG$
are equivalent.
Then every $V \in C(\xi,\, H)$ is stable under contraction induced by
the action $\mu_j$ $(j = 1,\dots,\ell)$. That is,
$\mu_j(x,\, \lambda_j) \in V$ holds for any $x \in V$ and any $0 < \lambda_j \le 1$.
\end{cor}
\begin{proof}
Let $(f,v) \in H$ and $1 \le j_0 \le \ell$. For simplicity, we set
$f(z) := f(|z^{(*)}|)$.
To show the corollary, it suffices to prove the following claim:
\begin{enumerate}
\item If $v \ne 0$, then $f(\mu_{j_0}(z,\lambda_{j_0})) = f(z)$.
\item If $v = 0$, there exist a non-negative rational number $\kappa$
such that $f(\mu_{j_0}(z,\lambda_{j_0})) = \lambda_{j_0}^\kappa f(z)$.
\end{enumerate}
Since $[H]$ is equivalent to $\widehat{\mathcal{G}}$ also,
by the definition, there exist $N \in \mathbb{Z}$, $\alpha_k \in \mathbb{Z}$ and
$\beta_j \in \mathbb{Z}_{\ge 0}$ satisfying
$$
f(\tau)^N = \underset{1 \le k \le m}{\prod} \left(\tau_k/\varphi(\lambda)\right)^{\alpha_k}
\underset{1 \le j \le \ell}{\prod} \lambda_j^{\beta_j}.
$$
Put
$
\tau = |(\mu_{j_0}(z,\,\lambda_{j_0}))^{(*)}|
$
and $\lambda = (1,\dots,1,\,\lambda_{j_0},\,1,\dots,1)$ into the above equation,
we have
$$
f(\mu_{j_0}(z,\,\lambda_{j_0}))^N = \lambda_{j_0}^{\beta_{j_0}} f(z)^N.
$$
By noticing that $v \ne 0$ implies $\beta_j = 0$ for any $j$, we
have obtained the claim.
\end{proof}

Now assume that $p = (x^{(0)};\,\xi) \in S_\chi$ is outside fixed points and that the action $\mu_\chi$ is non-degenerate.
Then, by Corollary \ref{cor:outside-fixed-points-coordinate}, we can take coordinates
blocks of a local model so that the conditions 1.~and 2.~of the corollary hold.
Furthermore, under this coordinates blocks, since $L = \ell$ and $q = 0$ hold,
we have $F^q = G$. Hence we have obtained the following corollary:
Let us recall the definition of $n_k(\xi)$ given in Remark \ref{oss:general_G}
(note that $n_k(\xi)$ is just $|\xi^{(k)}|$ when $\xi$ is normalized).
\begin{cor}
Under the situation described above,
the following family of subsets gives a cofinal family of $C(\xi)$:
$$
\left\{
(x^{(0)},\, x^{(1)},\dots,x^{(m)});\,
\begin{aligned}
&x^{(k)} \in W_k \qquad (k=1,\dots,m), \\
&\varphi^{-1}_k(|x^{(*)}|') < \epsilon \qquad (k \le \ell), \\
&n_k(\xi) - \epsilon < \dfrac{|x^{(k)}|}
{\varphi_k(\varphi^{-1}(|x^{(*)}|'))} < n_k(\xi) + \epsilon \quad
(k > \ell) \\
\end{aligned}
\right\},
$$
where $\epsilon > 0$, and $W_k$ runs
through open conic cones in $\mathbb{R}^{n_k}\setminus \{0\}$
containing the direction $\xi^{(k)}$ if $\xi^{(k)} \ne 0$ and
it is $\mathbb{R}^{n_k}$ if $\xi^{(k)}=0$.
\end{cor}

\subsection{A restriction condition}\label{subsection:A restriction condition}

Let $A=(\alpha_{i,j})_{\substack{1 \leq i \leq \ell \\ 1 \leq j \leq m}}$
be an $\ell \times m$ matrix with entries of
non-negative rationals. Set $L = \operatorname{rank}A$ ($1 \le L \le \min\{\ell,m\}$).
Assume that the first $L \times L$ sub-matrix $A'$ of $A$ is invertible.
Hereafter $\alpha_i$ denotes the $i$-th row of $A$ for $1 \leq i \leq \ell$.
%Set
%$$
%\lambda' = (\lambda_1,\dots,\lambda_L),\qquad
%\lambda'' = (\lambda_{L+1},\dots,\lambda_\ell)
%$$
%and
%$$
%\tau' = (\tau_1,\dots,\tau_L),\qquad
%\tau'' = (\tau_{L+1},\dots,\tau_m).
%$$

%Consider the equation of $\lambda'$
%$$
%\tau' = \varphi'(\lambda',\lambda''),
%$$
%where $\varphi = (\varphi', \varphi'')$.
%$$
%\tau_k = \varphi_k(\lambda',\lambda'')  \qquad (1 \le k \le L).
%$$
%Then we find the solution
%$$
%\lambda' = \varphi^{-1}(\tau',\,\lambda''),
%$$
%where $\varphi^{-1}$ consists of $L$ monomials of $\tau'$ and $\lambda''$, i.e.,
%$$
%\varphi^{-1}(\tau',\lambda'') = (\varphi_1^{-1}(\tau', \lambda''),\, \dots,\,
%\varphi_L^{-1}(\tau', \lambda'')).
%$$
%We also set
%$$
%\psi_k(\tau,\,\lambda'') = \dfrac{\tau_k}{\varphi_k(\varphi^{-1}(\tau',\lambda''),\, \lambda'')} \qquad (L < k \le m).
%$$

\begin{lem}\label{lem0} For $1 \leq i \leq \ell$, set $e^{\alpha_i}=(e^{\alpha_{i,1}},\dots,e^{\alpha_{i,m}})$.
\begin{itemize}
\item[(i)] Let $1 \le i \le L$. We have $\log\varphi_j^{-1}(e^{\alpha_i})|_{\lambda''=1} \geq 0$
for $j=1,\dots,L$.
\item[(ii)] Let $1 \le i \le \ell$. We have  $\log\psi_k(e^{\alpha_i})=0$ for $k=L+1,\dots,m$.
\end{itemize}
\end{lem}
\begin{proof}
Let $\beta_j=(\beta_{j,1},\dots,\beta_{j,L})$, where $(\beta_{j,1},\dots,\beta_{j,m})$ is the $j$-th column of $A'{}^{-1}$ and $\beta_{j,k}=0$, $k=L+1,\dots,m$.
\begin{itemize}
\item[(i)] Let $\varphi_j^{-1}(\tau')$, $j=1,\dots,L$. Then $\log\varphi_j^{-1}(e^{\alpha_i})|_{\lambda''=1}$ is nothing but the scalar product $\alpha_i \cdot \beta_j = \delta_{ij}$, where $\delta_{ij}$ is the Kronecker's delta.
\item[(ii)] First suppose that $1 \leq i \leq L$. Let $\psi_k(\tau')$, $k=L+1,\dots,m$. Then $\log\psi_k(e^{\alpha_i})$
is nothing but the scalar product
$\alpha_i \cdot \left(e_k - \sum_{j=1}^m\alpha_{jk}\beta_j\right) = \alpha_{ik} - \alpha_{ik}=0$. Moreover $\log(1/\psi_k(e^{\alpha_i}))=0$ as well.

Now suppose that $L+1 \leq i \leq \ell$. Then $\alpha_i=\sum_{j=1}^L c_j\alpha_j$, so $\log\psi_k(e^{\alpha_i})=\sum_{j=1}^L c_j\log\psi_k(e^{\alpha_j})=0$. Moreover $\log(1/\psi_k(e^{\alpha_i}))=0$ as well.
\end{itemize}
\end{proof}

\begin{cor}\label{cor1} Let $(f,v) \in \semiG$ and $1 \le j \le L$. Then
\begin{equation}\label{1}
\begin{array}{cl}
 \log f(e^{\alpha_j}) = 0 & \mbox{if} \ f=\prod_{k=L+1}^m\psi_k, \\
 \log f(e^{\alpha_j}) \geq 0 & \mbox{otherwise}.
\end{array}
\end{equation}
In particular
\begin{equation}%\label{1}
\begin{array}{c}
v=0 \Rightarrow \log f(e^{\alpha_j}) \geq 0, \\
v \neq 0 \Rightarrow \log f(e^{\alpha_j})=0.
\end{array}
\end{equation}
\end{cor}
\begin{proof} This follows immediately from Lemma \ref{lem0} and the fact that
	$[G]_=$
and $[G]_> * (1,0)$ generate $\semiG$.

%This follows immediately from Lemma \ref{lem1} and the fact that the radical of $[\mathcal{F}^q]$ is $\semiG$.
\end{proof}

Let us consider a $\ell \times m$ matrix $A$ and a $(\ell+1) \times m$ matrix $B$ obtained by adding to $A$ a $(\ell+1)$-th row $(\beta_{1},\dots,\beta_{m})$ ($\beta_{i} \in \mathbb{Q}_{\ge 0}$). We will use $(\cdot)_A$ (resp. $(\cdot)_B$) to indicate elements (sets, family of sets) related to the normal deformation associated to $A$ (resp. $B$). There are two possibilities: $\operatorname{rank}B=L$ or $\operatorname{rank}B=L+1$.\\

Let us consider the case $\operatorname{rank}B=L$. Set
\begin{equation}
b_j :=\log\imin{\varphi_{jA}}(e^{\beta})|_{\lambda''=1} \qquad (j=1,\dots,L).
\end{equation}

\begin{lem}\label{lem2.0} Let us assume that
\begin{equation}\label{2.0}
\begin{array}{c}
v=0 \Rightarrow \log f(e^{\beta}) \geq 0 \\
\end{array}
\end{equation}
for $(f,v) \in \semiG_A$. Then $\semiG_A$ and $\semiG_B$ have the same radical.
\end{lem}
\begin{proof}
For $j=1,\dots,L$, we have
$$
\tau_j=\varphi_{jB}(\lambda',\lambda'',\lambda_{\ell+1})=\varphi_{jA}(\lambda',\lambda'')\cdot\lambda_{\ell+1}^{\beta_j}.
$$
After a simple computation we obtain
$$
\lambda_j=\imin{\varphi_{jA}}(\tau,\lambda'')\lambda_{\ell+1}^{-b_j}=\imin{\varphi_{jB}}(\tau,\lambda'',\lambda_{\ell+1}).
$$
Hence, for $k \geq L+1$, we have
$$
\psi_{kB}(\tau)=\psi_{kA}(\tau).
$$
Set $J = \{1,\dots,L\}$, $K = \{L+1,\dots,m\}$ and $\tilde{J} = \{L+1,\dots,\ell\}$.
If $v \neq 0$, then $f=\displaystyle\prod_{k \in K}\psi_k^{\delta_k}$
for some $\delta_k \in \mathbb{Z}$, and hence,
$(f,v) \in \semiG_A$ if and only if $(f,v) \in \semiG_B$. \\

Assume now $v=0$. Let $\displaystyle{f=\prod_{j \in J}(\imin{\varphi_{jA}})^{a_j}\cdot\prod_{k \in K}(\psi_{kA})^{a_k}}$. Then
\begin{eqnarray*}
\log f(e^{\beta})|_{\lambda''=1} & = & \sum_{j \in J}a_j\log\imin{\varphi_{jA}}(e^{\beta})|_{\lambda''=1}+\sum_{k \in K}a_k\log\psi_{kA}(e^{\beta}) \\
& = &  \sum_{j \in J}a_j\log\imin{\varphi_{jA}}(e^{\beta})|_{\lambda''=1}+\sum_{k \in K}a_k\log\psi_{kB}(e^{\beta}) \\
& = & \sum_{j \in J}a_jb_j
\end{eqnarray*}
where the last equality follows since $\log\psi_{kB}(e^{\beta})=0$ by Lemma \ref{lem0}. Here we used the fact that $\operatorname{rank}A = \operatorname{rank}B$.

Now assume that $(f,0) \in \semiG_A$. Set $b_f=\log f(e^{\beta})|_{\lambda''=1}$.
%Being a linear combination of $\alpha_j$, then $b_f=\sum_j c_j b_j$, $j=1,\dots,L$.
Then
$$
\begin{aligned}
f& =\prod_{j \in J}(\imin{\varphi_{jA}})^{a_j}\cdot\prod_{k \in K}(\psi_{kA})^{a_k}
\cdot\prod_{j \in \tilde{J}}\lambda_j^{\delta_j} \\
&=\prod_{j \in J}(\imin{\varphi_{jB}})^{a_j}\cdot\lambda_{\ell+1}^{b_f}\cdot\prod_{k \in K}(\psi_{kB})^{a_k}
\cdot\prod_{j \in \tilde{J}}\lambda_j^{\delta_j},
\end{aligned}
$$
and $(f,0) \in \semiG_B$ if $b_f \geq 0$. Conversely, assume that $(f,0) \in \semiG_B$. Then
$$
\begin{aligned}
f&=\prod_{j \in J}(\imin{\varphi_{jB}})^{a_j}\cdot\lambda_{\ell+1}^{b_f}\cdot\prod_{k \in K}(\psi_{kB})^{a_k}
\cdot\prod_{j \in \tilde{J}}\lambda_j^{\delta_j} \\
&=\prod_{j \in J}(\imin{\varphi_{jA}})^{a_j}\cdot\prod_{k \in K}(\psi_{kA})^{a_k}
\cdot\prod_{j \in \tilde{J}}\lambda_j^{\delta_j},
\end{aligned}
$$
where the term $\lambda_{\ell+1}^{b_f}$ appears since the degree in $\lambda_{\ell+1}$ must be 0. Then $(f,0) \in \semiG_A$ and the result follows.
\end{proof}

\begin{cor}\label{cor2.0} Suppose that $\beta=\displaystyle\sum_{j=1}^L c_j\alpha_j$ with $c_j \geq 0$, $j=1,\dots,L$. Then $\semiG_A$ and $\semiG_B$ have the same radical.
\end{cor}
\begin{proof} By Lemma \ref{lem0} we obtain $b_f := \log f(e^\beta) \geq 0$ for each $(f,v) \in \semiG_A$ and the result follows by Lemma \ref{lem2.0}.
\end{proof}

Thanks to Proposition \ref{prop:radical-of-S} it is enough to check the conditions on the elements of $F_A^q$, $q := \#\left(J_Z(\xi) \cap \{1,2,\dots,L\}\right) \ge 0$, and we have
\begin{prop}\label{prop2,5.0} Let us assume that $\operatorname{rank} B = \operatorname{rank} A$ and
\begin{equation}\label{2,5.0}
\begin{array}{c}
v=0 \Rightarrow \log f(e^{\beta}) \geq 0
\end{array}
\end{equation}
for $(f,v) \in F^q_A$. Then $\semiG_A$ and $\semiG_B$ have the same radical.
\end{prop}

From Proposition \ref{prop2,5.0} we can deduce a condition for the restriction of multi-normal deformations. We omit the coordinate $x^{(0)}$ to lighten notations.
In what follows, given $\xi=(\xi^{(1)},\dots,\xi^{(m)}) \in \R^n$ we will call $p_A \in S_A$, the point $p_A=(\xi^{(1)},\dots,\xi^{(m)};\,0 \dots,0) \in \R^{n+\ell}$, and $p_B \in S_B$, the point $p_B=(\xi^{(1)},\dots,\xi^{(m)};\,0 \dots,0) \in \R^{n+\ell+1}$.

\begin{prop}\label{prop3.0}
Let $p_A=(\xi^{(1)},\dots,\xi^{(m)};\,0,\dots,0) \in S_A$ and
assume that $\operatorname{rank}{A} = \operatorname{rank}{B}$ and
\begin{equation}\label{3.0}
\begin{array}{c}
v=0 \Rightarrow \log f(e^{\beta}) \geq 0 \\
\end{array}
\end{equation}
for $(f,v) \in F^q_A$. Then $p_B \in C_{\chi_B}(Z)$
if and only if, for each $V \in C(\xi,\, \semiG_A)$ and each open neighborhood
$U$ of the origin,
we have $Z \cap V \cap U \neq \emptyset$ $($i.e., if and only if $p_A \in C_{\chi_A}(Z)$$)$.
\end{prop}

Thanks to Corollary \ref{cor2.0} we obtain a simpler sufficient condition.
Recall $L := \operatorname{rank} A$.
\begin{cor}\label{cor3.0} Let $p_A=(\xi^{(1)},\dots,\xi^{(m)};\,0,\dots,0) \in S_A$
and assume that
$\beta=\displaystyle\sum_{j=1}^L c_j\alpha_j$ with $c_j \geq 0$,
$j=1,\dots,L$.
Then $p_B \in C_{\chi_B}(Z)$
if and only if $p_A \in C_{\chi_A}(Z)$.
\end{cor}

\begin{es} In $\R^3$ let us consider the families $\chi_A=\{M_1,M_2,M_3\}$, $\chi_B=\{M_1,M_2,M_3,M_4\}$ with
$M_1=\{x_1=0\}$, $M_2=\{x_2=x_3=0\}$, $M_3=\{x_3=0\}$, $M_4=\{0\}$. If we consider the matrices
$$
A=\left(
\begin{array}{ccc}
1 & 0 & 0 \\
0 & 1 & 1 \\
0 & 0 & 1
\end{array}
\right)
\ \ \ \
B=\left(
\begin{array}{ccc}
1 & 0 & 0 \\
0 & 1 & 1 \\
0 & 0 & 1 \\
1 & 1 & 1
\end{array}
\right),
$$
the hypothesis of Corollary \ref{cor3.0} is satisfied.
\end{es}

In the following example, since $v$ of all the pairs $(f,v) \in F^q$ are zero,
we simply write a pair by $f$ instead of $(f,v)$.
%Also we omit the unit $(1,1)$ when describe the elements
%of $F^q$ for simplicity.
\begin{es} In $\R^3$ let us consider the families $\chi_A=\{M_1,M_2,M_3\}$, $\chi_B=\{M_1,M_2,M_3,M_4\}$ with
$M_1=\{x_1=x_3=0\}$, $M_2=\{x_2=x_3=0\}$, $M_3=\{x_3=0\}$, $M_4=\{0\}$. Set $\beta=(1,1,1)$. Consider the matrices
$$
A=\left(
\begin{array}{ccc}
1 & 0 & 1 \\
0 & 1 & 1 \\
0 & 0 & 1
\end{array}
\right)
\ \ \ \
B=\left(
\begin{array}{ccc}
1 & 0 & 1 \\
0 & 1 & 1 \\
0 & 0 & 1 \\
1 & 1 & 1
\end{array}
\right).
$$
\begin{enumerate}
\item Let $\xi_1,\xi_2,\xi_3 \neq 0$. Then $F^0_A=\displaystyle{\left\{\tau_1,\tau_2,{\tau_3 \over \tau_1\tau_2}\right\}}$. We have $\log{e^{\beta_1-\beta_2-\beta_3}}=1-1-1<0$ and \eqref{2,5.0} is not satisfied. Indeed one can check that $F^0_B=\displaystyle{\left\{\tau_1,\tau_2,{\tau_3 \over \tau_1},{\tau_3 \over \tau_2}\right\}}$ and $\displaystyle{\tau_3 \over \tau_1\tau_2} \notin \semiG_B$.
\item Let $\xi_3=0$, $\xi_1,\xi_2 \neq 0$. Then $F^1_A=\displaystyle{\left\{\tau_1,\tau_2,{\tau_3 \over \tau_1\tau_2}\right\}}$. We have $\log{e^{\beta_1-\beta_2-\beta_3}}=1-1-1<0$ and \eqref{2,5.0} is not satisfied. Indeed one can check that $F^1_B=\displaystyle{\left\{\tau_1,\tau_2,{\tau_3 \over \tau_1},{\tau_3 \over \tau_2}\right\}}$ and $\displaystyle{\tau_3 \over \tau_1\tau_2} \notin \semiG_B$.
\item Let $\xi_2=0$, $\xi_1,\xi_3 \neq 0$. Then $F^1_A=\displaystyle{\left\{\tau_1,\tau_2,{\tau_3 \over \tau_1}\right\}}$ and condition \eqref{2,5.0} is satisfied.
\item Let $\xi_1=0$, $\xi_2,\xi_3 \neq 0$. Then $F^1_A=\displaystyle{\left\{\tau_1,\tau_2,{\tau_3 \over \tau_2}\right\}}$ and condition \eqref{2,5.0} is satisfied. \\
\end{enumerate}
\end{es}

%Let us consider the case $\operatorname{rank}B=L+1$. If $L<\ell$, we exchange the $(L+1)$-th line with the $(\ell+1)$-th line and make a permutation of the columns in order to assume that the first $(L+1) \times (L+1)$ sub-matrix $B'$ of $B$ is invertible.

%Let $\imin{\varphi_{jA}}(\tau,\lambda'')$, $j=1,\dots,L$ and $\psi_{kA}(\tau)$, $k=L+1,\dots,m$ be the monomials defining $\semiG_A$. We are going to compute $\imin{\varphi_{jB}}(\tau,\lambda'')$, $j=1,\dots,L+1$ and $\psi_{kB}(\tau)$, $k=L+2,\dots,m$ using $\imin{\varphi_{jA}}(\tau,\lambda'')$ and $\psi_{kA}(\tau)$. Set
%$$
%b_j=
%\begin{cases}
%\log\imin{\varphi_{jA}}(e^{\beta})|_{\lambda''=1} & j=1,\dots,L, \\
%\log\psi_{jA}(e^{\beta}) & j=L+1,\dots,m.
%\end{cases}
%$$

%Assume for simplicity that all the columns of $B$ are non zero. We may also assume that $\exists\, k \in \{L+1,\dots,m\}$ such that $b_k \neq 0$. Otherwise $\log\psi_{kA}(e^{\alpha_j})=0$ for $j=1,\dots,L+1$ by Lemma \ref{lem0}. This implies that the rank of $B'$ is the same as the rank of $A'$. So, up to take a permutation of $\{1,\dots,m\}$, we may assume that $b_{L+1} \neq 0$. \\

Let us consider the case $\operatorname{rank}B=L+1$. If $L<\ell$, we exchange the $(L+1)$-th line with the $(\ell+1)$-th line.
% and make a permutation of the columns in order to assume that the first $(L+1) \times (L+1)$ sub-matrix $B'$ of $B$ is invertible.

Let $\imin{\varphi_{jA}}(\tau,\lambda'')$, $j=1,\dots,L$ and $\psi_{kA}(\tau)$, $k=L+1,\dots,m$ be the monomials defining $\semiG_A$. We are going to compute $\imin{\varphi_{jB}}(\tau,\lambda'')$, $j=1,\dots,L+1$ and $\psi_{kB}(\tau)$, $k=L+2,\dots,m$ using $\imin{\varphi_{jA}}(\tau,\lambda'')$ and $\psi_{kA}(\tau)$. Set
$$
b_j=
\begin{cases}
\log\imin{\varphi_{jA}}(e^{\beta})|_{\lambda''=1} & j=1,\dots,L, \\
\log\psi_{jA}(e^{\beta}) & j=L+1,\dots,m.
\end{cases}
$$

Assume for simplicity that all the columns of $B$ are non zero. We may also assume that $\exists\, k \in \{L+1,\dots,m\}$ such that $b_k \neq 0$. Otherwise $\log\psi_{kA}(e^{\alpha_j})=0$ for $j=1,\dots,L+1$ by Lemma \ref{lem0}. This implies that the rank of $B$ is the same as the rank of $A$. So, up to take a permutation of $\{1,\dots,m\}$, we may assume that $b_{L+1} \neq 0$. If $b_{L+1} \neq 0$, then the first $(L+1) \times (L+1)$ sub-matrix $B'$ of $B$ is invertible. Otherwise $\operatorname{rank}B'=\operatorname{rank}A'$ and, as in the proof of Lemma \ref{lem2.0}, $\log\psi_{(L+1)B'}(e^{\beta})=\log\psi_{(L+1)A'}(e^{\beta})$. We have
$$
\log\psi_{(L+1)B'}(e^{\beta})=\log\psi_{(L+1)A'}(e^{\beta})=\log\psi_{(L+1)A}(e^{\beta})=b_{L+1}
$$
and $\log\psi_{(L+1)B'}(e^{\beta})=0$ by Lemma \ref{lem0}.

Let $\tau$ denote $m$-variables $(\tau_1,\,\dots,\,\tau_m)$.
Consider the $(L+1)$-equations
$$
\tau_1 = \varphi_{1B}(\lambda), \dots, \tau_{L+1} = \varphi_{L+1B}(\lambda).
$$
It can be also written as
$$
\tau_1 = \varphi_{1A}(\lambda)\lambda_{L+1}^{\beta_{1}}, \dots, \tau_{L+1} = \varphi_{L+1A}(\lambda)\lambda_{L+1}^{\beta_{L+1}}.
$$
Then, by the assumption, the system of the equations can be solved
for the variables
$\lambda = (\lambda_1, \dots, \lambda_{L+1})$. After a computation, we can find
\begin{eqnarray*}
\imin{\varphi_{jB}}(\tau,\lambda'') & = & {\imin{\varphi_{jA}}(\tau,\lambda'') \over (\psi_{L+1A}(\tau))^{b_j \over b_{L+1}}} \ \ j=1,\dots,L, \\
\imin{\varphi_{L+1B}}(\tau) & = & (\psi_{L+1A}(\tau))^{1 \over b_{L+1}}, \\
\psi_{kB}(\tau) & = & {\psi_{kA}(\tau) \over (\psi_{L+1A}(\tau))^{b_k \over b_{L+1}}} \ \ k=L+2,\dots,m, \\
\end{eqnarray*}
and
\begin{eqnarray*}
\imin{\varphi_{jA}}(\tau,\lambda'') & = & \imin{\varphi_{jB}}(\tau,\lambda'') \cdot (\imin{\varphi_{L+1B}}(\tau))^{b_j} \ \ j=1,\dots,L, \\
\psi_{L+1A}(\tau) & = & (\imin{\varphi_{L+1B}}(\tau))^{b_{L+1}}, \\
\psi_{kA}(\tau) & = & \psi_{kB}(\tau) \cdot (\imin{\varphi_{L+1B}}(\tau))^{b_k} \ \ k=L+2,\dots,m. \\
\end{eqnarray*}
In what follows, we omit the variables $\tau,\lambda''$ to lighten notations.
\begin{lem}\label{lem2} Let us assume that
\begin{equation}\label{2}
\begin{array}{c}
v=0 \Rightarrow \log f(e^{\beta}) \geq 0, \\
v \neq 0 \Rightarrow \log f(e^{\beta}) = 0, \\
%b_k \geq 0 \ \ \forall\,k \in \{\ell+2,\dots,m\} \cap J_Z(\xi) \\
{\rm if} \; m \geq L+2, \ \ b_k = 0 \ \ \forall\,k \in \{L+2,\dots,m\} \setminus J_Z(\xi)
\end{array}
\end{equation}
for $(f,v) \in \semiG_A$. Then $\semiG_A$ and $\semiG_B$ have the same radical.
\end{lem}
\begin{proof}
%Assume for simplicity that all the columns of $B$ are non zero. We may also assume that $\exists j \in \{1,\dots,m\}$ such that $\xi_j=0$. Otherwise $(\psi_{kA},\xi_k) \in \semiG_A$ and $\log\psi_{kA}(e^{\alpha_j})=0$ for $j=1,\dots,\ell+1$ by \eqref{2} and Corollary \ref{cor1}. This implies that the rank of $B$ is the same as the rank of $A$. So, up to take a permutation of $\{1,\dots,m\}$, we may assume that $\xi_{\ell+1}=0$. \\
First remark that $b_{L+1} \neq 0$ and \eqref{2} imply $((\psi_{L+1A})^{a_{L+1}},v) \in \semiG_A$, $a_{L+1} \neq 0$, only if $v = 0$. Since the variable $\lambda''$ does not play a role in this proof, we assume $\ell=L$.\\

(i) Let $(f_A,v) \in \semiG_A$. Then
\begin{eqnarray*}
f_A & = & \prod_{j \in J}(\imin{\varphi_{jA}})^{a_j} \cdot (\psi_{\ell+1A})^{a_{\ell+1}} \cdot \prod_{k \in K}(\psi_{kA})^{a_k} \\
& = & \prod_{j \in J}(\imin{\varphi_{jB}})^{a_j} \cdot (\imin{\varphi_{\ell+1B}})^{c} \cdot \prod_{k \in K}(\psi_{kB})^{a_k},
\end{eqnarray*}
$J \subseteq \{1,\dots,\ell\}$, $K \subseteq \{\ell+2,\dots,m\}$ and
$$
c={\sum_{j \in J}a_jb_j + a_{\ell+1}b_{\ell+1}} + \sum_{k \in K}a_kb_k=\log f_A(e^{\beta}).
$$
\begin{itemize}
\item Suppose that $v=0$. It follows from \eqref{2} that $c \geq 0$.
\begin{itemize}
\item Suppose that $c>0$. Then $(f_A,0) \in \semiG_B$ since $((\imin{\varphi_{\ell+1B}})^c,0) \in [G_{B}]$ and $(\imin{\varphi_{\ell+1B}},v) \in G_B$ only if $v=0$.
\item Suppose that $c=0$. Then there are two possibilities:
\begin{itemize}
\item $a_j \neq 0$, $\exists\, j \in J$ or $a_k \neq 0$, $\exists\, k \in K \cap J_Z(\xi)$. In this case $(f_A,0) \in \semiG_B$ since $((\imin{\varphi_{jB}})^{a_j},0),((\psi_{kB})^{a_k},0) \in [G_{B}]$ and $(\imin{\varphi_{jB}},v),(\psi_{kB},v) \in G_B$ only if $v=0$.
\item $a_j=0$, $\forall\, j \in J$ and $a_k=0$, $\forall\, k \in K \cap J_Z(\xi)$. We may assume $K \cap J_Z(\xi)= \emptyset$ and by \eqref{2} we have
$$
f_A = \prod_{k \in K}{(\psi_{kA})^{a_k} \over (\psi_{\ell+1A})^{a_kb_k \over b_{\ell+1}}}=\prod_{k \in K}(\psi_{kA})^{a_k}=\prod_{k \in K}(\psi_{kB})^{a_k}.
$$
This is because by \eqref{2} we have $b_k=0$, $k \in K$.
We have $(f_A,v) \in [G_A]_> * [(1,0)]$. Hence $(f_A,0) \in [G_B]_> * [(1,0)]$. % since $f_A$ does not depend on $\tau_{\ell+1}$.
\end{itemize}
\end{itemize}

\item Suppose that $v \neq 0$. Then $(f_A,v) \in [G_A]_=$. It follows from \eqref{2} that $c=0$. Moreover $a_j=0$, $j \in J$ and $a_k=0$, $k \in K \cap J_Z(\xi)$ since $(\imin{\varphi_{jA}},v),(\psi_{kA},v) \in G_A$ only if $v=0$. Then
$$
f_A = \prod_{k \in K \setminus J_Z(\xi)}{(\psi_{kA})^{a_k} \over (\psi_{\ell+1A})^{a_kb_k \over b_{\ell+1}}}=\prod_{k \in K \setminus J_Z(\xi)}(\psi_{kA})^{a_k}=\prod_{k \in K \setminus J_Z(\xi)}(\psi_{kB})^{a_k}.
$$
This is because by \eqref{2} we have $b_k=0$ if $k \notin J_Z(\xi)$.
Then $(f_A, v) \in [G_A]_=$ implies $(f_A,v) \in [G_B]_=$. % since $f_A$ does not depend on $\tau_{\ell+1}$.
\end{itemize}

(ii) Let $(f_B,v) \in \semiG_B$. Then
\begin{eqnarray*}
f_B & = & \prod_{j \in J}(\imin{\varphi_{jB}})^{a_j} \cdot (\imin{\varphi_{\ell+1B}})^{a_{\ell+1}} \cdot \prod_{k \in K}(\psi_{kB})^{a_k} \\
& =& {\prod_{j \in J}(\imin{\varphi_{jA}})^{a_j} \cdot (\psi_{\ell+1A})^c} \cdot \prod_{k \in K}(\psi_{kA})^{a_k},
\end{eqnarray*}
$J \subseteq \{1,\dots,\ell\}$, $K \subseteq \{\ell+2,\dots,m\}$ and
$$
c={a_{\ell+1} - \sum_{j \in J}a_jb_j - \sum_{k \in K}a_kb_k \over b_{\ell+1}}.
$$
\begin{itemize}
\item Suppose that $v=0$. Then $((f_B)^{|b_{\ell+1}|},w) \in \semiG_A$
for some $w \geq 0$. In fact, it clearly holds if $\xi^{(\ell+1)} \ne 0$.
Otherwise, as $(f_B,v) \in \mathcal{G}_B$ implies $\nu_{\ell+1}(f_B) \ge 0$
when $\xi^{(\ell+1)} = 0$, we conclude $c \ge 0$ since
we see $c = \nu_{\ell+1}(f_B)$ by comparing the degree of the variable $\tau_{\ell+1}$.
\begin{itemize}
\item Suppose  $a_j \neq 0$, $\exists\, j \in J$ or $a_k \neq 0$, $\exists\, k \in K \cap J_Z(\xi)$ . In this case $((f_B)^{|b_{\ell+1}|},0) \in \semiG_A$ since $((\imin{\varphi_{jA}})^{a_j},0),((\psi_{kA})^{a_k},0) \in [G_{A}]$ and $(\imin{\varphi_{jA}},v),(\psi_{kA},v) \in G_A$ only if $v=0$.
\item Suppose that $a_j=0$, $\forall\, j \in J$ and $a_k=0$, $\forall\, k \in K \cap J_Z(\xi)$. We may assume $K \cap J_Z(\xi)= \emptyset$.
Let us argue by contradiction and suppose that $w \neq 0$. Then $((f_B)^{|b_{\ell+1}|},w) \in [G_A]_=$. Corollary \ref{cor1} implies that $a_{\ell+1}=\log f_B(e^{\beta})=0$ and by \eqref{2} we have
$$
f_B=\prod_{k \in K}(\psi_{kB})^{a_k} = \prod_{k \in K}{(\psi_{kA})^{a_k} \over (\psi_{\ell+1A})^{a_kb_k \over b_{\ell+1}}} = \prod_{k \in K}(\psi_{kA})^{a_k}.
$$
This implies that $((f_B)^{|b_{\ell+1}|},w) \in [G_B]_=$, %since $f_B$ does not depend on $\tau_{\ell+1}$,
a contradiction.
Hence we have $w = 0$ and $( (f_B)^{|b_{\ell+1}|}, 0) \in \mathcal{G}_A$ follows.
\end{itemize}
\item Suppose that $v \neq 0$. Then $(f_B,v) \in [G_B]_=$. In this case $a_j=0$, $j \in J \cup \{\ell+1\}$ and $a_k=0$, $k \in K \cap J_Z(\xi)$ since $(\imin{\varphi_{jB}},v),(\psi_{kB},v) \in G_B$ only if $v=0$. Then
$$
f_B = \prod_{k \in K \setminus J_Z(\xi)}(\psi_{kB})^{a_k}=\prod_{k \in K \setminus J_Z(\xi)}{(\psi_{kA})^{a_k} \over (\psi_{\ell+1A})^{a_kb_k \over b_{\ell+1}}}=\prod_{k \in K \setminus J_Z(\xi)}(\psi_{kA})^{a_k}.
$$
This is because by \eqref{2} we have $b_k=0$ if $k \notin J_Z(\xi)$.
Then $(f_B, v) \in [G_B]_=$ implies $(f_B,v) \in [G_A]_=$.
\end{itemize}

\end{proof}

Thanks to Proposition \ref{prop:radical-of-S} it is enough to check the conditions on the elements of $F_A^q$, $q := \#\left(J_Z(\xi) \cap \{1,2,\dots,L\}\right) \ge 0$, and we have
\begin{prop}\label{prop2,5} Let us assume that
\begin{equation}\label{2,5}
\begin{array}{c}
v=0 \Rightarrow \log f(e^{\beta}) \geq 0, \\
v \neq 0 \Rightarrow \log f(e^{\beta}) = 0, \\
%b_k \geq 0 \ \ \forall\,k \in \{\ell+2,\dots,m\} \cap J_Z(\xi) \\
{\rm if} \; m \geq L+2, \ \ b_k = 0 \ \ \forall\,k \in \{L+2,\dots,m\} \setminus J_Z(\xi)
\end{array}
\end{equation}
for $(f,v) \in F^q_A$. Then $\semiG_A$ and $\semiG_B$ have the same radical.
\end{prop}

From Proposition \ref{prop2,5} we can deduce a condition for the restriction of multi-normal deformations. We omit the coordinate $x^{(0)}$ to lighten notations. In what follows, given $\xi=(\xi^{(1)},\dots,\xi^{(m)}) \in \R^n$ we will call $p_A \in S_A$, the point $p_A=(\xi^{(1)},\dots,\xi^{(m)};\, 0 \dots,0) \in \R^{n+\ell}$, and $p_B \in S_B$, the point $p_B=(\xi^{(1)},\dots,\xi^{(m)};\,0 \dots,0) \in \R^{n+\ell+1}$.

\begin{prop}\label{prop3} Let $p_A=(\xi^{(1)},\dots,\xi^{(m)};\, 0,\dots,0) \in S_A$.
Assume that $\operatorname{rank} A + 1 = \operatorname{rank} B$ and
\begin{equation}\label{3}
\begin{array}{c}
v=0 \Rightarrow \log f(e^{\beta}) \geq 0, \\
v \neq 0 \Rightarrow \log f(e^{\beta}) = 0, \\
%b_k \geq 0 \ \ \forall\,k \in \{\ell+2,\dots,m\} \cap J_Z(\xi) \\
{\rm if} \; m \geq L+2,  \ \ b_k = 0 \ \ \forall\,k \in \{L+2,\dots,m\} \setminus J_Z(\xi)
\end{array}
\end{equation}
for $(f,v) \in F^q_A$. Then $p_B \in C_{\chi_B}(Z)$
if and only if, for each $V \in C(\xi,\, \semiG_A)$ and each open neighborhood $U$
	of the origin, we have $Z \cap V \cap U \neq \emptyset$
$($i.e., if and only if $p_A \in C_{\chi_A}(Z)$$)$. \\
\end{prop}

When $\chi_B$ is of normal type, i.e., $m=\ell+1=L+1$, Proposition \ref{prop3} can be reformulated as follows.

\begin{cor}\label{cor4} Let $m=\ell+1$.
Let $p_A=(\xi^{(1)},\dots,\xi^{(\ell+1)};\,0,\dots,0) \in S_A$. Assume that
$m = \operatorname{rank} B = \operatorname{rank} A + 1$ and
\begin{equation}\label{4}
\begin{array}{c}
v=0 \Rightarrow \log f(e^{\beta}) \geq 0, \\
((\psi_{\ell+1A})^{a_{\ell+1}},v) \in \semiG_A \Rightarrow v=0
\end{array}
\end{equation}
for $(f,v) \in F^q_A$, $a_{\ell+1} \neq 0$. Then $p_B \in C_{\chi_B}(Z)$
if and only if, for each $V \in C(\xi,\, \semiG_A)$ and each open neighborhood $U$
	of the origin,
we have $Z \cap V \cap U \neq \emptyset$ $($i.e., if and only if $p_A \in C_{\chi_A}(Z)$$)$.
\end{cor}

Corollary \ref{cor4} is a generalization of Corollary 4.3 of \cite{HP13}. Consider the family $\chi_B = \{M_1, \dots ,M_{\ell+1}\}$
satisfying condition H2 in \cite{HP13}, i.e.
\begin{itemize}
\item[H2] either $I_j \subset I_k$, $I_k \subset I_j$ or $I_j \cap I_k = \emptyset$ holds for $j \neq k$.
\end{itemize}
Consider a sub-family $\chi_A := \{M_{j_1} ,\dots ,M_{j_k}\}$ of $\chi_B$. We denote by $S_{\chi_B}$ the zero section
$\underset{X,1\leq j\leq\ell+1}{\times}T_{M_j}\iota_{\chi_B}(M_j)$ for the family $\chi_B$.
We also introduce the set
\begin{equation}{\label{def:S-chi-chi}}
S_{\chi_B/\chi_A} :=
\left(\underset{X,1\leq \alpha \leq k}{\times}T_{M_{j_\alpha}}\iota_{\chi_B}(M_{j_\alpha})
\right)
\underset{X}{\times} M,
\end{equation}
where $M=\bigcap_{j=1}^{\ell+1} M_j$.
We emphasize that, in the above definition, we use $\iota_{\chi_B}$ (not $\iota_{\chi_A}$).
Then we have the natural embedding
$$
S_{\chi_B} \hookleftarrow S_{\chi_B/\chi_A} \hookrightarrow S_{\chi_A}.
$$

\begin{cor} \label{cor: restriction multi-cones}
Assume that $\chi_B$ satisfies the condition H2.
Let $k \le \ell$ and $\{j_1,\dots,j_k\}$ be a subset of $\{1,2,\dots,\ell+1\}$. Set $\chi_A=\{M_{j_1},\ldots,M_{j_k}\}$.
Let $Z$ be a subset of $X$. Then we have
$$
C_{\chi_B} (Z) \cap S_{\chi_B/\chi_A}
= C_{\chi_A} (Z) \cap S_{\chi_B/\chi_A}.
$$
\end{cor}
\begin{proof}
It is enough to prove the result for $\sharp\chi_A=\ell$ and then apply it repeatedly. We can assume without loss of generality that $\chi_A=\{M_1,\ldots,M_\ell\}$. We use the coordinates $x=(x^{(0)},x^{(1)},\dots,x^{(\ell+1)})$ defined as follows: $x^{(j)}=(x_i)_{i \in \widehat{I}_j}$ with
\begin{eqnarray*}
\widehat{I}_0 & = & \{1,\dots,\ell+1\} \setminus \bigcup_{i=1}^{\ell+1}I_i, \\
\widehat{I}_j & = & I_j \setminus \bigcup_{I_i \subsetneqq I_j}I_i, \ \ j=1,\dots,\ell+1.
\end{eqnarray*}
Since the coordinate $x^{(0)}$ is irrelevant in the rest of the proof, we will assume $\widehat{I}_0=\emptyset$.
Then $A$ is a $\ell \times (\ell+1)$ matrix with entries $\alpha_{i,j}$, $1\leq i \leq \ell$, $1\leq j \leq \ell+1$, $B$ is a $(\ell+1) \times (\ell+1)$ matrix obtained adding the line $(\alpha_{\ell+1,1},\dots,\alpha_{\ell+1,\ell+1})$. Up to take a permutation of the coordinates we may assume that the first $\ell$ columns of $A$ are linearly independent. The entries of the matrices are given by
$\alpha_{ij}=1$ if $\widehat{I}_j \subseteq I_i$, $\alpha_{ij}=0$ otherwise.
%Moreover, if the $(\ell+1)$-th column is equal to the $j$-th column, $1 \leq j \leq \ell$, up to exchange the columns we assume $\alpha_{\ell+1,\ell+1}=1$.
The restriction maps $S_{\chi_B/\chi_A} \to S_{\chi_A}$,
$S_{\chi_B/\chi_A} \to S_{\chi_B}$ are given by $\xi^{(\ell+1)}=0$.\\

Thanks to the condition H2 it is easy to compute the inverse matrix $\imin B$, whose columns define the monomials $\imin{\varphi_{jB}}$, $j=1,\dots,\ell+1$ which are given by
\begin{eqnarray*}
\imin{\varphi_{jB}}(\tau) & = & \tau_j \mbox{ if $I_j$ is maximal}, \\
\imin{\varphi_{jB}}(\tau) & = & {\tau_j \over \tau_{k_j}} \mbox{ otherwise, where }I_{k_j}=\bigcap_{I_k \supsetneqq I_j}I_k.
\end{eqnarray*}
Remark that to define $\imin\varphi_{jB}$, $j=1,\dots,\ell+1$ we are working with the family of indices $I_j$, $j=1,\dots,\ell+1$.
Thanks to the condition H2 it is also easy to compute the inverse matrix $A'{}^{-1}$ of $A'$, the matrix defined by the first $\ell$ columns of $A$, whose columns define the monomials $\imin{\varphi_{jA}}$, $j=1,\dots,\ell$ which are given by
\begin{eqnarray*}
\imin{\varphi_{jA}}(\tau) & = & \tau_j \mbox{ if $I_j$ is maximal}, \\
\imin{\varphi_{jA}}(\tau) & = & {\tau_j \over \tau_{k_j}} \mbox{ otherwise, where }I_{k_j}=\bigcap_{I_k \supsetneqq I_j}I_k.
\end{eqnarray*}
Remark that to define $\imin\varphi_{jA}$, $j=1,\dots,\ell$ we are working with the family of indices $I_j$, $j=1,\dots,\ell$. Then one can also compute the monomial $\psi_{\ell+1A}$, which is given by
\begin{eqnarray*}
\psi_{\ell+1A}(\tau) & = & \tau_{\ell+1} \mbox{ if $I_{\ell+1}$ is maximal}, \\
\psi_{\ell+1A}(\tau) & = & {\tau_{\ell+1} \over \tau_{k_{\ell+1}}} \mbox{ otherwise, where }I_{k_{\ell+1}}=\bigcap_{I_k \supsetneqq I_{\ell+1}}I_k.
\end{eqnarray*}
It is easy to check that $\psi_{\ell+1A}(e^{\beta})=1$ and $\imin{\varphi_{jA}}(e^{\beta}) \in \{0,1\}$. Hence the conditions \eqref{4} are satisfied and the result follows from Corollary \ref{cor4}.
\end{proof}

\begin{es} Let us consider some examples in $\R^3$ of Corollary \ref{cor: restriction multi-cones}.
\begin{enumerate}
\item (Majima) Let $X={\mathbb R}^3$ with coordinates $(x_1, x_2, x_3)$ and let $\chi_A=\{M_1,M_2\}$, $\chi_B=\{M_1,M_2,M_3\}$ with
$M_i = \{x_i = 0\}$, $i=1,2,3$. Then $F^0_A=F^0_B=Q(\{(\tau_1, 0),
(\tau_2,0), (\tau_3, |\xi_3|)\})$. If we consider $\xi_3=0$ and the matrices
$$
A=\left(
\begin{array}{ccc}
1 & 0 & 0 \\
0 & 1 & 0
\end{array}
\right)
\ \ \ \
B=\left(
\begin{array}{ccc}
1 & 0 & 0 \\
0 & 1 & 0 \\
0 & 0 & 1
\end{array}
\right),
$$
condition \eqref{2,5} is satisfied.
\item (Takeuchi) Let $X={\mathbb R}^3$ with coordinates $(x_1, x_2, x_3)$. There are several configurations.  In the subsequent examples,
all the couples $(f,v)$ are with $v=0$, and we write it by $f$ instead of $(f,v)$ for short.
\begin{enumerate}
\item  let $\chi_A=\{M_1,M_2\}$, $\chi_B=\{M_1,M_2,M_3\}$ with
$M_1=\{0\}$, $M_2=\{x_2=x_3=0\}$, $M_3=\{x_3=0\}$. Let $\xi_3=0$, $\xi_1,\xi_2 \neq 0$. Then $F^1_A=\displaystyle{\left\{\tau_1,{\tau_2 \over \tau_1},{\tau_3 \over \tau_2}\right\}}$. If we consider the matrices
$$
A=\left(
\begin{array}{ccc}
1 & 1 & 1 \\
0 & 1 & 1
\end{array}
\right)
\ \ \ \
B=\left(
\begin{array}{ccc}
1 & 1 & 1 \\
0 & 1 & 1 \\
0 & 0 & 1
\end{array}
\right),
$$
condition \eqref{2,5} is satisfied.
\item  let $\chi_A=\{M_1,M_2\}$, $\chi_B=\{M_1,M_2,M_3\}$ with
$M_1=\{x_2=x_3=0\}$, $M_2=\{x_2=0\}$, $M_3=\{0\}$. Let $\xi_3=0$, $\xi_1,\xi_2 \neq 0$.  Then $F^1_A=\displaystyle{\left\{\tau_1,{\tau_2 \over \tau_1},\tau_3\right\}}$. If we consider the matrices
$$
A=\left(
\begin{array}{ccc}
1 & 1 & 0 \\
0 & 1 & 0
\end{array}
\right)
\ \ \ \
B=\left(
\begin{array}{ccc}
1 & 1 & 0 \\
0 & 1 & 0 \\
1 & 1 & 1
\end{array}
\right),
$$
condition \eqref{2,5} is satisfied.
\item  let $\chi_A=\{M_1,M_2\}$, $\chi_B=\{M_1,M_2,M_3\}$ with
$M_1=\{0\}$, $M_2=\{x_2=0\}$, $M_3=\{x_2=x_3=0\}$. Let $\xi_3=0$, $\xi_1,\xi_2 \neq 0$.  Then $F^1_A=\displaystyle{\left\{\tau_1,{\tau_2 \over \tau_1},{\tau_3 \over \tau_1}\right\}}$. If we consider the matrices
$$
A=\left(
\begin{array}{ccc}
1 & 1 & 1 \\
0 & 1 & 0
\end{array}
\right)
\ \ \ \
B=\left(
\begin{array}{ccc}
1 & 1 & 1 \\
0 & 1 & 0 \\
0 & 1 & 1
\end{array}
\right),
$$
condition \eqref{2,5} is satisfied.

\end{enumerate}
\item (Mixed) Let $X={\mathbb R}^3$ with coordinates $(x_1, x_2, x_3)$. There are several configurations.
\begin{enumerate}
\item  let $\chi_A=\{M_1,M_2\}$, $\chi_B=\{M_1,M_2,M_3\}$ with
$M_1=\{0\}$, $M_2=\{x_2=0\}$, $M_3=\{x_3=0\}$. Let $\xi_3=0$, $\xi_1,\xi_2 \neq 0$. Then $F^1_A=\displaystyle{\left\{\tau_1,{\tau_2 \over \tau_1},{\tau_3 \over \tau_1}\right\}}$. If we consider the matrices
$$
A=\left(
\begin{array}{ccc}
1 & 1 & 1 \\
0 & 1 & 0
\end{array}
\right)
\ \ \ \
B=\left(
\begin{array}{ccc}
1 & 1 & 1 \\
0 & 1 & 0 \\
0 & 0 & 1
\end{array}
\right),
$$
condition \eqref{2,5} is satisfied.
\item  let $\chi_A=\{M_1,M_2\}$, $\chi_B=\{M_1,M_2,M_3\}$ with
$M_1=\{x_1=0\}$, $M_2=\{x_2=0\}$, $M_3=\{0\}$. Let $\xi_3=0$, $\xi_1,\xi_2 \neq 0$. Then $F^1_A=\displaystyle{\{\tau_1,\tau_2,\tau_3\}}$. If we consider the matrices
$$
A=\left(
\begin{array}{ccc}
1 & 0 & 0 \\
0 & 1 & 0
\end{array}
\right)
\ \ \ \
B=\left(
\begin{array}{ccc}
1 & 0 & 0 \\
0 & 1 & 0 \\
1 & 1 & 1
\end{array}
\right),
$$
condition \eqref{2,5} is satisfied.
\end{enumerate}
\end{enumerate}
\end{es}

\begin{es} In $\R^3$ let us consider the family $\chi_A=\{M_1,M_2\}$, with $M_1=\{x_1=x_2=0\}$, $M_2=\{x_2=x_3=0\}$. The matrix $A$ is
$$
\left(
\begin{array}{ccc}
1 & 1 & 0 \\
0 & 1 & 1
\end{array}
\right)
$$
and $\imin{\varphi_{1A}}(\tau)=\tau_1$, $\imin{\varphi_{2A}}(\tau)=\displaystyle{\tau_2 \over \tau_1}$, $\imin{\varphi_{3A}}(\tau)=\displaystyle{\tau_1\tau_3 \over \tau_2}$. We have
$$
F^0_A=\left\{(\tau_1,0),\left({\tau_2 \over \tau_1},0\right),\left({\tau_1\tau_3 \over \tau_2},|\xi_3|\right),\left({\tau_2 \over \tau_1\tau_3},{1 \over |\xi_3|}\right)\right\}.
$$
In the following, all the couples $(f,v)$ are with $v=0$, so we omit $v$ for short.
\begin{itemize}
\item If $\xi_1=0$, $\xi_2,\xi_3 \neq 0$, $\displaystyle{F^1_A=\left\{\tau_1,\tau_2,\tau_3,{\tau_1\tau_3 \over \tau_2},{\tau_2 \over \tau_3}\right\}}$.
\item If $\xi_2=0$, $\xi_1,\xi_3 \neq 0$, $\displaystyle{F^1_A=\left\{\tau_1,{\tau_2 \over \tau_1},\tau_3,{\tau_2 \over \tau_1\tau_3}\right\}}$.
\item If $\xi_3=0$, $\xi_1,\xi_2 \neq 0$, $\displaystyle{F^1_A=\left\{\tau_1,{\tau_2 \over \tau_1},{\tau_1\tau_3 \over \tau_2}\right\}}$.
\item If $\xi_2 \neq 0$, $\xi_1,\xi_3=0$, $\displaystyle{F^2_A=\left\{\tau_1,\tau_2,\tau_3,{\tau_1\tau_3 \over \tau_2}\right\}}$.
\end{itemize}
We study the restriction when we consider $\chi_B=\{M_1,M_2,M_3\}$ obtained by adding a manifold $M_3$ to $\chi_A$, i.e. adding a line $\alpha=(\alpha_1,\alpha_2,\alpha_3)$ to $A$.
\begin{itemize}
\item Let $M_3=\{x_1=0\}$. We consider the matrix $B$ obtained by adding the line $\alpha=(1,0,0)$ to $A$.
\begin{itemize}
\item If $\xi_1=0$ condition \eqref{2,5} is satisfied.
\item If $\xi_2=0$, $\xi_1,\xi_3 \neq 0$, then $\log{e^{\alpha_2-\alpha_1}}=0-1<0$ and \eqref{2,5} is not satisfied. Indeed one can check that $\displaystyle{\tau_2 \over \tau_1} \notin \semiG_B$.
\item If $\xi_3=0$, $\xi_1,\xi_2 \neq 0$, then $\log{e^{\alpha_2-\alpha_1}}=0-1<0$ and \eqref{2,5} is not satisfied. Indeed one can check that $\displaystyle{\tau_2 \over \tau_1} \notin \semiG_B$.
\end{itemize}
\item Let $M_3=\{x_2=0\}$. We consider the matrix $B$ obtained by adding the line $\alpha=(0,1,0)$ to $A$.
\begin{itemize}
\item If $\xi_1=0$, $\xi_2,\xi_3 \neq 0$, then $\log{e^{\alpha_1+\alpha_3-\alpha_2}}=0+0-1<0$ and \eqref{2,5} is not satisfied. Indeed one can check that $\displaystyle{\tau_1\tau_3 \over \tau_2} \notin \semiG_B$.
\item If $\xi_2=0$ condition \eqref{2,5} is satisfied.
\item If $\xi_3=0$, $\xi_1,\xi_2 \neq 0$, then $\log{e^{\alpha_1+\alpha_3-\alpha_2}}=0+0-1<0$ and \eqref{2,5} is not satisfied. Indeed one can check that $\displaystyle{\tau_1\tau_3 \over \tau_2} \notin \semiG_B$.
\end{itemize}
\item Let $M_3=\{x_3=0\}$. We consider the matrix $B$ obtained by adding the line $\alpha=(0,0,1)$ to $A$.
\begin{itemize}
\item If $\xi_1=0$, $\xi_2,\xi_3 \neq 0$, then $\log{e^{\alpha_2-\alpha_3}}=0-1<0$ and \eqref{2,5} is not satisfied. Indeed one can check that $\displaystyle{\tau_2 \over \tau_3} \notin \semiG_B$.
\item If $\xi_2=0$, $\xi_1,\xi_3 \neq 0$, then $\log{e^{\alpha_2-\alpha_3}}=0-1<0$ and \eqref{2,5} is not satisfied. Indeed one can check that $\displaystyle{\tau_2 \over \tau_3} \notin \semiG_B$.
\item If $\xi_3=0$ condition \eqref{2,5} is satisfied.
\end{itemize}
\item Let $M_3=\{x_1=x_3=0\}$. We consider the matrix $B$ obtained by adding the line $\alpha=(1,0,1)$ to $A$.
\begin{itemize}
\item If $\xi_1=0$, $\xi_2,\xi_3 \neq 0$, then $\log{e^{\alpha_2-\alpha_3}}=0-1<0$ and \eqref{2,5} is not satisfied. Indeed one can check that $\displaystyle{\tau_2 \over \tau_3} \notin \semiG_B$.
\item If $\xi_2=0$, $\xi_1,\xi_3 \neq 0$, then $\log{e^{\alpha_2-\alpha_1-\alpha_3}}=0-1-1<0$ and \eqref{2,5} is not satisfied. Indeed one can check that $\displaystyle{\tau_2 \over \tau_1\tau_3} \notin \semiG_B$.
\item If $\xi_3=0$, $\xi_1,\xi_2 \neq 0$, then $\log{e^{\alpha_2-\alpha_1}}=0-1<0$ and \eqref{2,5} is not satisfied. Indeed one can check that $\displaystyle{\tau_2 \over \tau_1} \notin \semiG_B$.
\item If $\xi_1=\xi_3=0$ condition \eqref{2,5} is satisfied.
\end{itemize}
\item Let $M_3=\{x_1=x_2=x_3=0\}$. We consider the matrix $B$ obtained by adding the line $\alpha=(1,1,1)$ to $A$.
\begin{itemize}
\item If $\xi_1=0$ condition \eqref{2,5} is satisfied.
\item If $\xi_2=0$, $\xi_1,\xi_3 \neq 0$, then $\log{e^{\alpha_2-\alpha_1-\alpha_3}}=1-1-1<0$ and \eqref{2,5} is not satisfied. Indeed one can check that $\displaystyle{\tau_2 \over \tau_1\tau_3} \notin \semiG_B$.
\item If $\xi_3=0$ condition \eqref{2,5} is satisfied.
\end{itemize}
\end{itemize}
\end{es}

\begin{es}  In $\R^4$, let us consider the family $\chi_A=\{M_1,M_2\}$, with $M_1=\{x_1=x_2=x_4=0\}$, $M_2=\{x_2=x_3=0\}$. The matrix $A$ is
$$
\left(
\begin{array}{cccc}
1 & 1 & 0 & 1\\
0 & 1 & 1 & 0
\end{array}
\right)
$$
If $\xi_3=0$, $\xi_1,\xi_2,\xi_4 \neq 0$, we have
$$\displaystyle{F^1_A=\left\{(\tau_1,0),\left({\tau_2 \over \tau_1},0\right),\left({\tau_1\tau_3 \over \tau_2},0\right),\left({\tau_4 \over \tau_1},|\xi_4|\right),\left({\tau_1 \over \tau_4},{1 \over |\xi_4|}\right)\right\}}.$$
We study the restriction when we consider $\chi_B=\{M_1,M_2,M_3\}$ obtained by adding a manifold $M_3$ to $\chi_A$, i.e. adding a line $\alpha=(\alpha_1,\alpha_2,\alpha_3,\alpha_4)$ to $A$.
\begin{itemize}
\item Let $M_3=\{x_2=0\}$. We consider the matrix $B$ obtained by adding the line $\alpha=(0,1,0,0)$ to $A$. Then $\log{e^{\alpha_1+\alpha_3-\alpha_2}}=0+0-1<0$ and \eqref{2,5} is not satisfied. Indeed one can check that $\left(\displaystyle{\tau_1\tau_3 \over \tau_2},0\right) \notin \semiG_B$.
\item Let $M_3=\{x_2=x_4=0\}$. We consider the matrix $B$ obtained by adding the line $\alpha=(0,1,0,1)$ to $A$. Then $\log{e^{\alpha_1-\alpha_4}}=0-1<0$ and \eqref{2,5} is not satisfied. Indeed one can check that $\left(\displaystyle{\tau_1 \over \tau_4},{1 \over |\xi_4|}\right) \notin \semiG_B$.
\item Let $M_3=\{x_1=x_2=x_3=0\}$. We consider the matrix $B$ obtained by adding the line $\alpha=(1,1,1,0)$ to $A$. Then $\log{e^{\alpha_1-\alpha_4}}=1-0>0$ and \eqref{2,5} is not satisfied. Indeed one can check that $\left(\displaystyle{\tau_1 \over \tau_4},0\right) \in \semiG_B$.
\item Let $M_3=\{x_1=x_2=x_3=x_4=0\}$. We consider the matrix $B$ obtained by adding the line $\alpha=(1,1,1,1)$ to $A$. Then \eqref{2,5} is satisfied.
\end{itemize}
\end{es}

\section{Multi-asymptotic expansions}\label{section:Multi-asymptotic expansions}
\

In this section we extend the notion of asymptotic expansion associated to $\chi$ developed in \cite{HP13}. We refer to \cite{Ma84} for the classical theory of strongly asymptotically developable functions.

\subsection{Geometry of multi-asymptotic expansions}\label{subsection:Geometry of multi-asymptotic expansions}
Let $X = \mathbb{C}^n$ with coordinates $(z)=(z_1, \dots, z_n)$,
and let $\chi = \{Z_j\}_{j=1}^\ell$ be a family of closed complex submanifolds in $X$.
We divide complex coordinates $(z)$ into $(m+1)$-coordinates blocks
\begin{equation}
(z^{(0)},\,z^{(1)},\,\dots, z^{(m)}),
\end{equation}
for which each $Z_j$ is assumed to be linearized, that is,
there exists the subset $K_j$ in $\{1,2,\dots,m\}$ such that
each $Z_j$ is defined by the equations $z^{(k)} = 0$ with $k \in K_j$, i.e.,
\begin{equation}
Z_j = \left\{z^{(k)} = 0\,\,(k \in K_j)\right\} \qquad (j = 1,2,\dots,\ell).
\end{equation}

Let us consider a matrix $A_\chi = (a_{jk})$
of size $\ell \times m$ whose entries
are non-negative rational numbers. Set $L := \operatorname{Rank} A_\chi$.
We always assume the condition \eqref{eq:L_L_invertibility}, that is,
\begin{equation}
\text{the first $L \times L$ submatrix in $A_\chi$ is invertible.}
\end{equation}
As in Subsection \ref{subsection:Semigroups of generators},
we define the associated $m$-monomials of $\lambda = (\lambda_1,\dots,\lambda_\ell)$ by
$$
\varphi_k(\lambda) = \lambda_1^{a_{1k}}\lambda_2^{a_{2k}}\cdots \lambda_\ell^{a_{\ell k}}
\quad (k=1,\dots,m).
$$
We also define the action $\mu_j: X \times \mathbb{R}^+ \to X$
($j=1,\dots,\ell$) by
\begin{equation}
\mu_j(z,\,t) = (z^{(0)},\, t^{a_{j1}}z^{(1)},
t^{a_{j2}}z^{(2)},\,\dots,\, t^{a_{jm}}z^{(m)})
\end{equation}
and the action $\mu:X \times (\mathbb{R}^+)^\ell \to X$ by
\begin{equation}
\mu(z,\,\lambda) = (z^{(0)},\, \varphi_1(\lambda)z^{(1)},
\varphi_2(\lambda)z^{(2)},\,\dots,\, \varphi_m(\lambda)z^{(m)}).
\end{equation}
%It follows from the definition that we have
%$$
%\mu(z,\, \lambda)
%= \mu_1(\mu_2( \dots (\mu_\ell(z, \lambda_\ell),\, \dots),\, \lambda_2),\, \lambda_1).
%$$
For these actions,
we assume Condition A2.~introduced in
Subsection \ref{subsection:The local model}, that is,
\begin{equation}
\text{each $Z_j$ coincides with the set of fixed points of the action $\mu_j$ on $X$.}
\end{equation}
Note that this condition is equivalently saying that $A_\chi$ satisfies the
condition
\begin{equation}{\label{eq:cond-a2-alt}}
\text{$a_{jk} \ne 0$ if and only if $k \in K_j$} \qquad (j = 1,2,\dots,\ell).
\end{equation}
According to the system of coordinates $(z^{(0)},\, \dots,\, z^{(m)})$,
the space $X$ can be identified with
$$
\mathbb{C}_{z^{(0)}}^{n_0} \times
\mathbb{C}_{z^{(1)}}^{n_1} \times \dots \times
\mathbb{C}_{z^{(m)}}^{n_m},
$$
where $n_k$ denotes the number of coordinate variables in the coordinates block $z^{(k)}$.
\begin{oss}{\label{oss:complex-norm}}
Through the section, for a complex vector $a = (a_1,\,\dots,a_k) \in \mathbb{C}^k$,
its norm is defined by
\begin{equation}
|a| := \max\{|a_1|,\,|a_2|,\,\dots,|a_k|\}.
\end{equation}
\end{oss}
Let $p$ be a point
$
(z^{(0)};\, \zeta) = (0;\, \zeta^{(1)}, \dots, \zeta^{(m)}) \in S_\chi
$
and set
\begin{equation}
J_Z(\zeta) = \{k \in \{1,\dots,m\};\, \zeta^{(k)} = 0\}
\end{equation}
as usual. Note that we consider the problem near the
origin, i.e., $z^{(0)} = 0$.
Further, for the point $p \in S_\chi$, we assume the condition
{\eqref{eq:xi_0_for_zero_vector}, that is,
\begin{equation}
\text{if the $k$-th column of $A_\chi$ is a zero vector, then $\zeta^{(k)} = 0$}.
\end{equation}
%and
%\begin{equation}
%\text{$|\zeta_k| = 1$ for $1 \le k \le L$ with $k \notin J_Z(\zeta)$}.
%\end{equation}

We denote by $V(\zeta)$ the set of $\mathbb{R}^+$-cones in the form
\begin{equation}
\mathbb{C}^{n_0} \times V_1 \times \dots \times V_m \subset X,
\end{equation}
where each $V_k$ satisfies the following conditions:
\begin{itemize}
\item If $\zeta^{(k)} \ne 0$, then
$V_k$ is an open convex proper cone in
$\mathbb{C}^{n_k}$ containing the point $\zeta^{(k)}$.
\item If $\zeta^{(k)} = 0$, then $V_k = \mathbb{C}^{n_k}$.
\end{itemize}
%is an $\RP$-conic open convex subset in
%$\mathbb{C}_{z^{(k)}}^{n_k}$ containing the point $\zeta^{(k)}$.
%Note that, when $\zeta^{(k)} = 0$, an allowable $V_k$
%is only the whole space $\mathbb{C}_{z^{(k)}}^{n_k}$ itself.

In what follows, we constantly use notions introduced in
Section \ref{section:Multi-cones for a general case}.
Remember that the semigroup $\mathcal{G}$ was defined in
Subsection \ref{subsection:Semigroups of generators}, and
that $F^q$ (resp. $G$) is the finite subset of $\RRpair$ (resp. $\RRTpair$)
defined in the same subsection, where
we adopt the definition of $G$ given in Remark \ref{oss:general_G}.
%Since the point $p \in S_\chi$
%is outside fixed points and {\eqref{eq:normalized_condition_asymp}} is assumed,
%we have $q = 0$ and
That is,
$$
\begin{aligned}
F^0 &:= (\mathcal{L}^\lambda_\ell \circ \dots \circ \mathcal{L}^\lambda_{L+1})(G), \\
F^q &:= (\mathcal{L}_{k_q} \circ \dots \circ \mathcal{L}_{k_1})(F^0),
\end{aligned}
$$
where $\{k_1, \dots, k_q\} = J_Z(\zeta) \cap \{1,2,\dots,L\}$
$(1 \le k_1 < \dots < k_q \le L)$.
Note that,
by Proposition \ref{prop:radical-of-S},
the radical of $[F^q]$ is $\semiG$.

\begin{oss}
The most important and interesting cases are those where
the corresponding action $\mu_\chi$ is non-degenerate and $p$ is
outside the fixed points. In this case,
by adopting the coordinates blocks described
in Corollary \ref{cor:outside-fixed-points-coordinate},
we have $q = 0$ and $L = \operatorname{Rank} A_\chi = \ell$, and hence,
\begin{equation}
F^q = F^0 = G.
\end{equation}
\end{oss}
Set
\begin{equation}
|z^{(*)}| = (|z^{(1)}|,\,\dots,\,|z^{(L)}|,\,\dots,\,|z^{(m)}|)
\end{equation}
and
\begin{equation}
|z^{(*)}|' = (|z^{(1)}|,\,\dots,|z^{(L)}|).
\end{equation}
Let $H$ be a finite subset in $\RRpair$ satisfying that
$[H]$ and $\semiG$ defined in \eqref{eq:def_S} are equivalent, and
let $\{\epsilon\}$ be the set of positive real numbers
consisting of $\epsilon_0 > 0$, $\epsilon_{f,+} > 0$ and $\epsilon_{f,-} > 0$ for
each $(f,\,v) \in H$.
Then, for $V \in V(\zeta)$, we define the open subset
\begin{equation}{\label{def:multi-cone-complex}}
\Ssec
:= \left\{z \in X;\,
\begin{aligned}
&z \in V,\,\, |z^{(0)}| < \epsilon_0, \\
&v - \epsilon_{f,-} < f\left(|z^{(*)}|\right) < v+\epsilon_{f,+} \\
&\text{for any $(f,\,v) \in H$}.
\end{aligned}
\right\},
\end{equation}
which is often called a {\it{multicone}} or {\it{multisector}}.
Since $[H]$ and $\semiG$ are equivalent and
$(\tau_k^N,\, 0) \in \semiG$ ($k=1,2,\dots,m$)
holds for some $N \in \mathbb{N}$,
any open neighborhood of the origin contains the set $\Ssec$
if the positive real numbers in $\{\epsilon\}$ are sufficiently small.

Furthermore, for finite subsets $H$ and $H'$ in $\RRpair$
such that both the $[H]$ and $[H']$ are
equivalent, the families
$\left\{S_H(V,\,\{\epsilon\})\right\}_{\{\epsilon\}}$
and
$\left\{S_{H'}(V,\,\{\epsilon'\})\right\}_{\{\epsilon'\}}$
are equivalent with respect to inclusion of sets. That is,
for any $S_H(V,\,\{\epsilon\})$ (resp. $S_{H'}(V,\,\{\epsilon'\})$), there
exists $\{\epsilon'\}$ (resp. $\{\epsilon\}$) such that
$$
S_H(V,\{\epsilon\}) \subset S_{H'}(V,\,\{\epsilon'\}) \quad
\left(\text{resp. }\,\, S_{H'}(V,\{\epsilon'\}) \subset S_{H}(V,\,\{\epsilon\})\right)
$$
holds.
\begin{lem}{\label{lem:power_equivalent}}
Let $\semiG'$  be a semigroup in $\RRpair$ which is equivalent
to $\semiG$. Then, the family
$\left\{S_H(V,\{\epsilon\})\right\}_{H,\,\{\epsilon\}}$
and the one $\left\{S_{H'}(V,\,\{\epsilon'\})\right\}_{H',\,\{\epsilon'\}}$
are exactly same where
the former $H$ runs through the finite subset
of $\semiG$ with $[H]$ and $\semiG$ being equivalent
and the latter $H'$ runs through the finite subset of $\semiG'$
with $[H']$ and $\semiG'$ being equivalent.
That is, for any
$S_H(V,\,\{\epsilon\})$
$($resp. $S_{H'}(V,\,\{\epsilon'\}$$)$, there exists
a finite subset $H' \subset \semiG'$
$($resp. $H \subset \semiG$$)$ and $\{\epsilon'\}$ $($resp. $\{\epsilon\}$$)$
such that $S_H(V,\{\epsilon\}) = S_{H'}(V,\{\epsilon'\})$ holds.
\end{lem}
\begin{proof}
Let $H$ be a finite subset in $\semiG$ satisfying that $[H]$ and $\semiG$
are equivalent.
As $\semiG'$ is equivalent to $\semiG$, there exists a natural number $N$
such that $\langle H \rangle_N \subset \semiG'$ holds,
where $\langle H \rangle_N$ is the finite subset
$
\{(f,\,v)^N;\, (f,\,v) \in H\}.
$
Note that $[\langle H \rangle_N]$ is equivalent to $\semiG'$.

On the other hand, it follows from the definition that,
for any $N \in \mathbb{N}$ and $\{\epsilon\}$,
there exists a suitable $\{\epsilon'\}$ for which we have
$$
S_H(V,\,\{\epsilon\}) = S_{\langle H \rangle_N}(V,\,\{\epsilon'\}).
$$
Hence the claim of the lemma follows from this.
\end{proof}
\begin{oss}
When $Q(H) = H$ holds, for any $(f,\,v) \in H$ with $v \ne 0$, we have
$(1/f,\, 1/v) \in H$. Hence, by suitably choosing $\{\epsilon\}$ again,
the multicone defined by \eqref{def:multi-cone-complex} coincides with
the subset
\begin{equation}{\label{def:alt-multi-cone-complex}}
\left\{z \in X;\,
\begin{aligned}
&z \in V,\,\, |z^{(0)}| < \epsilon_0, \\
&f\left(|z^{(*)}|\right) < v+\epsilon_{f,+}\,\,\,
\text{for any $(f,\,v) \in H$}.
\end{aligned}
\right\}.
\end{equation}
\end{oss}
\begin{oss}
When $H = F^q \subset \semiG$ and all the positive real numbers
in $\{\epsilon\}$ are the same $\epsilon > 0$,
we write $S(V,\,\epsilon)$ instead of $\Ssec$.
\end{oss}
\begin{es}
Let $X = \mathbb{C}^n$ with coordinates blocks
$(z^{(0)},\,z^{(1)},\, z^{(2)},\, z^{(3)})$.
The submanifolds $Z_j$ ($j=1,2,3$) are defined by
$$
\begin{aligned}
Z_1 &= \{z^{(2)} = z^{(3)} = 0\}, \\
Z_2 &= \{z^{(1)} = z^{(3)} = 0\}, \\
Z_3 &= \{z^{(3)} = 0\}.
\end{aligned}
$$
Let us consider the $3 \times 3$ matrix $A_\chi$ associated with
the family $\chi = \{Z_1,\,Z_2,\,Z_3\}$.
Then clearly we have $\operatorname{Rank} A_\chi = 3$ and $\chi$ is
of normal type. In this case, the polynomials $\varphi_j(\lambda)$ is given by
$$
\varphi_1(\lambda) = \lambda_2,\,\,
\varphi_2(\lambda) = \lambda_1,\,\,
\varphi_3(\lambda) = \lambda_1 \lambda_2 \lambda_3,
$$
and hence, we have
$$
\varphi^{-1}_1(\tau) = \tau_1,\,\,
\varphi^{-1}_2(\tau) = \tau_2,\,\,
\varphi^{-1}_3(\tau) = \dfrac{\tau_3}{\tau_1 \tau_2}.
$$
Let $p = (0;\,\zeta) = (0;\,\zeta^{(1)}, \zeta^{(2)}, \zeta^{(3)}) \in S_\chi$
be a point outside fixed points.
Since $G$ is given by
$$
\{(\varphi^{-1}_1,\,0),\, (\varphi^{-1}_2,\,0),\, (\varphi^{-1}_3,\,0)\}
$$
and $F^q = G$ holds in this case, the cone $S(V,\epsilon)$ is defined by
$$
\left\{(z^{(0)},\,z^{(1)},\, z^{(2)},\, z^{(3)}) \in X;\,
\begin{aligned}
& z^{(k)} \in V_k\,(k=1,2,3),\,\, |z^{(j)}| < \epsilon\,\,(j=0,1,2,3)\\
&|z^{(3)}| < \epsilon |z^{(1)}| |z^{(2)}|\,\,
\end{aligned}
\right\},
$$
where $V_j$ ($j=1,2,3$) is a proper open cone in $\mathbb{C}^{n_j}$ containing
the vector $\zeta^{(j)}$.
\end{es}
\begin{es}
Let $X = \mathbb{C}^n$ with coordinates blocks
$(z^{(0)},\, z^{(1)},\, z^{(2)},\, z^{(3)},\, z^{(4)})$.
The submanifolds $Z_j$ ($j=1,2,3,4$) are defined by
$$
\begin{aligned}
Z_1 &= \{z^{(1)} = z^{(4)} = 0\}, \\
Z_2 &= \{z^{(2)} = z^{(4)} = 0\}, \\
Z_3 &= \{z^{(3)} = z^{(4)} = 0\}, \\
Z_4 &= \{z^{(1)} = z^{(2)} = z^{(3)} = z^{(4)} = 0\}.
\end{aligned}
$$
Let us consider the $4 \times 4$ matrix $A_\chi$ associated with
the family $\chi = \{Z_1,\,Z_2,\,Z_3,\,Z_4\}$.
Then $\operatorname{Rank} A_\chi = 4$ holds and $\chi$ is
of normal type.
In this case, the polynomials $\varphi_j(\lambda)$ is given by
$$
\varphi_1(\lambda) = \lambda_1 \lambda_4,\,\,
\varphi_2(\lambda) = \lambda_2 \lambda_4,\,\,
\varphi_3(\lambda) = \lambda_3 \lambda_4,\,\,
\varphi_4(\lambda) = \lambda_1 \lambda_2 \lambda_3 \lambda_4,
$$
and hence
$$
\varphi^{-1}_1(\tau) = \sqrt{\dfrac{\tau_1 \tau_4}{\tau_2 \tau_3}},\,\,
\varphi^{-1}_2(\tau) = \sqrt{\dfrac{\tau_2 \tau_4}{\tau_1 \tau_3}},\,\,
\varphi^{-1}_3(\tau) = \sqrt{\dfrac{\tau_3 \tau_4}{\tau_1 \tau_2}},\,\,
\varphi^{-1}_4(\tau) = \sqrt{\dfrac{\tau_1 \tau_2 \tau_3}{\tau_4}}.
$$
Let $p = (0;\,\zeta) = (0;\, \zeta^{(1)}, \zeta^{(2)}, \zeta^{(3)}, \zeta^{(4)})
\in S_\chi$ be a point outside fixed points.
By the same reasoning as that in the previous example,
the cone $S(V,\epsilon)$ is defined by
$$
\left\{(z^{(0)},\,z^{(1)},\, z^{(2)},\, z^{(3)},\, z^{(4)}) \in X;\,
\begin{aligned}
& z^{(k)} \in V_k\,\,(k=1,2,3,4), \,\, |z^{(0)}| < \epsilon\\
&|z^{(3)}||z^{(4)}| < \epsilon^2 |z^{(1)}||z^{(2)}| \\
&|z^{(2)}||z^{(4)}| < \epsilon^2 |z^{(1)}||z^{(3)}| \\
&|z^{(1)}||z^{(4)}| < \epsilon^2 |z^{(2)}||z^{(3)}| \\
&|z^{(1)}||z^{(2)}||z^{(3)}| < \epsilon^2 |z^{(4)}|
\end{aligned}
\right\},
$$
where $V_j$ ($j=1,2,3,4$) is a proper open cone in $\mathbb{C}^{n_j}$ containing
the vector $\zeta^{(j)}$.
\end{es}

%$\check{F}^q$ denotes the subset of $F^q$ consisting
%of pairs $(f,\,v)\in F^q$ for which $f(\tau)$ contains a negative power of
%some $\tau_k$, that is, there exists some $k \in \{1,\dots,m\}$ such that $\nu_k(f(\tau)) < 0$.
%
\begin{lem}{\label{lem:1-regular}}
Assume the positive real numbers in $\{\epsilon\}$ to be sufficiently small. Then
$\Ssec$ is $1$-regular, that is, there exists a constant $C>0$
satisfying that, for any points $p$ and $q$ in $\Ssec$, there
exists a rectifiable curve in $\Ssec$ which joins $p$ and $q$ and
whose length is less than or equal to $C|p -q|$.
\end{lem}
\begin{proof}
We first show the fact that the subset $D$ in $\mathbb{R}^{m}$ defined by
$$
\left\{(\tau_1,\dots,\tau_m) \in \mathbb{R}^m;\,
\begin{aligned}
&\tau_k > 0\,\, (k=1,\dots,m),\\
&v - \epsilon_{f,-} < f(\tau) < v+\epsilon_{f,+} \\
&\text{for any $(f,\,v) \in H$}.
\end{aligned}
\right\}
$$
is $1$-regular. Take points
$$
a = (a_1,\, \dots,\, a_m) \in D
$$
and
$$
b = (b_1,\, \dots,\, b_m) \in D.
$$
Consider the closed curve parameterized by $\theta$ ($0 \le \theta \le 1$)
$$
(a_1^\theta b_1^{(1-\theta)},\,\dots,\, a_m^\theta b_m^{(1-\theta)}).
$$
Since $f \in H$ is a monomial of the $\tau$ variables,
the curve is contained in $D$. Furthermore, as $a_k^\theta b_k^{(1-\theta)}$
is a monotonic function of $\theta$, the length of the curve is dominated by
$$
\sum_{k=1}^m|a_k - b_k|  \le m |a-b|.
$$
Hence $D$ is $1$-regular.  Now define
$$
Y = \{z \in \Ssec;\, |z^{(1)}| |z^{(2)}| \cdots |z^{(m)}| = 0\}.
$$
Since $\Ssec \setminus Y$ is dense in $\Ssec$,
we may assume $p$ and $q$ in the lemma to be outside $Y$.
Let
$$
p = (\hat{z}^{(0)},\, \hat{z}^{(1)},\, \dots,\,\hat{z}^{(m)}) \in \Ssec \setminus Y
$$
and
$$
q = (\hat{w}^{(0)},\, \hat{w}^{(1)},\, \dots,\,\hat{w}^{(m)}) \in \Ssec \setminus Y.
$$
Then we take points
$$
p_1 = (\hat{w}^{(0)},\, \hat{z}^{(1)},\, \dots,\,\hat{z}^{(m)})
\in S_H(V,\, \epsilon) \setminus Y
$$
and
$$
p_2 = \left(\hat{w}^{(0)},\,  \dfrac{|\hat{w}^{(1)}|}{|\hat{z}^{(1)}|}
\hat{z}^{(1)},\, \dots,\,
\dfrac{|\hat{w}^{(m)}|}{|\hat{z}^{(m)}|} \hat{z}^{(m)}\right)
\in \Ssec \setminus Y.
$$
We join these points successively:
First $p$ and $p_1$ are joined by the straight segment.
Let $P$ be the product of the spherical surfaces defined by
$$
\{(z^{(0)},\,z^{(1)},\,\dots,\,z^{(m)});\,
|z^{(k)}| = |\hat{w}^{(k)}|\,\,\, (k=1,\dots,m)\}.
$$
Then the point $p_2$ and $q$ are connected by the shortest arc
in the surface $P$. Since $D$ is $1$-regular, there exists a curve
$(\tau_1(\theta), \dots, \tau_m(\theta))$ $(0 \le \theta \le 1)$ in $D$ which
joins the points
$
a = (|z^{(1)}|, \dots, |z^{(m)}|)
$
and
$
b = (|w^{(1)}|, \dots, |w^{(m)}|)
$
and whose length is bounded by $m|a - b|$.
Then $p_1$ and $p_2$ are jointed by the
closed curve $(0 \le \theta \le 1)$
$$
\left(\hat{w}^{(0)},\,
\dfrac{\tau_1(\theta)}{|\hat{z}^{(1)}|}
\hat{z}^{(1)},\, \dots,\,
\dfrac{\tau_m(\theta)}{|\hat{z}^{(m)}|} \hat{z}^{(m)}\right).
$$
Clearly these curves are contained in $\Ssec$ and
the total length of these curves is bounded by a constant multiple of $|p-q|$.
\end{proof}

By the same method as above, we can also prove the following corollary:

\begin{cor} \label{cor: l.c.t. cone}
Assume the positive real numbers in $\{\epsilon\}$ to be sufficiently small. Then,
for any point $p \in X$, there exists a family $\{U_\kappa\}_{\kappa > 0}$
of open subanalytic neighborhoods of $p$ such that $U_\kappa \cap \Ssec$
is contractible. In particular, $\Ssec$ is locally cohomologically trivial.
\end{cor}
\begin{proof}
We may assume $I_0$ to be empty.
Then every point has a fundamental open neighborhood in the form
$$
\{(z^{(1)},\,\dots,z^{(m)});\,
a_k < |z^{(k)}| < b_k,\, z^{(k)} \in T_k \,\,(k=1,\dots,m)\},
$$
where $a_k$, $b_k$ and $T_k$ satisfy the one of the following conditions:
\begin{enumerate}
\item $0 < a_k < b_k$ and $T_k$ is an open convex proper cone
in $\mathbb{C}^{n_k}_{z^{(k)}}$.
\item $a_k < 0 < b_k$ and $T_k = \mathbb{C}^{n_k}_{z^{(k)}}$.
\end{enumerate}
Let $U$ be such an open subset.

We first suppose $\zeta^{(k)} \ne 0$ for any $k \in \{1,2,\dots,m\}$
and $U \cap \Ssec$
is non-empty. Take a point $q$ in $U \cap \Ssec$ and fix it. Then, by the
same argument as that for the proof of the previous lemma,
every point $p \in U \cap \Ssec$ can be joined with $q$ by an
explicitly specified path in $U \cap \Ssec$. Hence the family of these
paths give a contraction map.

Now let us consider the general case. Set
$$
K_Z := \{k \in \{1,\dots,m\};\, \zeta^{(k)} = 0,\,\, a_k < 0\}.
$$
For any $p = (z_*^{(1)},\,\dots,\, z_*^{(m)}) \in U \cap \Ssec$,
we determine the point
$\tilde{p} = (\tilde{z}_*^{(1)}, \, \dots,\, \tilde{z}_*^{(m)})$ by
$\tilde{z}_*^{(k)} = 0$ if $k \in K_Z$ and
$\tilde{z}_*^{(k)} = z_*^{(k)}$ if $k \notin K_Z$.
Then the straight segment joining $p$ and $\tilde{p}$ is completely contained in
$U \cap \Ssec$ because, for $(f(\tau),\,v) \in H$,
the variable $\tau_k$ with $k \in K_Z$  appears only in the numerator of $f$,
and further, $v = 0$ when $\tau_k$ appears in the numerator of $f$.
Therefore we first change every point $p$ in $U \cap \Ssec$ to $\tilde{p}$,
and then, apply the previous argument to the only coordinates
$z^{(k)}$ ($k \notin K_Z$), i.e., $z^{(k)}$ with $k \in K_Z$ is fixed to be the origin
in $\mathbb{C}^{n_k}$.
\end{proof}

Hereafter we always assume that the positive real numbers in $\{\epsilon\}$ are sufficiently small so that
the above lemma and corollary hold for $\Ssec$.
%We assume that, in what follows,
%\begin{equation}
%\text{$S(V,\,\epsilon)$ is 1-regular for any $V \in V(\zeta)$ and
%any sufficient small $\epsilon > 0$.}
%\end{equation}

We denote by $\subsetl$ the set of all the subsets of $\{1,\,2,\,\dots,\,\ell\}$
except for the empty set. For any $J \in \subsetl$, we set
\begin{equation}
K_J := \bigcup_{j \in J} K_j,\qquad Z_J := \bigcap_{j \in J} Z_j.
\end{equation}

For any subset $K \subset \{1,\,2,\,\dots,\,m\}$,
we denote by $z^{(K)}$ the set of the coordinates blocks $z^{(k)}$'s with $k \in K$.
%and by $z^{(K)}_{\mathrm{C}}$ the set of the coordinates
%which do not belong to $z^{(K)}$, i.e.,
%$z^{(k)}$'s with $k \in \{0,1,\dots,m\} \setminus K$.
We also set
$$
\compl K := \{0,1,\dots,m\} \setminus K.
$$
Note that the coordinates of $Z_J$ are given by $z^{(\compl K_J)}$.
Let $\pi_J$ denote the canonical projection from $X$ to $Z_J$ defined
by $z \to z^{(\compl K_J)}$.

\begin{lem}{\label{lemma:multi-cone-geometric}}
Let $V \in V(\zeta)$, and let $\{\epsilon\}$ be the set of positive real numbers.
Then we have the followings:
\begin{enumerate}
\item For a $\lambda = (\lambda_1, \dots, \lambda_\ell)$ with $0 < \lambda_j \le 1$
$(j = 1,2,\dots,\ell)$, we have
$$
\mu(\Ssec,\, \lambda) \subset \Ssec.
$$
\item For a $J \in \subsetl$, we have
$$
\overline{\pi_J(\Ssec)} = Z_J \cap \overline{\Ssec}
$$
\end{enumerate}
\end{lem}
\begin{proof}
The first part is already proved in Section \ref{section:Multi-cones for a general case}, Corollary \ref{cor: stability under contraction}.

Let us show the second claim: Set $S := \Ssec$ for short.  Clearly
$$
%\overline{\pi_J(S)} \supset \pi_J\left(\overline{S}\right)
%\cap Z_J = \overline{S} \cap Z_J
\overline{\pi_J(S)} \supset \pi_J\left(\overline{S}\right)
\supset \overline{S} \cap Z_J
$$
hold. Now we will show the converse inclusion.
For $\lambda = (\lambda_1, \dots, \lambda_\ell)$
and $J = \{j_1, \dots, j_r\}$ ($j_1 < \dots < j_r$), set
$
\lambda_J := (\lambda_{j_1}, \dots, \lambda_{j_r})
$
and define $\mu_J: X \times (\mathbb{R}^+)^{\#J} \to X$
by
$$
\mu_J(z, \lambda_J)
= \mu_{j_1}(\mu_{j_2}( \dots (\mu_{j_r}(z, \lambda_{j_r}),\, \dots),\,
\lambda_{j_2}),\, \lambda_{j_1}).
$$
Set
$$
\tilde{\mu}_J(z, t) = \mu_J(z, t,\,\dots,\,t) \quad (t \in \mathbb{R}).
$$
Suppose $z_\infty \in \overline{\pi_J(S)}$.
Then we can find a sequence $\{z_m\}_{m=1}^\infty$
such that $z_m \in S$ and $\pi_J(z_m) \to z_\infty$.
Since $\pi_J(z_m) = \tilde{\mu}_J(z_m, 0)$ holds, the sequence
$\tilde{z}_m := \tilde{\mu}_J(z_m, 1/m)$ tends to $z_\infty$.
By the claim 1.~we have $\tilde{z}_m \in S$, from which $z_\infty \in \overline{S}$ follows.
This completes the proof.
\end{proof}
For a $J \in \subsetl$ and for an $S := \Ssec$, we set
$S_J := \pi_J(S)$, and by the above lemma, we have
\begin{equation}
	\overline{S_J} = \overline{S} \cap Z_J.
\end{equation}

$S_J$ enjoy several good properties as those in $S$.
Let $\mu|_{Z_J}$ (resp. $\mu_j|_{Z_J}$)
be the restriction of the action $\mu$ (resp. $\mu_j$) to $Z_J$.
\begin{prop}{\label{prop:multi-cone-S_J}}
Assume $Q(H) = H$. Then, for $S:=\Ssec$ and $J \in \subsetl$,
the subset $S_J$ is a multicone in $Z_J$ with respect the action $\mu|_{Z_J}$.
In particular, we have:
\begin{enumerate}
\item
$S_J$ is stable under contraction by the action $\mu_j|_{Z_J}$,
that is, $\mu_j|_{Z_J}(z,t) \in S_J$ if
$z \in S_J$ and $0 < t \le 1$.
\item $S_J$ is $1$-regular and locally cohomologically trivial if $\{\epsilon\}$
consists of sufficiently small positive numbers.
\end{enumerate}
\end{prop}
\begin{proof}
It follows from Proposition \ref{prop:equivalence_G_hat_G} that
$\semiG$ and $\widehat{\semiG}$ are equivalent. Hence,
by Lemma {\ref{lem:power_equivalent}}, we can assume that
$H \subset \widehat{\semiG}$ and $[H]$ and $\widehat{\semiG}$ are
equivalent.
Furthermore, since $Q(H) = H$ holds, we also assume that $\Ssec$ is defined by
the special form \eqref{def:alt-multi-cone-complex}.
Suppose that $K_J$ is non-empty and
$$
K_J = \{m'+1,\,m'+2,\,\dots,\,m\}
$$
for some $0 \le m' < m$. Set
$$
\tau_J = (\tau_1,\,\dots,\,\tau_{m'})
$$
and denote by $\mathcal{H}_{\tau_J}$ the set of monomials of rational powers
of the variables $\tau_J$ with coefficients $1$.

Recall that $\widehat{\semiG}$ is given by
$$
(\RRpair)\, \bigcap\, [\,\widehat{G}\,],
$$
where $\widehat{G}$ is
$$
Q\left(\left(\dfrac{\tau_1}{\varphi_1(\lambda)},\,|\zeta^{(1)}|\right),\,\cdots,\,
\left(\dfrac{\tau_m}{\varphi_m(\lambda)},\,|\zeta^{(m)}|\right),\,
(\lambda_1,\,0),\,\cdots,\, (\lambda_\ell,\,0)\right).
$$
The corresponding semigroup $\widehat{\semiG}_J$ with respect to
the action $\mu|_{Z_J}$ in $Z_J$ is given by
$$
(\mathcal{H}_{\tau_J} \times \mathbb{R}_{\ge 0})\,\bigcap\, [\,\widehat{G}_J\,],
$$
where $\widehat{G}_J$ is
$$
Q\left(\left(\dfrac{\tau_1}{\varphi_1(\lambda)},\,|\zeta^{(1)}|\right),\,\cdots,\,
\left(\dfrac{\tau_{m'}}{\varphi_{m'}(\lambda)},\,|\zeta^{(m')}|\right),\,
(\lambda_1,\,0),\,\cdots,\, (\lambda_\ell,\,0)\right).
$$

By noticing the fact that the variable $\tau_k$ appears only in the term
$\tau_k/\varphi_k(\lambda)$, we can easily see
$$
(\mathcal{H}_{\tau_J} \times \mathbb{R}_{\ge 0})\, \bigcap\, \widehat{\semiG}
= \widehat{\semiG}_J.
$$
Recall the operation $\mathcal{L}_k$ defined in
Subsection \ref{subsection:Families of cones}. Let us define
the slightly modified operation $\widehat{\mathcal{L}}_k$ which consists
of the only F1 and F3 in the definition of $\mathcal{L}_k$.
Then, as the same reasoning as that for the claim 2.~in
Proposition \ref{prop:radical-of-S}, the finite subset defined by
$$
H_J:=(\widehat{\mathcal{L}}_m \circ \dots \circ \widehat{\mathcal{L}}_{m'+1})(H)
$$
is the radical of $\widehat{\semiG}_J$. Now, by repeated application of
 Lemma \ref{lem:proj_W} below,
we can conclude that $S_J = \pi_J(S)$ is a multicone $S_{H_J}(\pi_J(V),\{\epsilon'\})$
in $Z_J$ with suitably chosen $\{\epsilon'\}$.
This completes the proof.
\end{proof}

Before stating the lemma, we need some preparations:
Let $T$ be the product space $T_1 \times T_2 \times \cdots \times T_m$
with the coordinates $\tau = (\tau_1,\dots,\tau_m)$ where each $T_k$
is either $\RP$ or $\mathbb{R}_{\ge 0}$.
Set
$
T' := T_2 \times \cdots \times T_m
$
with the coordinates $\tau' = (\tau_2,\dots,\tau_m)$ and
$
\pi_1: T \to T'
$
is the canonical projection $\tau \to \tau'$.

Let $\{g_s\}_{s \in \Lambda}$ be a finite subset in $\mathcal{H}_\tau$. We assume that
each $g_s$ is well defined on $T$, that is, $\nu_k(g_s) \ge 0$
for any $1 \le k \le m$ with $T_k = \mathbb{R}_{\ge 0}$.
Let $\{a_s\}_{s \in \Lambda}$ be a set of positive real numbers and define
the subset
$$
W:=\{\tau \in T;\, g_s(\tau) < a_s \,\,(s \in \Lambda)\}.
$$
Then, under this situation, we have the lemma below:
\begin{lem}\label{lem:proj_W}
The subset $\pi_1(W) \subset T'$ is defined by the following family of inequalities:
\begin{equation}{\label{eq:projection_type_F1}}
g_s < a_s
\end{equation}
for any $g_s$ $(s \in \Lambda)$ with $\nu_1(g_s) = 0$ and
\begin{equation}{\label{eq:projection_type_F3}}
g_s^{a} g_{s'}^{b}
< a_s^{a} a_{s'}^{b}
\end{equation}
for any pair $g_s$ and $g_{s'}$
$(s,s' \in \Lambda)$
with $\nu_1(g_s) < 0$ and $\nu_1(g_{s'}) > 0$, and natural numbers $a$ and $b$ are
taken to be prime to each other with $a |\nu_1(g_s)| = b |\nu_1(g_{s'})|$.
\end{lem}
\begin{oss}
The inequalities {\eqref{eq:projection_type_F1}} and
{\eqref{eq:projection_type_F3}}
are counterparts of F1 and F3 of the operation $\mathcal{L}_1$, respectively.
\end{oss}
\begin{proof}
If $T_1 = \mathbb{R}_{\ge 0}$, then the lemma trivially holds. Hence we assume
$T_1 = \RP$ and, for simplicity,
we may also assume $\nu_1(g_s) \ne 0$ for any $s \in \Lambda$.
In this case, $\tau' \in \pi_1(W)$ if and only if we have
$$
\underset{s \in \Lambda,\,\nu_1(g_s) < 0}{\max}
\operatorname{sol}(g_s,a_s)(\tau')
\,\,\,<\,\,\,
\underset{s \in \Lambda,\,\nu_1(g_s) > 0}{\min} \operatorname{sol}(g_s,a_s)(\tau')
$$
where $\operatorname{sol}(g_s,a_s)(\tau')$ is the solution of the equation
$g_s(\tau_1,\,\tau') = a_s$ when we regarded it as an equation of $\tau_1$.
The latter condition is clearly equivalent to {\eqref{eq:projection_type_F3}}
for any pair $g_s$ and $g_{s'}$ $(s,s' \in \Lambda)$
with $\nu_1(g_s) < 0$ and $\nu_1(g_{s'}) > 0$.
\end{proof}
\begin{oss}
Some $\mu_j|_{Z_J}$ may be the identity action on $Z_J$.
Since the geometrical situation remains unchanged
if the corresponding deformation parameter $t_j$ is removed.
Hence we may
consider the situation where
such an action is ignored
and the corresponding deformation parameter is removed.
Furthermore, if there is a pair $\mu_j$ and $\mu_{j'}$ $(j \ne j')$
which gives the same induced actions $\mu_j|_{Z_J}$ = $\mu_{j'}|_{Z_J}$,
we can also eliminate such a duplication, for example,
by removing $\mu_j|_{Z_J}$ and the corresponding deformation parameter $t_j$.
\end{oss}

For a multicone, we can introduce the notion of a proper sub-multicone as follows:
Recall that,
for $V = \mathbb{C}^{n_0} \times V_1 \times \cdots \times V_m$
and
$V' = \mathbb{C}^{n_0} \times V'_1 \times \cdots \times V'_m$
in $V(\zeta)$, we say that $V$ is {\it{properly contained}} in $V'$
if $\overline{V_k}\setminus\{0\} \subset V'_k$ holds for any $k=1,2,\dots,m$.

\begin{df}\label{definition:proper sub-cone}
We say that $S=\Ssec$ is properly contained in
$S'=S_H(V',\,\{\epsilon'\})$ if $V$ is properly contained in $V'$ and
if the positive real number in $\{\epsilon\}$ is strictly smaller than
the corresponding one in $\{\epsilon'\}$, i.e.,
$\epsilon_0 < \epsilon'_0$,
$\epsilon_{f,-} < \epsilon'_{f,-}$ and
$\epsilon_{f,+} < \epsilon'_{f,+}$ hold for any $(f,v) \in H$.
\end{df}

\subsection{Asymptotic polynomials}\label{subsection:Asymptotic polynomials}
Let $J \in \subsetl$. We denote by $\mathbb{Z}^{J}_{\ge 0}$
the subset of $\mathbb{Z}^\ell_{\ge 0}$
\begin{equation}
\{(\lambda_1,\dots,\lambda_\ell) \in \mathbb{Z}^\ell_{\ge 0};\, \lambda_j = 0\,\,(j \notin J)\}
\end{equation}
which is sometimes identified with $\mathbb{Z}^{\# J}_{\ge 0}$ also.
Recall that the $(\ell \times m)$-matrix $A_\chi = (a_{jk})$
$(a_{jk} \in \mathbb{Q}_{\ge 0})$ defines the action
$\mu$ on $X$.
Let $\comD$ be a positive rational number so that
$a_{jk} \in \mathbb{Z}$ and all the $a_{jk}$ have no common divisors, and let us consider the
action on $X$ whose associated matrix is $\intA$ which is denoted by
$\mu_{\intA}$ hereafter. Note that $\mu_{\intA}(z,\lambda) = \mu(z,\,\lambda^\comD)$ holds.

Then we define,
for $\beta \in \mathbb{Z}^{J}_{\ge 0}$,
\begin{equation}{\label{eq:def_T_J_B}}
T_{J,\,\beta}\left(z^{(K_J)}\right) := \left.\dfrac{1}{\beta!}
\exp(-\mu_{\intA}(z,\, \lambda))\dfrac{\partial^\beta}{\partial \lambda^\beta}
\exp(\mu_{\intA}(z,\,\lambda))\right|_{\lambda = \check{e}_{J}},
\end{equation}
where $\check{e}_{J}\in \mathbb{C}^\ell$ denotes the point $(e_1, \dots, e_\ell)$ with
$e_j = 1$ if $j \notin J$ and $e_j = 0$ if $j \in J$. Let
$N = (n_1, \dots, n_\ell) \in \mathbb{Z}^\ell_{\ge 0}$ and $\beta \in \mathbb{Z}^{J}_{\ge 0}$.
Then we write
$
\beta <_J N
$
if and only if $\beta_j < n_j$ holds for any $j \in J$.

Now we define
\begin{equation}
T^{<N}_J\left(z^{(K_J)}\right) :=
\sum_{\beta \in \mathbb{Z}^{J}_{\ge 0},\, \beta <_J N}
T_{J,\,\beta}\left(z^{(K_J)}\right).
\end{equation}
Here, if the set $\{\beta \in \mathbb{Z}^{J}_{\ge 0};\, \beta <_J N\}$ of indices
is empty, we set $T^{<N}_J\left(z^{(K_J)}\right) := 0$ as usual convention.

For a multi-index $\alpha = (\alpha_1, \dots, \alpha_n) \in \mathbb{Z}^n_{\ge 0}$,
we denote by $\alpha^{(k)}$ ($k=0,\dots,m$) the part of a multi-index $\alpha$
which corresponds to the coordinates block $z^{(k)}$
%consisting of $\alpha_i$'s
%with $i \in \widehat{I}_k$
and call it an {\it{index block}}.
That is, $\alpha^{(k)}$ consists of $\alpha_i$'s
with $z_i \in z^{(k)}$,
where $z_i \in z^{(k)}$ means that the coordinate variable $z_i$ belongs
to the coordinates block $z^{(k)}$.
%which is a counterpart of the coordinate block $z^{(k)}$.
We denote by $|\alpha^{(k)}|$ the length of an index block $\alpha^{(k)}$, i.e.,
\begin{equation}
|\alpha^{(k)}| = \sum_{z_i \in z^{(k)}} \alpha_i.
\end{equation}
In subsequent arguments, $(\alpha_1,\, \dots,\, \alpha_n) \in \mathbb{Z}^n_{\ge 0}$
is often written in the form $(\alpha^{(1)}, \dots, \alpha^{(m)})$ of index blocks.
For a subset $K \subset \{1,\,2,\,\dots,\,m\}$,
we also define the subset of $\mathbb{Z}^n_{\ge 0}$ by
\begin{equation}
\mathbb{Z}^{(K)}_{\ge 0} := \left\{(\alpha^{(1)},\dots,\alpha^{(m)}) \in
\mathbb{Z}^n_{\ge 0};\, \alpha^{(k)} = 0\,\,(k \notin K) \right\}
\subset \mathbb{Z}^n_{\ge 0}.
\end{equation}

%The following lemma, in particular, shows that the polynomial $T^{<N}_J$ thus defined
%does not depend on the choice of $\comD \in \mathbb{N}$.
\begin{lem}{\label{eq:formula_development}}
For an $N = (n_1, \dots, n_\ell) \in \mathbb{Z}^\ell_{\ge 0}$, we have
\begin{equation}{\label{eq:lem-formula-polynomial}}
T^{<N}_J\left(z^{(K_J)}\right) = \sum_{\alpha \in A_J(N)} \dfrac{1}{\alpha!} z^\alpha,
\end{equation}
where $A_J(N) \subset \mathbb{Z}_{\ge 0}^{(K_J)}\, (\, \subset \mathbb{Z}_{\ge 0}^n\,)$ is defined by
\begin{equation}
A_J(N) = \left\{\alpha \in \mathbb{Z}^{(K_J)}_{\ge 0};\,
\sum_{k \in K_J} a_{jk}|\alpha^{(k)}| < n_j/\comD \,\,\, \text{for any $j \in J$}\right\}.
\end{equation}
\end{lem}
\begin{proof}
For $\beta \in \mathbb{Z}^J_{\ge 0}$ and
for any $k \in K_J$ and any $z_i$ in $z^{(k)}$,
we have
$$
\begin{aligned}
&\dfrac{\partial}{\partial z_i} T_{J,\,\beta}\left(z^{(K_J)}\right) =
\dfrac{\partial}{\partial z_i} \left(\left.\dfrac{1}{\beta!}
\exp(-\mu_{\intA}(z,\, \lambda)) {\partial_\lambda^\beta}
\exp(\mu_{\intA}(z,\,\lambda))\right|_{\lambda = \check{e}_{J}}\right) \\
\qquad&
= \left.\dfrac{1}{\beta!}
\exp(-\mu_{\intA}(z,\, \lambda)) {\partial_\lambda^\beta}
(\varphi_k(\lambda)^{\comD} \exp(\mu_{\intA}(z,\,\lambda)))\right|_{\lambda = \check{e}_{J}} \\
&
= \left.\dfrac{1}{(\beta \,-_J\, \comD A_\chi^{*,k})!}
\exp(-\mu_{\intA}(z,\, \lambda))
\partial_\lambda^{\beta\,-_J\,\comD A_\chi^{*,k}}
\exp(\mu_{\intA}(z,\,\lambda))\right|_{\lambda = \check{e}_{J}} \\
&
= T_{J,\,\left(\beta\,-_J\,\comD A_\chi^{*,k}\right)}\left(z^{(K_J)}\right).
\end{aligned}
$$
Here $A_\chi^{*,k}$ denotes the $k$-th column vector of the matrix $A_\chi$ and
we regarded $\sigma_A A_\chi^{*,k}$
as a multi-index in $\mathbb{Z}_{\ge 0}^\ell$. Furthermore the operation
$\alpha\, -_J\, \beta$ means that subtraction is performed only in $j$-th components with $j \in J$, that is,
for $\alpha = (\alpha_j)$, $\beta = (\beta_j)$ and $\gamma = (\gamma_j)$,
the equation $\gamma = \alpha\,-_J\, \beta$ holds if and only if
$\gamma_j = \alpha_j - \beta_j$ for $j \in J$ and $\gamma_j = \alpha_j$ for $j \notin J$.
We also note that,
when some entry of the vector $\beta\,-_J\,\comD A_\chi^{*,k}$ becomes negative, we set
$T_{J, \left(\beta\,-_J\, \comD A_\chi^{*,k}\right)} = 0$ as usual convention.
Therefore the argument goes in the same way as Lemma 7.6 in \cite{HP13}. However,
for reader's convenience, we repeat its argument here:
Let us show the lemma by  induction on $n = n_1 + \dots + n_\ell$.
If $n = 0$, as both sides of the equation
\eqref{eq:lem-formula-polynomial} are zero by definition,
the lemma is true.

Now we prove the lemma for a general $n > 0$.
Let us consider the system of partial differential equations
of an unknown function
$u(z^{(K_J)})$ defined by
\begin{equation}{\label{eq:partial-diff}}
	\frac{\partial}{\partial {z_i}} u =
	T^{<N -_J \sigma_A A_{\chi}^{*,k}}_J \quad
	(k \in K_J,\, z_i \text{ in } z^{(k)}).
\end{equation}
Then, by the induction hypothesis, the right hand side of the above equation
is given by
$$
T^{<N -_J \sigma_A A_{\chi}^{*,k}}_J =
\displaystyle\sum_{\alpha \in A_J(N-_J \sigma_A A_{\chi}^{*,k})} \frac{1}{\alpha!}z^\alpha.
$$
Then we can confirm that both $u=T^{<N}_{J}(z^{(K_J)})$ and
$
u=\displaystyle\sum_{\alpha \in A_J(N)} \frac{1}{\alpha!}z^\alpha
$
satisfy the same equation (\ref{eq:partial-diff}).
Note that the solution of (\ref{eq:partial-diff}) is uniquely determined
if the initial value of $u$ at $z^{(K_J)} = 0$ is given.
It is easy to see that $T^{<N}_{J}(0) = 1$ if $A_J(N) \ne \emptyset$
and $T^{<N}_{J}(0) = 0$ if $A_J(N) = \emptyset$.
Hence we have obtained \eqref{eq:lem-formula-polynomial} for $n$.
This completes the proof.
\end{proof}
\begin{df}{\label{df:total_coefficients}}
Let $S := \Ssec$ with $V \in V(\zeta)$ and $\{\epsilon\}$
of positive real numbers.
We say that $\{F_J\}_{J \in \subsetl}$ is a total family of coefficients of
multi-asymptotic expansion along $\chi$ on $S$ if each $F_J$
consists of a family
$\left\{f_{J,\,\alpha}\right\}_{\alpha \in
\mathbb{Z}^{(K_J)}_{\ge 0}}$
of holomorphic functions on $S_J$.
\end{df}
Let $\Gamma$ be a finite subset of $\mathbb{Z}^{(K_J)}_{\ge 0}$ and
$p = \sum_{\alpha \in \Gamma} c_\alpha z^\alpha$ a polynomial of the
variables $z^{(K_J)}$. We define
$$
\tau_{F,J}(p)(z) := \sum_{\alpha \in \Gamma} c_\alpha
f_{J,\,\alpha}\left(z^{(\compl K_J)}\right) z^\alpha.
$$
Then we set
\begin{equation}
T_{J,\,\beta}(F;\,z) := \tau_{F,J}(T_{J,\,\beta}),
\end{equation}
\begin{equation}
	T^{<N}_J(F;\,z) := \sum_{\beta <_J N} T_{J,\,\beta}(F;\,z),
\end{equation}
and
\begin{equation}
\operatorname{App}^{<N}(F;\,z)
:= \sum_{J \in \subsetl}(-1)^{(\#J + 1)} T_J^{<N}(F;\,z).
\end{equation}

\

Let $f(z)$ be a holomorphic function on $\overline{S}$, and let us define
the family $\{f_J\}_{J \in \subsetl}$ by
\begin{equation}
f_J = \left\{ f_{J,\,\alpha}\left(z^{(\compl K_J)}\right) := \left.\dfrac{\partial^\alpha f}{\partial z^\alpha}\right|_{Z_J}
\right\}_{\alpha \in \mathbb{Z}^{(K_J)}_{\ge 0}}
\end{equation}
for each $J \in \subsetl$.
\begin{lem}{\label{lem:asymptotic_f}}
We have
\begin{equation}
\begin{aligned}
T_J^{<N}(f_J;\, z)
&=
\sum_{\alpha \in A_J(N)}
\dfrac{f_{J,\,\alpha}\left(z^{(\compl K_J)}\right)}{\alpha!} z^\alpha\\
&=
\sum_{\beta \in \mathbb{Z}^{J}_{\ge 0},\, \beta <_J N}
\left.\dfrac{1}{\beta!}
\dfrac{\partial^\beta}{\partial \lambda^\beta}
f(\mu_{\intA}(z,\,\lambda))\right|_{\lambda = \check{e}_{J}}.
\end{aligned}
\end{equation}
\end{lem}
\begin{proof}
The lemma can be proved in the same way as
in the proof of Lemma \ref{eq:formula_development}.
Set
$$
\widetilde{T}_J^{<N}(f;\,z) :=
\sum_{\beta \in \mathbb{Z}^{J}_{\ge 0},\, \beta <_J N}
\left.\dfrac{1}{\beta!}
\dfrac{\partial^\beta}{\partial \lambda^\beta}
f(\mu_{\intA}(z,\,\lambda))\right|_{\lambda = \check{e}_{J}},
$$
where, if the set $\{\beta \in \mathbb{Z}^{J}_{\ge 0};\, \beta <_J N\}$
is empty, we set $\widetilde{T}_J^{<N}(f;\,z) := 0$ as usual.

Let $k \in K_J$ and $z_i$ in $z^{(k)}$.
For $\beta \in \mathbb{Z}^J_{\ge 0}$,
by the same computation as in the proof of Lemma {\ref{eq:formula_development}},
we have
$$
\begin{aligned}
& \dfrac{\partial}{\partial z_i}
\left(
\left.
\dfrac{1}{\beta!}
\partial_\lambda^\beta
f(\mu_{\intA}(z,\,\lambda))\right|_{\lambda = \check{e}_{J}}
\right)\\
&\qquad =
\left.
\dfrac{1}{(\beta -_J \sigma_A A_\chi^{*,k})!}
\left(
\partial_\lambda^{\beta -_J \sigma_A A_\chi^{*,k}}
\dfrac{\partial f}{\partial z_i}(\mu_{\intA}(z,\,\lambda))\right)\right|_{\lambda = \check{e}_{J}},
\end{aligned}
$$
where we set $\partial_\lambda^{\beta -_J \sigma_A A_\chi^{*,k}}g = 0$ if
some entries of $\beta -_J \sigma_A A_\chi^{*,k}$ take negative values,
from which we get
$$
\dfrac{\partial}{\partial z_i} \widetilde{T}_J^{<N}(f;\,z)
= \widetilde{T}_J^{<N-_J \sigma_A A_\chi^{*,k}}
\left(\dfrac{\partial f}{\partial z_i};\,z\right).
$$
By direct computations, we also have
$$
\dfrac{\partial}{\partial z_i} T_J^{<N}(f_J;\,z)
= T_J^{<N-_J \sigma_A A_\chi^{*,k}}
\left( \left(\dfrac{\partial f}{\partial z_i}\right)_J;\,z\right)
$$
and
$$
\left.\widetilde{T}_J^{<N}(f;\,z)\right|_{Z_J} = \left.T_J^{<N}(f_J;\,z)\right|_{Z_J}
=
\left\{
\begin{array}{ll}
f|_{Z_J}
\quad &\text{if $\{\beta \in \mathbb{Z}^{J}_{\ge 0};\, \beta <_J N\} \ne \emptyset$}, \\
\\
0
\quad &\text{if $\{\beta \in \mathbb{Z}^{J}_{\ge 0};\, \beta <_J N\} = \emptyset$}.
\end{array}
\right.
$$

Therefore the claim $T_J^{<N} = \widetilde{T}_J^{<N}$ can be shown by
the induction on the length of $N$.
\end{proof}
\subsection{A family of level functions of multi-asymptotic expansions along $\chi$}\label{subsection:Level functions}

Let
$
(0;\, \zeta) = (0;\, \zeta^{(1)}, \dots, \zeta^{(m)}) \in S_\chi.
$
Hereafter we always assume
\begin{equation}
\text{$p$ is outside fixed points.}
\end{equation}
Hence, by Corollary
\ref{cor:outside-fixed-points-coordinate}, we may assume
\begin{equation}{\label{eq:normalized_condition_asymp}}
\text{$|\zeta^{(k)}| \ne 0$ for $1 \le k \le L$.}
\end{equation}
%Note that we have $J_Z(\zeta) \cap \{1,2,\dots,L\} = \emptyset$.

In this subsection, we introduce level functions which are needed to give an estimate of
the remainder term of an asymptotic expansion. When the associated action $\mu_\chi$ is
non-degenerate, the level functions $\rho_{\Lambda,j}$'s are nothing but
$\varphi_j^{-1}$'s.  Hence the readers who are interested only in
a non-degenerate case may skip this subsection.

Through this subsection, we use the same notations as those introduced in Subsection
\ref{subsection:Semigroups of generators}.
For $f(\tau,\,\lambda) \in \mathcal{H}_{\tau,\lambda}$ with $a := \nu^\lambda_j(f) < 0$,
we set
\begin{equation}
\operatorname{sol}^\lambda_j(f)(\tau,\lambda) :=
\lambda_j\,f^{1/|a|},
\end{equation}
that is, we have $f(\tau,\lambda) = 1$ by putting
$\lambda_j = \operatorname{sol}^\lambda_j(f)$ into $f$.
Note that, by  definition,  $\operatorname{sol}^\lambda_j(f)$ is independent of
the variable $\lambda_j$.  Recall that, for $L < j \le \ell$, we set
$$
F^{0,j} := (\mathcal{L}^\lambda_j \circ \dots \circ \mathcal{L}^\lambda_{L+1})(G)
$$
and, for convenience, we also set
$$
F^{0,L} = G.
$$
Now we define
\begin{equation}
\rho_{j+1}(\tau',\lambda''):=
\underset{(f,\,v) \in F^{0,j}_<}{\max}\,\, \operatorname{sol}^\lambda_{j+1}(f)(\tau',\,\lambda'')
\qquad (L \le j < \ell),
\end{equation}
where we set
$$
F^{0,j}_< := \{(f,v) \in F^{0,j};\, \nu^\lambda_{j+1}(f) < 0\}
\qquad (L \le j < \ell).
$$
Note that, when $F^{0,j}_<$ is an empty set, we set
$\rho_{j+1}(\tau',\,\lambda'') = 1$.
\begin{lem}{\label{lem:non-empty-negative-set}}
Let $L \le j < \ell$.
If $\mu_{j+1}$ is not the identity action on $X$, then $F^{0,j}_<$ is non-empty.
\end{lem}
\begin{proof}
By the assumption, the $(j+1)$-th row is non-zero. Hence there exists $1 \le k \le L$
such that $a_{j+1,k} \ne 0$. In fact, if $a_{j+1,k} = 0$ holds for any $1 \le k \le L$,
since $(j+1)$-th row is non-zero, we have $\operatorname{Rank} A_\chi > L$, which
contradicts $\operatorname{Rank} A_\chi = L$. As
$$
\tau_k = \varphi_k(\varphi^{-1}(\tau', \lambda''),\, \lambda'')
$$
holds, we get, for some $w \in \mathbb{R}_{\ge 0}$ and some integer $N$,
$$
\left( \left(\tau_k/\lambda_{j+1}^{a_{j+1,k}}\right)^N,\,w\right) \in
T:=\{(f,v) \in [G];\, \nu^\lambda_{L+1}(f) = \dots = \nu^\lambda_j(f) = 0 \}.
$$
By the same reasoning as that in the proof of Proposition \ref{prop:radical-of-S},
the radical of $[F^{0,j}]$ coincides with $T$, which implies
$F^{0,j}_< \ne \emptyset$.
\end{proof}
\begin{oss}
As $\psi_k(\tau)$ ($L < k \le m$) does not contain the variables $\lambda$,
for $(f,\,v) \in F^{0,j}_<$,  $f$ is independent of the variables $\tau''$.
Furthermore, since it follows from the definition of $F^{0,j}$
that  $f$ contains the only variables $\tau'$ and
$\lambda_{j+1}$, $\dots$, $\lambda_\ell$,
 $\rho_j$ is a function
of the variables $\tau'$ and $\lambda_{j+1}$, $\dots$, $\lambda_\ell$.
That is,
$$
\begin{aligned}
\rho_{L+1}(\tau',\lambda'') &= \rho_{L+1}(\tau',\,\lambda_{L+2},\,\lambda_{L+3},\,\dots,\,\lambda_\ell), \\
\rho_{L+2}(\tau',\lambda'') &= \rho_{L+2}(\tau',\,\lambda_{L+3},\dots,\,\lambda_\ell), \\
	& \dots \\
\rho_{\ell-1}(\tau',\lambda'') &= \rho_{\ell-1}(\tau',\,\lambda_\ell), \\
\rho_{\ell}(\tau',\lambda'') &= \rho_{\ell}(\tau').
\end{aligned}
$$
\end{oss}
Now we define the closed subset $\Lambda \subset (\RP)^{m+\ell}$ by
\begin{equation}
\left\{(\tau,\,\lambda) \in (\RP)^{m+\ell};\,
\lambda_j = \rho_j(\tau',\,\lambda'') \qquad (L < j \le \ell)
\right\}.
\end{equation}
Let $\pi_{\tau,\lambda'}: (\RP)^{m+\ell} \to (\RP)^{m+L}$ be the canonical projection
by $(\tau,\,\lambda',\,\lambda'') \to (\tau,\lambda')$.
Then, by the above remark,
 $\pi_{\tau,\lambda'}|_\Lambda$ gives an isomorphism between $\Lambda$ and
$(\RP)^{m+L}$. Hence, for any $f \in \mathcal{H}_{\tau,\lambda}$, the restriction
$f|_\Lambda$ can be regarded as a function of the variables $\tau$ and $\lambda'$.
In particular, if $f$ is independent of the variables $\lambda'$, then $f|_\Lambda$
becomes a function of $\tau$.
%Hence, for $(f,\,v) \in F^{0,j}$, since $f$ is independent of
%the variables $\tau''$ and $\lambda'$, the $f|_\Lambda$ can be regarded as a function of $\tau'$ only.
\begin{oss}{\label{rem:successive_composition}}
The inverse of $\pi_{\tau,\lambda'}|_\Lambda$ is concretely
given by the following procedures:
First define
\begin{equation}{\label{eq:lambda_subst_step1}}
\lambda_\ell = \rho_\ell|_\Lambda(\tau') = \rho_\ell(\tau').
\end{equation}
Then we eliminate the variable $\lambda_\ell$ in $\rho_{\ell-1}(\tau', \lambda_\ell)$
by \eqref{eq:lambda_subst_step1} and define
\begin{equation}{\label{eq:lambda_subst_step2}}
\lambda_{\ell-1} =
\rho_{\ell-1}|_\Lambda(\tau') = \rho_{\ell-1}(\tau',\, \rho_\ell(\tau')).
\end{equation}
Again we eliminate the variables $\lambda_{\ell-1}$ and $\lambda_{\ell}$
in $\rho_{\ell-2}(\tau',\, \lambda_{\ell-1},\,\lambda_{\ell})$
by \eqref{eq:lambda_subst_step1} and \eqref{eq:lambda_subst_step2}, we define
$$
\lambda_{\ell-2} =
\rho_{\ell-2}|_\Lambda(\tau') =
\rho_{\ell-2}(\tau',\, \rho_{\ell-1}|_\Lambda(\tau'),\, \rho_{\ell}|_\Lambda(\tau')).
$$
We successively apply the same procedures to $\rho_j$'s, we finally obtain
$$
\lambda_{L+1} =
\rho_{L+1}|_\Lambda(\tau') =
\rho_{L+1}(\tau',\, \rho_{L+2}|_\Lambda(\tau'),\,\dots,\,\rho_{\ell}|_\Lambda(\tau')).
$$
Then, the inverse is given by the correspondence
$$
(\tau,\, \lambda') \to
(\tau,\,\lambda'\,\,\lambda'') = (\tau,\, \lambda',\, \rho_{L+1}|_\Lambda(\tau'),\,\dots,\, \rho_{\ell}|_\Lambda(\tau')).
$$
\end{oss}
In what follows, we also denote by $\rho_{\Lambda,j}$ the restriction of $\rho_j$ to $\Lambda$, i.e.,
$$
\rho_{\Lambda,\,j}(\tau') := \rho_j|_\Lambda(\tau') \qquad (L < j \le \ell).
$$
Set
\begin{equation}
\begin{aligned}
	\rho''(\tau',\lambda'') &:= (\rho_{L+1}(\tau',\lambda''),\,\dots,\,
	\rho_{\ell}(\tau',\lambda'')),\\
	\rho''_\Lambda(\tau') &:= (\rho_{\Lambda,L+1}(\tau'),\,\dots,\,
	\rho_{\Lambda,\ell}(\tau')).
\end{aligned}
\end{equation}
Now we also set
\begin{equation}
\begin{aligned}
	\rho'(\tau',\lambda'') &:= \varphi^{-1}(\tau',\,\lambda'') \\
	\rho'_\Lambda(\tau') &:= \varphi^{-1}|_\Lambda(\tau')
	= \varphi^{-1}(\tau',\,\rho''_\Lambda(\tau')).
\end{aligned}
\end{equation}
We often write
$\varphi^{-1}_{\Lambda,j}$ instead of $\varphi_j^{-1}|_\Lambda$ hereafter.
Finally we define
$$
\rho = (\rho',\,\rho'')\quad \text{and} \quad
\rho_\Lambda = (\rho'_\Lambda,\,\rho''_\Lambda).
$$
%Note that, for $f(\tau', \lambda'') \in F^{0,j}$, its restriction
%$f|_\Lambda(\tau')$ is realized by the substitution
%$\lambda''= \rho''_\Lambda(\tau')$, i.e., it is
%given by $f(\tau',\, \rho''_\Lambda(\tau'))$.
We here emphasize that $\rho_\Lambda$ are functions
of the variables $\tau'$, i.e., they are independent of
the variables $\tau''$.

\begin{oss}
When $\mu_\chi$ is non-degenerate, since $\ell = L$ holds, in particular,
no $\lambda''$ and $\rho''$ appear in this case,
we have
$$
\rho = \rho_\Lambda = \varphi^{-1}(\tau').
$$
\end{oss}
\begin{lem}{\label{lem:inverse_rho}}
We have
$$
\tau_k = \varphi_k(\rho_{\Lambda}(\tau')) \qquad (1 \le k \le L).
$$
\end{lem}
\begin{proof}
As we have
$$
\varphi^{-1}_j|_\Lambda = \rho_{\Lambda,j}\,\,\, (1 \le j \le L),\quad
\lambda_j|_\Lambda = \rho_{\Lambda,j}\,\,\, (L < j \le \ell),
$$
the equality is obtained by the restriction of the following usual one
to the subset $\Lambda$:
$$
\tau_k = \varphi_k(\varphi^{-1}, \lambda'') \qquad (1 \le k \le L).
$$
\end{proof}

\begin{es}
	Let us consider $X$ with coordinates blocks
$(z^{(0)},\,z^{(1)},\, z^{(2)},\, z^{(3)})$ and
$3$-closed submanifolds
$$
Z_1 = \{z^{(1)} =z^{(3)} = 0\}, \quad
Z_2 = \{z^{(2)} = z^{(3)} = 0\}, \quad
Z_3 = \{z^{(3)} = 0\}.
$$
In this case, we have $m=\ell=3$. Let $A_\chi$ be the matrix associated with the family
$\chi = \{Z_1,\,Z_2,\,Z_3\}$. Clearly $L = \operatorname{Rank} A_\chi = 3$, and thus,
$\chi$ is of normal type.

Then
$$
\varphi_1(\lambda) = \lambda_1, \quad
\varphi_2(\lambda) = \lambda_2, \quad
\varphi_3(\lambda) = \lambda_1 \lambda_2 \lambda_3,
$$
and
$$
\varphi^{-1}_1(\tau) = \tau_1, \quad
\varphi^{-1}_2(\tau) = \tau_2, \quad
\varphi^{-1}_3(\tau) = \tau_3/(\tau_1 \tau_2).
$$
Because of $\ell = L = 3$, we get
$$
F^q = G = \{(\varphi_1^{-1},\,0),\,
(\varphi_2^{-1},\,0),\,
(\varphi_3^{-1},\,0) \}.
$$
As $\chi$ is of normal type, by the above remark, we have
$$
\rho_1(\tau) = \tau_1, \quad
\rho_2(\tau) = \tau_2, \quad
\rho_3(\tau) = \tau_3/(\tau_1 \tau_2).
$$
\end{es}
\begin{es}{\label{ex:transitive_non_degenerate}}
Let us consider $X$ with coordinates blocks
$(z^{(0)},\,z^{(1)},\, z^{(2)},\, z^{(3)})$ and
$4$-closed submanifolds
$$
\begin{aligned}
&Z_1 = \{z^{(1)} =z^{(3)} = 0\}, \quad
Z_2 = \{z^{(2)} = z^{(3)} = 0\}, \quad
Z_3 = \{z^{(3)} = 0\}, \\
&Z_4 = \{z^{(1)}=z^{(2)}=z^{(3)} = 0\}.
\end{aligned}
$$
In this case,
we have $m=3$ and $\ell=4$. Let $A_\chi$ be the $4 \times 3$ matrix
associated with the family $\chi = \{Z_1,\,Z_2,\,Z_3,\,Z_4\}$.
Clearly $L = \operatorname{Rank} A_\chi = 3$, and thus,
$\chi$ is of transitive type.
Then
$$
\varphi_1(\lambda) = \lambda_1 \lambda_4, \quad
\varphi_2(\lambda) = \lambda_2 \lambda_4, \quad
\varphi_3(\lambda) = \lambda_1 \lambda_2 \lambda_3 \lambda_4,
$$
and
$$
\varphi^{-1}_1(\tau,\,\lambda_4) = \tau_1/\lambda_4, \quad
\varphi^{-1}_2(\tau,\,\lambda_4) = \tau_2/\lambda_4, \quad
\varphi^{-1}_3(\tau,\,\lambda_4) = \tau_3 \lambda_4/(\tau_1 \tau_2).
$$
By the definition, we have
$$
F^{0,3} = G =
\{(\varphi_1^{-1},\,0),\,
(\varphi_2^{-1},\,0),\,
(\varphi_3^{-1},\,0),\,
(\lambda_4,\,0)\}.
$$
By applying the operation $\mathcal{L}^\lambda_4$ to $F^{0,3}$, we obtain
$$
F^q = F^{0,4} =
\{(\tau_1,\,0),\,
(\tau_2,\,0),\,
(\tau_3/\tau_2,\,0),\,
(\tau_3/\tau_1,\,0)\}.
$$
Since we easily see
$$
F^{0,3}_< = \{(\tau_1/\lambda_4,\,0),\, (\tau_2/\lambda_4,\,0)\},
$$
we get
$$
\rho_4 = \max\,\{\tau_1, \tau_2\}
$$
and
$$
\Lambda = \{ \lambda_4 = \rho_4(\tau) \}.
$$
Hence we finally obtain
$$
\begin{aligned}
\rho_{\Lambda,1}(\tau) &= \tau_1/\max\{\tau_1,\,\tau_2\}, \quad
\rho_{\Lambda,2}(\tau) = \tau_2/\max\{\tau_1,\,\tau_2\},\\
\rho_{\Lambda,3}(\tau) &= \tau_3 \max\{\tau_1,\,\tau_2\} /(\tau_1 \tau_2), \quad
\rho_{\Lambda,4}(\tau) = \max\{\tau_1,\,\tau_2\}.
\end{aligned}
$$
\end{es}
\begin{es}{\label{ex:transitive_non-degenerate}}
Let us consider $X$ with coordinates blocks
$(z^{(0)},\,z^{(1)},\, z^{(2)},\, z^{(3)})$ and
$5$-closed submanifolds
$$
\begin{aligned}
&Z_1 = \{z^{(1)} = 0\}, \quad
Z_2 = \{z^{(2)} = 0\}, \quad
Z_3 = \{z^{(3)} = 0\}, \\
&Z_4 = \{z^{(1)}=z^{(3)}= 0\}, \quad
 Z_5 = \{z^{(2)}=z^{(3)} = 0\}.
\end{aligned}
$$
In this case,
we have $m=3$ and $\ell=5$. Let $A_\chi$ be the $5 \times 3$ matrix
associated with the family $\chi = \{Z_1,\,Z_2,\,Z_3,\,Z_4,\,Z_5\}$.
Clearly $L = \operatorname{Rank} A_\chi = 3$, and thus,
$\chi$ is of transitive type.
Then
$$
\varphi_1(\lambda) = \lambda_1 \lambda_4, \quad
\varphi_2(\lambda) = \lambda_2 \lambda_5, \quad
\varphi_3(\lambda) = \lambda_3 \lambda_4 \lambda_5,
$$
and
$$
\varphi^{-1}_1(\tau,\,\lambda'') = \tau_1 /\lambda_4, \quad
\varphi^{-1}_2(\tau,\,\lambda'') = \tau_2/\lambda_5, \quad
\varphi^{-1}_3(\tau,\,\lambda'') = \tau_3/\lambda_4\lambda_5,
$$
where $\lambda'' = (\lambda_4,\,\lambda_5)$.
By the definition, we have
$$
F^{0,3} = G =
\{(\varphi_1^{-1},\,0),\,
(\varphi_2^{-1},\,0),\,
(\varphi_3^{-1},\,0),\,
(\lambda_4,\,0),\,
(\lambda_5,\,0)\}.
$$
By applying the operation $\mathcal{L}^\lambda_4$ to $F^{0,3}$, we obtain
$$
F^{0,4} =
\{(\tau_1,\,0),\,
(\tau_2/\lambda_5,\,0),\,
(\tau_3/\lambda_5,\,0),\,
(\lambda_5,\,0)\},
$$
and, by applying the operation $\mathcal{L}^\lambda_5$ to $F^{0,4}$, we get
$$
F^q = F^{0,5} =
\{(\tau_1,\,0),\,
(\tau_2,\,0),\,
(\tau_3,\,0)\}.
$$
Since we see
$$
F^{0,3}_< = \{(\tau_1/\lambda_4,\,0),\,
(\tau_3/ \lambda_4\lambda_5,\,0)\},
$$
and thus,
$$
\rho_4(\tau,\,\lambda_5) = \max\,\{\tau_1,\, \tau_3/\lambda_5\}.
$$
In the same way, we have
$$
F^{0,4}_< = \{(\tau_2/\lambda_5,\,0), (\tau_3/\lambda_5,\,0)\},
\,\,\text{ and thus, }\quad
\rho_5(\tau) = \max\,\{\tau_2,\,\tau_3\}.
$$
Hence we have
$$
\Lambda = \{ \lambda_4 = \rho_4(\tau, \lambda_5),\, \lambda_5 = \rho_5(\tau) \}
$$
and
$$
\begin{aligned}
\rho_{\Lambda,1}(\tau) &= \tau_1/\max\,\{\tau_1,\,\tau_3/\max\,\{\tau_2,\,\tau_3\}\},
\quad
\rho_{\Lambda,2}(\tau) = \tau_2/\max\,\{\tau_2,\, \tau_3\},\\
\rho_{\Lambda,3}(\tau) &= \tau_3/\max\,\{\tau_1 \max\,\{\tau_2,\,\tau_3\},\, \tau_3\}, \\
\rho_{\Lambda,4}(\tau) &= \max\,\{\tau_1,\, \tau_3/\max\,\{\tau_2,\,\tau_3\}\},\\
\rho_{\Lambda,5}(\tau) &= \max\,\{\tau_2,\,\tau_3\}.
\end{aligned}
$$
\end{es}
For families of the level functions obtained so far,
each action $\mu_j$ becomes strict
in the sense of Subsection \ref{subsection:Uniqueness of coefficients}.
However the following family gives a non-strict action.
\begin{es}{\label{ex:transitive_degenerate}}
Let us consider $X$ with coordinates blocks
$(z^{(0)},\,z^{(1)},\, z^{(2)},\, z^{(3)})$ and
$5$-closed submanifolds
$$
\begin{aligned}
&Z_1 = \{z^{(1)} = 0\}, \quad
Z_2 = \{z^{(2)} = 0\}, \quad
Z_3 = \{z^{(1)} = z^{(2)} = z^{(3)} = 0\}, \\
&Z_4 = \{z^{(1)}=z^{(2)}= 0\}, \quad
 Z_5 = \{z^{(2)}=z^{(3)} = 0\}.
\end{aligned}
$$
In this case,
we have $m=3$ and $\ell=5$. Let $A_\chi$ be the $5 \times 3$ matrix
associated with the family $\chi = \{Z_1,\,Z_2,\,Z_3,\,Z_4,\,Z_5\}$.
Clearly $L = \operatorname{Rank} A_\chi = 3$, and thus,
$\chi$ is of transitive type.
Then
$$
\varphi_1(\lambda) = \lambda_1 \lambda_3 \lambda_4, \quad
\varphi_2(\lambda) = \lambda_2 \lambda_3 \lambda_4 \lambda_5, \quad
\varphi_3(\lambda) = \lambda_3 \lambda_5,
$$
and
$$
\varphi^{-1}_1(\tau,\,\lambda'') = \tau_1 \lambda_5/(\tau_3 \lambda_4), \quad
\varphi^{-1}_2(\tau,\,\lambda'') = \tau_2/(\tau_3 \lambda_4), \quad
\varphi^{-1}_3(\tau,\,\lambda'') = \tau_3/\lambda_5,
$$
where $\lambda'' = (\lambda_4,\,\lambda_5)$.
By the definition, we have
$$
F^{0,3} = G =
\{(\varphi_1^{-1},\,0),\,
(\varphi_2^{-1},\,0),\,
(\varphi_3^{-1},\,0),\,
(\lambda_4,\,0),\,
(\lambda_5,\,0)\}.
$$
By applying the operation $\mathcal{L}^\lambda_4$ to $F^{0,3}$, we obtain
$$
F^{0,4} =
\{(\tau_2/\tau_3,\,0),\,
(\tau_1 \lambda_5/\tau_3,\,0),\,
(\tau_3/\lambda_5,\,0),\,
(\lambda_5,\,0)\},
$$
and, by applying the operation $\mathcal{L}^\lambda_5$ to $F^{0,4}$, we get
$$
F^q = F^{0,5} =
\{(\tau_2/\tau_3,\,0),\,
(\tau_3,\,0),\,
(\tau_1,\,0)\}.
$$
Since we see
$$
F^{0,3}_< = \{(\tau_2/(\tau_3\lambda_4),\,0),\,
(\tau_1 \lambda_5/(\tau_3 \lambda_4),\,0)\},
$$
and thus,
$$
\rho_4(\tau,\,\lambda_5) = \max\,\{\tau_2/\tau_3,\, \tau_1 \lambda_5/\tau_3\}.
$$
In the same way, we have
$$
F^{0,4}_< = \{(\tau_3/\lambda_5,\,0)\},
\,\,\text{ and thus, }\quad
\rho_5(\tau) = \tau_3.
$$
Hence we have
$$
\Lambda = \{ \lambda_4 = \rho_4(\tau, \lambda_5),\, \lambda_5 = \rho_5(\tau) \}
$$
and
$$
\begin{aligned}
\rho_{\Lambda,1}(\tau) &= 1/\max\{1,\,\tau_2/(\tau_1 \tau_3)\}, \quad
\rho_{\Lambda,2}(\tau) = 1/\max\{1,\,\tau_1 \tau_3/\tau_2\},\\
\rho_{\Lambda,3}(\tau) &= 1, \quad
\rho_{\Lambda,4}(\tau) = \max\{\tau_1,\,\tau_2/\tau_3\},\quad
\rho_{\Lambda,5}(\tau) = \tau_3.
\end{aligned}
$$
\end{es}

The following proposition plays a key role in the theory of multi-asymptotic expansions.
\begin{prop}{\label{prop:prop_basic_estimate}}
Let $S$ be a subset in $(\RP)^m$ with coordinates $\tau = (\tau', \tau'')$.
Assume that, for any $(f,\,v) \in F^q=F^{0,\ell}$,
 $f(\tau)$ is bounded on $S$. Then we have:
\begin{enumerate}
	\item For any $(f,\,v) \in G = F^{0,L}$,
 $f|_\Lambda(\tau)$ is bounded on
$S$. That is, there exists a constant $C>0$ for which we have
\begin{equation}\label{eq:level-cond-2.1}
	\rho_{\Lambda,j}(\tau') = \varphi^{-1}_{\Lambda,j}(\tau') \le C
	\quad (1 \le j \le L,\, \tau \in S),
\end{equation}
\begin{equation}{\label{eq:level-cond-1.1}}
\begin{array}{ll}
	\psi_k(\tau) \le C
	&\quad (L < k \le m,\, \tau \in S), \\
	C^{-1} \le \psi_k(\tau)
	&\quad (L < k \le m,\,k \notin J_Z(\zeta),\, \tau \in S)
\end{array}
\end{equation}
and
\begin{equation}{\label{eq:level-cond-2.2}}
\rho_{\Lambda,j}(\tau') = \lambda_j|_\Lambda \le C
\qquad \quad (L < j \le \ell,\, \tau \in S).
\end{equation}

\item
There exist $C > 0$ and $N > 0$ such that
\begin{equation}{\label{eq:level-cond-1.2}}
\begin{aligned}
&C^{-1}t^{-N} \rho_{\Lambda,j}(\tau')
\le
\rho_{\Lambda,j}\left(\tau_1, \dots,
\overset{\text{$k$-th}}{t \tau_k}, \dots, \tau_L\right)
\le
Ct^{N} \rho_{\Lambda,j}(\tau') \\
&\qquad\qquad\qquad\qquad\qquad\qquad
(t \ge 1,\,\, \tau' \in (\RP)^L)
\end{aligned}
\end{equation}
hold for any $1 \le j \le \ell$ and $1 \le k \le L$.
\end{enumerate}
\end{prop}
\begin{proof}
We may assume that all the actions $\mu_j$ ($j=1,2,\dots, \ell$)
are different from the identity action on $X$.
We prove the following claim by induction on $j = \ell,\dots,L$:
\begin{center}
For any $(f,\,v) \in F^{0,j}$,  $f|_\Lambda(\tau)$ is bounded
on $S$.
\end{center}
When $j = \ell$, the claim follows from the assumption.

Now assume that the claim is true for $j+1$ and we will prove the claim for $j$.
We consider the following $3$ cases:
Let $(f,\,v) \in F^{0,j}$ with $\nu^\lambda_{j+1}(f) = 0$. In this case,
$(f,\,v)$ belongs to $F^{0,j+1}$ also, and thus,
$f|_\Lambda(\tau)$ is bounded on $S$ by the induction hypothesis.

Let $(f,\,v) \in F^{0,j}$ with $\nu^\lambda_{j+1}(f) < 0$. By the definition of
$\rho_{j+1}$, we have
$$
\operatorname{sol}^\lambda_{j+1}(f)|_\Lambda(\tau') \le
\rho_{\Lambda,j+1}(\tau').
$$
By noticing $\nu^\lambda_{j+1}(f) < 0$, we get
$$
\begin{aligned}
f|_\Lambda(\tau')
&= f(\tau',\, \rho_{\Lambda,j+1}(\tau'),\, \rho_{\Lambda,j+2}(\tau'),
\,\dots,\,\rho_{\Lambda,\ell}(\tau')) \\
&\le f(\tau',\, \operatorname{sol}^\lambda_{j+1}(f)|_\Lambda(\tau'),\,
\rho_{\Lambda,j+2}(\tau'),\,
\dots,\,\rho_{\Lambda,\ell}(\tau')) \\
& = f(\tau',\, \operatorname{sol}^\lambda_{j+1}(f),\,
\lambda_{j+2},\, \dots,\,\lambda_\ell)\big|_\Lambda
= 1.
\end{aligned}
$$
Hence we have obtained the boundness of $f|_\Lambda(\tau')$.

Let $(f,\,v) \in F^{0,j}$ with $\nu^\lambda_{j+1}(f) > 0$. Since
$F^{0,j}_<$ is non-empty by Lemma \ref{lem:non-empty-negative-set},
for each fixed $\tau \in S$,
we can find $(g_\tau,v_\tau) \in F^{0,j}$ with $\nu^\lambda_{j+1}(g_\tau) < 0$ such that
\begin{equation}{\label{eq:lambda_identity_one}}
\operatorname{sol}^\lambda_{j+1}(g_\tau)|_\Lambda(\tau') =
\rho_{\Lambda,j+1}(\tau').
\end{equation}
For such a $g_\tau$, we see, by F3 of the operation $\mathcal{L}^\lambda_{j+1}$, we have
$$
h := g_\tau^{\nu^\lambda_{j+1}(f)} f^{|\nu^\lambda_{j+1}(g_\tau)|} \in F^{0,j+1}.
$$
Since $h|_\Lambda(\tau')$ is bounded on $S$ by the induction hypothesis and
since $g_\tau|_\Lambda(\tau') = 1$ holds by \eqref{eq:lambda_identity_one},
we obtain
$$
h|_\Lambda(\tau') = (f|_\Lambda(\tau'))^{|\nu^\lambda_{j+1}(g_\tau)|},
$$
which implies the boundness of $f|_\Lambda(\tau')$ on $S$.

Hence we can show the claim for $j$, and thus, the claim is true for $j = L$ which gives
the claim 1.~of the proposition.

The claim 2.~of the proposition easily follows from the fact that $\rho_{\Lambda,j}(\tau)$ is
obtained by successive compositions of functions defined by maximal of
several monomials as shown in Remark \ref{rem:successive_composition}.
\end{proof}

\subsection{Definition of multi-asymptotic expansions along $\chi$ and
their basic properties}\label{subsection:Definition of multi-asymptotic}

We first define the notion of multi-asymptotical developability on a multisector $S$
along $\chi$ using level functions $\rho_{\Lambda,j}(\tau')$ $(j=1,2,\dots,\ell)$
introduced in the previous subsection.
As in the previous subsection, we continue to assume that
\begin{equation}
\text{$p = (z^{(0)};\,\zeta) \in S_\chi$ is outside the fixed points.}
\end{equation}
Note that, through the rest of subsetctions in Section
\ref{section:Multi-asymptotic expansions}, we continue to assume the point $p$ to be
outside the fixed points.
\begin{oss}
We here give a short summary for a major case where
the associated action $\mu_\chi$ is non-degenerate, which is the most important and
interesting case: Let $\ell$ be the number of manifolds in $\chi$, and let $m+1$ be the
number of coordinates blocks $(z^{(0)},\,z^{(1)},\, \dots,\,z^{(m)})$.

Then, by the non-degenerate condition, we have $L=\operatorname{Rank} A_\chi = \ell \le m$ and
$$
F^q = F^{0,\ell} = G,
$$
where $G$ is the finite subset of $\RRTpair$
defined in Remark \ref{oss:general_G} of
Subsection \ref{subsection:Semigroups of generators}.

The level functions of this case are given by
$$
\rho_{\Lambda,j}(\tau) = \varphi^{-1}_j(\tau') \qquad (1 \le j \le \ell),
$$
for which Proposition \ref{prop:prop_basic_estimate} in the previous subsection
clearly holds because of $F^q=G$.
%Note that $\tau = \tau'$ and no variables $\tau''$ appear in this case.
\end{oss}
Let $H$ be a finite subset in $\RRpair$ such that
$[H]$ and $\semiG$ defined in \eqref{eq:def_S} are equivalent, and
let $V \in V(\zeta)$ and $\{\epsilon\}$ a set of positive real numbers.
Let $f$ be a holomorphic function on $S:= \Ssec$ and $\{F_J\}_{J\in\subsetl}$ a
total family of coefficients of multi-asymptotic expansion along $\chi$ on $S$.
Here $\Ssec$ was defined in \eqref{def:multi-cone-complex} and
the definition of a total family of coefficients was given
in Definition {\ref{df:total_coefficients}}.
\begin{oss}
In what follows, to make the notation light, we set
$$
\rho_{\Lambda,j}(z) := \rho_{\Lambda,j}(|z^{(*)}|')
\,\,\text{ and }\,\,
\rho_{\Lambda}(z) := \rho_{\Lambda}(|z^{(*)}|')
$$
for $z \in \Ssec$.
\end{oss}

Recall that $\comD$ is a positive rational number so that
$a_{jk} \in \mathbb{Z}$ and all the $a_{jk}$ have no common divisors.
\begin{df}{\label{def:multi-asymptotics}}
We say that $f$ is multi-asymptotically developable to $F=\{F_J\}$ along
$\chi$ on $S$ if and only if for any cone
$S'=S_H(V',\{\epsilon'\})$ properly contained in $S$ and for any $N=(n_1,\dots,n_\ell)
\in \mathbb{Z}^\ell_{\ge 0}$, there exists a constant $C_{S',N}$ such that
\begin{equation}
\left|
f(z) - \operatorname{App}^{<N}(F;\, z)
\right|
\le C_{S',N}
\underset{1 \le j \le \ell}{\prod}
\rho_{\Lambda,j}(z)^{n_j/\comD} \qquad (z \in S').
\end{equation}
\end{df}
Then the most fundamental result is:
\begin{prop}{\label{prop:stability-differentiale}}
If $f(z)$ is multi-asymptotically developable along $\chi$ on $S$, then
any derivative of $f$ is also multi-asymptotically developable along $\chi$ on $S$.
\end{prop}
For $z_i \in z^{(0)}$, it easily follows from the definition that
$\dfrac{\partial}{\partial z_i}f$ becomes multi-asymptotically developable.

Let $k \in \{1, \dots,m\}$ and $z_i$ in $z^{(k)}$. It suffices to show
that $\dfrac{\partial}{\partial z_i}f$ is also multi-asymptotically developable
along $\chi$ on $S$. Assume that $f$ is multi-asymptotically developable to
the total family $F=\{F_J\}_{J \in \subsetl}$ of coefficients
%of multi-asymptotic expansion
with $F_J = \{f_{J,\alpha}\}_{\alpha \in \mathbb{Z}^{(K_J)}_{\ge 0}}$.
The following lemma can be easily proved in the same way as that of Lemma 7.12
in \cite{HP13}.
\begin{lem}
Let $N=(n_1,\dots,n_\ell) \in \mathbb{Z}^\ell_{\ge 0}$.
Set $N^+= N + \comD\, A_\chi^{*,k}$ where $A_\chi^{*,k}$ denotes
the $k$-th column of the matrix $A_\chi$
and define the total family $F' = \{F'_J\}$ by
$$
F'_J := \left\{\dfrac{\partial f_{J,\alpha}}{\partial z_i}\right\}_{\alpha \in
\mathbb{Z}^{(K_J)}_{\ge 0}}
\quad \text{for $J \in \subsetl$ with $k \notin K_J$}
$$
and
$$
F'_J := \left\{f_{J,\alpha + e_i}\right\}_{\alpha \in \mathbb{Z}^{(K_J)}_{\ge 0}}
\quad \text{for $J \in \subsetl$ with $k \in K_J$},
$$
where $e_i$ is the unit vector in $\mathbb{Z}^n$ with the $i$-th element being $1$.
Then we have
$$
\dfrac{\partial}{\partial z_i} \operatorname{App}^{<N^+}(F;\,z) =
\operatorname{App}^{<N}(F';\,z).
$$
\end{lem}
\begin{proof}
Set $N^+ = (n'_1,\dots,n'_\ell)$.
For $J$ with $k \notin K_J$, as the polynomial $T_J^{<N}$ does not contain
the variable $z_i$ and $n'_j = n_j$ for $j \in J$
by the condition \eqref{eq:cond-a2-alt}, we have
$$
\frac{\partial}{\partial {z_i}} T^{<N_+}_J(F;\, z)
=
\frac{\partial}{\partial {z_i}} \tau_{F, J}\left(T^{<N_+}_J\right)
=
\tau_{F', J}\left(T^{<N_+}_J\right)
=
T^{<N}_J(F';\, z).
$$

Let us consider the case $k \in K_J$.
Then, by the equality given at the beginning of the proof for Lemma \ref{eq:formula_development}, we have
$$
\frac{\partial}{\partial {z_i}} T^{<N_+}_J = T^{<N}_J.
$$
For a polynomial $p(z^{(K_J)})$ of the variables $z^{(K_J)}$, we have
$$
\displaystyle\frac{\partial}{\partial z_i}\tau_{F,J}(p) =
\tau_{F',J}\left(\displaystyle\frac{\partial p}{\partial z_i}\right).
$$
Hence we have obtained
$$
\begin{aligned}
\frac{\partial}{\partial {z_i}} T^{<N_+}_J(F;\, z) &=
\frac{\partial}{\partial {z_i}} \tau_{F, J}(T^{<N_+}_J) =
\tau_{F', J}\left(\frac{\partial}{\partial {z_i}} T^{<N_+}_J\right) \\
&= \tau_{F', J}\left(T^{<N}_J\right)
= T^{<N}_J(F';\, z).
\end{aligned}
$$
This completes the proof.
\end{proof}
We also need the following lemma.
\begin{lem}
For $(f,\,v) \in F^q = F^{0,\ell}$,  $f(|z^{(*)}|)$ is bounded on $\Ssec$.
\end{lem}
\begin{proof}
This is a consequence of the fact that $[H]$, $[F^q]$ and $\semiG$
are all equivalent.
\end{proof}

Now we will prove the proposition.
\begin{proof}
Take $S''=S_H(V'',\,\{\epsilon''\})$
so that $S \supset\supset S'' \supset\supset S'$ where
$S \supset\supset S''$ means that $S''$ is properly contained in $S$.
By the lemma, we have
$$
\left| \dfrac{\partial}{\partial z_i} f - \operatorname{App}^{<N}(F';\,z)\right|
=
\left|
\dfrac{\partial}{\partial z_i} \left(f - \operatorname{App}^{<N^+}(F;\,z)\right)\right|.
$$
On $S''$, we get
$$
\left|g := f - \operatorname{App}^{<N^+}(F;\,z)\right|
\le C \underset{1\le j \le \ell}{\prod} \rho_{\Lambda,j}(z)^{n_j/\comD + a_{jk}}.
$$
%We first note that, by the condition \eqref{eq:level-cond-1.1}, the
%\eqref{eq:level-cond-1.2} always holds for not only $r \in J_Z(\zeta)$ but also
%any $r \in \{1,2,\dots,m\}$.
We first note that the condition
$\eqref{eq:level-cond-1.2}$ is equivalently saying that existence of
$N \in \mathbb{N}$ and $C>0$ such that,
for $1 \le j \le \ell$ and $1 \le r \le m$,
\begin{equation}{\label{eq:level-cond-1.2_alt}}
\begin{aligned}
&Ct^{N} \rho_{\Lambda,j}(\tau')
\le
\rho_j(\tau_1, \dots, \tau_{r-1}, t \tau_r, \tau_{r+1},\dots \tau_L)
\le
C^{-1}t^{-N} \rho_{\Lambda,j}(\tau') \\
&\qquad \qquad (0 < t \le 1,\,\, \tau' \in (\RP)^L).
\end{aligned}
\end{equation}
%that is, the similar estimate holds for $0 < \lambda \le 1$ also.

Let us evaluate the remainder term $g$ on $S'$.
%We first note that we may consider the problem under the assumptions
%$z^{(k)} \ne 0$ ($k \in J_Z(\zeta)$)
%because the estimate
%at any point in $S'$ follows by letting $z^{(k)} \to 0$ ($k \in J_Z(\zeta)$).
First consider the case $k \notin J_Z(\zeta)$.
Then we can find a sufficiently small $\epsilon > 0$ such that
$$
B(z) := z +
\left\{w \in \mathbb{C}^n;\, |w_i| \le \epsilon |z^{(k)}|, \,
w_r = 0\,\,(r \ne i)\right\}
\subset S''
$$
holds for any $z \in S'$, which is a consequence of
the following easy lemma.
\begin{lem}
Let $Y=\mathbb{C}^m$. Set, for $\delta > 0$,
$$
D(\delta) := \{z \in Y;\, |z| < \delta\},
$$
where $|z|$ denotes $\max\,\{|z_1|,\dots,|z_m|\}$ for $z = (z_1,\dots,z_m) \in Y$.
\begin{enumerate}
\item We have, for $z \in Y$ and $w \in z + D(\epsilon|z|)$,
\begin{equation}
(1 - \epsilon)|z| < |w| < (1+\epsilon) |z|.
\end{equation}
\item Let $V' \subset V$ be open cones in $Y$ with
$\overline{V'}\setminus\{0\} \subset V$. Then there exists $\kappa > 0$
such that
\begin{equation}
z + D(\kappa|z|) \subset V\qquad (z \in V').
\end{equation}
\end{enumerate}

\end{lem}
By noticing
$$
B(z) \subset
\left\{w \in \mathbb{C}^n;\,
(1-\epsilon)|z^{(k)}| \le |w^{(k)}| \le (1+\epsilon)|z^{(k)}|,\,
w^{(j)}=z^{(j)}\,\,(j \ne k)\right\},
$$
it follows from
\eqref{eq:level-cond-1.2} and
\eqref{eq:level-cond-1.2_alt}
that there exists a constant $C> 0$ satisfying that, for any $z \in S'$,
$$
\underset{1\le j \le \ell}{\prod} \rho_{\Lambda,j}(w)^{n_j/\comD + a_{jk}}
\le C
\underset{1\le j \le \ell}{\prod} \rho_{\Lambda,j}(z)^{n_j/\comD + a_{jk}}
\qquad (w \in B(z)).
$$
Let $\gamma$ be the circle in $\mathbb{C}$ of center $z_i$ and radius $\epsilon |z^{(k)}|$
with the anti-clockwise direction.
Then we obtain
$$
\begin{aligned}
\left|\dfrac{\partial g}{\partial z_i}(z) \right|
&=
\left|\dfrac{1}{2\pi\sqrt{-1}}
\int_\gamma \dfrac{g(z_1,\dots,z_{i-1}, t, z_{i+1}, \dots, z_n)}{(t - z_i)^2} dt \right| \\
&\le
\dfrac{C}{\epsilon|z^{(k)}|}
\underset{1\le j \le \ell}{\prod} \rho_{\Lambda,j}(z)^{n_j/\comD + a_{jk}}
\qquad(z \in S').
%\\
%%
%&\le
%\dfrac{CC_2}{\epsilon \rho_k(z)}
%\underset{1\le j \le \ell}{\prod} \left(\varphi_j^{-1}(\rho(z))\right)^{n_j + a_{jk}}
%\quad (z \in S').
\end{aligned}
$$
If $1 \le k \le L$, then we have, by Lemma \ref{lem:inverse_rho},
$$
|z^{(k)}| =
\varphi_k(\rho_{\Lambda}(z)) =
\underset{1\le j \le \ell}{\prod} \rho_{\Lambda,j}(z)^{a_{jk}}.
$$
If $k > L$, then it follows from \eqref{eq:level-cond-1.1} that we have
$$
C |z^{(k)}| \ge
\varphi_k(\rho_{\Lambda}(z)) =
\underset{1\le j \le \ell}{\prod} \rho_{\Lambda,j}(z)^{a_{jk}}.
$$
Hence we have obtained, for some $C'>0$,
$$
\left|\dfrac{\partial g}{\partial z_i}(z) \right|
\le
C'\underset{1\le j \le \ell}{\prod} \rho_{\Lambda,j}(z)^{n_j/\comD}\quad
(z \in S'),
$$
which concludes the claim of the proposition when $k \notin J_Z(\zeta)$.

Now let us consider the case $k \in J_Z(\zeta)$.
Since $p$ is outside fixed points, we have $k > L$. Note that
$\rho_{\Lambda}$ does not depend on the variables $\tau''$
and only $\psi_k(\tau)$ contains the variable $\tau_k$.
Set for $\epsilon > 0$ and $z \in S'$,
$$
B(z) :=
\left\{w \in \mathbb{C}^n;\, |w_i| = |z_i| +  \epsilon h(z),\,\,
w_r = z_r\,\,(r \ne i)\right\},
$$
where $h(z) :=  \varphi_k(\rho_{\Lambda}(z))$.
Note that $h(z)$ is independent of the variables $z^{(k)}$.
\begin{lem}
If $\epsilon > 0$ is sufficiently small, then we have $B(z) \subset S''$ for any
$z \in S'$.
\end{lem}
\begin{proof}
Since $\psi_k$ does not contain the variables $\lambda$, it belongs to $F^q$.
Note that, in $F^q$,  only $\psi_k$ contains the variable $\tau_k$.
Note also that, as $k \in J_Z(\zeta)$, there is no $(f,v) \in H$ with $\nu_k(f) < 0$.
Let $(f,\,v) \in H$ with $\alpha := \nu_k(f) > 0$, and thus, we can set
$f = \tau_k^\alpha/g(\tau)$ with $\nu_k(g) = 0$.
Then, as the radical of $[F^q]$ is $\semiG$, there exist
$\beta, \gamma \in \mathbb{Q}$ and $r \in [F^q]$ with $\nu_k(r) = 0$ such that
$$
f(\tau) = \tau_k^\alpha/g(\tau) = \psi_k(\tau)^\beta r(\tau)^\gamma.
$$
Then, by comparing the order of the variable $\tau_k$ in both sides, we have $\alpha = \beta$, and
hence, we obtain
$$
r(\tau)^{\gamma/\alpha} g(\tau)^{1/\alpha} =
\varphi_k(\varphi^{-1}(\tau', \lambda''), \lambda'').
$$
Since $r(|z^{(*)}|)$ is bounded on $S$, by restricting the above equality to $\Lambda$,
we can find $C > 0$ such that
$$
C g(|z^{(*)}|)^{1/\alpha} \ge \varphi_k(\rho_{\Lambda}(z)) = h(z) \qquad (z \in S).
$$
By noticing that $f(|z^{(*)}|) < \epsilon_{f,+}$
determines the region in the $z_i$-space
$$
\left\{z_i \in \mathbb{C};\,
|z_i| < \left(\epsilon_{f,+}\right)^{1/\alpha}
g(|z^{(*)}|)^{1/\alpha}\right\},
$$
we conclude $B(z) \subset S''$ for any $z \in S'$
by taking $\epsilon$ sufficiently small.
\end{proof}

We continue the proof.
As $\rho_{\Lambda}(\tau)$ does not depend on the variables $\tau''$,
we have
$$
\underset{1\le j \le \ell}{\prod} \rho_{\Lambda,j}(w)^{n_j/\comD + a_{jk}}
=
\underset{1\le j \le \ell}{\prod} \rho_{\Lambda,j}(z)^{n_j/\comD + a_{jk}}
\qquad (w \in B(z)).
$$
Let $\gamma$ be the closed circle in $\mathbb{C}$ of
center the origin and radius $|z^{(k)}| + \epsilon h(z)$ with the anti-clockwise
direction.
Then we have
$$
\begin{aligned}
\left|\dfrac{\partial g}{\partial z_i}(z) \right|
&=
\left|\dfrac{1}{2\pi\sqrt{-1}}
\int_\gamma \dfrac{g(z_1,\dots,z_{i-1}, t, z_{i+1}, \dots, z_n)}{(t - z_i)^2} dt \right| \\
&\le
\dfrac{C' (|z^{(k)}| + \epsilon h(z))}{(\epsilon h(z))^2}
\underset{1\le j \le \ell}{\prod} \rho_{\Lambda,j}(z)^{n_j/\comD + a_{jk}}
\quad (z \in S').
%\\
%&\le
%\dfrac{C(\epsilon + 1/C_1)}{\epsilon^2 \rho_k(z)}
%\underset{1\le j \le \ell}{\prod} \left(\varphi_j^{-1}(\rho(z))\right)^{n_j + a_{jk}}
%\quad (z \in S').
\end{aligned}
$$
Hence, by \eqref{eq:level-cond-1.1}, we get
$$
\left|\dfrac{\partial g}{\partial z_i}(z) \right|
\le \dfrac{C'(\epsilon + C)}{\epsilon^2 \varphi_k(\rho_{\Lambda}(z))}
\underset{1\le j \le \ell}{\prod} \rho_{\Lambda,j}(z)^{n_j/\comD + a_{jk}}
\quad (z \in S').
$$

Therefore we obtained the desired estimate in this case,
and the proof of the proposition has been completed.
\end{proof}

\begin{teo}{\label{theorem:multi-asymptotic}}
Let $V \in V(\zeta)$ and $\{\epsilon\}$ be the set of positive real numbers.
Let $f$ be a holomorphic function on $S=\Ssec$.
Then the following conditions are equivalent.
\begin{enumerate}
\item $f$ is multi-asymptotically developable along $\chi$ on $S$.
\item  For any multicone $S':=S_H(V', \{\epsilon'\})$
properly contained in $S$ and for any $\alpha \in \Z_{\ge 0}^n$,
$\left|\displaystyle\frac{\partial^\alpha f}{\partial z^\alpha}\right|$ is
bounded on $S'$.
\item For any multicone $S':=S_H(V', \{\epsilon'\})$
properly contained in $S$,
the holomorphic function $f\vert_{S'}$ on $S'$ can be extended to a $\C^\infty$-function
on $X_{\mathbb R}$ ($X_{\mathbb R}$ denotes the underlying real analytic manifold of $X$).
\end{enumerate}
\end{teo}
\begin{proof}
We first show 1.~implies 2.
it follows from Definition {\ref{def:multi-asymptotics}} with $N=(0,\dots,0)$ that
$f$ is bounded on $S'$. Since each higher derivative of $f$
is still multi-asymptotically developable thanks to
Proposition {\ref{prop:stability-differentiale}},
$\displaystyle\frac{\partial^\alpha f}{\partial z^\alpha}$
is also bounded for any $\alpha \in \Z^n_{\ge 0}$ on $S'$.

\medskip
As $S'$ is 1-regular by Lemma \ref{lem:1-regular}
and as $f$ is holomorphic (i.e.
$\displaystyle\frac{\partial^\alpha f}{\partial \bar{z}^\alpha} = 0$),
the claim 3.~follows from 2.~by the result of Whitney in \cite{Wh34}.
Clearly 3.~implies 2. Hence the claim 2.~and 3.~are equivalent.

\medskip

\medskip

Now we will show 3.~implies 1.
Assume that $f$ satisfies 3. In particular, any derivative of $f$
extends to $\overline{S'}$ and is bounded on $\overline{S'}$.
It follows from Lemma {\ref{lemma:multi-cone-geometric}} that
for $z \in S'$ and $0 \le \lambda_j \le 1$ ($j \in \{1,2,\dots, \ell\}$),
we get $\mu(z,\lambda) \in \overline{S'}$, which also implies
$\mu_{\comD A_\chi}(z,\,\lambda) = \mu(z,\,\lambda^\comD) \in \overline{S'}$.
Therefore,
for $N = (n_1, \dots, n_\ell) \in \Z^\ell_{\ge 0}$,
$$
\varphi_N(f;z) :=
\int_{0}^{\infty}\dots\int_{0}^{\infty}
\frac{\partial^N f(\mu_{\comD A_\chi}(z,\lambda))}{\partial \lambda^N}
\underset{1 \le j \le \ell}{\prod}
K_{n_j - 1}(1 - \lambda_j) d\lambda_1 \dots d\lambda_\ell
$$
where
$$
K_n(t) :=
\left\{
\begin{aligned}
\delta(t)\qquad &(n = -1), \\
\frac{t_+^n}{n!}      \qquad &(n \ge 0)
\end{aligned}
\right.
$$
is well-defined on $S'$. Here $\delta$ denotes the Dirac delta function and $t_+= t$ if $t \ge 0$, $t_+=0$ if $t < 0$.
Then, by integration by parts, we have
$$
\varphi_N(f;z) = f(z) - \sum_{J \in \Powl}
(-1)^{\#J + 1}\sum_{\beta \in \Z^{J}, \beta<_J N}
\left.
\frac{1}{\beta!}\partial^\beta_{\lambda, J}f(\mu_{\comD A_\chi}(z,\,\lambda))
\right\vert_{\lambda = \check{e}_J}.
$$
By taking Lemma \ref{lem:asymptotic_f} into account,
we will evaluate the remainder term $\varphi_N(f; z)$.
%We may assume $z^{(k)} \ne 0$ for $k \in J_Z(\zeta)$ in
%the subsequent argument
%because the result at an arbitrary point $z \in S'$ follows by letting $z^{(k)} \to 0$
%($k \in J_Z(\zeta)$).
Let us consider a coordinate transformation of $\lambda$ defined by
$$
\tilde{\lambda}_j = (\rho_{\Lambda,j}(z))^{1/\comD} \lambda_j
\qquad (j\in \{1,2,\dots,\ell\}).
$$
Note that $\rho_{\Lambda,j}(z)$ is different from $0$.
Then,
by Lemma \ref{lem:inverse_rho} and the fact that $\varphi_k(\varphi^{-1}(\tau',\lambda''),\lambda'')$ ($L < k \le m$) is independent of the variables $\lambda''$
due to Lemma \ref{lem:independent_psi_lambda}, we have
$$
\varphi_k(\lambda^\comD) =
\dfrac{\varphi_k(\tilde{\lambda})^\comD}{|z^{(k)}|}
\qquad (1 \le k \le L)
$$
and
$$
\varphi_k(\lambda^\comD) =
\varphi_{k}(\tilde{\lambda})^\comD
\dfrac{\psi_k(|z^{(*)}|)}{|z^{(k)}|} \qquad (L < k \le m).
$$
Hence we get, by noticing $\mu_{\intA}(z,\,\lambda) = \mu(z,\,\lambda^\comD)$,
$$
\begin{aligned}
&\dfrac{\partial^N} {\partial \lambda^N} f(\mu_{\comD A_\chi}(z,\lambda)) \\
&\qquad= \underset{1\le j \le \ell}{\prod}
\rho_{\Lambda,j}(z)^{n_j/\comD}
\frac{\partial^N }{\partial \tilde{\lambda}^N}
f\bigg(z^{(0)},\, \varphi_1(\tilde{\lambda})^\comD \dfrac{z^{(1)}}{|z^{(1)}|},\,
\dots,\, \varphi_L(\tilde{\lambda})^\comD \dfrac{z^{(L)}}{ | z^{(L)}|},\, \\
&\qquad\qquad
\varphi_{L+1}(\tilde{\lambda})^\comD
\psi_{L+1}(|z^{(*)}|)
\dfrac{z^{(L+1)}}{|z^{(L+1)}|},\,
\dots,\,
\varphi_{m}(\tilde{\lambda})^\comD
\psi_m(|z^{(*)}|)
\dfrac{z^{(\ell)}}{|z^{(\ell)}|}\bigg).
\end{aligned}
$$
%where $\psi_k(z) := \dfrac{|z^{(k)}|}{\rho_k(z)}$ for $k \le \ell$ and
%$\psi_k(z) := \dfrac{|z^{(k)}|}{\varphi_k(\varphi^{-1}(\rho(z)))}$ for $k > \ell$.

It follows from the condition
\eqref{eq:level-cond-2.1}
and
\eqref{eq:level-cond-2.2}
that
$\tilde{\lambda}_j$ is bounded if $0 \le \lambda \le 1$ and $z \in S'$.
Each $\psi_k(|z^{(*)}|)$ is also bounded on $z \in S'$ by the
conditions \eqref{eq:level-cond-1.1}.
Hence there exists $C > 0$ such that
$$
\sup_{\lambda \in [0,1]^\ell }
\left| \frac{\partial^N f(\mu_{\intA}(z,\lambda))}{\partial \lambda^N} \right| \le
C \underset{1\le j \le \ell}{\prod} \rho_{\Lambda,j}(z)^{n_j/\comD}
\quad (z \in S').
$$
This gives the desired estimate of the remainder term.
The proof has been completed.
\end{proof}

Assume $f$ is multi-asymptotically developable along $\chi$ on $S:=\Ssec$.
Then, for any multicone $S'$ properly contained in $S$ and for any multi-index $\alpha$,
it follows from the theorem that
$\dfrac{\partial^\alpha f}{\partial z^\alpha}$ can uniquely extend to
$\overline{S'}$, and then, by restricting the extension of $f$ to
$\pi_J(S') \subset \overline{S'}$, we have the holomorphic function on $\pi_J(S')$.
Since the union of $\pi_J(S')$'s where $S'$ runs through multicones properly
contained in $S$ is $S_J$, the family of these holomorphic functions
determines the holomorphic function on $S_J$ and we denote it by
$\left.\dfrac{\partial^\alpha f}{\partial z^\alpha}\right|_{Z_J}$.
Now, as an immediate consequence of the above proof, we have:
\begin{cor}{\label{cor:asymptptotic_F_J}}
If $f$ is multi-asymptotically developable along $\chi$ on $S$, then
it is multi-asymptotically developable to the particular $F=\{F_J\}$
where $F_J$ is defined by
\begin{equation}{\label{eq:consistent_data_from_f}}
F_J = \left\{ f_{J,\,\alpha} = \left.\dfrac{\partial^\alpha f}{\partial z^\alpha}\right|_{Z_J}
\right\}_{\alpha \in \mathbb{Z}^{(K_J)}_{\ge 0}}.
\end{equation}
\end{cor}
Note that, in general, a total family of coefficients
of multi-asymptotic expansion is not uniquely determined by $f$. We study details
on this problem in the subsequent subsections.

\subsection{Uniqueness of coefficients I}\label{subsection:Uniqueness of coefficients}

Let $S:=\Ssec$ and $\rho_\Lambda = (\rho_{\Lambda,1},\,\dots,\,\rho_{\Lambda,\ell})$
be a family of level functions of $\chi$.
Recall that, for $J, J' \in \subsetl$, the action $\mu_J$ is the one generated
by the actions $\mu_j$'s $(j \in J)$ and $\mu_J|_{Z_{J'}}$ denotes the restriction
of the action $\mu_J$ to the submanifold $Z_{J'}$.

\begin{df}
We say that the action $\mu_j$ is strict on $S$ with respect to
the family $\rho_\Lambda$ of level functions
if there exists an open dense subset $\widetilde{S}$
of $S$ such that
$$
\lim_{t \to 0+0} \rho_{\Lambda,j}(\mu_j(z,\, t)) = 0 \qquad (z \in \widetilde{S}).
$$
Furthermore, for $J \in \subsetl$, the action $\mu_J$ (resp. $\mu$) is said to be strict
if each $\mu_j$ with $j \in J$ (resp. any $j$)
is strict on $S$ with respect to
the family of level functions.
\end{df}

\begin{es}
Let us see Examples \ref{ex:transitive_non_degenerate}
and \ref{ex:transitive_non-degenerate}
given in
Subsection \ref{subsection:Level functions}.
Here all the actions $\mu_j$ are strict with respect to
the level functions $\rho_{\Lambda}$.
Note that the action $\mu$ is transitive and degenerate.
On the other hand, in Example \ref{ex:transitive_degenerate}, the
action $\mu_3$ is not strict with respect to
the level functions $\rho_{\Lambda}$ of this example.
\end{es}

We have the following lemma for the strictness of actions.
\begin{lem}{\label{lem:strictness}}
We have:
\begin{enumerate}
\item If the action $\mu$ is non-degenerate, each action $\mu_j$ $(1 \le j \le \ell)$
is strict with respect to $\rho_\Lambda$.
\item $\mu_\ell$ is always strict with respect to $\rho_\Lambda$.
\end{enumerate}
\end{lem}
\begin{proof}
The claim 1.~follows from the equality
$$
\rho_{\Lambda,j}(\mu_j(z,\, t)) =
%\rho_{\Lambda,j}(\varphi'(1,\dots,1,t,1,\dots,1))
%\rho_{\Lambda,j}(|z^{(*)}|') =
t\,\, \rho_{\Lambda,j}(z).
$$

Now we prove the second claim. We follow notations used in
Section \ref{subsection:Level functions}.
We first show that, for any $f$ in $(f,v) \in F^{0,j}_{<}$,
$$
f(\mu_\ell(z,\,t),\, t \lambda'') = f(z,\, \lambda'')
$$
holds for any $t > 0$, where
$
f(z,\,\lambda'')
$
denotes $f(|z^{(*)}|',\,\lambda'')$ for simplicity and
$$
t \lambda'' = (\lambda_{L+1}, \dots, \lambda_{\ell-1}, t \lambda_\ell).
$$
Since $\psi_k$ does not contain
the variables $\lambda''$, any $f$ with $(f,v) \in F^{0,j}_{<}$ is a
product of a rational power of $\varphi_j^{-1}(\tau',\,\lambda'')$'s.
Hence it suffices to show the above claim for $f = \varphi_j^{-1}(\tau',\lambda'')$.
Then, we have, by definitions,
$$
\tau' = \varphi'(\lambda',\, \lambda''),\quad
\mu_\ell(\tau',\, t) = \varphi'(\lambda', t \lambda''),
$$
from which the claim
$$
(\,\lambda' = \,)\,
\varphi_j^{-1}(z,\lambda'')
=
\varphi_j^{-1}(\mu_\ell(z,\,t), t \lambda'').
$$
follows.

Let $(f,v) \in F^{0,\ell-1}_<$. Then, by noticing that
$\operatorname{sol}^\lambda_{\ell}(f)$ is a function of the variable $\tau'$ only,
we have
$$
\operatorname{sol}^\lambda_{\ell}(f)(\mu_\ell(z,\,t))
= t \lambda_\ell f(\mu_\ell(z,\,t),\, t \lambda'')^{1/|a|}
= t \lambda_\ell f(z,\, \lambda'')^{1/|a|}
$$
with $a := \nu_\ell^\lambda(f) < 0$.
Hence we have obtained
$$
\operatorname{sol}^\lambda_{\ell}(f)(\mu_\ell(z,\,t)) \to 0
\quad (t \to 0+0),
$$
which immediately implies
$
\rho_{\Lambda,\ell}(\mu_\ell(z,\,t)) \to 0
$.
This completes the proof.
\end{proof}

%
%\begin{df}
%We say that the action $\mu$ is covered by strict actions on $S$ with respect
%to the family of level functions if, for any $J \in \subsetl$, there exists
%$1 \le j  \le \ell$ such that $Z_J \subset Z_j$ and $\mu_j$ is
%strict on $S$ with respect to the family of level functions.
%\end{df}
%\begin{es}
%As we have already seen, the action $\mu_3$ in Example \ref{ex:transitive_degenerate}
%is not strict, however, the action $\mu$ in this example itself is covered
%by strict actions.
%\end{es}
Let $F = \{F_J\}_{J \in \subsetl}$
with $F_J = \{f_{J,\,\alpha}\}_{\alpha \in \mathbb{Z}^{(K_J)}_{\ge 0}}$
be a total family of
coefficients of multi-asymptotic expansion along $\chi$ on $S$, and
$f$ a holomorphic function on $S$ which is strongly asymptotically developable
to $F$ along $\chi$ on $S$. Under this situation, we have:
%Then, by Theorem \ref{theorem:multi-asymptotic},
%$f$ can be uniquely extended to $\overline{S'}$ as a $C^\infty$ function for any
%multicone
%$S' = S_H(V',\, \{\epsilon'\})$ properly contained in $S$.
%We can regard $f$ as a function on $\overline{S'}$, and hence, for any $J \in \subsetl$,
%the function
%$
%\left.\dfrac{\partial^\alpha f}{\partial z^\alpha}\right|_{Z_J}
%$
%is well-defined as a function on $S_J$ by considering all the $S'$ properly
%contained in $S$.
\begin{teo}{\label{th:uniqueness_of_coe}}
Let $\widetilde{J}$ be the subset of $\{1,2,\dots,\ell\}$
which consists of an index $j \in \{1,2,\dots,\ell\}$
such that $\mu_j$ is strict on $S$ with respect to the family $\rho_\Lambda$ of level
functions.  Then we have the following.
\begin{enumerate}
	\item For any non-empty subset $J$ of $\{1,2,\dots,\ell\}$ and
$\alpha \in \mathbb{Z}^{(K_J)}_{\ge 0}$,
any (higher) derivative of $f_{J,\alpha}$ is bounded on $S'_J$
for every $S'=S_H(V',\, \{\epsilon'\})$ properly contained in $S$.
\item For any non-empty subset $J$ of $\widetilde{J}$, the family
$F_J = \{f_{J,\,\alpha}\}_{\alpha \in \mathbb{Z}^{(K_J)}_{\ge 0}}$
is uniquely determined
as follows.
$$
f_{J,\,\alpha} = \left.\dfrac{\partial^\alpha f}{\partial z^\alpha}\right|_{Z_J}
\quad (\alpha \in \mathbb{Z}^{(K_J)}_{\ge 0}).
$$
\end{enumerate}
\end{teo}
\begin{proof}
We only prove the claim 2.~because the proof of 1.~goes in the same way as that
of 2.
It suffices to show the claim on $S'$ properly contained in $S$. Hence,
from the beginning, we may assume that $f$ can be extended to $\overline{S}$ as
a $C^\infty$ function. Furthermore, by Proposition \ref{prop:stability-differentiale},
it suffices to show the claim for $\alpha = 0$, i.e.,
$
f_{J,\,0} =
f|_{Z_J}
$
for any non-empty subset $J$ of $\widetilde{J}$.

We will show the claim by induction of the number of elements in $J$.
Assume $J = \{k\}$, i.e. $\# J = 1$. Let $N=(n_1,\dots, n_\ell)$ with $n_j = 0$
($j \ne k$) and $n_k = 1$, and consider the asymptotic expansion for this $N$.
Then we have
$$
|f(z) - f_k(\pi_J(z))| \le C \rho_{\Lambda,k}(z)^{1/\comD}.
$$
Let $z' \in S_J$ and take $\hat{z}$ with $\pi_J(\hat{z}) = z'$.
We may assume $\hat{z} \in \widetilde{S}$. Otherwise we take an
sequence $\{\hat{z}_r\}_{r=1}^\infty$ of
$\widetilde{S}$ with $\hat{z}_r \to \hat{z}$ ($r \to \infty$)
and apply the following argument to each point $\pi_J(\hat{z}_r)$.
Then, by taking $r \to \infty$,
we finally obtain the result at $z'$.

By putting $z = \mu_k(\hat{z}, \lambda) \in S$ into the asymptotic expansion
and by letting $\lambda$ to
the zero, we have
$$
\rho_{\Lambda,k}(\mu_k(\hat{z},\, \lambda)) \to 0
$$
as $\mu_k$ is strict. Hence we obtain $f(z') = f_k(z')$.

Assume that the theorem were true for $\#J = \kappa \ge 1$. We will show the theorem
for $\#J = \kappa + 1$. Take $j_0 \in J$ and fix it. We denote by $\widetilde{S}$
the exceptional set with respect to the action $\mu_{j_0}$
which appears in the definition.
Let $z' \in S_J$ and take $\hat{z}$ with $\pi_J(\hat{z}) = z'$.
Let $\{\lambda_r\}_{r=1}^\infty$ be a
sequence of real numbers in $(0,1]$ which tends to $0$ $(r \to \infty)$.
Set $\hat{z}_r = \mu_J(\hat{z}, \lambda_r)$. If $\hat{z}_r \notin \widetilde{S}$,
we replace $\hat{z}_r$ with a point $\hat{z}'_r$ in $\widetilde{S}$ sufficiently
close to the $\hat{z}_r$ with $|\hat{z}_r - \hat{z}'_r| \le 1/r$.
Then we determine a sequence $\{\delta_r\}_{r=1}^\infty$
of real numbers in $(0,1]$ so that
$$
\rho_{\Lambda,j_0}(\mu_{j_0}(\hat{z}_r,\, \delta_r)) \le 1/r
\qquad r=1,2,\dots.
$$
This is possible because $\mu_{j_0}$ is strict.
Set $N=(n_1,\dots, n_\ell)$ with $n_j = 1$ if $j \in J$ and $n_j=0$ if $j \notin J$,
and consider the asymptotic expansion for this $N$. Then we obtain
$$
\left|f(z) - \sum_{J'}
(-1)^{\#J'+1} f_{J'}(\pi_{J'}(z))
- (-1)^{\#J+1} f_J(z)\right|
\le C\prod_{j \in J} \rho_{\Lambda,j}(z)^{1/\comD},
$$
where the first sum is taken over all the subsets $J'$ in $J$ except for $J$ itself
and the empty set. Set
$$
z_r = \mu_{j_0}(\hat{z}_r,\, \delta_r) \in S.
$$
By putting $z_r$ into the expansion and by letting $r \to \infty$,
the right-hand side of the expansion tends to $0$, and $f(z)$ and $f_{J'}(\pi_{J'}(z))$
in the left-hand side tends to $f(z')$ by the induction hypothesis.
Note that each $f_{J'}$ ($J' \ne J$) extends to $\overline{S_{J'}}$ as a $C^\infty$ function because $f$ extends to $\overline{S}$ by the assumption
and $\overline{S_{J'}} = Z_{J'} \cap \overline{S}$ holds by Lemma \ref{lemma:multi-cone-geometric}.
Hence we have $f(z') = f_J(z')$, which completes the proof.
\end{proof}

Let $J \in \subsetl$ and $S:=\Ssec$ be a multicone, and let
$\{f_\alpha\}_{\alpha \in \mathbb{Z}_{\ge 0}^{(K_J)}}$ be a family of
holomorphic functions on $S_J$. Note that $S_J$ is also a multicone
in $Z_J$ with respect to the action $\mu|_{Z_J}$ by Proposition \ref{prop:multi-cone-S_J}.
We assume:
For any
$\alpha \in \mathbb{Z}^{(K_J)}_{\ge 0}$,
$\beta \in \mathbb{Z}^{(\compl K_J)}_{\ge 0}$ and for any multicone $S_J' \subset Z_J$
properly contained in $S_J$,
\begin{equation}{\label{eq:cond-bound-S_J}}
\text{$\left|\dfrac{\partial^\beta f_\alpha}{\partial z^\beta}\right|$
is bounded on $S_J'$.}
\end{equation}

Let $R \in \subsetl$ with $J \subset R$. We denote by $\pi_{R,J}: Z_J \to Z_{R}$
the canonical projection, i.e.,
$$
\left(z^{(\compl K_J)}\right) \in Z_J \to \left(z^{(\compl K_R)}\right) \in Z_R.
$$
Then it follows from the definition that we have
$$
\pi_{R} = \pi_{R,J} \circ \pi_J
$$
and, by Lemma \ref{lemma:multi-cone-geometric} and
Proposition \ref{prop:multi-cone-S_J},
$$
S_{R} = \pi_{R,J}(S_J) \quad\text{and}\quad
\overline{S_{R}} = \overline{S_J} \cap Z_{R}.
$$

For a multicone $S_J' \subset Z_J$ properly contained in $S_J$, it follows from
the condition \eqref{eq:cond-bound-S_J} and the $1$-regularity of $S_J'$ that
$\dfrac{\partial^\beta f_\alpha}{\partial z^\beta}$ can uniquely extend to
$\overline{S'_J}$ and its restriction to $\pi_{R,J}(S'_J) \subset \overline{S'_J}$
is a holomorphic function
on $\pi_{R,J}(S_J')$. Since
the union of $\pi_{R,J}(S_J')$ is $S_{R}$
when $S_J'$ runs through multicones in $Z_J$ is properly contained in $S_J$,
we have obtained a holomorphic function on $S_{R}=\pi_{R,J}(S_J)$ from these restrictions.
Hereafter we denote it by
$
\left.\dfrac{\partial^\beta f_\alpha}{\partial z^\beta}\right|_{Z_{R}}.
%\quad\text{or}\quad \varpi_{J',J}\left(\dfrac{\partial^\beta f_\alpha}{\partial z^\beta}\right).
$
%Note that the holomorphic function
%$\left.\dfrac{\partial^\beta f_\alpha}{\partial z^\beta}\right|_{Z_{R}}$
%satisfies the
%condition \eqref{eq:cond-bound-S_J} with replace of $J$ with $R$ also.

Now by using this restriction, we define
a family $\{g_\gamma\}_{\gamma \in \mathbb{Z}^{(K_{R})}_{\ge 0}}$ of
holomorphic functions on $S_{R}$ by
\begin{equation}
	g_\gamma =
	\left.\dfrac{\partial^\beta f_\alpha}{\partial z^\beta}\right|_{Z_{R}},
\end{equation}
where $\alpha \in \mathbb{Z}^{(K_J)}_{\ge 0}$ and
$\beta \in \mathbb{Z}^{(\compl K_J)}_{\ge 0}$ are determined by the formula
$$
\gamma = \alpha + \beta.
$$
We denote by $\varpi_{R,J}$
the map from a family of holomorphic functions on $S_J$
to the one on $S_{R}$ defined by the above correspondence, that is,
\begin{equation}
\varpi_{R,J}\left(\{f_\alpha\}_{\alpha \in \mathbb{Z}_{\ge 0}^{(K_J)}}\right) :=
\{g_\gamma\}_{\gamma \in \mathbb{Z}^{(K_{R})}_{\ge 0}}.
\end{equation}

\begin{df}
Let $S := \Ssec$ and let $\{F_J\}_{J \in \subsetl}$
with $F_J = \left\{f_{J,\,\alpha}\right\}_{\alpha}$ be a total family of coefficients of
multi-asymptotic expansion along $\chi$ on $S$. Then $\{F_J\}$ is said to be
consistent if it satisfies the following conditions:
\begin{enumerate}
\item For any $J$ and any $\alpha$, the holomorphic function $f_{J,\,\alpha}$
satisfies the condition \eqref{eq:cond-bound-S_J} on $S_J$.
\item For any pair $J \subset R$ $(J,\, R \in \subsetl)$, we have
$$
\varpi_{R,J}\left(\{f_{J,\,\alpha}\}\right) = \left(\{f_{R,\,\alpha}\}\right).
$$
\end{enumerate}
\end{df}

The following theorem easily follows from Corollary \ref{cor:asymptptotic_F_J} and
Theorem {\ref{th:uniqueness_of_coe}}.
\begin{teo}
Assume $f$ to be multi-asymptotically developable along $\chi$ on $S$.
Then there exists
a consistent total family of coefficients of
multi-asymptotic expansion such that $f$
is multi-asymptotically developable to this consistent family.
Furthermore, if all the actions $\mu_j$'s are strict on $S$ with respect to
the family $\rho_\Lambda$ of level functions, then a total family of
coefficients to which a holomorphic function on $S$
is multi-asymptotically developable is consistent.
\end{teo}

Note that, in the next subsection, we study the case
where some $\mu_j$ is not strict on $S$.
%the family $\rho_\Lambda$ of level functions, then a total family of
%If the action $\mu$ is not strict, we can no longer expect uniqueness of
%a total family of coefficients to which a holomorphic function is multi-asymptotically
%developable, however, if we restrict a family to the consistent one, then
%we can obtain the following weak uniqueness if $\mu$ is covered by strict actions.
%\begin{prop}
%Assume that $\mu$ is covered by strict actions, and assume that
%a holomorphic function $f$ on $S$ is multi-asymptotically developable
%to both the consistent families $\{F_J\}$ and $\{F'_J\}$. Then
%we have $\{F_J\} = \{F'_J\}$.
%\end{prop}
We also have the following theorem.
\begin{teo} \label{teo: consistent family Whitney}
Let $S := \Ssec$ and let $\{F_J\}_{J \in \subsetl}$
with $F_J = \left\{f_{J,\,\alpha}\right\}_{\alpha}$ be a total family of coefficients of
multi-asymptotic expansion along $\chi$ on $S$.
Then the following conditions are equivalent:
\begin{enumerate}
\item $\{F_J\}$ is consistent.
\item For any multicone $S'$ properly contained in $S$,
there exists a $\C^\infty$-function $g(z)$ on $X_{\R}$ such that for any
$J \in \subsetl$ and $\alpha \in \Z^{(K_J)}_{\ge 0}$ we have
$$
f_{J, \alpha} =
\left. \displaystyle\frac{\partial^\alpha g}
{\partial z^\alpha}\right|_{S'_J} \quad (z \in S'_J).
$$
\end{enumerate}
\end{teo}
\begin{proof}
Clearly the claim 2.~implies the claim 1. Hence we will show
the implication from 1.~to 2.
For each $J \in \subsetl$,
since $S'_J$ is $1$-regular by Proposition \ref{prop:multi-cone-S_J} and
since the condition \eqref{eq:cond-bound-S_J} holds on $S_J$,
the family $\{F_J\}$ defines a Whitney jet on $\overline{S'_J}$
(see \cite{M66} for Malgrange's definition of a $\C^\infty$-Whitney jet).
Note that, by the claim 2.~of Lemma \ref{lemma:multi-cone-geometric}, we have
$$
\overline{S'} \cap (Z_1 \cup \dots \cup Z_\ell)
= \underset{J \in \subsetl}{\bigcup}\,\, (\overline{S'} \cap Z_J)
= \underset{J \in \subsetl}{\bigcup}\,\, \overline{S'_J}.
$$
Note also that, for $J$ and $J'$ in $\subsetl$,
the Whitey jets $\{F_J\}$ and $\{F_{J'}\}$ coincide with
$\{F_{J \cup J'}\}$ on $S'_{J \cup J'}$ by the second condition in the definition
of a consistent family,
and hence, they coincide on
$$
\overline{S'_J} \cap \overline{S'_{J'}} =
(\overline{S'} \cap Z_J) \cap (\overline{S'} \cap Z_{J'})
= \overline{S'} \cap (Z_J \cap Z_{J'})
= \overline{S'_{J \cup J'}}.
$$
Hence it follows from Theorem 5.5 of \cite{M66} that we obtain the Whitney
jet defined on $\overline{S'} \cap (Z_1 \cup \dots \cup Z_\ell)$, which
implies the claim 2.
\end{proof}

\subsection{Uniqueness of coefficients II}\label{subsection:Uniqueness of coefficients 2}

In the previous subsection, we see that, if the action $\mu$ is degenerate,
then uniqueness of coefficients of multi-asymptotics is not necessarily guaranteed.
Furthermore Lemma \ref{lem:strictness} indicates that the construction
of level functions $\rho_\Lambda$ depends on the order of actions $\mu_1$, $\dots$,
$\mu_\ell$, and in particular
the last action $\mu_\ell$ is always strict.
Hence we can recover the uniqueness of coefficients
by considering all the possible order of actions as we will see.

Let $\Theta_\mu$ be the set of permutations $\theta$ on $\{1,2,\dots,\ell\}$ such that
the $L$-rows in $A_\chi$ corresponding to the actions $\mu_{\theta(1)}$, $\dots$,
$\mu_{\theta(L)}$ are linearly independent.
For $\theta \in \Theta_\mu$, we construct the family of level
functions $\rho^\theta_\Lambda$ with respect to the sequence of actions
\begin{equation}
\mu^\theta = \{\mu_{\theta(1)},\, \dots,\, \mu_{\theta(\ell)}\}
\end{equation}
as in Subsection \ref{subsection:Level functions} where we keep the order of
actions specified in the sequence $\mu^\theta$. Then we define a generalized family of
level functions
$\hat{\rho}_\Lambda = (\hat{\rho}_{\Lambda,1},\, \dots,\, \hat{\rho}_{\Lambda,\ell})$ by
\begin{equation}
	\hat{\rho}_{\Lambda,j}(z) := \min_{\theta \in \Theta_\mu}\, \rho^\theta_{\Lambda, \theta^{-1}(j)}(z) \qquad j =1,2,\dots,\ell.
\end{equation}
By using $\hat{\rho}_{\Lambda}$, we can define the multi-asymptotic developability
along $\chi$ with level functions $\hat{\rho}_{\Lambda}$.
That is:
We say that $f$ is multi-asymptotically developable to $F=\{F_J\}$ along
$\chi$ on $S$ with level functions $\hat{\rho}_{\Lambda}$ if for any cone
$S'=S_H(V',\{\epsilon'\})$ properly contained in $S$ and for any $N=(n_1,\dots,n_\ell)
\in \mathbb{Z}^\ell_{\ge 0}$, there exists a constant $C_{S',N}$ such that
$$
\left|
f(z) - \operatorname{App}^{<N}(F;\, z)
\right|
\le C_{S',N}
\underset{1 \le j \le \ell}{\prod}
\hat{\rho}_{\Lambda,j}(z)^{n_j/\comD} \qquad (z \in S').
$$

\

Clearly, if $f$ is multi-asymptotically developable to $F=\{F_J\}$
with level functions $\hat{\rho}_{\Lambda}$, then
$f$ is multi-asymptotically developable to $F=\{F_J\}$ with level functions
$\rho^\theta_{\Lambda}$ for any $\theta \in \Theta_\mu$. Conversely, if
$f$ is multi-asymptotically developable to $F=\{F_J\}$
with level functions $\rho^\theta_{\Lambda}$ for some $\theta \in \Theta_\mu$,
it follows from Theorem \ref{theorem:multi-asymptotic} that
$f$ is multi-asymptotically developable to $F=\{F_J\}$
with level functions $\hat{\rho}_{\Lambda}$.
Hence any level functions either $\rho^\theta_{\Lambda}$ or
$\hat{\rho}_{\Lambda}$ give the same asymptotical developability.
An advantage of $\hat{\rho}_{\Lambda}$ is:
\begin{prop}
Each action $\mu_j$ $(1 \le j \le \ell)$ is strict
with respect to $\hat{\rho}_{\Lambda}$,
that is, for any $1 \le j \le \ell$,
$$
\lim_{t \to 0+0} \hat{\rho}_{\Lambda,j}(\mu_j(z,\, t)) = 0 \qquad (z \in S).
$$
\end{prop}
\begin{proof}
Take $1 \le i_0 \le L$.

First assume that any $j$-th row $(L+1 \le j \le \ell)$
is a linear combination of $i$-th rows $(1 \le i \le L,\, i \ne i_0)$.
Then it is easy to see that $\varphi^{-1}_{i_0}(\tau', \lambda'')$ is
independent of $\lambda''$, and hence, we have
$$
\rho_{\Lambda,i_0}(\tau') = \varphi^{-1}_{i_0}(\tau').
$$
Therefore we get
$$
\rho_{\Lambda,i_0}(\mu_{i_0}(z,t)) =
\varphi^{-1}_{i_0}(\mu_{i_0}(z,t)) =
t\varphi^{-1}_{i_0}(z)  \to 0 \quad (t \to 0+0),
$$
which implies that $\mu_{i_0}$ is strict with respect to
$\hat{\rho}_{\Lambda}$ in this case.

Next, assume that some $j_0$-th row $(L+1 \le j_0 \le \ell)$ is not
given by a linear combination of $i$-th rows $(1 \le i \le L, i\ne i_0)$.
Then there exists a permutation in $\Theta_\mu$ which exchanges
$j_0$ and $i_0$. Hence we can find a permutation $\theta \in \Theta_\mu$
with $\theta(\ell) = i_0$, that is,
the $L$-rows in $A_\chi$
corresponding to the actions $\mu_{\theta(1)}$, $\dots$, $\mu_{\theta(L)}$
are linearly independent and the last action $\mu_{\theta{(\ell)}}$ of $\mu^\theta$ is $\mu_{i_0}$.
Then, by Lemma \ref{lem:strictness},
the $\mu_{i_0} = \mu_{\theta(\ell)}$ is strict with respect to $\rho^\theta_{\Lambda}$.

From these observations,
we see that, for any $i_0 \le L$, the action $\mu_{i_0}$
is strict with respect to $\hat{\rho}_{\Lambda}$.
By the same argument, we can show that $\mu_{i_0}$ for
$i_0 > L$ is also strict with respect to $\hat{\rho}_{\Lambda}$. This completes the proof.
\end{proof}
By this proposition, we finally established uniqueness of coefficients.
\begin{teo}{\label{th:general_uniqueness_of_coe}}
Let $F = \{F_J\}_{J \in \subsetl}$
with $F_J = \{f_{J,\,\alpha}\}_{\alpha \in \mathbb{Z}^{(K_J)}_{\ge 0}}$
be a total family of
coefficients of multi-asymptotic expansion along $\chi$ on $S$, and let
$f$ be a holomorphic function on $S$ which is strongly asymptotically developable
to $F$ along $\chi$ on $S$
with level functions $\hat{\rho}_{\Lambda}$. Then
$F_J$
is uniquely determined by
$$
f_{J,\,\alpha} = \left.\dfrac{\partial^\alpha f}{\partial z^\alpha}\right|_{Z_J}
\quad (\alpha \in \mathbb{Z}^{(K_J)}_{\ge 0}).
$$
In particular, $F=\{F_J\}_{J \in \subsetl}$ is consistent.
\end{teo}

In what follows, a family of level functions
 for which
every $\mu_j$ $(1 \le j \le \ell)$ is strict is assumed to be chosen.
Note that this assumption is satisfied if we choose
\begin{enumerate}
\item the usual $\rho_{\Lambda}$
when the action $\mu$ is non-degenerate.
\item the generalized family $\hat{\rho}_{\Lambda}$ of level functions.
\end{enumerate}
Then we can define the notion of flatness along $\chi$
for a multi-asymptotically developable function as in the usual asymptotics.
\begin{df}
Let $S := \Ssec$ and let $f$ be a holomorphic function on $S$.
We say that $f$ is flat on $S$ along $\chi$ if
$f$ is multi-asymptotically developable on $S$ along $\chi$ to
the total family $\{0_J\}_{J \in \subsetl}$.
\end{df}

Then we immediately obtain the following theorem.
\begin{teo}{\label{teo:flat-holomorphic}}
Let $S := \Ssec$ and let $f$ be a holomorphic function on $S$.
Then $f$ is flat on $S$ along $\chi$ if and only if,
for any multicone $S':=S_H(V', \{\epsilon'\})$ properly contained in $S$,
the holomorphic function $f|_{S'}$ on $S'$ extends to the whole space as a
$\mathcal{C}^\infty$-function and this extension satisfies
\begin{equation}
\left. \displaystyle\frac{\partial^\alpha f}{\partial z^\alpha}\right|_{
\overline{S'}\, \cap\, \left(\bigcup_{j=1}^\ell Z_j\right)} = 0
\qquad (\alpha \in \Z_{\ge 0}^n).
\end{equation}

\end{teo}

One important and interesting fact is that holomorphic functions
appearing in a consistent family themselves are also strongly asymptotically developable,
which we will explain from now.

As was explained in Subsection \ref{subsection:The local model}, the local model
is determined by the base space $X$ and the action $\mu$ on $X$. Hence, in what follows,
we sometimes denote by $\widetilde{(X,\,\mu)}$ the local model associated with $X$ and
$\mu$.%, whose zero section is also denoted by $S_{\widetilde{(X,\mu)},\,\chi}$.

Let $J \in \subsetl$. Define $J^* \subset \{1,2,\dots,\ell\}$ by
\begin{equation}
J^* := \{j \in \{1,2,\dots,\ell\};\,
\text{$\mu_j|_{Z_J}$ is not the identity action on $Z_J$}\}.
\end{equation}
Note that $J^* \subset \{1,2,\dots,\ell\}\setminus J$ hold.
\begin{oss}
If there exists a pair $j \ne j' \in J^*$ which gives the same
induced action $\mu_j|_{Z_J} = \mu_{j'}|_{Z_J}$ on $Z_J$, we may
remove one of either $j$ or $j'$ in $J^*$, and we can
resolve such a duplication.
\end{oss}
We denote by $\mu_{J^*}$ the action on $X$ generated by $\mu_j$'s with
$j \in J^*$, i.e.,
\begin{equation}
	\mu_{J^*} = \underset{j \in J^*}{\circ} \mu_{j},
\end{equation}
and we also denote by $\chi_{J^*}$ the family of closed submanifolds in $Z_J$
which consists of $Z_j \cap Z_J$ with $j \in J^*$.
That is,
\begin{equation}
	\chi_{J^*} = \{Z_j \cap Z_J\}_{j \in J^*}.
\end{equation}
Now set
\begin{equation}
\widetilde{Z_J} := \widetilde{(Z_J,\,\mu_{J^*}|_{Z_J})}.
\end{equation}
Then we have $\widetilde{Z_J} = Z_J \times \mathbb{R}^{\# J^*}$, and
the coordinates of $\widetilde{Z_J}$
and those of the zero section $S_{\chi_{J^*}}$ of $\widetilde{Z_J}$
are given by
$$
\left(z^{(\compl K_J)};\, \{t_j\}_{j \in J^*}\right) \quad\text{and}\quad
\left(z^{(0)};\, \zeta^{(\compl K_J \setminus \{0\})}\right),
$$
respectively.
For $p = (z^{(0)};\,\zeta^{(1)},\dots,\zeta^{(m)}) \in S_\chi$,
we define the map $(\pi_J)_*: S_\chi \to S_{\chi_{J^*}}$ by
$$
(z^{(0)};\,\zeta^{(1)},\dots,\zeta^{(m)}) \to
\left(z^{(0)};\, \zeta^{(\compl K_J \setminus \{0\})}\right).
$$
Then, for simplicity, we introduce the following assumption for the point $p$ under
consideration:
\begin{equation}{\label{eq:cond-position}}
\text{For any $J \in \subsetl$ with $J^* \ne \emptyset$,
$(\pi_J)_*(p)$ is outside the fixed points of $S_{\chi_{J^*}}$.}
%and the action $\mu_{J^*}|_{Z_J}$
%is strict on $S_J$ with respect to the family of level functions in $Z_J$.
\end{equation}

Let $F:=\{F_J\}$ be a total family of coefficients of multi-asymptotic expansion,
and let $J \in \subsetl$ and $\alpha \in \mathbb{Z}^{(K_J)}_{\ge 0}$.
For any non-empty subset $R \subset J^*$, we define a family
$G_R := \{g_{R,\,\beta}\}_{\beta \in \mathbb{Z}^{(K_R \setminus K_J)}_{\ge 0}}$ of holomorphic functions
on $\pi_{R \cup J,J}(S_J) = \pi_{R\cup J}(S)$ by
\begin{equation}
g_{R,\,\beta} = f_{R \cup J,\, \beta + \alpha}
\qquad \left(\beta \in \mathbb{Z}^{(K_R \setminus K_J)}_{\ge 0}\right).
\end{equation}
Hence, for any $J \in \subsetl$ and for any $\alpha \in \mathbb{Z}^{(K_J)}_{\ge 0}$,
we can obtain, in the space $Z_J$, a total family $\{G_R\}_{\emptyset \ne R \subset J^*}$
of coefficients of multi-asymptotic expansion on $S_J$ along $\chi_{J^*}$, and
set
\begin{equation}
\varpi_{J,\alpha}(F) := \{G_R\}_{\emptyset \ne R \subset J^*}.
\end{equation}
Then, since $S_J$ is a multicone in $Z_J$ with respect to the
action $\mu_{J^*}|_{Z_J}$ by Proposition \ref{prop:multi-cone-S_J},
we can apply Theorem \ref{theorem:multi-asymptotic} to each $f_{J,\alpha}$ and
obtain the following theorem.
\begin{teo}
Let $F:=\{F_J\}$ be a total family of coefficients of multi-asymptotic expansion
on $S:=\Ssec$ along $\chi$, and let $p$ be a point in $S_\chi$ satisfying
the condition \eqref{eq:cond-position}.
Then the following two conditions are equivalent.
\begin{enumerate}
\item $F$ is consistent.
\item $F$ satisfies the conditions C1.~and C2.~below.
\begin{enumerate}
\item[C1.] $F_J = F_{J'}$ holds for $J$ and $J'$ in $\subsetl$ with $Z_J = Z_{J'}$,
\item[C2.] For any $J \in \subsetl$ with $J^* \ne \emptyset$ and for any
$\alpha \in \mathbb{Z}_{\ge 0}^{(K_J)}$,  $f_{J,\alpha}$ is,
in the space $Z_J$, strongly asymptotically developable to $\varpi_{J,\alpha}(F)$ on $S_J$ along $\chi_{J^*}$.
\end{enumerate}
\end{enumerate}
\end{teo}
\begin{oss}
The condition 2.~of the theorem was adopted
as the definition of a consistent family in our previous paper \cite{HP13} or in
Majima's original definition for the case of a normal crossing divisor \cite{Ma84}.
As a matter of fact, both  cases satisfy the condition
\eqref{eq:cond-position}.
\end{oss}

\subsection{The classical sheaves of multi-asymptotically developable functions}{\label{subsec:classical_sheaves}}

As we have seen in the previous subsections, the multi-asymptotical developability is
a local notion with respect to a point in $S_\chi$. Hence, by a usual argument,
we can construct the classical sheaves on $S_\chi^\circ$ which reflect this local
nature. For reader's convenience, we quickly review their
constructions.

We denote by $S_\chi^\circ$ the set of points in $S_\chi$ outside the fixed points
hereafter.
For $p = (z^{(0)};\,\zeta) \in S_\chi^\circ$, we define
\begin{equation}
S(z^{(0)},\, V,\, \epsilon) := (z^{(0)},0,\dots,0) + S(V,\,\epsilon),
\end{equation}
for $V \in V(\zeta)$ and $\epsilon > 0$. Note that
it follows from the construction of $F_q$ that $v$ of $(f,v) \in F^q$
can be considered as a real analytic function of $|\zeta|^{(k)}$'s
$(k = 1,\dots,m)$ on $S_\chi^\circ$.

Let $W$ be an open subset in $S_\chi$ and $U$ an open subset in $X$.
Recall that $F = \{F_J\}_{J \in \subsetl}$
with $F_J = \{f_{J,\,\alpha}\}_{\alpha \in \mathbb{Z}^{(K_J)}_{\ge 0}}$ is said to
be a total family of coefficients along $\chi$ on $U$ if
each $f_{J,\,\alpha}$ is a holomorphic function on $\pi_J(U)$, where
$\pi_J$ is the canonical projection from $X$ to $Z_J$.

Let $\{S_{p}\}_{p \in W}$ be a family of multicones indexed by
the points in $W$, where
each $S_p$ for $p = (z^{(0)};\,\zeta)$ is $S(z^{(0)}, V, \epsilon)$
for some $V \in V(\zeta)$ and $\epsilon > 0$. Let $\mathcal{S}(W)$
be the set of such families $\{S_p\}_{p \in W}$. Note that $\mathcal{S}(W)$
becomes filtrant by the partial order
$$
\{S_p\}_{p \in W} \prec \{S'_p\}_{p \in W}
\iff
S_p \supset S'_p \qquad (p \in W).
$$
Now we define the presheaves on $S_\chi^\circ$:
\begin{equation}
\widetilde{\mathcal{A}}_\chi(W)
:= \varinjlim_{\{S_p\}_p \in \mathcal{S}(W)}
\left\{f \in \mathcal{O}_X(\displaystyle\bigcup_{p \in W} S_p);\,
\begin{aligned}
&\text{$f|_{S_p}$ is multi-asymptotically} \\
&\text{developable along $\chi$ on each $S_p$} \\
&\text{for $p \in W$}
\end{aligned}
\right\},
\end{equation}
\begin{equation}
	\widetilde{\mathcal{A}}^{<0}_\chi(W)
:= \varinjlim_{\{S_p\}_p \in \mathcal{S}(W)}
\left\{f \in \mathcal{O}_X(\displaystyle\bigcup_{p \in W} S_p);\,
\begin{aligned}
&\text{$f|_{S_p}$ is flat along $\chi$ on each $S_p$} \\
&\text{for $p \in W$}
\end{aligned}
\right\},
\end{equation}
\begin{equation}
\widetilde{\mathcal{A}}^{CF}_\chi(W)
:= \varinjlim_{\{S_p\}_p \in \mathcal{S}(W)}
\left\{F=\{F_J\};\,
\begin{aligned}
&\text{$F$ is a total family on $\displaystyle\bigcup_{p \in W} S_p$, and} \\
&\text{$F|_{S_p}$ is consistent on each $S_p$ }\\
&\text{for $p \in W$}
\end{aligned}
\right\}.
\end{equation}

Then, by sheafication, we have obtained the sheaves
$\mathcal{A}_\chi$, $\mathcal{A}^{<0}_\chi$ and
$\mathcal{A}^{CF}_\chi$ on $S_\chi^\circ$ associated with the presheaves
$\widetilde{\mathcal{A}}_\chi$, $\widetilde{\mathcal{A}}^{<0}_\chi$ and
$\widetilde{\mathcal{A}}^{CF}_\chi$, respectively. Note that
we have, for $p = (z^{(0)};\,\zeta) \in S_\chi^\circ$,
\begin{equation}
\mathcal{A}_{\chi,p}
:= \varinjlim_{V,\epsilon}
\left\{f \in \mathcal{O}_X(S(z^{(0)},V,\epsilon));\,
\begin{aligned}
&\text{$f$ is multi-asymptotically} \\
&\text{developable along $\chi$ on $S(z^{(0)},V,\epsilon)$}
\end{aligned}
\right\},
\end{equation}
\begin{equation}
\mathcal{A}^{<0}_{\chi,p}
:= \varinjlim_{V,\epsilon}
\left\{f \in \mathcal{O}_X(S(z^{(0)},V,\epsilon));\,
\text{$f$ is flat along $\chi$ on $S(z^{(0)},V,\epsilon)$}
\right\},
\end{equation}
\begin{equation}
\mathcal{A}^{CF}_{\chi,p}
:= \varinjlim_{V,\epsilon}
\left\{F=\{F_J\};\,
\begin{aligned}
&\text{$F$ is a consistent family of coefficients}\\
&\text{along $\chi$ on $S(z^{(0)},V,\epsilon)$}
\end{aligned}
\right\},
\end{equation}
where these $V$ and $\epsilon$ runs through $V(\zeta)$ and $\epsilon > 0$.
Hence, in particular, these shaves are $(\RP)^\ell$-conic.

\section{The functor of multi-specialization}\label{section:Multi-specialization}

In this section we extend the notion of multi-specialization of \cite{HP13} to our new setting. References are made to \cite{KS90} for the classical theory of specialization along a submanifold.

\subsection{A decomposition result on the normal deformation}\label{subsection:A decomposition result on the normal deformation}

\begin{df} Denote by $\op(\widetilde{X})$ the category of open subsets of $\widetilde{X}$, and let $Z$ be a subset of $\widetilde{X}$.

(i) %Let $Z$ be a subset of $X$.
We set $\RP_j Z=\mu_j(Z,\RP).$
If $U \in \op(\widetilde{X})$, then $\RP_j U \in \op(\widetilde{X})$ since $\mu_j$ is open for each $j=1,\ldots,\ell$.

(ii)  %Let $Z$ be a subset of $X$.
Let $J = \{j_1,\dots,j_k\} \subset \{1,\ldots,\ell\}$. We set $$\RP_JZ=\RP_{j_1}\cdots\RP_{j_k}Z=\mu_{j_1}(\cdots\mu_{j_k}(Z,\RP),\dots,\RP).$$
We set $(\RP)^\ell Z=\RP_{\{1,\ldots,\ell\}}Z=\mu(Z,(\RP)^\ell).$
If $U \in \op(\widetilde{X})$, then $\RP_J U \in \op(\widetilde{X})$ since $\mu_j$ is open for each $j \in \{1,\ldots,\ell\}$.

(iii) %Let $S$ be a subset of $X$.
We say that $Z$ is $(\RP)^\ell$-conic ($\ell$-conic for short) if
$Z=(\RP)^\ell Z$. In other words, $Z$ is invariant by the action of
$\mu_j$, $j=1,\dots,\ell$.

\end{df}

\begin{df}
(i)  We say that a subset $Z$ of $\widetilde{X}$ is $\RP_j$-connected if
$Z \cap \RP_j x$ is connected for each $x \in Z$.

(ii) We say that a subset $Z$ of $\widetilde{X}$ is $(\RP)^\ell$-connected if
there exists a permutation $\sigma:\{1,\ldots,\ell\}\to\{1,\ldots,\ell\}$ such that
\begin{equation*}
 \begin{cases}
 \text{$Z$ is $\RP_{\sigma(1)}$-connected,}\\
 \text{$\RP_{\sigma(1)}Z$ is $\RP_{\sigma(2)}$-connected,}\\
 \text{\ \ $\vdots$}\\
 \text{$\RP_{\sigma(1)}\cdots\RP_{\sigma(\ell-1)}Z$ is $\RP_{\sigma(\ell)}$-connected}.
  \end{cases}
\end{equation*}
\end{df}

The proof of the following is almost the same as that of Proposition 4.1.3 of \cite{KS90}, and we shall not repeat it.

\begin{prop}\label{m4.1.3} Let $V$ be an $(\RP)^\ell$-conic open subset of the zero section $S_\chi$.
\begin{itemize}
\item[(i)] Let $W$ be an open neighborhood of $V$ in $\widetilde{X}$, and let $U=\widetilde{p}(W \cap \Omega)$. Then $V \cap C_\chi(X \setminus U)=\emptyset$.
\item[(ii)] Let $U$ be an open subset of $X$ such that $V \cap C_\chi(X \setminus U)=\emptyset$. Then $\imin {\widetilde{p}}(U) \cup V$ is an open neighborhood of $V$ in $\overline{\Omega}$.
\end{itemize}
\end{prop}

%{\bf WE NEED THE CORRESPONDING VERSION OF LEMMA 4.1.3 OF SHEAVES OF MANIFOLDS}

Set, for $j=1,\dots,\ell$,
$$
S_{\chi,j} :=
\{(x^{(0)},\xi^{(1)},\dots,\xi^{(\ell)}) \in S_\chi ;\, \xi^{(j)} = 0\} \subset S_\chi
$$
and
$
\pi_j: S_\chi \to S_{\chi,j}
$
be the canonical projection. Let $V$ be an $(\RP)^\ell$-conic
subanalytic subset of $S_\chi$.
We introduce the following conditions Va.~and Vb.~of $V$ for each $j$.
\begin{enumerate}
\item[{\bf{Va.}}] $V$ does not intersect $S_{\chi,j}$.
\item[{\bf{Vb.}}] $\pi_j(V) \subset \pi_j(V \cap S_{\chi,j})$.
\end{enumerate}

Then we have the following proposition whose proof will be given in Appendix \ref{C}.
\begin{prop}\label{m4.1.4}
Let $V$ be an $(\RP)^\ell$-conic subanalytic subset of $S_\chi$ which
satisfies the condition either Va.~or Vb.~for each $j$.
Assume that $X=\mathbb{R}^n$ with coordinates $(x^{(0)},x^{(1)},\dots,x^{(m)})$ defined by $\chi$.
Then any subanalytic neighborhood $W$
of $V$ in $\widetilde{X}$ contains $\widetilde{W}$ open and subanalytic in
$\widetilde{X}$ such that:
\begin{equation} \label{mpitilde}
\begin{cases}
%\text{(i) the fibers of the map $\widetilde{p}:\widetilde{W} \cap \Omega \to X$ are connected,} \\
\text{(i) $\widetilde{W} \cap \Omega$ is $(\RP)^\ell$-connected,} \\
%\text{(ii) $\widetilde{p}(\widetilde{W} \cap \Omega)$ is subanalytic in $X$.}
\text{(ii) $\RP_1 \cdots \RP_{\ell}(\widetilde{W} \cap \Omega)=\imin {\widetilde{p}}(\widetilde{p}(W \cap \Omega))$ is subanalytic in $\widetilde{X}$.}
\end{cases}
\end{equation}
\end{prop}

Let $Z \subset S_\chi$ be closed conic and subanalytic.
%Let $Z$ be a closed sub-bundle of the vector bundle $S$ over a closed submanifold of $M$ which is $(\RP)^\ell$-conic.
%which is also $(\RP)^\ell$-conic,
%and let $V$ be an $(\RP)^\ell$-conic open subanalytic subset in $Z$.
%Then we can easily find a finite family $\{V^{(k)}\}$ of
%$(\RP)^\ell$-conic open subanalytic subsets in $Z$ with
%$\cup_k V^{(k)} = V$ such that every $V^{(k)}$ satisfies the condition
%either Va.~or Vb.~for each $j$.
We have the following corollary.
\begin{cor}\label{cor: Rpl neighborhood}
Let $V$ be an $(\RP)^\ell$-conic open subanalytic subset of $Z$.
Then there exists a locally finite family $\{V^{(\alpha)}\}_\alpha$
of $(\RP)^\ell$-conic open subanalytic subsets in $Z$
which satisfies the following conditions.
\begin{enumerate}
\item $V = \bigcup_{\alpha} V^{(\alpha)}$ and  each $V^{(\alpha)}$ satisfies the condition either Va.~or Vb.~for each $j$.
\item For any subanalytic open neighborhood $W$ of $V$ in $\widetilde{X}$,
	there exists a subanalytic open neighborhood $\widetilde{W}^{(\alpha)} \subset W$ of $V^{(\alpha)}$ for which
	(\ref{mpitilde}) of Proposition \ref{m4.1.4} holds.
	Furthermore $\{\widetilde{p}(\widetilde{W}^{(\alpha)} \cap \Omega)\}_\alpha$
	is a locally finite family of subanalytic open subsets of $X$.
\end{enumerate}
\end{cor}

\subsection{Multi-specialization}\label{subsection:Multi-specialization}

Let $k$ be a field and denote by $\mod(k_{X_{sa}})$ (resp. $D^b(k_{X_{sa}})$) the category (resp. bounded derived category) of sheaves on the subanalytic site $X_{sa}$.
For the theory of sheaves on subanalytic sites we refer to \cite{KS01,Pr08}. %For the theory of multi-specialization we refer to \cite{HP13}.
For classical sheaf theory we refer to \cite{KS90}. For an exposition on $(\RP)^\ell$-conic sheaves see the appendix.
%Most of the proofs of the results concerning the functor of multi-specialization can be obtained adapting the ones of \cite{HP13}. We will  omit them and refer to \cite{HP13}.

\begin{df}\label{def1.11} The multi-specialization along $\chi$ is the functor
\begin{eqnarray*}
\nu^{sa}_\chi \colon D^b(k_{X_{sa}})  \to  D^b(k_{S_{\chi sa}}) , &&
F \mapsto  \imin sR\Gamma_{\Omega_\chi}\imin p F.
\end{eqnarray*}
\end{df}
%\begin{oss}\label{stalk-nu}
It follows from the definition (see also the proof of Theorem 6.3 (i) of \cite{HP13}) that $\nu^{sa}_\chi F$ is conic, i.e., it belongs to $D^b_{(\RP)^\ell}(k_{S_{\chi sa}})$. The following result can be obtained adapting the proof of Lemma 6.1 of \cite{HP13}.

\begin{lem} Let $F \in D^b(k_{X_{sa}})$. There is a natural isomorphism
$$
\imin s R\Gamma_\Omega \imin p F \simeq s^!(\imin p F)_\Omega.
$$
\end{lem}

Thanks to Proposition \ref{m4.1.4} we can give a description
of the sections of the multi-specialization of
$F \in D^b(k_{X_{sa}})$. %The proof is similar to the one of Theorem 6.3 of \cite{HP13}.

\begin{teo}\label{sections of multispec} Let $V$ be a conic subanalytic open subset of
$S_\chi$ satisfying either condition Va. or Vb. for each $j$. Then:
\[H^j(V;\nu^{sa}_\chi F) \simeq \lind U H^j(U;F),\]
where $U$ ranges through the family of open subanalytic subsets of $X$ such that
$C_\chi(X \setminus U) \cap V=\emptyset$.
\end{teo}
\begin{proof}
Let $U \in \op(X_{sa})$ such that $V \cap C_\chi(X \setminus
U)=\emptyset$. We have the chain of morphisms
\begin{eqnarray*}
R\Gamma(U;F) & \to & R\Gamma(\imin p(U);\imin p F) \\
& \to & R\Gamma(\imin p(U) \cap \Omega;\imin p F) \\
& \to & R\Gamma(\imin {\widetilde{p}}(U) \cup V;R\Gamma_\Omega\imin p F) \\
& \to & R\Gamma(V;\nu^{sa}_\chi F)
\end{eqnarray*}
where the third arrow exists since $\imin {\widetilde{p}}(U) \cup
V$ is a neighborhood of $V$ in $\overline{\Omega}$ by Proposition
\ref{m4.1.3} (ii).

Let us show that it is an isomorphism.
%By Corollary \ref{cor: Rpl neighborhood} we may assume that $V$ satisfies the condition either Va. or Vb. for each $j$.
We have
\begin{eqnarray*}
H^k(V;\nu^{sa}_\chi F) & \simeq & \lind W H^k(W;R\Gamma_\Omega\imin p F) \\
& \simeq & \lind W H^k(W \cap \Omega;\imin p F),
\end{eqnarray*}
where $W$ ranges through the family of subanalytic open
neighborhoods of $V$ in $\widetilde{X}$. By Proposition \ref{m4.1.4},
we may assume that $W$ satisfies \eqref{mpitilde}.
%($\RP$-connectedness is equivalent to \eqref{mpitilde} (i)).
Since $\imin p F$ is $(\RP)^\ell$-conic, we have
\begin{eqnarray*}
H^k(W \cap \Omega;\imin p F) & \simeq & H^k(\imin p(p(W \cap
\Omega));\imin p F) \\
& \simeq & H^k(p(W \cap \Omega) \times \{(1)_\ell\};\imin p F) \\
& \simeq & H^k(p(W \cap \Omega);F),
\end{eqnarray*}
where $(1)_\ell=(1,\dots,1) \in \R^\ell$. The second isomorphism follows since every subanalytic
neighborhood of $p(W \cap \Omega) \times \{(1)_\ell\}$ contains an
$(\RP)^\ell$-connected subanalytic neighborhood (the proof is similar to
that of Proposition \ref{m4.1.4}). By Proposition \ref{m4.1.3} (i) we have that
$p(W \cap \Omega)$ ranges through the family of subanalytic open
subsets $U$ of $X$ such that $V \cap C_\chi(X \setminus
U)=\emptyset$ and we obtain the result.
\end{proof}

Remember that a sheaf $F \in \mod(k_{X_{sa}})$ is said to be quasi-injective if the restriction morphism $\Gamma(U;F) \to \Gamma(V;F)$ is surjective for each $U,V \in \op^c(X_{sa})$ with $U \supseteq V$. Since quasi-injective sheaves are $\Gamma(U;\cdot)$-acyclic for each $U \in \op(X_{sa})$ we get the following result.

\begin{cor} Let $F \in \mod(k_{X_{sa}})$ be quasi-injective. Then $F$ is $\nu^{sa}_\chi$-acyclic.
\end{cor}
\begin{proof}
It is enough to prove that $H^k(V;\nu^{sa}_\chi F)=0$, $k \neq 0$, on a basis for the topology of $S_{\chi sa}$. Being $\nu^{sa}_\chi F$ conic, we may assume that $V$ is conic as well. By Corollary \ref{cor: Rpl neighborhood} we may assume that $V$ satisfies the condition either Va. or Vb. for each $j$. Thanks to Theorem \ref{sections of multispec} it is enough to show that $H^k(U;F)=0$, $k \neq 0$ when $U \in \op(X_{sa})$. Since quasi-injective sheaves are $\Gamma(U;\cdot)$-acyclic for each $U \in \op(X_{sa})$ we get the result.
\end{proof}

Let $p=(x^{(0)};\,\xi^{(1)},\dots,\xi^{(\ell)}) \in S_\chi$,
let $\pi:S_\chi \to \bigcap_{j=1}^\ell M_j$ be the projection. Namely, $\pi(x^{(0)};\,\xi^{(1)},\dots,\xi^{(\ell)})=(x^{(0)},0^{(1)},\dots,0^{(\ell)})$.
Let $B_\epsilon$ be an open ball of radius $\epsilon>0$ with its center at
$\pi(p)$ and set
\[\operatorname{Cone}_\chi(p,\, \epsilon):=\{V \cap B_\epsilon;\, V \in C(\xi,F^q)\}.\] Applying the functor $\imin \rho \colon D^b(k_{S_{\chi sa}}) \to D^b(k_{S_\chi})$ (see \cite{Pr08} for details) and Theorem \ref{sections of multispec} (see also Corollary 6.5 of \cite{HP13}) we can calculate the fibers at $p \in S_\chi$ which are given by
\begin{equation}\label{eq:fibers}
(H^j\imin \rho \nu^{sa}_\chi F)_p \simeq \lind {W} H^j(W;F),
\end{equation}
where $W$ ranges through the family $\operatorname{Cone}_\chi(p,\,\epsilon)$
 for $\epsilon > 0$.\\

From Proposition \ref{prop3} we can deduce a condition for the restriction of multi-specialization. Set $\chi_A=\{M_1,\dots,M_\ell\}$, $\chi_B=\chi_A \cup \{M_{\ell+1}\}$. We follow the notations of Section \ref{section:Multi-cones for a general case}. We assume that the first $\ell$ (resp. $\ell+1$) columns of the matrix $A$ (resp. $B$) are linearly independent. We omit the coordinate $x^{(0)}$ to lighten notations. Given $\xi=(\xi^{(1)},\dots,\xi^{(m)}) \in \R^n$
we will call $p_A \in S_A$, the point $p_A=(\xi^{(1)},\dots,\xi^{(m)};\,0 \dots,0) \in \R^{n+\ell}$, and $p_B \in S_B$, the point $p_B=(\xi^{(1)},\dots,\xi^{(m)};\,0 \dots,0) \in \R^{n+\ell +1}$.

\begin{cor} Let $F \in D^b(k_{X_{sa}})$. Let $K \subseteq \{1,\dots,m\}$. Suppose that on $S_K=\{\xi^{(k)}=0,\,k \in K\}$ conditions \eqref{3} are satisfied.
Then we have
$$
(\nu^{sa}_{\chi_A} F)|_{S_K}
\simeq (\nu^{sa}_{\chi_B}F)|_{S_K}.
$$
\end{cor}

\subsection{Operations}

Let $X$ and $Y$ be real analytic manifolds
and let
$\chi^M=\{M_j\}_{j=1}^\ell$, $\chi^N=\{N_j\}_{j=1}^\ell$ be
families of smooth closed submanifolds of $X$ and $Y$ %\eqref{eq:geometric-condition}
respectively.
Let $f:X \to Y$ be a morphism of real analytic manifolds such that
$f(M_j) \subseteq N_j$, $j=1,\ldots,\ell$.

\begin{oss} In order to study more general situations, we will assume $X=Y=\R^n$ and consider the local model.
Even if we work locally, the theory below completely covers the manifold case.
\end{oss}

\begin{oss} Let $\comD$ be a positive rational number so that $a_{jk} \in \mathbb{Z}$ and all the $a_{jk}$ have no common divisors. Let $\mu_{\intA}$ be the action on $X$ associated with the matrix $\intA$. Since multicones associated with the action $\mu$ and those with the action $\mu_{\intA}$ determine equivalent families, we may assume from the beginning that all the entries of $A_\chi$ are non-negative integers and $\comD = 1$.
\end{oss}

We want to extend $f$ to a morphism $\widetilde{f}:\widetilde{X} \to \widetilde{Y}$. This is done by repeatedly employing the usual construction
of a morphism between normal deformations,
i.e. we extend $f$ to  $\widetilde{f}_1:\widetilde{X}_{M_1} \to \widetilde{Y}_{N_1}$,
then we extend $\widetilde{f}_1$ to $\widetilde{f}_{1,2}:\widetilde{X}_{M_1,M_2} \to \widetilde{Y}_{N_1,N_2}$.
Then we can define recursively $\widetilde{f}:\widetilde{X} \to \widetilde{Y}$
by extending the morphism $\widetilde{f}_{1,\ldots,\ell-1}:\widetilde{X}_{M_1,\ldots,M_{\ell-1}} \to \widetilde{Y}_{N_1,\ldots,N_{\ell-1}}$
to the normal deformations with respect to $M_\ell$ and $N_\ell$ respectively. We also denote by $S^M$ (resp. $S^N$) the zero section of
$\widetilde{X}$ (resp. $\widetilde{Y}$). In this general case we have to put some conditions to assure that $\widetilde{f}$ is well defined on the zero section.

In a local coordinate system set
\begin{eqnarray*}
x & = & (x^{(0)},x^{(1)},\ldots,x^{(m_X)}), \\
y & = & (y^{(0)},y^{(1)},\ldots,y^{(m_Y)}), \\
f(x) & = & (f^{(0)}(x),f^{(1)}(x),\ldots,f^{(m_Y)}(x)).
\end{eqnarray*}
Hereafter we omit the coordinates blocks $x^{(0)}$ and $y^{(0)}$ for simplicity.
Outside $\{t_1t_2\ldots t_\ell = 0\}$ in $\widetilde{X}$,
the morphism $\widetilde{f}:\widetilde{X}\to\widetilde{Y}$ $$
(x_1,\dots,x_n,t_1,\dots,t_\ell) \to
(y_1,\dots,y_m,\hat{t}_1,\dots,\hat{t}_\ell)
$$
is explicitly written in the form
$$
\left\{
\begin{aligned}
&\varphi_k^N(\hat{t})y^{(k)}=f^{(k)}(\varphi_1^M(t)x^{(1)},\ldots,\varphi_{m_X}^M(t)x^{(m_X)})
\qquad (k = 1,\dots,m_Y),\\
&\hat{t}_j = t_j \qquad (j = 1,\dots,\ell),
\end{aligned}
\right.
$$
which naturally extends to the whole $\widetilde{X}$ if the limit
$$
\lim_{t \to 0}{f^{(k)}(\varphi_1^M(t)x^{(1)},\ldots,\varphi_{m_X}^M(t)x^{(m_X)}) \over \varphi_k^N(t)}
$$
exists for each $k=1,\dots,m_Y$. In order to understand the restriction of $\widetilde{f}$ to the whole $\widetilde{X}$,
we set
$$
I=\bigcup_{j \in J^N_k} I^M_j,\,\, M_I=\bigcap_{j \in J^N_k} M_j,
$$
where $J^N_k$ denotes the set $\{j \in\{1,\dots,\ell\},\, \hat{I}^N_k \subset I^N_j\}$.
By expanding $f^{(k)}(x)$ along the submanifold $M_I$,
we obtain
$$
f^{(k)}(x) =
\sum_{i \in I}\left.\frac{\partial f^{(k)}}{\partial x^{(i)}}\right|_{M_I} x^{(i)}
+ \frac{1}{2}\sum_{i_1,i_2 \in I}
\left.
\frac{\partial^2 f^{(k)}}{\partial x^{(i_1)}\partial x^{(i_2)}}\right|_{M_I}
x^{(i_1)}x^{(i_2)}
+ \dots,
$$
as $f_k\vert_{M_I} = 0$ holds.
Then we get, on $t_1\dots t_\ell \ne 0$,
\begin{equation*}{\label{eq:morphism-expansions}}
y_k =
\sum_{i \in I}\left.\frac{\partial f^{(k)}}{\partial x^{(i)}}\right|_{M_I}
\frac{\varphi^M_i(t)}{\varphi^N_k(t)}x^{(i)}
+ \frac{1}{2}\sum_{i_1,i_2 \in I}
\left.
\frac{\partial^2 f^{(k)}}{\partial x^{(i_1)}\partial x^{(i_2)}}\right|_{M_I}
\frac{\varphi^M_{i_1}(t)\varphi^M_{i_2}(t)}{\varphi^N_k(t)}x^{(i_1)}x^{(i_2)}
+ \dots.
\end{equation*}
Note that the limit for $t \to 0$ exists under the condition
\begin{equation} \label{eq:cond existence}
\left.{\partial f^{(k)} \over \partial x^{(i_1)} \cdots \partial x^{(i_p)}}\right|_{M_I}=0 \ \ \ \ \mbox{if } \varphi^N_k \nmid \varphi^M_{i_1}\cdots\varphi^M_{i_p}.
\end{equation}

\begin{oss}  %Let $J^N_k=\{j \in\{1,\dots,\ell\},\, \hat{I}^N_k \subset I^N_j\}$.
If there exists $j \in J^N_k$ such that
$\{i_1,\ldots,i_p\} \cap I^M_j=\emptyset$, as
$f^{(k)} \vert_{M_j} = 0$, we have
$$
\left.{\partial f^{(k)} \over \partial x^{(i_1)} \cdots \partial x^{(i_p)}}\right|_{M_I}=0.
$$
If the matrix $A$ is associated with $\chi$, i.e., its
entries are either $0$ or $1$,
\begin{equation*}
\mbox{$\varphi^N_k \nmid \varphi^M_{i_1}\cdots\varphi^M_{i_p}$ $\Rightarrow$ $\exists\,j \in J^N_k$ such that
$\{i_1,\ldots,i_p\} \cap I^M_j=\emptyset$. }
\end{equation*}
Then \eqref{eq:cond existence} always holds.
\end{oss}

Suppose that $\varphi^N_k \mid \varphi^M_{i_1}\cdots\varphi^M_{i_p}$.
\begin{itemize}
\item[(i)] If there exists $j \in J^N_k$ such that
$\{i_1,\ldots,i_p\} \cap I^M_j=\emptyset$, we have
$$
\left.{\partial f^{(k)} \over \partial x^{(i_1)} \cdots \partial x^{(i_p)}}\right|_{M_I}=0.
$$
\item[(ii)] If $\alpha^N_{jk}<\alpha^M_{ji_1}+\cdots+\alpha^M_{ji_p}$, for some $j \in \{1,\dots,\ell\}$, we have
$$
\lim_{t\to 0}\frac{\varphi^M_{i_1}(t)\cdots\varphi^M_{i_p}(t)}{\varphi^N_k(t)} \to 0.
$$
\item[(iii)] If $\alpha^N_{jk}=\alpha^M_{ji_1}+\cdots+\alpha^M_{ji_p}$, for all $j \in \{1,\dots,\ell\}$, we have
$$
\lim_{t\to 0}\frac{\varphi^M_{i_1}(t)\cdots\varphi^M_{i_p}(t)}{\varphi^N_k(t)} \to 1.
$$
\end{itemize}
By the above observations, assuming \eqref{eq:cond existence}, in order to compute the restriction to the zero section we are interested in sequences $P=\{i_1,\dots,i_{\sharp P}\}$ ($i_p \in \{1,\dots,m_X\}$ ($p=1,\dots,\sharp P$)) of indices of finite length, such that
\begin{equation}\label{eq:P}
\alpha^N_{jk}=\sum_{p = 1}^{\sharp P}\alpha^M_{ji_p}, \qquad j=1,\dots,\ell.
\end{equation}

Let us consider the family $\mathcal{F}_k$ of sequences $P$ of indices of finite length
satisfying \eqref{eq:P}.
%of subsets of $\{1,\dots,m_X\}$ defined as follows:
%$$
%\mathcal{F}_k=\left\{P=\{i_1,\dots,i_{\sharp P}\},\, i_p \in \{1,\dots,m_X\},\, p=1,\dots,\sharp P,\, \alpha^N_{jk}=\sum_{i \in P}\alpha^M_{ji},\, j=1,\dots,\ell\right\}.
%$$
%From the above observations,
The morphism $\widetilde{f}$ is
described by, on $S \subset \widetilde{X}$,
         $$
	 y_k = \sum_{P \in \mathcal{F}_k} \left.\frac{1}{\sharp P!}{\partial^{\sharp P} f^{(k)} \over \partial x^{(i_1)} \cdots \partial x^{(i_{\sharp P})}}\right|_{M}x^{(i_1)} \cdots x^{(i_{\sharp P})}
 \qquad (k=1,\dots,m_Y).
         $$
Here $M:=M_1 \cap \dots \cap M_\ell$.

\begin{oss}\label{oss:order}
The order is important in $\mathcal{F}_k$
as it is the set of sequences.
For example, $\{1,2\} \neq \{2,1\}$.
\end{oss}

\begin{oss} When $m_X=\ell$, i.e. $\chi^M$ is of normal type, the family $\mathcal{F}_k$ consists of at most one element. This follows from the invertibility of the matrix associated to $\chi^M$.
\end{oss}

\begin{es} Let $X=Y=\R^3$, $M_1=\{x_1=x_2=0\}$, $M_2=\{x_2=x_3=0\}$, $N_1=0$, $N_2=\{x_2=x_3=0\}$, $f(x_1,x_2,x_3)=(x_1,x_1x_3+x_2,x_1x_3)$. Then
\begin{eqnarray*}
\varphi^M_1=t_1, & \varphi^M_2=t_1t_2, & \varphi^M_3=t_2, \\
\varphi^N_1=t_1, & \varphi^N_2=t_1t_2, & \varphi^N_3=t_1t_2.
\end{eqnarray*}
Condition \eqref{eq:cond existence} is satisfied and we have
\begin{eqnarray*}
y_1 & = & \left.{\partial f_1 \over \partial x_1}\right|_{0}x_1=x_1, \\
y_2 & = & \left.{\partial f_2 \over \partial x_2}\right|_{0}x_2+{1 \over 2}\left.{\partial^2 f_2 \over \partial x_1\partial x_3}\right|_{0}x_1x_3+{1 \over 2}\left.{\partial^2 f_2 \over \partial x_3\partial x_1}\right|_{0}x_3x_1=x_2+x_1x_3, \\
y_3 & = & {1 \over 2}\left.{\partial^2 f_3 \over \partial x_1\partial x_3}\right|_{0}x_1x_3+{1 \over 2}\left.{\partial^2 f_2 \over \partial x_3\partial x_1}\right|_{0}x_3x_1=x_1x_3.
\end{eqnarray*}
Note that, as stressed in Remark \ref{oss:order}, the order of the elements $P$ in $\mathcal{F}_2$ and $\mathcal{F}_3$ is important.
\end{es}

\begin{es} \label{es: cusp morphism} Let $X=Y=\R^2$, $M_1=\{x_1=0\}$, $M_2=\{x_2=0\}$, $N_1=N_2=\{0\}$,
$$
A_{\chi^M}=
\left(
\begin{array}{cc}
1 & 0 \\
0 & 1
\end{array}
\right),
\ \ \ \
A_{\chi^N}=
\left(
\begin{array}{cc}
3 & 2 \\
1 & 1
\end{array}
\right)
$$
$f(x_1,x_2)=(x_1^3x_2,x_1^2x_2)$. Then
\begin{eqnarray*}
\varphi^M_1=t_1, & \varphi^M_2=t_2, \\
\varphi^N_1=t_1^3t_2, & \varphi^N_2=t_1^2t_2. \\
\end{eqnarray*}
Condition \eqref{eq:cond existence} is satisfied and we have
\begin{eqnarray*}
y_1 & = & {1 \over 4!} \cdot 4 \cdot \left.{\partial^4 f_1 \over \partial x_1^3\partial x_2}\right|_{0}x_1^3x_2=x_1^3x_2, \\
y_2 & = & {1 \over 3!} \cdot 3 \cdot \left.{\partial^3 f_1 \over \partial x_1^2\partial x_2}\right|_{0}x_1^2x_2=x_1^2x_2.
\end{eqnarray*}
Note that we multiplied the first (resp. second) line by 4 (resp. 3) taking in account the order of the elements $P$ in $\mathcal{F}_1$ (resp. $\mathcal{F}_2$).
\end{es}

%Let $f:X \to Y$ be a morphism of real analytic manifolds, $\chi^M=\{M_1,\ldots,M_\ell\}$, $\chi^N=\{N_1,\ldots,N_\ell\}$ two families of closed analytic submanifolds of $X$ and $Y$ respectively (in the local model, set $X=Y=\R^n$).
Suppose that $f$ satisfies \eqref{eq:cond existence} and $f(M_i) \subseteq N_i$, $i=1,\ldots,\ell$.
 We call $T_\chi f$ the induced map on the zero section.
Let $J = \{j_1,\dots,j_k\}$ be a non-empty subset in $\{1,\ldots,\ell\}$
and set $\chi^M_J=\{M_{j_1},\ldots,M_{j_k}\}$. The map $T_{\chi_J}f$ denotes the restriction of $\widetilde{f}$ to $\{t_{j_k}=0,\,j_k\in J\}$. In the following we will denote with the same symbol $C_{\chi^M_J}(Z)$ the normal cone with respect to $\chi^M_J$ and its inverse image via the map $\widetilde{X} \to \widetilde{X}_{M_{j_1},\dots,M_{j_k}}$.

\begin{prop} \label{prop: direct image} Let $F \in D^b(k_{X_{sa}})$.

\begin{itemize}
\item[(i)] There exists a commutative diagram of canonical
morphisms
$$
\xymatrix{R(T_\chi f)_{!!}\nu^{sa}_{\chi^M} F  \ar[d] \ar[r] & \nu^{sa}_{\chi^N}Rf_{!!}F \ar[d] \\
R(T_\chi f)_*\nu^{sa}_{\chi^M} F & \nu^{sa}_{\chi^N}Rf_*F. \ar[l]}
$$
\item[(ii)] Moreover if $f$ is proper over $\supp F$ and
	$T_{\chi_J}f$ is proper over $C_{\chi^M_J}(\supp F)$
	for each (non-empty) $J = \{j_1,\dots,j_k\}$,
	and if $\supp F \cap \imin f(N_j) \subseteq M_j$, $j \in \{1,\ldots,\ell\}$, then
%$C_{M_{i_1}\dots M_{i_k}}(\supp F) \to \oplus_{1\leq j\leq k} T_{M'_{i_j}}\iota(M'_{i_j})$
%for each $J = \{j_1,\dots,j_k\}$ are proper,
the above morphisms are isomorphisms.
\end{itemize}
\end{prop}
\begin{proof} (i) The existence of the arrows is done as in \cite{KS90} Proposition 4.2.4.

	(ii) If $\widetilde{f}$ is proper over $\overline{\widetilde{p}^{-1}(\supp F)}$, then all the morphisms are isomorphisms. We have to prove that for a closed subset $Z$ of $X$, the restriction of $\widetilde{f}$ to $\overline{\imin {\widetilde{p}}(Z)}$ is proper, if $Z \to Y$ and $C_{\chi^M_J}(Z) \to
Y$
%$C_{M_{j_1}\dots M_{j_k}}(Z) \to \oplus_{1\leq j\leq k}
%{M'_{i_j} \neq \iota(M'_{i_j})} T_{M'_{i_j}}\iota(M'_{i_j})$
for $k\leq \ell$ are proper, and if $Z \cap \imin f(N_j) \subseteq M_j$, $j \in \{1,\ldots,\ell\}$. We argue by induction on $\sharp \chi^M$. If $\sharp \chi^M=0$ the proof is trivial. Suppose it is true for $\sharp\chi^M\leq\ell-1$. It follows from the hypothesis that the fibers of $\widetilde{f}$ restricted to $\overline{\widetilde{p}^{-1}(Z)}$ are compact
(if $t_{j_1}=\cdots=t_{j_k}=0$, this is a consequence of the fact that $C_{\chi^M_J}(Z) \to
Y$
%$C_{M_{j_1}\dots M_{j_k}}(Z) \to \oplus_{1\leq j\leq k}
%{M'_{i_j} \neq \iota(M'_{i_j})} T_{M'_{i_j}}\iota(M'_{i_j})$
 is proper). Then it remains to prove that it is a closed map. Let $\{u_n\}_{n\in\N}$ be a sequence in $\overline{\imin{\widetilde{p}}(Z)}$ and suppose that $\{\widetilde{f}(u_n)\}_{n\in\N}$ converges. We shall find a convergent subsequence of $\{u_n\}_{n\in\N}$.
 We may also assume that $\{\widetilde{p}(u_n)\}_{n\in\N}$ converges as $f$ is proper over $Z$. The map $\imin{\widetilde{p}}(Z) \setminus \{t_1=\cdots=t_\ell=0\} \to \widetilde{Y} \setminus \{t_1=\cdots=t_\ell=0\}$ is proper. Indeed, let $K$ be a compact subset of $\widetilde{Y} \setminus \{t_1=\cdots=t_\ell=0\}$, and
%suppose without loss of generality
reduce to the case that $K$ is contained in $\{c \leq t_j \leq d\} \subset \widetilde{Y}$, $c,d>0$, for some $j \in \{1,\dots,\ell\}$. Suppose without loss of generality that $j=1$. Then $K_1:=\widetilde{p}_{N_1}(K)$ is a compact subset of $\widetilde{Y}_{N_2\ldots N_\ell}$. Let us identify $K_1$ with $\imin{\widetilde{p}_{N_1}}(K_1) \cap \{t_1=1\}$. Then $\imin{\widetilde{f}}(K_1)$ is compact by the induction hypothesis. Hence $\imin{\widetilde{f}}(K) \subseteq \mu_1(\imin{\widetilde{f}}(K_1),[c,d])$ is compact since it is closed and contained in a compact subset. %(here we identified $\imin{\widetilde{f}}(K_1)$ with $\imin{\widetilde{p}_{M_1}}(\imin{\widetilde{f}}(K_1)) \cap \{t_1=1\}$).
We may assume that $\{\widetilde{f}(u_n)\}_{n\in\N}$ converges to a point of $\{t_1=\cdots=t_\ell=0\}$. Then $\{\widetilde{p}(u_n)\}_{n\in\N}$ converges to a point of $Z \cap \imin f(N_1 \cap \cdots \cap N_\ell) \subseteq M_1 \cap \cdots \cap M_\ell$.

\

Taking local coordinates systems of $X$ and $Y$, let
$$
u_n=(x_n^{(1)},\ldots,x_n^{(m)};\, t_{1n},\ldots,t_{\ell n}),  \quad t_{jn}>0,\,\, j=1,\dots,\ell.
$$
Then $t_{jn} \to 0$, $j=1,\dots,\ell$ and $\varphi^M_i(t_n)x_n^{(i)}\to 0$, $i=1,\ldots,m$. It is enough to show that $\{|x_n^{(i)}|\}_{n\in\N}$ is bounded for each $i=1,\dots,m$. We argue by contradiction. Up to make a permutation of $\{1,\dots,m\}$, we may suppose without loss of generality that $|x_n^{(i)}|\to+\infty$ for $0 < i \leq m_0 \leq m$.\\

Suppose we can find a non infinitesimal bounded sequence $\{\widetilde{u}_n\}_{n\in\N} \in \imin{\widetilde{p}}(Z)$ such that $\{\widetilde{f}(\widetilde{u}_n)\}_{n\in\N}$ converges to 0. Up to extract a subsequence we may assume that $\{u_n\}_{n\in\N}$ converges to a non zero vector $v$. Then $\widetilde{f}(w)=0$ for each $w \in (\RP)^\ell v$ which contradicts the fact that the fibers of $\overline{\imin{\widetilde{p}}(Z)} \to \widetilde{Y}$ are compact. \\

Suppose that there exist sequences $\{c_{jn}\}_{n \in \N}$, $j=1,\dots,\ell$, satisfying the following properties:

\begin{itemize}
\item[(i)] $c_{jn} \to +\infty$, $j=1,\dots,\ell,$
\item[(ii)] $\displaystyle{\prod_{j=1}^\ell c_{jn}^{\alpha_{m_0j}}=|x^{(m_0)}_n|}$,
\item[(iii)] $t_{jn}c_{jn} \to 0$, $j=1,\dots,\ell$.
\end{itemize}
Let us construct the sequence $\{\overline{u}_n\}_{n\in\N}$, $\overline{u}_n=(\overline{x}^{(1)}_n,\dots,\overline{x}^{(m)}_n,\overline{t}_{1n},\dots,\overline{t}_{\ell n})$ as follows:
$$
\overline{x}^{(i)}_n={x^{(i)}_n \over \prod_{j=1}^\ell c_{jn}^{\alpha_{ij}}}, \ \ \ \ \overline{t}_{jn}=t_{jn}c_{jn}.
$$
By (i) $\overline{x}^{(i)}_n$ is bounded if $x^{(i)}_n$ is bounded, $i=1,\dots,m$, by (ii) $\overline{x}^{(m_0)}_n$ is bounded non infinitesimal and by (iii) $\overline{t}_{jn} \to 0$, $j=1,\dots,\ell$. Moreover $\varphi^M_i(\overline{t}_n)\overline{x}_n^{(i)}=\varphi^M_i(t_n)x_n^{(i)}$, $i=1,\dots,m$. Applying this process (at most) $m_0$ times we may find a non infinitesimal bounded sequence $\{\widetilde{u}_n\}_{n\in\N} \in \imin{\widetilde{p}}(Z)$, $\widetilde{u}_n=(\widetilde{x}^{(1)}_n,\dots,\widetilde{x}^{(m)}_n,\widetilde{t}_{1n},\dots,\widetilde{t}_{\ell n})$  such that
$$
{f_k(\varphi^M_1(\widetilde{t}_n)\widetilde{x}^{(1)}_n,\dots,\varphi^M_1(\widetilde{t}_n)\widetilde{x}^{(m)}_n) \over \varphi^N_k(\widetilde{t}_n)} = {f_k(\varphi^M_1(t_n)x_n^{(1)},\dots,\varphi^M_n(t_n)x_n^{(n)}) \over \varphi^N_k(t_n)d_n} \to 0.
$$
This is because $\{\widetilde{f}(u_n)\}_{n\in\N}$ converges and, by construction of $\widetilde{t}_n$,
$$
d_n:={\varphi^N_k(\widetilde{t}_n) \over \varphi^N_k(t_n)} \to +\infty.
$$

So we are reduced to find sequences $\{c_{jn}\}_{n\in\N}$, $j=1,\dots,\ell$ satisfying (i)-(iii), given $u_n=(x_n^{(1)},\ldots,x_n^{(m)},t_{1n},\ldots,t_{\ell n})$, $t_{jn}>0$, $j=1,\dots,\ell$ with $t_{jn} \to 0$ ($j=1,\dots,\ell$), $\varphi^M_i(t_n)x_n^{(i)}\to 0$ ($i=1,\ldots,m$) and $|x_n^{(m_0)}|\to+\infty$ for $1 \leq m_0 \leq m$. Assume without loss of generality that $m_0=1$. There exists $j_0 \in \{1,\dots,\ell\}$ such that
$$
\prod_{j \leq j_0}t_{jn}^{\alpha_{1j}}|x^{(1)}_n| \to 0, \ \ \ \ \prod_{j < j_0}t_{jn}^{\alpha_{1j}}|x^{(1)}_n| \to +\infty.
$$
Let us construct the required sequence as follows:
\begin{eqnarray*}
c_{jn} & = & {1 \over t_{jn}\epsilon_n^{1 \over \alpha_{1j}}} \ \ j<j_0, \\
c_{j_0n}& = & \left(\prod_{j<j_0}t_{jn}^{\alpha_{1j}}|x^{(1)}_n|\epsilon_n^{2j_0-m-1}\right)^{1 \over \alpha_{1j_0}}, \\
c_{jn} & = & \epsilon_n^{1 \over \alpha_{1j}} \ \ j>j_0.
\end{eqnarray*}
Here $\{\epsilon_n\}_{n\in\N}$ is a sequence, $\epsilon_n \to +\infty$ slowly enough to satisfy
$c_{jn} \to +\infty$ if $j \leq j_0$ and $t_{jn}c_{jn} \to 0$ if $j \geq j_0$. One checks easily that (i)-(iii) are satisfied.
\end{proof}

\begin{prop} \label{prop: inverse image} Let $F \in D^b(k_{Y_{sa}})$.

\begin{itemize}
\item[(i)] There exists a commutative diagram of canonical
morphisms
$$
\xymatrix{\omega_{S_X/S_Y}\otimes\imin {(T_\chi f)}\nu^{sa}_{\chi^N} F  \ar[d] \ar[r] & \nu^{sa}_{\chi^M}(\omega_{S_X/S_Y}\otimes \imin f F) \ar[d] \\
T_{\chi} f^!\nu^{sa}_{\chi^N} F & \nu^{sa}_{\chi^M}f^!F. \ar[l]}
$$
\item[(ii)]
%Let $A \subseteq \{1,\ldots,\ell\}$ and set $\chi_A=\{M_{j_1},\ldots,M_{j_k}\}$, $j_k \in A$. Let $T_{\chi_A}f$ denote the restriction of $\widetilde{f}$ to $\{t_{j_k}=0,\,j_k\in A\}$.
The above morphisms are isomorphisms on the open set where $T_\chi f$ is smooth. %$T_{\chi_J} f$ is smooth for each $J \subseteq \{1,\ldots,\ell\}$.
\end{itemize}
\end{prop}
\begin{proof}
(i) The existence of the arrows is done as in \cite{KS90} Proposition 4.2.5. When (ii) is satisfied the function $\widetilde{f}$ is smooth at any point of the %boundary of $\Omega$
zero section and all the above morphisms become isomorphisms.
\end{proof}

\subsection{Multi-specialization and asymptotic expansions}\label{subsection:Multi-specialization and asymptotic expansions}

We give a functorial construction of the sheaf of multi-asymptotically  developable
functions using the sheaf of Whitney $\C^\infty$-functions of \cite{KS01}. Let $X$ be a real analytic manifold ($X=\R^n$ in the local model).
As usual, given $F \in D^b(\CC_{X})$ we set $D'F=\rh(F,\CC_X)$. Remember that an open subset $U$ of $X$ is locally cohomologically trivial (l.c.t. for short) if $D'\CC_U \simeq \CC_{\overline{U}}$. Let $\wtens$ denote the Whitney tensor product of \cite{KS96}.

\begin{df} Let $F \in \mod_{\rc}(\CC_X)$ and let $U \in \op(X_{sa})$. We define the presheaf $\CW_{X|F}$ as follows:
$$U \mapsto \Gamma(X;H^0D'\CC_U \otimes F \wtens \C^\infty_X).$$
\end{df}
Let $U,V \in \op(X_{sa})$, and consider the exact sequence
$$\exs{\CC_{U\cap V}}{\CC_U \oplus \CC_V}{\CC_{U\cup V}},$$
applying the functor $\ho(\cdot,\CC_X)=H^0D'(\cdot)$ we obtain
$$\lexs{H^0D'\CC_{U\cup V}}{H^0D'\CC_U \oplus H^0D'\CC_V}{H^0D'\CC_{U\cap V}},$$
applying the exact functors $\cdot \otimes F$, $\cdot\wtens \C^\infty_X$ and taking global sections
we obtain
$$\lexs{\CW_{X|F}(U\cup V)}{\CW_{X|F}(U) \oplus \CW_{X|F}(V)}{\CW_{X|F}(U\cap V)}.$$
This implies that $\CW_{X|F}$ is a sheaf on $X_{sa}$. Moreover if $U \in \op(X_{sa})$ is l.c.t., the morphism $\Gamma(X;\CW_{X|F}) \to \Gamma(U;\CW_{X|F})$ is surjective and $R\Gamma(U;\CW_{X|F})$ is concentrated in degree zero. Let $\exs{F}{G}{H}$ be an exact sequence in $\mod_{\rc}(\CC_X)$, we obtain an exact sequence in $\mod(\CC_{X_{sa}})$
\begin{equation}\label{exsFGH}
\exs{\CW_{X|F}}{\CW_{X|G}}{\CW_{X|H}}.
\end{equation}

We can easily extend the sheaf $\CW_{X|F}$ to the case of $F \in D^b_{\rc}(\CC_X)$, taking a finite resolution of $F$ consisting of locally finite sums $\oplus \CC_V$ with $V$ l.c.t. in $\op^c(X_{sa})$. In fact, the sheaves $\CW_{X|\oplus \CC_V}$ form a complex quasi-isomorphic to $\CW_{X|F}$ consisting of acyclic objects with respect to $\Gamma(U;\cdot)$, where $U$ is l.c.t. in $\op^c(X_{sa})$.

As in the case of Whitney $\C^\infty$-functions one can prove that, if $G \in D^b_{\rc}(\CC_X)$ one has
$$\imin \rho \rh(G,\CW_{F|X}) \simeq D'G \otimes F \wtens \C^\infty_X.$$

%\begin{oss} The sheaf $\CW_{M|F}$ is different from the sheaf $\ho(F,\CWM)$. For example, let $p \in M$ and let $U \in %\op(M_{sa})$ l.c.t. with $p \in \partial U$. Then $\Gamma(U;\ho(\CC_p,\CWM)=0$ but %$\Gamma(U;\CW_{M|\CC_p})=\Gamma(M;\CC_p \wtens \C^\infty_M)$.
%\end{oss}

\begin{es} Setting $F=\CC_X$ we obtain the sheaf of Whitney $\C^\infty$-functions. Indeed, let $U$ be a l.c.t. subanalytic open subset of $X$. Then $\Gamma(U;\CW_{X_\R}) \simeq \Gamma(X; \CC_{\overline{U}} \wtens C^\infty_{X_\R})$ is nothing
but the set of $\C^\infty$-Whitney jets on $\overline{U}$. Let $Z$ be a closed subanalytic subset of $X$. Similarly one checks that $\CW_{X|\CC_{X\setminus Z}}$ is the sheaf of Whitney $\C^\infty$-functions vanishing on $Z$ with all their derivatives.
\end{es}

\begin{nt} Let $S$ be a locally closed subanalytic subset of $X$. We set for short $\CW_{X|S}$ instead of $\CW_{X|\CC_S}$.
\end{nt}

Let $\chi=\{M_1,\ldots,M_\ell\}$ be a family of closed analytic submanifolds of $X$, set $\ckM= \bigcup_{k=1}^\ell M_k$ and consider $\widetilde{X}$.
Set $F=\CC_{X\setminus \ckM}$, $G=\CC_X$, $H=\CC_\ckM$ in \eqref{exsFGH}. The exact sequence
$$\exs{\CW_{X|X\setminus \ckM}}{\CW_X}{\CW_{X|\ckM}}$$
induces an exact sequence
\begin{equation}\label{asym}
\exs{\imin \rho\nu^{sa}_\chi\CW_{X|X\setminus \ckM}}{\imin \rho\nu^{sa}_\chi\CW_X}{\imin \rho \nu^{sa}_\chi \CW_{X|\ckM}},
\end{equation}
where the surjective arrow is the map which associates to a function its asymptotic expansion. In fact, let $V$ be a l.c.t. conic subanalytic subset of the zero section of $\widetilde{X}$ and $U \in \op(X_{sa})$ such that
$C_\chi(X \setminus U) \cap V=\emptyset$, then we can find a l.c.t. $U'\subset U$ satisfying the same property. \\

Let $X$ be a complex manifold ($X=\CC^n$ in the local model) and let $X_\R$ denote the underlying real analytic manifold
of $X$. Let $\chi=\{Z_1,\ldots,Z_\ell\}$ be a family of complex submanifolds of $X$ and set $\ckZ= \bigcup_{k=1}^\ell Z_k$. Let $F \in D^b_{\rc}(\CC_X)$. We denote by $\OW_{X|F}$ the sheaf defined as follows:
$$
\OW_{X|F} := \rh_{\rho_!\D_{\overline{X}}}(\rho_!\OO_{\overline{X}},\CW_{X_\R|F}).
 $$
Let $\exs{F}{G}{H}$ be an exact sequence in $\mod_{\rc}(\CC_X)$. Then the exact sequence \eqref{exsFGH} gives rise to the distinguished triangle
\begin{equation}\label{dtFGH}
\dt{\OW_{X|F}}{\OW_{X|G}}{\OW_{X|H}}.
\end{equation}

Setting $F=\CC_{X\setminus \ckZ}$, $G=\CC_X$, $H=\CC_\ckZ$ in \eqref{dtFGH} and applying the functor of specialization, we have the distinguished triangle
\begin{equation}\label{dtColin}
\dt{\imin \rho \nu^{sa}_\chi\OW_{X|X \setminus \ckZ}}{\imin \rho \nu^{sa}_\chi\OWX}{\imin \rho \nu^{sa}_\chi \OW_{X|\ckZ}}.
\end{equation}

Let $U$ be an open l.c.t. subanalytic subset in $X$.
Then $H^0(U;\OW_X) \simeq H^0(X; \CC_{\overline{U}} \wtens \OO_X)$ consists of
$\C^\infty$-Whitney jets on $\overline{U}$ that satisfy the Cauchy-Riemann system.
Therefore, if a proper cone $S:=S(V, \epsilon)$
is subanalytic, in view of Theorem \ref{theorem:multi-asymptotic} the set of holomorphic functions on $S$
that are multi-asymptotically developable along
$\chi$ is equal to
$$
\varprojlim_{S'}
H^0(X;\CC_{\overline{S'}} \wtens \OO_X)
$$
where $S'$ runs through the family of open subanalytic proper cones
$S(V',\epsilon')$ properly contained in $S$. %Note that, here,
Let us consider the multi-specialization of Whitney holomorphic functions.
Let $W$ be an $(\RP)^\ell$-conic subanalytic open subset in $S_\chi$.
Remember that $(\RP)^\ell$-conic open subset $V'$ in $S_\chi$ is said to be
compactly generated in $W$ if there exists a compact subset $K$ in $W$ with
$V' \subset (\RP)^\ell K$.

We have
$$
\begin{aligned}
H^0(W;\imin \rho \nu^{sa}_\chi \OW_X) &= \lpro {V'}\,\lind {U'} H^0(U';\OW_X),
\end{aligned}
$$
where $V'$ ranges through the family of $(\RP)^\ell$-conic
subanalytic open subsets in $S_\chi$
which are compactly generated in $W$, $U'$ ranges through the family of $\op(X_{sa})$ such that
$C_\chi(X \setminus U') \cap V'=\emptyset$. One can check easily that, given $p=(0;\zeta) \in S_\chi$, then
\begin{eqnarray*}
(\imin\rho H^0\nu^{sa}_\chi\OW_X)_p & \simeq & \lind U H^0(U;\nu^{sa}_\chi\OW_X) \\
& \simeq & \lind {V,\epsilon} H^0(S(V,\epsilon);\OW_X) \\
& \simeq & \lind {V,\epsilon} H^0(X;\CC_{\overline{S(V,\epsilon)}} \wtens \OO_X),
\end{eqnarray*}
where $U$ ranges through the $(\RP)^\ell$-conic neighborhoods of $p$, $V$ through the family $V(\zeta)$ and $\epsilon>0$. The last isomorphism follows from Corollary \ref{cor: l.c.t. cone}.
Hence we can construct functorially the sheaf $\imin\rho H^0\nu^{sa}_\chi\OW_X$ and its stalk can be described by
multi-asymptotics. \\

Similarly
we have
\begin{eqnarray*}
(\imin\rho H^0\nu^{sa}_\chi\OW_{X|X \setminus \ckZ})_p & \simeq & \lind {V,\epsilon} H^0(X;\CC_{\overline{S(V,\epsilon)} \setminus \ckZ} \wtens \OO_X),
\end{eqnarray*}
where %$W$ ranges through the $(\RP)^\ell$-conic neighborhoods of $p$,
$V$ through the family $V(\zeta)$ and $\epsilon>0$.
Hence we can construct functorially the sheaf $\imin\rho H^0\nu^{sa}_\chi\OW_{X|X \setminus \ckZ}$ and its stalk can be described by
flat multi-asymptotics by Theorem \ref{teo:flat-holomorphic}.

\

Finally, thanks to Theorem \ref{teo: consistent family Whitney} one can check
\begin{eqnarray*}
(\imin\rho H^0\nu^{sa}_\chi\OW_{X| \ckZ})_p & \simeq & \lind {V,\epsilon} H^0(X;\CC_{\overline{S(V,\epsilon)} \cap \ckZ} \wtens \OO_X),
\end{eqnarray*}
where
%$W$ ranges through the $(\RP)^\ell$-conic neighborhoods of $p$,
$V$ through the family $V(\zeta)$ and $\epsilon>0$.
Hence we can construct functorially the sheaf $\imin\rho H^0\nu^{sa}_\chi\OW_{X| \ckZ}$ and its stalk can be described by
consistent families of coefficients. \\

\subsection{Vanishing theorems and Borel-Ritt exact sequence}\label{subsection:Borel-Ritt}

In this subsection, we establish several vanishing theorems related to
multi-specializations of the sheaves of Whitney holomorphic functions. As an application we obtain a Borel-Ritt exact sequence for multi-asymptotics. We continue to consider
the problems in the complex domain. \\

With the notations of Subsection \ref{subsection:Families of cones}, let $X=\CC^n$ with coordinates $z=(z_1,\dots,z_n)$.
 We consider the family $F^q$, $q = \#\left(J_Z(\xi) \cap \{1,2,\dots,L\}\right) \ge 0$, whose elements are couples $(f,v)$ with $f$ non-zero rational monomial and $v \in \R_{\geq 0}$.

In $F^q$ there is a multiplication $\star$ limited to the pairs $(f,v)$ and $(g,w)$ such that $\nu_k(f)>0$ and $\nu_k(g)<0$ for some $k \in \{1,\dots,L\}$. It is defined by $(f,g) \star (g,w) =(f^ag^b,v^aw^b)$, where the natural numbers $a$ and $b$ are  taken to be prime to each other and such that $a|\nu_k(f)|=b|\nu_k(g)|$.

Let $K$ be the set of monomials generated by the operation $\star$ in $F^q$. If $(f,v) \in K$, then  $f=\prod f_i^{a_i}$, $v=\prod v_i^{a_i}$ with $a_i$ natural numbers and $(f_i,v_i) \in F^q$. Remark that there are finitely many monomials belonging to $K$. Given $(f,v) \in K$, let
$$K_f=\left\{(f_i,v_i) \in F^q, \ (f,v)=\left(\prod f_i^{a_i}, \prod v_i^{a_i}\right)\right\}.$$

\begin{lem} \label{lem: closure cone} The closure of
$$
\bigcap_{(f,v) \in F^q}\left\{v-\epsilon < f(|z|) < v+\epsilon\right\}% \qquad f={f_n \over f_d}\right\},\ \ q \in 1,\dots,\ell
$$
is
$$
\bigcap_{(f,v) \in K}\left\{f_d(|z|)\prod_{(f_i,v_i)\in K_f}(v_i-\epsilon)^{a_i} \leq f_n(|z|) \leq f_d(|z|)\prod_{(f_i,v_i) \in K_f}(v_i+\epsilon)^{a_i}, \ \ f={f_n \over f_d}\right\}
$$
%where $f=\prod f_i^{a_i}$, $i \in K_f \subseteq F^q$, $f_i \in F^q$.
where $f_n,f_d$ are monomials.
%$$
%\bigcap_{\substack{(f,w) \in F^{q'} \\ q' \geq q}}\left\{f_d(|z|)\prod_{(f_i,v_i)\in K_f}(v_i-\epsilon)^{a_i} \leq f_n(|z|) \leq f_d(|z|)\prod_{(f_i,v_i) \in K_f}(v_i+\epsilon)^{a_i}, \ \ f={f_n \over f_d}\right\}
%$$
%where $f=\prod f_i^{a_i}$, $i \in K_f \subseteq F^q$, $f_i \in F^q$.
%where $f_n,f_d$ are monomials and
%$$K_f=\left\{(f_i,v_i) \in F^q, \ (f,v)=\left(\prod f_i^{a_i}, \prod v_i^{a_i}\right)\right\}.$$
%$$K_f=\left\{(f_i,v_i) \in F^q, \ f = \prod f_i^{a_i}\right\}.$$

\end{lem}
\begin{proof}
In what follows we write for short $f \in F^q$ if $(f,v) \in F^q$ and $i \in K_f$ if $(f_i,v_i) \in K_f$. It is enough to prove that the second one is contained in the closure of the first one. Indeed one can find from
$$
v-\epsilon < f(|z|) < v+\epsilon \qquad f={f_n \over f_d} \in F^q
$$
that
$$
f_d(|z|)\prod_{i\in K_f}(v_i-\epsilon)^{a_i} < f_n(|z|) < f_d(|z|)\prod_{i\in K_f}(v_i+\epsilon)^{a_i} \qquad f={f_n \over f_d} \in F_{q'},\ \ q' \geq q.
$$

It is enough to prove that, given $p_j>0$, $j=1,\dots,m$, with
\begin{equation}\label{eq:varphip}
f_d(p)\prod_{i\in K_f}(v_i-\epsilon)^{a_i} \leq f_n(p) \leq f_d(p)\prod_{i\in K_f}(v_i+\epsilon)^{a_i} \qquad f={f_n \over f_d} \in F_{q'},\ \ q' \geq q,
\end{equation}
%\begin{equation}\label{eq:varphip}
%\imin{\varphi_j}(p) \leq \epsilon \qquad j=1,\dots,\ell
%\end{equation}
%\begin{equation}\label{eq:p}
%p_j \leq \epsilon^{{\rm deg}\varphi_j} \qquad j=1,\dots,\ell
%\end{equation}
($p=(p_1,\dots,p_m)$) then for each $n \in \mathbb{N}$ there exist $|d_i|<\displaystyle{1\over n}$, $j=1,\dots,m$, such that
\begin{equation}\label{eq:varphip<}
f_d(|p-d|)\prod_{i\in K_f}(v_i-\epsilon)^{a_i} < f_n(|p-d|) < f_d(|p-d|)\prod_{i\in K_f}(v_i+\epsilon)^{a_i}  \qquad f={f_n \over f_d} \in F_{q'},\ \ q' \geq q.
\end{equation}
%\begin{equation}\label{eq:varphip<}
%\imin{\varphi_j}(|p-d|) < \epsilon \qquad j=1,\dots,\ell
%\end{equation}
%\begin{equation}\label{eq:p<}
%|p_j-d_j|< \epsilon^{{\rm deg}\varphi_j} \qquad j=1,\dots,\ell
%\end{equation}
Here $|p-d|$ means $(|p_i-d_i|)_i$. We may assume $0<\epsilon<1$. Suppose that $p_i \neq 0$ for $i=1,\dots,m$. We may assume that $d_i$ is such that $|p_i-d_i|=\epsilon^{\delta_i}p_i$. If $\delta_i$ is small enough, then $|d_i|<\displaystyle{1\over n}$. Let $\epsilon^\delta=(\epsilon^{\delta_1},\dots,\epsilon^{\delta_m})$. We have
$$
f_n(|p-d|)=f_n(p)f_n(\epsilon^\delta), \ \ f_d(|p-d|)=f_d(p)f_d(\epsilon^\delta).
$$
By \eqref{eq:varphip} we have
$$
f_d(p)\prod_{i\in K_f}(v_i-\epsilon)^{a_i} \leq f_n(p) \leq f_d(p)\prod_{i\in K_f}(v_i+\epsilon)^{a_i}.
$$
Suppose that
$$
f_d(p)\prod_{i\in K_f}(v_i-\epsilon)^{a_i} < f_n(p)
$$
(when the equality holds it means $v_i \neq 0$, in that case we consider the inequality
$$
{1 \over f_d(p)\prod_{i\in K_f}(v_i+\epsilon)^{a_i}} < {1 \over f_n(p)}
$$
and the proof is similar). We are reduced to prove that for each $n \in \N$ there exists $0<\delta<\displaystyle{1\over n}$ such that
$$
f_n(\epsilon^\delta)<f_d(\epsilon^\delta) \qquad f={f_n \over f_d} \in F_{q'},\ \ q' \geq q.
$$
Remark that if $\delta$ is small enough (and hence $\epsilon^\delta$ is close enough to 1) then
$$
f_d(p)f_d(\epsilon^\delta)\prod_{i\in K_f}(v_i-\epsilon)^{a_i} < f_n(p)f_n(\epsilon^\delta).
$$
Setting
$$
f_n=\prod_{i \in K_{f_n}}x^{a_i}_i, \ \ f_d=\prod_{i \in K_{f_d}}x^{a_i}_i
$$
and taking $\log_\epsilon$ we are reduced to prove that for each $n \in \N$ there exists $0<\delta<\displaystyle{1\over n}$ such that
$$
\sum_{i \in K_{f_n}}a_i\delta_i>\sum_{i \in K_{f_d}}a_i\delta_i \qquad f={f_n \over f_d} \in F_{q'},\ \ q' \geq q.
$$
We argue by induction on the number of variables $m$. If $m=1$ the result is easy to prove. We suppose that we proved it for $m' \leq m-1$. Hence we are reduced to solve systems
$$
\begin{cases}
\displaystyle{\sum_{i \in K_{f_n}}a_i\delta_i>\sum_{i \in K_{f_d}}a_i\delta_i} & \\
\displaystyle{\sum_{i \in K_{g_n}}a_i\delta_i>\sum_{i \in K_{g_d}}a_i\delta_i}
\end{cases}
$$
with the exponents $a_i>0$. We may assume that there exists $i_0 \in K_{f_n} \cap K_{g_d}$ (if not the system has always solutions). Arguing by induction again, we may assume that $i_0=m$.
$$
\begin{cases}
\displaystyle{\sum_{i \in K_{f_n} \setminus \{m\}}a_i\delta_i+a_n\delta_m>\sum_{i \in K_{f_d}}a_i\delta_i} & \\
\displaystyle{\sum_{i \in K_{g_n}}a_i\delta_i>\sum_{i \in K_{g_d} \setminus \{m\}}a_i\delta_i+a_d\delta_m}
\end{cases}
$$
So, up to normalize the coefficient of $\delta_{m}$ we may reduce to
$$
\begin{cases}
f'_n(\delta')+\delta_{m}>f'_d(\delta') \\
g'_n(\delta')>g'_d(\delta')+\delta_{m}
\end{cases}
$$
where $\delta'=(\delta_1,\dots,\delta_{m-1})$. This system has solutions if
$$
g'_n(\delta')-g'_d(\delta') > \delta_m > f'_d(\delta')-f'_n(\delta').
$$
So we can find a solution $\delta_m$ if
$(f'_n+g'_n)(\delta') > (f'_d+g'_d)(\delta')$ which follows from the induction hypothesis. Indeed, by definition, $\displaystyle{fg={f_ng_n \over f_dg_d}} \in F^{q'}$, for some $q' \geq q$. Such a solution belongs to the interval
$$
(f'_n+g'_n)(\delta') > \delta_m > (f'_d+g'_d)(\delta').
$$
Again by induction we may assume that $\delta'$ is small enough to imply $0<\delta_m<\displaystyle{1\over n}$ and the result follows.\\

Let us consider the case when some $p_i$ is zero. Up to take a permutation of $\{1,\dots,m\}$ we may assume that $p_m=0$. In this case we assume that $|p_i-d_i|=\epsilon^{\delta_i}$, with $0<\epsilon<1$. We argue by induction again. We have to prove that for each $N>0$ there exists $\delta_m>N$ such that
$$
f_d(\epsilon^\delta)\prod_{i\in K_f}(v_i-\epsilon)^{a_i} \leq f_n(\epsilon^\delta) \leq f_d(\epsilon^\delta)\prod_{i\in K_f}(v_i+\epsilon)^{a_i}.
$$
with $\delta_i$, $i=1,\dots,m-1$ is such that $|p_i-\epsilon^{\delta_i}|<\displaystyle{1 \over n}$.
Taking $\log_\epsilon$ as in the previous case we are reduced to solve systems
$$
\begin{cases}
\displaystyle{Q^-_f+\sum_{i \in K_{f_d}}a_i\delta_i > \sum_{i \in K_{f_n} \setminus \{m\}}a_i\delta_i+a_n\delta_m>Q^+_f+\sum_{i \in K_{f_d}}a_i\delta_i} & \\
\displaystyle{Q^-_g+\sum_{i \in K_{g_d} \setminus \{m\}}a_i\delta_i+a_d\delta_m>\sum_{i \in K_{g_n}}a_i\delta_i>Q^+_g+\sum_{i \in K_{g_d} \setminus \{m\}}a_i\delta_i+a_d\delta_m}
\end{cases}
$$
where
\begin{eqnarray*}
Q^+_f=\sum\log_\epsilon(v_i+\epsilon)^{a_i}, && Q^-_f=\sum_{i \in K_{f}}\log_\epsilon(v_i-\epsilon)^{a_i}, \\
Q^+_g=\sum\log_\epsilon(v_i+\epsilon)^{a_i}, && Q^-_g=\sum_{i \in K_{g}}\log_\epsilon(v_i-\epsilon)^{a_i},
\end{eqnarray*}
with the convention $\log_\epsilon(v_i-\epsilon)=+\infty$ if $v_i=0$. So we may reduce to systems of this kind
$$
\begin{cases}
f'_n(\delta')+\delta_m>a_f'+f'_d(\delta') \\
g'_n(\delta')>a'_g+g'_d(\delta')+\delta_m
\end{cases}
$$
where $\delta'=(\delta_1,\dots,\delta_{m-1})$. This system has solutions if $(f'_n+g'_n)(\delta') > a'_f+a'_g+ (f'_d+g'_d)(\delta')$ which follows from the induction hypothesis and a solution belongs to the interval
$$
(f'_n+g'_n)(\delta') > \delta_m > a'_f+a'_g+ (f'_d+g'_d)(\delta').
$$
Again by induction we may assume that $\delta'$ is big enough to imply $M<g'_n(\delta')$ and the result follows. Indeed there must be $i \in K_{g_n}$ such that $p_i=0$ and hence we may assume that $\delta_i$ is big enough for such an $i$.

\end{proof}

As a consequence of Lemma \ref{lem: closure cone} we get a description of the closure of the open sets $S(V,\, \epsilon)$.

\begin{prop}{\label{prop:closure_of_S}} Let $p=(0,\zeta)=(0,\zeta^{(1)},\dots,\zeta^{m}) \in S_\chi$. Then for $\epsilon>0$ and $V \in V(\zeta)$ the closure of
$$
S(V,\, \epsilon)
= \left\{z \in X;\,
\begin{aligned}
&z \in V,\,\, |z^{(0)}| < \epsilon,\,\,
v - \epsilon < f\left(|z^{(*)}|\right) < v+\epsilon \\
&\text{for any $(f,\,v) \in F^q$}.
\end{aligned}
\right\}
$$
is
$$
\overline{S(V,\, \epsilon)}
= \bigcap_{(f,v) \in K}\left\{z \in X;\,
\begin{aligned}
& z \in \overline{V},\,\, |z^{(0)}| \leq \epsilon, \\
& (v - \epsilon_f)f_d\left(|z^{(*)}|\right) \leq f_n\left(|z^{(*)}|\right) \leq (v+\epsilon_f)f_d\left(|z^{(*)}|\right)% \\
%&\text{for any $(f,\,v) \in F^{q'}$ and any $q' \geq q$}.
\end{aligned}
\right\}.
$$
Here $f=\displaystyle{f_n \over f_d}=\prod f_i^{a_i}$ ($f_n,f_d$ monomials) and $v\pm\epsilon_f=\prod (v_i\pm\epsilon)^{a_i}$, $(f_i,v_i) \in K_f \subseteq F^q$.% $f_i \in F^q$.
\end{prop}

\begin{cor}  Let $p=(0,\zeta)=(0,\zeta^{(1)},\dots,\zeta^{m}) \in S_\chi$. A cofinal family for
$$
\{\overline{S(V,\, \epsilon)},\, \epsilon>0,\, V \in V(\zeta)\}
$$
is given by the sets
$$
\bigcap_{(f,v) \in K}\left\{z \in X;\,
\begin{aligned}
&z \in \overline{V},\,\, |z^{(0)}| \leq \epsilon, \\
& (v - \epsilon)f_d\left(|z^{(*)}|\right) \leq f_n\left(|z^{(*)}|\right) \leq (v+\epsilon)f_d\left(|z^{(*)}|\right)% \\
%&\text{for any $(f,\,v) \in F^{q'}$ and any $q' \geq q$}.
\end{aligned}
\right\}%_{\substack{\epsilon >0 \\ V \in V(\zeta)}}
$$
where $f=\displaystyle{f_n \over f_d}$, $(f,v)=(\prod f_i^{a_i},\prod v_i^{a_i})$, $(f_i,v_i) \in K_f$, $\epsilon >0$ and $V \in V(\zeta)$.
\end{cor}

\begin{oss} In general,
the closure of $S(V,\epsilon)$ is strictly smaller than
the set defined by equalities where $<$ is replaced with $\leq$.
One  can observe this fact for a clean intersection case:
$Z_1 = \{z_1 = z_3 = 0\}$,
$Z_2 = \{z_2 = z_3 = 0\}$,
$Z_3 = \{z_1 = z_2 = z_3 = 0\}$.
Then $S(V,\epsilon)$ is defined by
\begin{eqnarray*}
|z_3|/|z_1| & < & \epsilon \\
|z_3|/|z_2| & < & \epsilon \\
|z_1 z_2|/ |z_3| & < & \epsilon
\end{eqnarray*}
and $z_k \in G_k$, where $G_k$ is an open cone.
Note that this set is relatively compact.
Then,
\begin{eqnarray*}
|z_3| & \leq & \epsilon |z_1| \\
|z_3| & \leq & \epsilon |z_2| \\
|z_1 z_2| & \leq & \epsilon |z_3|
\end{eqnarray*}
with $z_k \in \overline{G_k}$,
is not compact.
In fact, the unbounded set $z_1 = z_3=0$ and $z_2 \in G_2$
is contained in this set.
So the set is strictly bigger than the closure of $S(V,\epsilon)$.
To avoid this, we need to add some extra inequalities.
\end{oss}

\begin{teo} \label{teo:vanishing_asymptotic}
We have
\begin{equation}
H^k(\imin\rho\nu^{sa}_\chi\OW_{X}) = 0 \qquad (k \ne 0).
\end{equation}
\end{teo}
\begin{proof}
In the proof, we follow the notations used in Subsection
\ref{subsection:Geometry of multi-asymptotic expansions}, and we may also assume
from the beginning that all the entries of $A_\chi$ are non-negative integers and
$\comD = 1$.
Let $p = (0;\,\zeta) = (0;\zeta^{(1)},\dots,\zeta^{(m)}) \in S_\chi$.
Define $\widetilde{\zeta} := (\widetilde{\zeta}^{(1)},\dots, \widetilde{\zeta}^{(m)})$ by
$\widetilde{\zeta}^{(k)} = \zeta^{(k)}/|\zeta^{(k)}|$ ($k \notin J_Z(\zeta)$) and
$\widetilde{\zeta}^{(k)} = 0$ otherwise.
%It follows from the fiber formula and the cohomological
%triviality of $S(V,\epsilon)$ that we have
By the observations in this section, we have
$$
(\imin\rho H^k\nu^{sa}_\chi\OW_X)_p
 \simeq  \lind {V,\epsilon} H^k(X;\CC_{\overline{S(V,\epsilon)}} \wtens \OO_X),
$$
where $V$ runs through the family $V(\zeta)$ and $\epsilon>0$. As the
closed subset $\overline{S(V,\epsilon)}$ does not satisfy the ``property $(\lambda)$''
of \cite{Du79} in general,
we need to replace the family $\{\overline{S(V,\epsilon)}\}_{V,\epsilon}$ with an equivalent family of closed analytic polyhedra which, in particular, satisfy the ``property $(\lambda)$''.
Then the result immediately follows from the main theorem in \cite{Du79}.

Let us construct an equivalent family. We may assume, in what follows, the cone $V$
is always taken to be an analytic polyhedron.
Set
$$
\check{N} :=
\{ \alpha = (\alpha_1,\dots,\alpha_m) \in \mathbb{Z}_{\ge 0}^m;\,
\alpha_k =      0 \,\, (k \in J_Z(\zeta))\}
$$
and, for $\alpha = (\alpha_1,\dots,\alpha_m) \in \mathbb{Z}_{\ge 0}^m$,
$$
h_{\alpha}(z) := \underset{1 \le k \le m}{\prod}\, |z^{(k)}|^{\alpha_k}.
$$
We also denote by $\Theta$ the set of strictly increasing continuous functions
$\theta(t)$ on $(0,+\infty)$ satisfying $\displaystyle\lim_{t \to 0 + 0} \theta(t) = 0$.
By Proposition \ref{prop:closure_of_S}
on the closure of $S(V,\epsilon)$, we may assume that
each closed subset $\overline{S(V,\epsilon)}$ in the family
$\{\overline{S(V,\epsilon)}\}_{V,\epsilon}$ has the following form:
$$
\begin{aligned}
&\left(\underset{(\alpha, \beta, \theta) \in P}{\bigcap}
\{z \in \overline{V};\, h_{\alpha}(z) \le \theta(\epsilon) h_{\beta}(z)\}\right)\,\,\, \bigcap \\
&\qquad
\left(
\underset{(\gamma, \delta,v, \theta) \in Q}{\bigcap}
\{z \in \overline{V};\, h_{\gamma}(z) \le (v + \theta(\epsilon)) h_{\delta}(z)\}
\right)\, \bigcap\, \left\{|z^{(0)}| \le \epsilon \right\},
\end{aligned}
$$
where $P$ is a finite subset of $(\mathbb{Z}_{\ge 0}^m) \times \check{N} \times \Theta$
and
$Q$ is also a finite subset of
$\check{N} \times \check{N} \times \mathbb{R}_{> 0} \times \Theta$. Note that
$\alpha$ belongs to $\mathbb{Z}_{\ge 0}^m$, however, the other
$\beta$, $\delta$ and $\gamma$ belong to $\check{N}$, which
is crucial in the subsequent arguments.

Define the family $\{H(V,\epsilon)\}_{V,\epsilon > 0}$ by
$$
\{z \in \overline{V};\, h_{\alpha}(z) \le \theta(\epsilon) h_{\beta}(z)\}
$$
with a fixed $(\alpha, \beta, \theta) \in (\mathbb{Z}^m_{\ge0}) \times \check{N} \times \Theta$,
and define also the family $\{G(V,\epsilon)\}_{V,\epsilon > 0}$ by
$$
\{z \in \overline{V};\, h_{\gamma}(z) \le (v + \theta(\epsilon)) h_{\delta}(z)\}
$$
with a fixed
$(\gamma, \delta, v, \theta) \in
\check{N} \times \check{N} \times \mathbb{R}_{>0} \times \Theta$.

First we will make an equivalent
family $\{\widetilde{H}(\Xi,V,\epsilon)\}_{\Xi,V,\epsilon}$
to the family $\{H(V,\epsilon)\}_{V, \epsilon > 0}$
which consists of closed analytic polyhedra as follows:
For a vector
$$
\xi = (\xi^{(1)},\,\dots,\xi^{(m)}) \in \mathbb{C}^{n_1} \times \dots \times \mathbb{C}^{n_m}
$$
with $|\xi^{(k)}| \ne 0$ ($1 \le k \le m$),
and for $\alpha \in \mathbb{Z}_{\ge 0}^m$, set
$$
\widehat{h}_{\xi,\alpha}(z) :=
\operatorname{Re}\,\, \underset{1 \le k \le m}{\prod}\,
\left(\dfrac{\left|\xi^{(k)}\right|\, \left\langle z^{(k)},\,\, \overline{\xi^{(k)}} \right\rangle}
{\left\langle \xi^{(k)},\,\,\overline{\xi^{(k)}}\right\rangle}
\right)^{\alpha_k},
$$
where $\langle a,b \rangle := a_1 b_1 + \dots + a_r b_r$ for
$a = (a_1,\dots,a_r)$ and $b = (b_1,\dots, b_r)$ and $\overline{a}$ denotes
the complex conjugate of a complex vector $a$. Then we can easily confirm the
following facts: For any $z \in X$, we have
\begin{equation}{\label{eq:norm_est_x1}}
\widehat{h}_{\xi,\alpha}(z) \le (\sqrt{n})^{|\alpha|} h_\alpha(z).
\end{equation}
Furthermore, for $\beta \in \check{N}$ and $ V \in V(\zeta) = V(\widetilde{\zeta})$,
there exists constant $\kappa > 1$ satisfying
\begin{equation}{\label{eq:norm_est_x2}}
	\kappa^{-1} h_\beta (z) \le \widehat{h}_{\widetilde{\zeta},\,\beta}(z)
	\le \kappa h_\beta (z) \qquad (z \in \overline{V}).
\end{equation}
Note that, here, if $V$ becomes smaller as a cone,
then we can take the above $\kappa$ to be closer to $1$.

Let $\Xi = \{\Xi_k\}_{k \in J_Z(\zeta)}$ be a family of the sets of vectors
where each set $\Xi_k$ ($k \in J_Z(\zeta))$ consists
of finitely many vectors $\xi^{(k)}$ in $\mathbb{C}^{n_k}$ with $|\xi^{(k)}| = 1$.
Then we define the closed analytic polyhedra $\widetilde{H}(\Xi,V,\epsilon)$ by
$$
\underset{\xi}{\bigcap}
\left\{z \in \overline{V};\, \widehat{h}_{\xi,\alpha}(z) \le \epsilon \widehat{h}_{\widetilde{\zeta},\beta}(z)
\right\},
$$
where the vector $\xi = (\xi^{(1)},\dots,\xi^{(m)})$ runs through the ones satisfying
$\xi^{(k)} = \widetilde{\zeta}^{(k)}$ for $k \notin J_Z(\zeta)$ and
$\xi^{(k)} \in \Xi_k$ for $k \in J_Z(\zeta)$.
We can easily confirm that, by
\eqref{eq:norm_est_x1} and
\eqref{eq:norm_est_x2},  the family $\{\widetilde{H}(\Xi,V,\epsilon)\}_{\Xi,V,\epsilon}$ and
the one $\{H(V,\epsilon)\}_{V,\epsilon}$ are equivalent.

For an equivalent family to the family $\{G(V,\epsilon)\}_{V, \epsilon > 0}$,
by noticing the remark after \eqref{eq:norm_est_x2},
we can construct an equivalent family $\{\widetilde{G}(V,\epsilon)\}_{V,\epsilon}$
of $\{G(V,\epsilon)\}_{V,\epsilon > 0}$ consisting of closed analytic polyhedra by
%in the same way as that for $\{H(V,\epsilon)\}_{V,\epsilon}$. That is, each $\widetilde{G}(V,\epsilon)$ is given by
$$
\widetilde{G}(V,\epsilon) :=
\left\{z \in \overline{V};\, \widehat{h}_{\widetilde{\zeta}, \gamma}(z) \le (v+\epsilon) \widehat{h}_{\widetilde{\zeta},\delta}(z)
\right\}.
$$

Hence we have obtained a required equivalent family of
$\{\overline{S(V,\epsilon)}\}_{V,\epsilon}$ and the theorem follows.
\end{proof}
%The identity morphism induces a morphism of presheaves
%\begin{eqnarray*}
%\widetilde{\A}_\chi & \to & \imin \rho H^0\nu^{sa}_\chi \OW_X, \\
%\end{eqnarray*}
%The above morphism is isomorphisms in the stalks (i.e. in the limit of multicones containing a given direction
%$p \in S^\circ$). Let  $\A_\chi$ be the sheaf associated to $\widetilde{\A}_\chi$. We get
%\begin{eqnarray*}
%\A_\chi & \iso & \imin \rho H^0\nu^{sa}_\chi \OW_X, \\
%\end{eqnarray*}

%By Proposition \ref{prop: exact sequence}, on $S^\circ$ we have the exact sequence of sheaves
%\begin{equation}\label{eq:Borel-Ritt}
%\exs{\A^{<0}_{\chi}}{\A_{\chi}}{\A^{CF}_{\chi}}.
%\end{equation}
%When $Z$ forms a normal crossing divisor, these sheaves are nothing but the sheaves of strongly asymptotically developable functions defined by Majima in \cite{Ma84}.

%\begin{oss} The exact sequence \eqref{eq:Borel-Ritt} is nothing but a Borel-Ritt exact sequence for multi-asymptotically developable functions. This was already proven for formal specialization in the single divisor case in \cite{Co01} and for Majima's asymptotic in \cite{GS99,HS00}. Both results are based on the Borel-Ritt theorem in dimension one (for a proof, see \cite{Wa65}). Thanks to Proposition \ref{vanishing} we obtained a purely cohomological proof of the Borel-Ritt exact sequence.
%\end{oss}

\begin{teo} \label{teo: flat asymptotics in deg 0}
Let $p \in S_\chi$ be a point outside fixed points.
Assume that the associated action $\mu$ is non-degenerate.
Then we have
$$
H^k(\imin\rho\nu^{sa}_\chi\OW_{X|X\setminus \ckZ})_p = 0 \qquad (k \ne 0).
$$
\end{teo}

Let $p = (0;\,\zeta) = (0;\,\zeta^{(1)},\, \dots, \zeta^{(m)})$.
%Let $\comD$ be a positive rational number so that $a_{jk} \in \mathbb{Z}$ and all the $a_{jk}$ have no common divisors. Let $\mu_{\intA}$ be the action on $X$ associated with the matrix $\intA$. Since multicones associated with the action $\mu$ and those with the action $\mu_{\intA}$ determine equivalent families, we may assume from the beginning that all the entries of $A_\chi$ are non-negative integers and $\comD = 1$.
Let $\varphi_k$ ($1 \le k \le m$) be the monomials associated with the
action $\mu$. We may assume from the beginning that all the entries of $A_\chi$ are non-negative integers and $\comD = 1$. Since $\mu$ is non-degenerate and $p$ is outside fixed
points,  we have $\operatorname{Rank} A_\chi = \ell \le m$ and
we may assume that the first $\ell \times \ell$ sub-matrix in $A_\chi$
is invertible and $\zeta^{(k)} \ne 0$ holds for $1 \le k \le \ell$.

Now we prove this theorem after some preparations.
We regard
$$
w=(w_1,\dots,w_\ell) = \varphi(t)
=(\varphi_1(t),\,\dots,\,\varphi_\ell(t))
$$
as a holomorphic map from $\mathbb{C}_t^\ell$ to $\mathbb{C}_w^\ell$.
Set
$$
D := \{w_1w_2\cdots w_\ell = 0\} \subset \mathbb{C}^n_w \text{ and }
L := \{t_1 \cdots t_\ell = 0\} \subset \mathbb{C}^n_t.
$$
Then it is easy to see
$$
\varphi^{-1}(D) = L.
$$
Furthermore, as $w=\varphi(t)$ has the inverse $t = \varphi^{-1}(w)$ consisting of
monomials of rational powers, the map
$$
\varphi|_{\mathbb{C}^n_t \setminus L}: (\mathbb{C}^n_t \setminus L) \longrightarrow
(\mathbb{C}^n_w \setminus D)
$$
is a finite covering.

Let $\Xi$ be the set of open poly-sectors
$V = V_1 \times V_2 \times \dots \times V_\ell$
in $\mathbb{C}^\ell$ where each $V_k$ is an open proper sector in $\mathbb{C}$
and $1 \in V_k$.
We denote by $\hat{A}_\chi$ the sub-matrix of $A_\chi$ consisting
of the leftmost $\ell$ columns which is, by the assumption, invertible.
Then, as
\begin{equation}\label{eq:angle_relation}
\arg(w) = (\arg(w_1),\dots, \arg(w_\ell))
= \arg(t) {\hat{A}_\chi} \qquad (t \in \mathbb{C}^\ell \setminus L),
\end{equation}
holds, we have the following:
$$
\text{For any $V \in \Xi$, there exists $V', W \in \Xi$ such that $V' \subset \varphi(W)
\subset V$}.
$$

\

Remember that $n_k$ is dimension of the $k$-th coordinates block.
Now, by a linear coordinates transformation on each coordinate blocks,
we may assume that each direction $\zeta^{(k)}$ has a form
$$
(|\zeta^{(k)}|,\,0,\,\dots,\,0) \in \mathbb{C}^{n_k} \qquad (1 \le k \le m).
$$
Furthermore we may assume $|\zeta^{(k)}| = 1$ for $1 \le k \le \ell$ by considering
normalization of $\zeta$.
Set
$$
\begin{array}{cccccccc}
Y = &\mathbb{C}^{n_0}
&\times &\left(\underset{1 \le k \le \ell}{\times} \mathbb{C}^{n_k-1}\right)
&\times &\left(\underset{\ell < k \le m}{\times} \mathbb{C}^{n_k}\right)
&\times &\mathbb{C}^\ell. \\
\\
&\big( z^{(0)}
& , &{z}_*^{(1)},\dots,{z}_*^{(\ell)}
& , &z^{(\ell+1)},\cdots,z^{(m)}
& , &t \big).
\end{array}
$$
Then define the map $f: Y \to X$ by
$$
\begin{aligned}
&\big( z^{(0)},
\,\,z_*^{(1)},\dots,z_*^{(\ell)},
\,\,z^{(\ell+1)},\cdots,z^{(m)},
\,\,t \big)  \\
&\quad \longrightarrow
\big( z^{(0)},
\,\,\varphi_1(t)(1,z_*^{(1)}),\dots,\varphi_\ell(t)(1,z_*^{(\ell)}),
\,\,\varphi_{\ell+1}(t) z^{(\ell+1)},\cdots, \varphi_m(t)z^{(m)} \big)
\end{aligned}
$$
We denote by $z^{(k)}_1$  the first coordinate of the coordinates block $z^{(k)}$
and by $z^{(k)}_*$ the rest, i.e., $z^{(k)} = (z^{(k)}_1,\,z^{(k)}_*)$.
Define the subset in $X$ by
$$
T_X := \{z^{(1)}_1 z^{(2)}_1 \cdots z^{(\ell)}_1 = 0\} \subset X
$$
and the one in $Y$ by
$$
T_Y := \{t_1 t_2 \cdots t_\ell = 0\} \subset Y.
$$
By the previous observations and by noticing the fact $\varphi_k(t) \ne 0$
if $t_1 \dots t_\ell \ne 0$, we have
$$
f^{-1}(T_X) = T_Y
$$
and
$$
\text{the map $f|_{Y\setminus T_Y}$ is a finite covering over
$X \setminus T_X$.}
$$
Furthermore, we have
$$
Z_1\cup Z_2 \cup \dots \cup Z_\ell = f(T_Y) \subset T_X.
$$
Let $\epsilon_Y > 0$ and $V_Y \in \Xi$.
Then define the open subset $S_Y(V_Y,\,\epsilon_Y)$ in $Y$ by
$$
\left\{
\begin{aligned}
&|z^{(0)}| < \epsilon_Y,\,\,\, |z_*^{(k)}| < \epsilon_Y\,\,(1 \le k \le \ell)\\
&|z^{(k)} - \zeta^{(k)}| < \epsilon_Y\,\,(\ell < k \le m)\\
&t \in V_Y,\,\, |t_k| < \epsilon_Y\,\,(1 \le k \le \ell)
\end{aligned}
\right\}.
$$
We first show the following lemma:
\begin{lem}{\label{lem:set_equiv_f}}
For a multicone $S_X(V_X,\,\epsilon_X) \subset X$
with $\epsilon_X > 0$ and $V_X \in V(\zeta)$,
we can find some $S_Y(V_Y,\,\epsilon_Y)$
with $\epsilon_Y >0$ and $V_Y \in \Xi$ such that
$$
f(S_Y(V_Y,\,\epsilon_Y)) \subset S_X(V_X,\,\epsilon_X).
$$
\end{lem}
\begin{proof}
In fact, if $\epsilon_Y$ is sufficiently small and $V_Y$ is sufficient thin,
then clearly
$
f(S_Y(V_Y,\,\epsilon_Y)) \subset V_X
$
holds. As we have
$$
\varphi^{-1}_k(
\varphi_1(t)(1,z_*^{(1)}),\dots,\varphi_\ell(t)(1,z_*^{(\ell)}))
= t_k\varphi^{-1}_k((1,z_*^{(1)}),\dots,(1,z_*^{(\ell)})),
$$
and as $|\varphi^{-1}_k((1,z_*^{(1)}),\dots,(1,z_*^{(\ell)}))|$ is closed to 1,
we see that
$
|\varphi_k^{-1}| < \epsilon_X
$
holds at a point in $f(S_Y(V_Y,\,\epsilon_Y))$.  In the same way, we have
$$
\begin{aligned}
&\psi_k(\varphi_1(t)(1,z_*^{(1)}),\dots,\varphi_\ell(t)(1,z_*^{(\ell)}),
\varphi_{\ell+1}(t) z^{(\ell+1)},\cdots, \varphi_m(t)z^{(m)}) \\
&\qquad= z^{(k)}/\varphi_k(\varphi^{-1}((1,z_*^{(1)}),\dots,(1,z_*^{(\ell)})).
\end{aligned}
$$
Since $|\varphi_k(\varphi^{-1}((1,z_*^{(1)}),\dots,(1,z_*^{(\ell)}))|$ is very close to $1$ and
$|z^{(k)}|$ is nearly $|\zeta^{(k)}|$, we know that
$$
|\zeta^{(k)}| - \epsilon_X < |\psi_k| < |\zeta^{(k)}| + \epsilon_X
$$
holds at a point in $f(S_Y(V_Y,\,\epsilon_Y))$. Hence we have obtained the desired inclusion.
\end{proof}

On the other hand, we have the
inverse of $f$ defined on a fixed $S_X(V_X,\,\epsilon_X)$
($\epsilon_X > 0$ and $V_X \in V(\zeta)$),
which is explicitly given by
$$
\begin{aligned}
&\big( z^{(0)},
\,\,z^{(1)},\dots, z^{(m)}\big) \in X  \\
&\quad \longrightarrow
\big( z^{(0)},
\,\,z_*^{(1)}/z_1^{(1)},\dots,
z_*^{(\ell)}/z_1^{(\ell)},
\,\,z^{(\ell+1)}/\varphi_{\ell+1}(t) ,\cdots, z^{(m)}/\varphi_m(t),\,
t\big) \in Y,
\end{aligned}
$$
where $t$ is determined by
\begin{equation}\label{eq:determined_t}
t_k = \varphi^{-1}_k(z^{(1)}_1,z^{(2)}_1,\dots,z^{(\ell)}_1)\quad (1 \le k \le \ell).
\end{equation}
\begin{lem}{\label{lem:set_equiv_f_inv}}
For a given $S_Y(V_Y,\,\epsilon_Y)$
with $\epsilon_Y > 0$ and $V_Y \in \Xi$,
we can find $S_X(V'_X,\,\epsilon_X')$ with
$0 < \epsilon_X' < \epsilon_X$ and $V'_X \subset V_X$ in $V(\zeta)$
such that
$$
f^{-1}(S_X(V_X',\,\epsilon_X')) \subset S_Y(V_Y,\,\epsilon_Y).
$$
\end{lem}
\begin{proof}
This is shown by the same argument as the one in the previous lemma.
First we note that, for any $\delta>0$, if we take $V'_X \in V(\zeta)$
sufficiently thin, then we have
\begin{equation}\label{eq:bounded_other_by_top}
(1 - \delta) |z^{(k)}_1| \le |z^{(k)}| \le (1+\delta) |z^{(k)}_1|
\qquad (1 \le k \le m)
\end{equation}
hold on this $V'_X$.

It follows from \eqref{eq:angle_relation}
and \eqref{eq:determined_t}
that we have
$
t \in V_Y
$ if $\epsilon'_X$ is sufficiently small and $V'_X$ is thin.
Furthermore, $z \in S_X(V'_X,\,\epsilon'_X)$ implies
$
\varphi^{-1}_k(|z^{(*)}|) < \epsilon'_X
$
for $1 \le k \le \ell$.
Hence, by noticing \eqref{eq:bounded_other_by_top},
we get $|t_k| < \epsilon_Y$.

For $1 \le k \le \ell$,
as $z^{(k)}$ is a point in a cone in $\mathbb{C}^{n_k}$ with direction $(1,0,\dots,0)$,
by \eqref{eq:bounded_other_by_top}, we have
$$
|z_*^{(k)}/z_1^{(k)}| < \epsilon_Y \qquad	(1 \le k \le \ell).
$$

Since $|\zeta^{(k)}| - \epsilon'_X < \psi_k(|z^{(*)}|) < |\zeta^{(k)}| + \epsilon'_X$ ($\ell < k \le m$)
follows from $z \in S_X(V'_X,\,\epsilon'_X)$,
by taking \eqref{eq:bounded_other_by_top} into account,
for any $\delta' > 0$, we get
$$
(1-\delta')|\zeta^{(k)}| < |z^{(k)}/t_k| < (1+ \delta')|\zeta^{(k)}| \qquad (\ell <  k \le m)
$$
if $\epsilon'_X > 0$ is taken to be so small and $V'_X$ is so thin.
Hence, since $z \in V'_X$, in particular, $z^{(k)}$ belongs to a thin cone in $\mathbb{C}^{n_k}$ with
direction $\zeta^{(k)}$ and since $\arg (t_k)$ is nearly $0$, we conclude
$$
|z^{(k)}/t_k - \zeta^{(k)}| < \epsilon_Y.
$$
This completes the proof.
\end{proof}

\

For any $V_Y \in \Xi$, since $Z_1 \cup \dots \cup Z_\ell = f(T_Y) \subset T_X$
and $f$ is locally isomorphic outside $T_Y$, we have
$$
\overline{f\left( (\mathbb{C}^{n-\ell} \times V_Y) \cap U\right)} \cap T_X =
\overline{f\left( (\mathbb{C}^{n-\ell} \times V_Y) \cap U\right)}
\cap (Z_1 \cup Z_2 \cup \dots Z_\ell)
$$
for a bounded open subset $U$ in $Y$.
Therefore we have
$$
\begin{aligned}
f\left(\overline{S_Y(V_Y,\,\epsilon_Y)} \setminus T_Y\right)
&=f\left(\overline{S_Y(V_Y,\,\epsilon_Y)}\right) \setminus T_X \\
&= \overline{f\left(S_Y(V_Y,\,\epsilon_Y)\right)} \setminus
T_X  \\
&=\overline{f(S_Y(V_Y,\,\epsilon_Y))} \setminus (Z_1 \cup \dots \cup Z_\ell).
\end{aligned}
$$
Set $S_Y := S_Y(V_Y,\,\epsilon_Y)$ for short.
Then, as $f|_{\overline{S_Y} \setminus T_Y}$ is an isomorphism onto its image,
we have
$$
Rf_!\CC_{\overline{S_Y}\setminus T_Y} =
\CC_{f\left(\overline{S_Y}\setminus T_Y\right)} =
\CC_{\overline{f(S_Y)}\setminus(Z_1 \cup \cdots \cup Z_\ell)}.
$$
Hence it follows from Theorem 5.7 of {\cite{KS96}} that we have
\begin{equation}{\label{eq:direct-image-ks}}
R\Gamma(Y;\rh_{\D_Y}(\D_{Y \stackrel{f}{\to} X},\CC_{\overline{S_Y} \setminus T_Y}
\wtens \OO_Y) \simeq
R\Gamma(X;\CC_{\overline{f(S_Y)}\setminus (Z_1 \cup \cdots \cup Z_\ell)} \wtens \OO_X).
\end{equation}
Since $f$ is locally isomorphic outside $T_Y$, its Jacobian $J_f$ does not vanish
outside $T_Y$. Therefore, for any point $p \in Y$,
$J_f$ has the form $g\,\, t_1^{\alpha_1}t_2^{\alpha_2}\cdots t_\ell^{\alpha_\ell}$ for
some multi-index $\alpha = (\alpha_1,\dots,\alpha_\ell) \in \mathbb{Z}_{\ge 0}^\ell$ and
some holomorphic function $g$ with $g(p) \ne 0$.

Let $h:=t_1 t_2 \cdots t_\ell$, and let $\OO_{Y,h}$ denote
the sheaf of meromorphic functions whose poles are contained in $T_Y$.
(i.e., $\OO_{Y, h} = \OO_Y[t_1^{-1},\, \dots,\, t_{\ell}^{-1}]$).
We set $\D_{Y,h} := \OO_{Y, h} \underset{\OO_Y}{\otimes} \D_{Y}$.
Since $J_f$ is invertible in $\OO_{Y,h}$, we have
\begin{equation}\label{eq:frac_free}
\D_{Y,h} \underset{\D_Y}{\otimes}\D_{Y \stackrel{f}{\to}X}\simeq \D_{Y,h}.
\end{equation}
(See p. 479 in {\cite{HP13}} also).
As $\CC_{\overline{S_Y}\setminus T_Y} \wtens \OO_Y$ is a $\D_{Y,h}$-module and
$\D_{Y,h}$ is flat over $\D_Y$, we have
$$
\rh_{\D_Y}(\D_{Y \stackrel{f}{\to}X}, \CC_{\overline{S_Y}\setminus T_Y}
\wtens \OO_Y)
\simeq  \rh_{\D_{Y,h}}(\D_{Y,h} \underset{\D_Y}{\otimes}
\D_{Y \stackrel{f}{\to}X},\CC_{\overline{S_Y}\setminus T_Y} \wtens \OO_Y),
$$
and hence, by
\eqref{eq:direct-image-ks}
and
\eqref{eq:frac_free}, we finally obtain
\begin{equation}{\label{eq:equiv-X_Y-vanishing}}
R\Gamma(Y; \CC_{\overline{S_Y} \setminus T_Y} \wtens \OO_Y) \simeq
R\Gamma(X;\CC_{\overline{f(S_Y)}\setminus (Z_1 \cup \cdots \cup Z_\ell)} \wtens \OO_X).
\end{equation}

Now we are ready to give the proof of the theorem.
\begin{proof}
By the fiber formula, it suffices to show
$$
\begin{aligned}
&H^k(\imin\rho\nu^{sa}_\chi(\OW_{X|X\setminus \ckZ}))_p  \\
&\qquad =
\lind{\epsilon_X > 0,\, V_X \in V(\zeta)}
H^k(S_X(V_X,\,\epsilon_X);\, \OW_{X|X\setminus \ckZ})\\
&\qquad =
\lind{\epsilon_X > 0,\, V_X \in V(\zeta)}
H^k(X;\CC_{\overline{S_X(V_X,\,\epsilon_X)}\setminus
(Z_1 \cup \cdots \cup Z_\ell)} \wtens \OO_X) = 0 \quad (k \ne 0).
\end{aligned}
$$
By Lemmas {\ref{lem:set_equiv_f}} and {\ref{lem:set_equiv_f_inv}} and
\eqref{eq:equiv-X_Y-vanishing}, the claim is equivalent to
$$
\lind{\epsilon_Y > 0,\, V_Y \in \Xi}
H^k(Y;\CC_{\overline{S_Y(V_Y,\,\epsilon_Y)}\setminus T_Y} \wtens \OO_Y) = 0 \quad (k \ne 0).
$$
Recall that $V_Y$ has a form $V_1 \times \cdots \times V_\ell$ where
each $V_k$ is a proper open sector in $\mathbb{C}$ containing $1$.
Since we have
$$
\overline{S_Y(V_Y,\,\epsilon_Y)}\setminus T_Y =
B \times (\overline{V_1} \setminus\{0\}) \times \cdots \times
(\overline{V_\ell}\setminus\{0\}),
$$
where $B$ is a closed convex set in $\mathbb{C}^{n-\ell}$, it follows from
the topological tensor product formula in \cite{KS96} that we get
$$
\begin{aligned}
&R\Gamma(Y;\CC_{\overline{S_Y(V_Y,\,\epsilon_Y)}\setminus T_Y} \wtens \OO_Y)  \\
&\quad =
R\Gamma(\mathbb{C}^{n-\ell}; \CC_{B} \wtens \OO) \widehat{\boxtimes}
R\Gamma(\mathbb{C}; \CC_{\overline{V_1}\setminus\{0\}} \wtens \OO)
\widehat{\boxtimes} \cdots
\widehat{\boxtimes}
R\Gamma(\mathbb{C}; \CC_{\overline{V_\ell}\setminus\{0\}} \wtens \OO).
\end{aligned}
$$
As $B$ is a convex closed set,
$R\Gamma(\mathbb{C}^{n-\ell}; \CC_{B} \wtens \OO)$ is concentrated in degree $0$
and the $0$-th cohomology group is naturally equipped with {\bf{FN}} topology.
For $1 \le k \le \ell$,
$R\Gamma(\mathbb{C}; \CC_{\overline{V_k}\setminus\{0\}} \wtens \OO)$ is also known to
be concentrated in degree $0$ (see, for example, Lemma 8.5 \cite{HP13})
and the $0$-th cohomology group has {\bf{FN}} topology.
Hence we can conclude that
$R\Gamma(Y;\CC_{\overline{S_Y(V_Y,\,\epsilon_Y)}\setminus T_Y} \wtens \OO_Y)$
is concentrated in degree $0$.
This completes the proof.
\end{proof}

Thanks to Theorems \ref{teo:vanishing_asymptotic} and \ref{teo: flat asymptotics in deg 0} we obtain the following result

\begin{teo} \label{teo: Borel-Ritt}
Assume the associated action $\mu$ is non-degenerate. Then	
the distinguished triangle \eqref{dtColin} induces an exact sequence (outside the fixed points)
\begin{equation}\label{exspT}
\exs{\imin \rho H^0 \nu^{sa}_\chi\OW_{X|X \setminus \ckZ}}{\imin \rho H^0 \nu^{sa}_\chi\OWX}{\imin \rho H^0 \nu^{sa}_\chi \OW_{X|\ckZ}}.
\end{equation}
All the complexes $\imin \rho \nu^{sa}_\chi\OW_{X|X \setminus \ckZ}$, $\imin \rho \nu^{sa}_\chi\OWX$ and $\imin \rho \nu^{sa}_\chi \OW_{X|\ckZ}$ are concentrated in degree zero.
\end{teo}
Note that, by Theorem \ref{teo:vanishing_asymptotic}, the complex $\imin \rho \nu^{sa}_\chi\OWX$ is always concentrated in degree zero
even if the action $\mu$ is degenerate and the point $p$ is in the fixed points. \\

Let us summarize our previous results. We assume $X = \CC^n$.
Recall that the set of points in $S_\chi$ outside the fixed points is denoted by $S_\chi^\circ$.
Let $W$ be an $(\RP)^\ell$-conic subanalytic open subset in $S_\chi^\circ$.
Let us consider the multi-specialization of Whitney holomorphic functions. We have
$$
\begin{aligned}
\Gamma(W;\imin \rho H^0\nu^{sa}_\chi \OW_{X|X\setminus \ckZ}) & \simeq \lpro {V'}\,\lind {U'} \Gamma(U';\OW_{X|X\setminus \ckZ}),
\end{aligned}
$$
$$
\begin{aligned}
\Gamma(W;\imin \rho H^0\nu^{sa}_\chi \OW_X) & \simeq \lpro {V'}\,\lind {U'} \Gamma(U';\OW_X),
\end{aligned}
$$
$$
\begin{aligned}
\Gamma(W;\imin \rho H^0\nu^{sa}_\chi \OW_{X|\ckZ}) & \simeq \lpro {V'}\,\lind {U'} \Gamma(U';\OW_{X|\ckZ}),
\end{aligned}
$$
where $V'$ ranges through a family of $(\RP)^\ell$-conic subanalytic open subsets
which are compactly generated in $W$,
$U'$ ranges through the family of $\op(X_{sa})$ such that
$C_\chi(X \setminus U') \cap V'=\emptyset$.
Here we say that $V'$ is compactly generated in $W$ if
there exists a compact subset $K \subset W$ with $V' \subset (\RP)^\ell K$.

%The identity morphism induces morphisms of presheaves
%\begin{eqnarray*}
%\widetilde{\A}^{<0}_\chi & \to & \imin \rho H^0\nu^{sa}_\chi \OW_{X|X\setminus \ckZ}, \\
%\widetilde{\A}_\chi & \to & \imin \rho H^0\nu^{sa}_\chi \OW_X, \\
%\widetilde{\A}^{CF}_\chi & \to & \imin \rho H^0\nu^{sa}_\chi \OW_{X|\ckZ}.
%\end{eqnarray*}
%The above morphisms are isomorphisms in the stalks (i.e. in the limit of multicones containing a given direction $p$).
Recall that the classical sheaf $\A^{<0}_\chi$ (resp. $\A_\chi$, resp. $\A^{CF}_\chi$)
on $S_\chi^\circ$ was defined in Subsection {\ref{subsec:classical_sheaves}}.
For $\{S_p\}_{p \in W} \in \mathcal{S}(W)$ (for these definitions, see Subsection
{\ref{subsec:classical_sheaves}}) and a compactly generated $(\RP)^\ell$-conic subset
$V'$ in $W$, we can find finitely many points $P \subset W$ such that
$$
C_\chi\left(X \setminus (\displaystyle\bigcup_{p \in P} S_p)\right) \cap V' = \emptyset.
$$
Hence,
by Theorems {\ref{theorem:multi-asymptotic}},
\ref{teo: consistent family Whitney}
and \ref{teo:flat-holomorphic},
the identity morphism induces morphisms of presheaves on $S_\chi^\circ$
\begin{eqnarray*}
\widetilde{\A}^{<0}_\chi & \to & \imin \rho H^0\nu^{sa}_\chi \OW_{X|X\setminus \ckZ}, \\
\widetilde{\A}_\chi & \to & \imin \rho H^0\nu^{sa}_\chi \OW_X, \\
\widetilde{\A}^{CF}_\chi & \to & \imin \rho H^0\nu^{sa}_\chi \OW_{X|\ckZ}.
\end{eqnarray*}
Since the above morphisms are isomorphisms in the stalks,
we have obtained the isomorphisms of the sheaves on $S_\chi^\circ$
% associated to $\widetilde{\A}^{<0}_\chi$ (resp. $\widetilde{\A}_\chi$, resp. $\widetilde{\A}^{CF}_\chi$). We get ù
\begin{eqnarray*}
\A^{<0}_\chi & \simeq & \imin \rho H^0\nu^{sa}_\chi \OW_{X|X\setminus \ckZ}, \\
\A_\chi & \simeq & \imin \rho H^0\nu^{sa}_\chi \OW_X, \\
\A^{CF}_\chi & \simeq & \imin \rho H^0\nu^{sa}_\chi \OW_{X|\ckZ}.
\end{eqnarray*}
By Proposition \ref{teo: Borel-Ritt}, on $S_\chi^\circ$, we have the exact sequence of sheaves
\begin{equation}\label{eq:Borel-Ritt}
\exs{\A^{<0}_{\chi}}{\A_{\chi}}{\A^{CF}_{\chi}},
\end{equation}
when the associated action $\mu$ is non-degenerate.
The above exact sequence \eqref{eq:Borel-Ritt} is a generalization of the Borel-Ritt exact sequence for multi-asymptotically developable functions.

\subsection{Vanishing theorems for tempered holomorphic functions}\label{subsection:Vanishing theorems for tempered holomorphic functions}

We are now going to consider the vanishing of the cohomology for the multi-specialization  of the sheaf of  tempered holomorphic functions of \cite{KS01}. Let $X$ be a real analytic manifold ($X=\R^n$ in the local model). Let $\db_X$ denote the sheaf of distributions on $X$.

\begin{df} One denotes by $\dbt_X$ the presheaf of tempered
distributions on $X_{sa}$ defined as follows:
$$U \mapsto \Gamma(X;\db_X)/\Gamma_{M\setminus U}(X;\db_X).$$
\end{df}
 As a consequence of the \L ojasievicz's inequalities
\cite{Lo59}, for $U,V \in \op(X_{sa})$ the sequence
$$\exs{\dbt_X(U \cup
V)}{\dbt_X(U)\oplus\dbt_X(V)}{\dbt_X(U \cap V)}$$ is exact. This implies that
$\dbt_X$ is a sheaf on $X_{sa}$. Moreover it follows by definition
that $\dbt_X$ is quasi-injective.
%We have the following result %(see \cite{KS01}, Proposition 7.2.6)

Let $X$ be a complex manifold ($X=\CC^n$ in the local model) and let $X_\R$ denote the underlying real analytic manifold
of $X$. One denotes by $\ot_X$ the sheaf defined as follows:
$$
\ot_X := \rh_{\rho_!\D_{\overline{X}}}(\OO_{\overline{X}},\dbt_{X_\R}).
$$

Let $U$ be an open subset in X. Then $H^0(U;\ot_X)$ consists of tempered distributions on $U$ that satisfy the Cauchy-Riemann system. \\

 Let $\chi=\{Z_1,\ldots,Z_\ell\}$ be a family of complex submanifolds of $X$. Theorem \ref{teo:vanishing_asymptotic} holds for $\OO_X^t$ and $\OO_X$ (the sheaf of holomorphic functions).
\begin{teo} \label{teo: vanishing o}
Under the same geometrical situation as in Theorem \ref{teo:vanishing_asymptotic},
we have
\begin{equation}
H^k(\imin\rho\nu^{sa}_\chi\OO^t_{X}) = 0 \qquad (k \ne 0),
\end{equation}
and
\begin{equation}
H^k(\nu_\chi\OO_{X}) = 0 \qquad (k \ne 0).
\end{equation}
\end{teo}
\begin{proof}
We may assume that all the entries of $A_\chi$ are non-negative integers
and $\comD = 1$.
Let $p = (0;\, \zeta) \in S_\chi$ and $H$ a finite subset in $\RRpair$ for which
$[H]$ and $\mathcal{G}$ are equivalent.
It follows from the fiber formula that
$$
(\imin\rho H^k\nu^{sa}_\chi\OO^t_X)_p
\simeq  \lind {V,\{\epsilon\}} H^k(S_H(V,\{\epsilon\});\, \OO^t_X),
$$
where $V$ runs through $V(\zeta)$ and $\{\epsilon\}$ ranges through
families of positive real numbers.
We may additionally assume that $H$ satisfies $H = Q(H)$,
and thus, a multicone $S_H(V,\{\epsilon\})$ has the form
\begin{equation}
\left\{z \in X;\,
\begin{aligned}
&z \in V,\,\, |z^{(0)}| < \epsilon_0, \\
&f\left(|z^{(*)}|\right) < v+\epsilon_{f,+}\,\,\,
\text{for any $(f,\,v) \in H$}
\end{aligned}
\right\}
\end{equation}
for a convex cone $V \in V(\zeta)$ and a family $\{\epsilon\}$ of positive real numbers.
To prove the theorem, it suffices to show $S_H(V,\{\epsilon\})$ to be Stein.
Set
$$
Y = \{\eta = (\eta_1,\dots,\eta_m) \in \mathbb{C}^m;\,
\eta_k \ne 0\,\,(k \notin J_Z(\zeta))\}.
$$
Then, for each $(f,v) \in H$, $f(\eta)$ is a holomorphic function on $Y$
as $\nu_k(f) \ge 0$ holds for $k \in J_Z(\zeta)$. Let us
consider the set
$$
\begin{aligned}
Y_H &:= \{\eta \in Y;\, |f(\eta)| < v+\epsilon_{f,+}\,\,\,
\text{for any $(f,\,v) \in H$}\}\\
&=
\bigcap_{(f,\,v) \in H} f^{-1}(D_{v+\epsilon_{f,+}}).
\end{aligned}
$$
Here $D_\delta \subset \mathbb{C}$ denotes the open disk
with radius $\delta > 0$ and center at the origin, and each $f$ is regarded as a
holomorphic map from $Y$ to $\mathbb{C}$.
Then $Y_H$ becomes a Stein open subset because the inverse image of a Stein subset
by a holomorphic map on a Stein subset is again Stein
and $Y$ itself is a Stein subset in $\mathbb{C}^m$.
Let $\Theta$ be the set of projections
$$
p(z) = (z_{i_1},\dots, z_{i_m}) \in \mathbb{C}^m, \qquad z = (z_1,\dots,z_n) \in \mathbb{C}^n
$$
from $\mathbb{C}^n$ to $\mathbb{C}^m$ such that $z_{i_k} \in z^{(k)}$ ($k=1,2,\dots,m$).
Note that $\Theta$ is a finite set. Then, by noticing Remark \ref{oss:complex-norm},
we have
$$
S_H(V,\{\epsilon\}) =
\left\{z \in V;\, |z^{(0)}| < \epsilon_0 \right\}\,
\bigcap\,
\left(
\bigcap_{p \in \Theta}\,
p^{-1}(Y_H)
\right),
$$
which shows that $S_H(V,\{\epsilon\})$ is Stein.
This completes the proof.
\end{proof}

\section{Examples}\label{section:Examples}

In order to understand the previous constructions, in this section we construct several examples of multi-specialization and the associated multi-asymptotics.
We shall work in the local model, in particular we are going to consider examples in $\CC^n$ starting from the matrices defining the actions and the associated multi-normal deformations.

\subsection{Majima's asymptotics}\label{subsection:Majima's asymptotics}

Let us consider $X=\CC^n$ with coordinates $z=(z_1,\dots,z_n)$. We are going to check how our previous constructions permit to obtain strongly asymptotically developable functions in the sense of Majima \cite{Ma84}. Let $\chi=\{Z_1,\dots,Z_n\}$, $Z_i=\{z_i=0\}$, $Z=\bigcap_{i=1}^nZ_i$. Set $\widetilde{X}=\CC^n \times \R^n$ and consider the action
\begin{eqnarray*}
\mu_\chi:\CC^n \times \R^n & \to & \CC^n, \\
(z_1,\dots,z_n,t_1,\dots,t_n) & \mapsto & (z_1t_1,\dots,z_nt_n).
\end{eqnarray*}
The monomials $\varphi_i$, $i=1,\dots,n$, defining the action are $\varphi_i(t)=t_i$, $i=1,\dots,n$, and the associated matrix $A_\chi$ is nothing but the identity matrix. By Proposition \ref{prop:bundle-description} the zero section $S_\chi$ of $\widetilde{X}$ has a vector bundle structure
$$
S_\chi \simeq T_{Z_1}X \underset{X}{\times} \cdots \underset{X}{\times} T_{Z_n}X.
$$
Following the notations of Section \ref{section:Multi-cones for a general case} we are going to study the family of multicones $C(\xi)$ associated to a point $\xi=(\xi_1,\dots,\xi_n) \in S_\chi$, $\xi_i \neq 0$, $i=1,\dots,n$. First of all we easily see that the rational monomials $\imin{\varphi_j}$, $j=1,\dots,n$ are $\imin{\varphi_j}(\tau)=\tau_j$ (no monomials $\psi_k$ appear since $A_\chi$ is an invertible square matrix),
and hence $F^0=\{\tau_1,\dots,\tau_n\}$.
A cofinal family of $C(\xi)$ is $C(\xi,F^0)$ defined as follows
$$
\begin{aligned}
& \quad \left\{(z_1,\dots,z_n) \in X;
\begin{aligned}
&z_i \in W_i\quad(k=1,\dots,n), \\
&|z_i|< \epsilon
\end{aligned}
\right\},
\end{aligned}
$$
where $\epsilon > 0$ and each $W_i$ is an $\mathbb{R}^+$-conic open
subset in $T_{Z_i}X \simeq \CC$ containing the point $\xi_i$. It means that a cofinal family of $C(\xi)$ is given by the family of multisectors in $\CC^n$ containing the direction $\xi$. \\

We are now ready to compute the fibers of the specialization of Whitney holomorphic functions. Let $\OW_X$ be the subanalytic sheaf of Whitney holomorphic functions on $X$. Let $\xi \in S_\chi \simeq \CC^n$. Then
$$
H^0(\rho^{-1}\nu^{sa}_\chi\OW_X)_\xi \simeq \lind S H^0(S;\OW_X) \simeq \lind S \lpro {S'}H^0(S';\OW_X),
$$
where $S=S_1 \times \cdots \times S_n$ ranges through the family of multisectors containing $\xi$ and $S'=S'_1 \times \cdots \times S'_n$ ranges through the family of multisectors properly contained in $S$ (i.e. $\overline{S'_i} \setminus \{0\} \subset S_i$, $i=1,\dots,n$). By Theorem \ref{theorem:multi-asymptotic},
$$
\lpro {S'}H^0(S';\OW_X) \simeq \Gamma(S;\mathcal{A}_\chi),
$$
where $\mathcal{A}_\chi$ denotes the sheaf of strongly asymptotically developable functions in the sense of Majima \cite{Ma84}. Moreover $H^k(\imin\rho\nu^{sa}_\chi(\OW_{X})) = 0$, $k \neq 0$, in view of Theorem \ref{teo:vanishing_asymptotic}, hence we get the isomorphism of sheaves (outside the fixed points of $S_\chi$)
$$
\imin\rho\nu^{sa}_\chi\OW_X \simeq \A_\chi.
$$
Moreover, following the notations of Section \ref{section:Multi-specialization}, by Theorem \ref{teo: Borel-Ritt} the sequence
$$
\exs{\imin\rho \nu^{sa}_\chi\OW_{X|X \setminus \ckZ}}{\imin\rho\nu^{sa}_\chi\OW_X}{\imin\rho\nu^{sa}_\chi\OW_{X|\ckZ}}
$$
is exact. This corresponds to the sequence
$$
\exs{\A_\chi^{<0}}{\A_\chi}{\A^{CF}_\chi},
$$
where $\A_\chi^{<0}$ and $\A^{CF}_\chi$ denote the sheaves of flat asymptotics and consistent families of coefficients respectively. This is nothing but the Borel-Ritt exact sequence for strongly asymptotically developable holomorphic functions.

\subsection{Examples in $\CC^2$}

Let us consider some interesting cases in $\CC^2$ with variables $z=(z_1, z_2)$. We are going to consider the normal deformation (and the associated multi-specialization and asymptotics) for some interesting cases.
%\begin{itemize}
%%%%%%%%%%%%%%%%%%%%%%%%%%%%%%%%%%%%%%%%%%%
% MAJIMA C2
%%%%%%%%%%%%%%%%%%%%%%%%%%%%%%%%%%%%%%%%%%%
%\item

\

\noindent
$\checkmark\,${\bf{Majima asymptotics.}} In this case
$$
A_\chi=
\left(
\begin{array}{ccc}
1 & 0 \\
0 & 1
\end{array}
\right).
$$
We have $Z_i=\{z_i=0\}$, $i=1,2$. Then we can define a  normal deformation $\widetilde{X}=\CC^2 \times \R^2$ with the map $p(z;\,t)=(t_1z_1, t_2z_2)$. By Proposition \ref{prop:bundle-description} the zero section $S_\chi$ of $\widetilde{X}$ has a vector bundle structure
$$
S_\chi \simeq T_{Z_1}X \underset{X}{\times} T_{Z_2}X.
$$
The rational monomials $\imin{\varphi_j}$, $j=1,2$, are $\imin{\varphi_1}(\tau)=\tau_1$, $\imin{\varphi_2}(\tau)=\tau_2$. Let $\xi=(\xi_1,\xi_2)$, $\xi_i \neq 0$, $i=1,2$. A cofinal family of $C(\xi)$ is $C(\xi,F^0)$ consisting of elements $S(W,\epsilon)$, $W=W_1 \times W_2$, defined as follows
$$
S(W,\epsilon)=
\begin{aligned}
 \left\{z \in X;
\begin{aligned}
& z_k \in W_k \quad (k=1,2), \\
& |z_i| < \epsilon \quad (i=1,2)
\end{aligned}
\right\},
\end{aligned}
$$
where $\epsilon > 0$ and each $W_k$ is an $\mathbb{R}^+$-conic open subset in $\CC$ containing the point $\xi_k$. \\

We now construct the asymptotics.
%We have $K_{\{j\}}=K_j$, $j=1,2$, $K_{\{1,2\}}=K_1 \cup K_2$ and $Z_J=\{z^{(k)}=0,\,k \in K_J\}$ for $\emptyset \neq J \subseteq \{1,2\}$. For any subset $K \subset \{1,2\}$, we denote by $z^{(K)}$ the set of the coordinates $z^{(k)}$'s with $k \in K$ and by $z^{(K)}_{\mathrm{C}}$ the set of the coordinates which do not belong to $z^{(K)}$, i.e., $z^{(k)}$'s with $k \in \{0,1,\dots,m\} \setminus K$. Hence, the coordinates of $Z_J$ are given by $z^{(K_J)}_{\mathrm{C}}$. The $\pi_J$ denotes the canonical projection from $X$ to $Z_J$ defined by $z \to z^{(K_J)}_{\mathrm{C}}$. Given $S(W,\epsilon)$, set $S_J=\pi_J(S(W,\epsilon))$.
We have $S_{\{1,2\}}=\{0\}$ and, for $i \neq j \in \{1,2\}$
$$
S_{\{i\}}=
\begin{aligned}
 \left\{z \in Z_i;
\begin{aligned}
&z_j \in W_j, \\
& |z_j| < \epsilon
\end{aligned}
\right\}.
\end{aligned}
$$
A total family of coefficients of multi-asymptotic expansion is given by
\begin{eqnarray*}
F & = & \{F_{\{1\}},F_{\{2\}},F_{\{1,2\}}\} \\
& = & \left\{\{f_{\{1\},k}(z_2)\}_{k \in \Z_{\geq 0}},
\{f_{\{2\},k}(z_1)\}_{k \in
\Z_{\geq 0}},
\{f_{\{1,2\},\alpha}\}_{\alpha \in \Z_{\geq 0}^2}\right\},
\end{eqnarray*}
where $f_{\{1\},k}(z_2)$ (resp. $f_{\{2\},k}(z_1)$) is holomorphic in $S_{\{1\}}$ (resp. $S_{\{2\}}$) and $f_{\{1,2\},\alpha} \in \CC$. An asymptotic expansion $\operatorname{App}^{<N}(F;\,z)$, $N=(n_1,n_2) \in \Z_{\geq 0}^2$ is given by
$$
\begin{aligned}
	T_{\{1\}}^{<N}(F;\, z) &= \sum_{k < n_1} f_{\{1\},k}(z_2)\frac{z_1^k}{k!}, \qquad k \in \Z_{\geq 0}, \\
	T_{\{2\}}^{<N}(F;\, z) &= \sum_{k < n_2} f_{\{2\},k}(z_1)\frac{z_2^k}{k!}, \qquad k \in \Z_{\geq 0}, \\		
T_{\{1,2\}}^{<N}(F;\, z) &= \sum_{\substack{\alpha_1 < n_1 \\ \alpha_2 < n_2}} f_{\{1,2\},\alpha}\frac{z_1^{\alpha_1}z_2^{\alpha_2}}{\alpha_1!\alpha_2!}, \qquad \alpha_1,\alpha_2 \in \Z_{\geq 0}, \\
	\operatorname{App}^{<N}(F;\, z) &=
T_{\{1\}}^{<N}(F;\, z) + T_{\{2\}}^{<N}(F;\, z) - T_{\{1,2\}}^{<N}(F;\, z).
\end{aligned}
$$
We say that $f$ is multi-asymptotically developable to $F=\{F_{\{1\}},F_{\{2\}},F_{\{1,2\}}\}$ along
$\chi$ on $S=S(W,\epsilon)$ if and only if for any cone
$S'=S(W',\epsilon')$ properly contained in $S$ and for any $N=(n_1,n_2)
\in \mathbb{Z}^2_{\ge 0}$, there exists a constant $C_{S',N}$ such that
$$
\left|
f(z) - \operatorname{App}^{<N}(F;\, z)
\right|
\le C_{S',N}|z_1|^{n_1}|z_2|^{n_2}
 \qquad (z \in S').
$$
The family $F$ is consistent if
\begin{itemize}
\item $f_{\{1\},k}(z_2)$ is strongly asymptotically developable to
     $$
     \{f_{\{1,2\},(k,\alpha_2)}\}_{\alpha_2 \in \Z_{\geq 0}}
     $$
     on $S_{\{1\}}$ for each $k \in \Z_{\geq 0}$,
\item $f_{\{2\},k}(z_1)$ is strongly asymptotically developable to
    $$
    \{f_{\{1,2\},(\alpha_1,k)}\}_{\alpha_1 \in \Z_{\geq 0}}
    $$
    on $S_{\{2\}}$ for each $k \in \Z_{\geq 0}$.
\end{itemize}
These are the asymptotics of \cite{Ma84} in $\CC^2$.

%%%%%%%%%%%%%%%%%%%%%%%%%%%%%%%%%%%%%%%%%%%
% TAKEUCHI C2
%%%%%%%%%%%%%%%%%%%%%%%%%%%%%%%%%%%%%%%%%%%

\

\noindent
$\checkmark\,${\bf{Takeuchi asymptotics.}} In this case
$$
A_\chi=
\left(
\begin{array}{ccc}
1 & 1 \\
0 & 1
\end{array}
\right).
$$
We have $Z_1=\{0\}$, $Z_2=\{z_2=0\}$. Then we can define a  normal deformation $\widetilde{X}=\CC^2 \times \R^2$ with the map $p(z;\,t)=(t_1z_1, t_1t_2z_2)$. By Proposition \ref{prop:bundle-description} the zero section $S_\chi$ of $\widetilde{X}$ has a vector bundle structure
$$
S_\chi \simeq T_{Z_1}Z_2 \underset{X}{\times} T_{Z_2}X.
$$
The rational monomials $\imin{\varphi_j}$, $j=1,2$, are $\imin{\varphi_1}(\tau)=\tau_1$, $\imin{\varphi_2}(\tau)=\displaystyle{\tau_2 \over \tau_1}$. Let $\xi=(\xi_1,\xi_2)$, $\xi_i \neq 0$, $i=1,2$. A cofinal family of $C(\xi)$ is $C(\xi,F^0)$ consisting of elements $S(W,\epsilon)$, $W=W_1 \times W_2$, defined as follows
$$
S(W,\epsilon)=
\begin{aligned}
 \left\{z \in X;
\begin{aligned}
& z_k \in W_k \quad (k=1,2), \\
& |z_1| < \epsilon, \\
& |z_2| < \epsilon |z_1|
\end{aligned}
\right\},
\end{aligned}
$$
where $\epsilon > 0$ and each $W_k$ is an $\mathbb{R}^+$-conic open subset in $\CC$ containing the point $\xi_k$. \\

We now construct the asymptotics.
%We have $K_{\{j\}}=K_j$, $j=1,2$, $K_{\{1,2\}}=K_1 \cup K_2$ and $Z_J=\{z^{(k)}=0,\,k \in K_J\}$ for $\emptyset \neq J \subseteq \{1,2\}$. For any subset $K \subset \{1,2\}$, we denote by $z^{(K)}$ the set of the coordinates $z^{(k)}$'s with $k \in K$ and by $z^{(K)}_{\mathrm{C}}$ the set of the coordinates which do not belong to $z^{(K)}$, i.e., $z^{(k)}$'s with $k \in \{0,1,\dots,m\} \setminus K$. Hence, the coordinates of $Z_J$ are given by $z^{(K_J)}_{\mathrm{C}}$. The $\pi_J$ denotes the canonical projection from $X$ to $Z_J$ defined by $z \to z^{(K_J)}_{\mathrm{C}}$. Given $S(W,\epsilon)$, set $S_J=\pi_J(S(W,\epsilon))$.
We have $S_{\{1\}}=S_{\{1,2\}}=\{0\}$ and
$$
S_{\{2\}}=
\begin{aligned}
 \left\{z \in Z_2;
\begin{aligned}
&z_1 \in W_1, \\
& |z_1| < \epsilon
\end{aligned}
\right\}.
\end{aligned}
$$
A total family of coefficients of multi-asymptotic expansion is given by
\begin{eqnarray*}
F & = & \{F_{\{1\}},F_{\{2\}},F_{\{1,2\}}\} \\
& = & \left\{\{f_{\{1\},\alpha}\}_{\alpha \in \Z_{\geq 0}^2},
\{f_{\{2\},k}(z_1)\}_{k \in
\Z_{\geq 0}},
\{f_{\{1,2\},\alpha}\}_{\alpha \in \Z_{\geq 0}^2}\right\},
\end{eqnarray*}
where $f_{\{2\},k}(z_1)$ is holomorphic in $S_{\{2\}}$ and $f_{\{1\},\alpha},f_{\{1,2\},\alpha} \in \CC$. An asymptotic expansion $\operatorname{App}^{<N}(F;\,z)$, $N=(n_1,n_2) \in \Z_{\geq 0}^2$ is given by
$$
\begin{aligned}
	T_{\{1\}}^{<N}(F;\, z) &= \sum_{\alpha_1 + \alpha_2 < n_1} f_{\{1\},\alpha}\frac{z_1^{\alpha_1}z_2^{\alpha_2}}{\alpha_1!\alpha_2!}, \qquad \alpha=(\alpha_1,\alpha_2) \in \Z_{\geq 0}^2, \\
	T_{\{2\}}^{<N}(F;\, z) &= \sum_{k < n_2} f_{\{2\},k}(z_1)\frac{z_2^{k}}{k!}, \qquad k \in \Z_{\geq 0}, \\		
T_{\{1,2\}}^{<N}(F;\, z) &= \sum_{\substack{\alpha_1+\alpha_2 < n_1 \\ \alpha_2 < n_2}} f_{\{1,2\},\alpha}\frac{z_1^{\alpha_1}z_2^{\alpha_2}}{\alpha_1!\alpha_2!}, \qquad \alpha=(\alpha_1,\alpha_2) \in \Z_{\geq 0}^2, \\
	\operatorname{App}^{<N}(F;\, z) &=
T_{\{1\}}^{<N}(F;\, z) + T_{\{2\}}^{<N}(F;\, z) - T_{\{1,2\}}^{<N}(F;\, z).
\end{aligned}
$$
We say that $f$ is multi-asymptotically developable to $F=\{F_{\{1\}},F_{\{2\}},F_{\{1,2\}}\}$ along
$\chi$ on $S=S(W,\epsilon)$ if and only if for any cone
$S'=S(W',\epsilon')$ properly contained in $S$ and for any $N=(n_1,n_2)
\in \mathbb{Z}^2_{\ge 0}$, there exists a constant $C_{S',N}$ such that
$$
\left|
f(z) - \operatorname{App}^{<N}(F;\, z)
\right|
\le C_{S',N}|z_1|^{n_1-n_2}|z_2|^{n_2}
 \qquad (z \in S').
$$
The family $F$ is consistent if
\begin{itemize}
\item $f_{\{1\},\alpha}=f_{\{1,2\},\alpha}$ for each $\alpha \in \Z^2_{\geq 0}$,
\item $f_{\{2\},k}(z_1)$ is strongly asymptotically developable to
    $$
    \{f_{\{1,2\},(\alpha_1,k)}\}_{\alpha_1 \in \Z_{\geq 0}}
    $$
    on $S_{\{2\}}$ for each $k \in \Z_{\geq 0}$.
\end{itemize}
These are the asymptotics of \cite{HP13} associated to the bispecialization of \cite{ST94} in $\CC^2$.\\

Set $X=Y=\CC^2$, $X$ (resp. $Y$) with coordinates $z=(z_1,z_2)$ (resp. $w=(w_1,w_2)$). Let $\chi_Y=\{Z_{Y1},Z_{Y2}\}$, with $Z_{Y1}=\{0\}$, $Z_{Y2}=\{w_{2}=0\}$ and $\chi_X=\{Z_{X1},Z_{X2}\}$, with $Z_{X1}=\imin f(0)=\{z_{1}=0\}$, $Z_{X2}=\{z_{2}=0\}$. We have
$$
A_{\chi_X}=
\left(
\begin{array}{cc}
1 & 0 \\
0 & 1
\end{array}
\right),
\ \ \ \
A_{\chi_Y}=
\left(
\begin{array}{cc}
1 & 1 \\
0 & 1
\end{array}
\right).
$$
Consider the map
\begin{eqnarray*}
f:X & \to & Y \\
(z_{1},z_{2}) & \mapsto & (z_{1},z_{1}z_{2})
\end{eqnarray*}
(locally) the blow-up at the origin. Remark that condition \eqref{eq:cond existence} holds.
Let us consider $\widetilde{f}$. On the zero section (we keep the same coordinates for simplicity) the map $T_\chi f$ is defined as follows:
\begin{eqnarray*}
w_{1} & = & {\partial f \over \partial z_{1}}(0,0) z_{1} = z_{1}, \\
w_{2} & = & {\partial^2f \over \partial z_{1} \partial z_{2}}(0,0) z_{1}z_{2} = z_{1}z_{2}.
\end{eqnarray*}
It is a conic map with respect to the $(\RP)^2$-actions on $S_{\chi_X}$ and $S_{\chi_Y}$. By Propositions \ref{prop: direct image} and \ref{prop: inverse image}, $f$ makes a link between specializations with respect to $\chi_X$ and $\chi_Y$.
%%%%%%%%%%%%%%%%%%%%%%%%%%%%%%%%%%%%%%%%%%%%
% CUPS C2
%%%%%%%%%%%%%%%%%%%%%%%%%%%%%%%%%%%%%%%%%%%%

\

\noindent
$\checkmark\,${\bf{Cusp asymptotics.}} In this case
$$
A_\chi=
\left(
\begin{array}{ccc}
3 & 2 \\
1 & 1
\end{array}
\right).
$$
We have $Z_1=Z_2=\{0\}$. Then we can define a  normal deformation $\widetilde{X}=\CC^2 \times \R^2$ with the map $p(z;\,t)=(t_1^3t_2z_1, t_1^2t_2z_2)$.
The rational monomials $\imin{\varphi_j}$, $j=1,2$, are $\imin{\varphi_1}(\tau)=\displaystyle{\tau_1 \over \tau_2}$, $\imin{\varphi_2}(\tau)=\displaystyle{\tau_2^3 \over \tau_1^2}$. Let $\xi=(\xi_1,\xi_2)$, $\xi_i \neq 0$, $i=1,2$. A cofinal family of $C(\xi)$ is $C(\xi,F^0)$ consisting of elements $S(W,\epsilon)$, $W=W_1 \times W_2$, defined as follows
$$
S(W,\epsilon)=
\begin{aligned}
 \left\{z \in X;
\begin{aligned}
& z_k \in W_k \quad (k=1,2), \\
& |z_1| < \epsilon |z_2|, \\
& |z_2|^3 < \epsilon |z_1|^2
\end{aligned}
\right\},
\end{aligned}
$$
where $\epsilon > 0$ and each $W_k$ is an $\mathbb{R}^+$-conic open subset in $\CC$ containing the point $\xi_k$. \\

We now construct the asymptotics.
%We have $K_{\{j\}}=K_j$, $j=1,2$, $K_{\{1,2\}}=K_1 \cup K_2$ and $Z_J=\{z^{(k)}=0,\,k \in K_J\}$ for $\emptyset \neq J \subseteq \{1,2\}$. For any subset $K \subset \{1,2\}$, we denote by $z^{(K)}$ the set of the coordinates $z^{(k)}$'s with $k \in K$ and by $z^{(K)}_{\mathrm{C}}$ the set of the coordinates which do not belong to $z^{(K)}$, i.e., $z^{(k)}$'s with $k \in \{0,1,\dots,m\} \setminus K$. Hence, the coordinates of $Z_J$ are given by $z^{(K_J)}_{\mathrm{C}}$. The $\pi_J$ denotes the canonical projection from $X$ to $Z_J$ defined by $z \to z^{(K_J)}_{\mathrm{C}}$. Given $S(W,\epsilon)$, set $S_J=\pi_J(S(W,\epsilon))$.
We have $S_{\{1\}}=S_{\{2\}}=S_{\{1,2\}}=\{0\}$.
A total family of coefficients of multi-asymptotic expansion is given by
\begin{eqnarray*}
F & = & \{F_{\{1\}},F_{\{2\}},F_{\{1,2\}}\} \\
& = & \left\{\{f_{\{1\},\alpha}\}_{\alpha \in \Z_{\geq 0}^2},
\{f_{\{2\},\alpha}\}_{\alpha \in \Z_{\geq 0}^2},
\{f_{\{1,2\},\alpha}\}_{\alpha \in \Z_{\geq 0}^2}\right\},
\end{eqnarray*}
where $f_{\{1\},\alpha},f_{\{2\},\alpha},f_{\{1,2\},\alpha} \in \CC$. An asymptotic expansion $\operatorname{App}^{<N}(F;\,z)$, $N=(n_1,n_2) \in \Z_{\geq 0}^2$ is given by
$$
\begin{aligned}
	T_{\{1\}}^{<N}(F;\, z) &= \sum_{3\alpha_1 + 2\alpha_2 < n_1} f_{\{1\},\alpha}\frac{z_1^{\alpha_1}z_2^{\alpha_2}}{\alpha_1!\alpha_2!}, \qquad \alpha=(\alpha_1,\alpha_2) \in \Z_{\geq 0}^2, \\
	T_{\{2\}}^{<N}(F;\, z) &= \sum_{\alpha_1 + \alpha_2 < n_2} f_{\{2\},\alpha}\frac{z_1^{\alpha_1}z_2^{\alpha_2}}{\alpha_1!\alpha_2}, \qquad \alpha=(\alpha_1,\alpha_2) \in \Z_{\geq 0}^2, \\		
T_{\{1,2\}}^{<N}(F;\, z) &= \sum_{\substack{3\alpha_1+2\alpha_2 < n_1 \\ \alpha_1+\alpha_2 < n_2}} f_{\{1,2\},\alpha}\frac{z_1^{\alpha_1}z_2^{\alpha_2}}{\alpha_1!\alpha_2!}, \qquad \alpha=(\alpha_1,\alpha_2) \in \Z_{\geq 0}^2, \\
	\operatorname{App}^{<N}(F;\, z) &=
T_{\{1\}}^{<N}(F;\, z) + T_{\{2\}}^{<N}(F;\, z) - T_{\{1,2\}}^{<N}(F;\, z).
\end{aligned}
$$
We say that $f$ is multi-asymptotically developable to $F=\{F_{\{1\}},F_{\{2\}},F_{\{1,2\}}\}$ along
$\chi$ on $S=S(W,\epsilon)$ if and only if for any cone
$S'=S(W',\epsilon')$ properly contained in $S$ and for any $N=(n_1,n_2)
\in \mathbb{Z}^2_{\ge 0}$, there exists a constant $C_{S',N}$ such that
$$
\left|
f(z) - \operatorname{App}^{<N}(F;\, z)
\right|
\le C_{S',N}|z_1|^{n_1-2n_2}|z_2|^{3n_2-n_1}
 \qquad (z \in S'). \\
$$
The family $F$ is consistent if
\begin{itemize}
\item $f_{\{1\},\alpha}=f_{\{2\},\alpha}=f_{\{1,2\},\alpha}$  for each $\alpha \in \Z_{\geq 0}^2$.
\end{itemize}

\

Set $X=Y=\CC^2$, $X$ (resp. $Y$) with coordinates $z=(z_1,z_2)$ (resp. $w=(w_1,w_2)$). Let $\chi_Y=\{Z_{Y1},Z_{Y2}\}$, with $Z_{Y1}=Z_{Y2}=\{0\}$ and $\chi_X=\{Z_{X1},Z_{X2}\}$, with $Z_{X1}=\{z_1=0\}$, $Z_{X2}=\{z_2=0\}$. We have
$$
A_{\chi_X}=
\left(
\begin{array}{cc}
1 & 0 \\
0 & 1
\end{array}
\right),
\ \ \ \
A_{\chi_Y}=
\left(
\begin{array}{cc}
3 & 2 \\
1 & 1
\end{array}
\right).
$$
Consider the map
\begin{eqnarray*}
f:X & \to & Y \\
(z_{1},z_{2}) & \mapsto & (z_{1}^3z_2,z_{1}^2z_{2})
\end{eqnarray*}
(desingularization of a cusp at the origin). Remark that condition \eqref{eq:cond existence} holds.
Let us consider $\widetilde{f}$. On the zero section (we keep the same coordinates for simplicity) the map $T_\chi f$ is defined as follows:
\begin{eqnarray*}
w_1 & = & {1 \over 4!} \cdot 4 \cdot {\partial^4 f_1 \over \partial z_1^3\partial z_2}(0,0)z_1^3z_2={z_1^3z_2} \\
w_2 & = & {1 \over 3!} \cdot 3 \cdot {\partial^3 f_1 \over \partial z_1^2\partial z_2}(0,0)z_1^2z_2={z_1^2z_2}.
\end{eqnarray*}
It is a conic map with respect to the $(\RP)^2$-actions on $S_{\chi_X}$ and $S_{\chi_Y}$. By Propositions \ref{prop: direct image} and \ref{prop: inverse image}, $f$ makes a link between specializations with respect to $\chi_X$ and $\chi_Y$.

%\end{itemize}

%%%%%%%%%%%%%%%%%%%%%%%%%%%%%%%%%%%%%%%%%%%%
% END OF C2
%%%%%%%%%%%%%%%%%%%%%%%%%%%%%%%%%%%%%%%%%%%%

\subsection{Examples in $\CC^3$}

Let us consider some interesting cases in $\CC^3$ with variables $z=(z_1, z_2,z_3)$. We are going to consider the normal deformation (and the associated multi-specialization and asymptotics) for some interesting cases. We shall consider $2 \times 3$ and $3 \times 3$ matrices with entries $0,1$, avoiding the degenerate cases.
%\begin{itemize}
%%%%%%%%%%%%%%%%%%%%%%%%%%%%%%%%%%%%%%%%%%%%
%  CLEAN INTERSECTION C3
%%%%%%%%%%%%%%%%%%%%%%%%%%%%%%%%%%%%%%%%%%%%

\

\noindent
$\checkmark\,${\bf{Clean intersection asymptotics (2 lines).}} In this case
$$
A_\chi=
\left(
\begin{array}{ccc}
1 & 1 & 0 \\
0 & 1 & 1
\end{array}
\right).
$$
We have $Z_1=\{z_1=z_2=0\}$, $Z_2=\{z_2=z_3=0\}$. Then we can define a  normal deformation $\widetilde{X}=\CC^3 \times \R^2$ with the map $p(z;\,t)=(t_1z_1, t_1t_2z_2, t_2z_3)$.
The rational monomials $\imin{\varphi_j}$, $j=1,2$, are $\imin{\varphi_1}(\tau)=\tau_1$, $\imin{\varphi_2}(\tau)=\displaystyle{\tau_2 \over \tau_1}$ and $\psi_3(\tau)=\displaystyle{\tau_1\tau_3 \over \tau_2}$. Let $\xi=(\xi_1,\xi_2,\xi_3)$, $\xi_i \neq 0$, $i=1,2$. A cofinal family of $C(\xi)$ is $C(\xi,F^0)$ consisting of elements $S(W,\epsilon)$, $W=W_1 \times W_2 \times W_3$, defined as follows:
$$
S(W,\epsilon)=
\begin{aligned}
 \left\{z \in X;
\begin{aligned}
& z_k \in W_k \quad (k=1,2,3), \\
& |z_1| < \epsilon,  \\
& |z_2| < \epsilon |z_1|, \\
& (n_3(\xi)-\epsilon)|z_2| < |z_1||z_3| < (n_3(\xi)+\epsilon) |z_2|
\end{aligned}
\right\},
\end{aligned}
$$
where $n_3(\xi) := |\xi_1\xi_3/\xi_2|$, $\epsilon > 0$ and each $W_k$ is an $\mathbb{R}^+$-conic open subset in $\CC$
containing the point $\xi_k$. \\

We now construct the asymptotics.
%We have $K_{\{j\}}=K_j$, $j=1,2$, $K_{\{1,2\}}=K_1 \cup K_2$ and $Z_J=\{z^{(k)}=0,\,k \in K_J\}$ for $\emptyset \neq J \subseteq \{1,2\}$. For any subset $K \subset \{1,2\}$, we denote by $z^{(K)}$ the set of the coordinates $z^{(k)}$'s with $k \in K$ and by $z^{(K)}_{\mathrm{C}}$ the set of the coordinates which do not belong to $z^{(K)}$, i.e., $z^{(k)}$'s with $k \in \{0,1,\dots,m\} \setminus K$. Hence, the coordinates of $Z_J$ are given by $z^{(K_J)}_{\mathrm{C}}$. The $\pi_J$ denotes the canonical projection from $X$ to $Z_J$ defined by $z \to z^{(K_J)}_{\mathrm{C}}$. Given $S(W,\epsilon)$, set $S_J=\pi_J(S(W,\epsilon))$.
We have $S_{\{1,2\}}=\{0\}$ and
$$
S_{\{1\}}=
\begin{aligned}
 \left\{z \in Z_1;
\begin{aligned}
&z_3 \in W_3, \\
& |z_3| < \epsilon
\end{aligned}
\right\},
\end{aligned}
$$
$$
S_{\{2\}}=
\begin{aligned}
 \left\{z \in Z_2;
\begin{aligned}
&z_1 \in W_1, \\
& |z_1| < \epsilon
\end{aligned}
\right\}.
\end{aligned}
$$

A total family of coefficients of multi-asymptotic expansion is given by
\begin{eqnarray*}
F & = & \{F_{\{1\}},F_{\{2\}},F_{\{1,2\}}\} \\
& = & \Big\{\{f_{\{1\},\alpha}(z_3)\}_{\alpha \in \Z^2_{\geq 0}},
\{f_{\{2\},\alpha}(z_1)\}_{\alpha \in
\Z^2_{\geq 0}},\{f_{\{1,2\},\beta}\}_{\beta \in \Z_{\geq 0}^3}\Big\},
\end{eqnarray*}
where $f_{\{1\},\alpha}(z_3)$ (resp. $f_{\{2\},\alpha}(z_1)$) is holomorphic in $S_{\{1\}}$ (resp. $S_{\{2\}}$) and $f_{\{1,2\},\beta} \in \CC$. An asymptotic expansion $\operatorname{App}^{<N}(F;\,z)$, $N=(n_1,n_2) \in \Z_{\geq 0}^2$ is given by
$$
\begin{aligned}
	T_{\{1\}}^{<N}(F;\, z) &= \sum_{\alpha_1+\alpha_2 < n_1} f_{\{1\},\alpha}(z_3)\frac{z_1^{\alpha_1}z_2^{\alpha_2}}{\alpha_1!\alpha_2!}, \qquad \alpha_1,\alpha_2 \in \Z_{\geq 0}, \\
	T_{\{2\}}^{<N}(F;\, z) &= \sum_{\alpha_2+\alpha_3 < n_2} f_{\{2\},\alpha}(z_1)\frac{z_1^{\alpha_2}z_2^{\alpha_3}}{\alpha_2!\alpha_3!}, \qquad \alpha_2,\alpha_3 \in \Z_{\geq 0}, \\	
	T_{\{1,2\}}^{<N}(F;\, z) &= \sum_{\substack{\beta_1+\beta_2 < n_1 \\ \beta_2+\beta_3 < n_2}} f_{\{1,2\},\beta}\frac{z_1^{\beta_1}z_2^{\beta_2}z_3^{\beta_3}}{\beta_1!\beta_2!\beta_3!}, \qquad \beta_1,\beta_2,\beta_3 \in \Z_{\geq 0}, \\
	\operatorname{App}^{<N}(F;\, z) &=
T_{\{1\}}^{<N}(F;\, z) + T_{\{2\}}^{<N}(F;\, z)- T_{\{1,2\}}^{<N}(F;\, z).
\end{aligned}
$$
We say that $f$ is multi-asymptotically developable to $F$ along
$\chi$ on $S=S(W,\epsilon)$ if and only if for any cone
$S'=S(W',\epsilon')$ properly contained in $S$ and for any $N=(n_1,n_2)
\in \mathbb{Z}^2_{\ge 0}$, there exists a constant $C_{S',N}$ such that
$$
\left|
f(z) - \operatorname{App}^{<N}(F;\, z)
\right|
\le C_{S',N}|z_1|^{n_1-n_2}|z_2|^{n_2}
 \qquad (z \in S').
$$
The family $F$ is consistent if
\begin{itemize}
\item $f_{\{1\},(\alpha_1,\alpha_2)}(z_3)$ is strongly asymptotically developable to
    $$
    \{f_{\{1,2\},(\alpha_1,\alpha_2,\beta_3)}\}_{\beta_3 \in \Z_{\geq 0}}
    $$
    on $S_{\{1\}}$ for each $(\alpha_1,\alpha_2) \in \Z^2_{\geq 0}$,
\item $f_{\{2\},(\alpha_2,\alpha_3)}(z_1)$ is strongly asymptotically developable to
    $$
    \{f_{\{1,2\},(\beta_1,\alpha_2,\alpha_3)}\}_{\beta_1 \in \Z_{\geq 0}}
    $$
    on $S_{\{2\}}$ for each $(\alpha_2,\alpha_3) \in \Z^2_{\geq 0}$.
\end{itemize}
%%%%%%%%%%%%%%%%%%%%%%%%%%%%%%%%%%%%%%
% MAJIMA ASYMPTOTIC C3
%%%%%%%%%%%%%%%%%%%%%%%%%%%%%%%%%%%%%%

\

\noindent
$\checkmark\,${\bf{Majima asymptotics.}} In this case
$$
A_\chi=
\left(
\begin{array}{ccc}
1 & 0 & 0 \\
0 & 1 & 0 \\
0 & 0 & 1
\end{array}
\right).
$$
We have $Z_i=\{z_i=0\}$, $i=1,2,3$. Then we can define a  normal deformation $\widetilde{X}=\CC^3 \times \R^3$ with the map $p(z;\,t)=(t_1z_1, t_2z_2, t_3z_3)$. By Proposition \ref{prop:bundle-description} the zero section $S_\chi$ of $\widetilde{X}$ has a vector bundle structure
$$
S_\chi \simeq T_{Z_1}X \underset{X}{\times} T_{Z_2}X \underset{X}{\times} T_{Z_3}X.
$$
The rational monomials $\imin{\varphi_j}$, $j=1,2,3$, are $\imin{\varphi_1}(\tau)=\tau_1$, $\imin{\varphi_2}(\tau)=\tau_2$, $\imin{\varphi_3}(\tau)=\tau_3$. Let $\xi=(\xi_1,\xi_2,\xi_3)$, $\xi_i \neq 0$, $i=1,2,3$. A cofinal family of $C(\xi)$ is $C(\xi,F^0)$ consisting of elements $S(W,\epsilon)$, $W=W_1 \times W_2 \times W_3$, defined as follows
$$
S(W,\epsilon)=
\begin{aligned}
 \left\{z \in X;
\begin{aligned}
& z_k \in W_k \quad (k=1,2,3), \\
& |z_i| < \epsilon \quad (i=1,2,3)
\end{aligned}
\right\},
\end{aligned}
$$
where $\epsilon > 0$ and each $W_k$ is an $\mathbb{R}^+$-conic open subset in $\CC$ containing the point $\xi_k$. \\

We now construct the asymptotics.
%We have $K_{\{j\}}=K_j$, $j=1,2$, $K_{\{1,2\}}=K_1 \cup K_2$ and $Z_J=\{z^{(k)}=0,\,k \in K_J\}$ for $\emptyset \neq J \subseteq \{1,2\}$. For any subset $K \subset \{1,2\}$, we denote by $z^{(K)}$ the set of the coordinates $z^{(k)}$'s with $k \in K$ and by $z^{(K)}_{\mathrm{C}}$ the set of the coordinates which do not belong to $z^{(K)}$, i.e., $z^{(k)}$'s with $k \in \{0,1,\dots,m\} \setminus K$. Hence, the coordinates of $Z_J$ are given by $z^{(K_J)}_{\mathrm{C}}$. The $\pi_J$ denotes the canonical projection from $X$ to $Z_J$ defined by $z \to z^{(K_J)}_{\mathrm{C}}$. Given $S(W,\epsilon)$, set $S_J=\pi_J(S(W,\epsilon))$.
We have $S_{\{1,2,3\}}=\{0\}$ and, for $i \neq j \neq k \in \{1,2,3\}$
$$
S_{\{i,j\}}=
\begin{aligned}
 \left\{z \in Z_i \cap Z_j;
\begin{aligned}
&z_k \in W_k, \\
& |z_k| < \epsilon
\end{aligned}
\right\}.
\end{aligned}
$$
$$
S_{\{i\}}=
\begin{aligned}
 \left\{z \in Z_i;
\begin{aligned}
&z_j \in W_j, \\
&z_k \in W_k, \\
& |z_j|,|z_k| < \epsilon
\end{aligned}
\right\}.
\end{aligned}
$$

A total family of coefficients of multi-asymptotic expansion is given by
\begin{eqnarray*}
F & = & \{F_{\{1\}},F_{\{2\}},F_{\{3\}},F_{\{1,2\}},F_{\{1,3\}},F_{\{2,3\}},F_{\{1,2,3\}}\} \\
& = & \Big\{\{f_{\{1\},k}(z_2,z_3)\}_{k \in \Z_{\geq 0}},
\{f_{\{2\},k}(z_1,z_3)\}_{k \in
\Z_{\geq 0}},\{f_{\{3\},k}(z_1,z_2)\}_{k \in
\Z_{\geq 0}},  \\
&& \{f_{\{1,2\},\alpha}(z_3)\}_{\alpha \in \Z_{\geq 0}^2},\{f_{\{1,3\},\alpha}(z_2)\}_{\alpha \in \Z_{\geq 0}^2},\{f_{\{2,3\},\alpha}(z_1)\}_{\alpha \in \Z_{\geq 0}^2}, \\
&& \{f_{\{1,2,3\},\beta}\}_{\beta \in \Z_{\geq 0}^3}\Big\},
\end{eqnarray*}
where $f_{\{1\},k}(z_2,z_3)$ (resp. $f_{\{2\},k}(z_1,z_3)$, resp. $f_{\{3\},k}(z_1,z_2)$) is holomorphic in $S_{\{1\}}$ (resp. $S_{\{2\}}$, resp. $S_{\{3\}}$), $f_{\{1,2\},\alpha}(z_3)$ (resp. $f_{\{1,3\},\alpha}(z_2)$, resp. $f_{\{2,3\},\alpha}(z_1)$) is holomorphic in $S_{\{1,2\}}$ (resp. $S_{\{1,3\}}$, resp. $S_{\{2,3\}}$) and $f_{\{1,2,3\},\beta} \in \CC$. An asymptotic expansion $\operatorname{App}^{<N}(F;\,z)$, $N=(n_1,n_2,n_3) \in \Z_{\geq 0}^3$ is given by
$$
\begin{aligned}
	T_{\{1\}}^{<N}(F;\, z) &= \sum_{k < n_1} f_{\{1\},k}(z_2,z_3)\frac{z_1^k}{k!}, \qquad k \in \Z_{\geq 0}, \\
	T_{\{2\}}^{<N}(F;\, z) &= \sum_{k < n_2} f_{\{2\},k}(z_1,z_3)\frac{z_2^k}{k!}, \qquad k \in \Z_{\geq 0}, \\	
	T_{\{3\}}^{<N}(F;\, z) &= \sum_{k < n_3} f_{\{3\},k}(z_1,z_2)\frac{z_3^k}{k!}, \qquad k \in \Z_{\geq 0}, \\	
T_{\{1,2\}}^{<N}(F;\, z) &= \sum_{\substack{\alpha_1 < n_1 \\ \alpha_2 < n_2}} f_{\{1,2\},\alpha}(z_3)\frac{z_1^{\alpha_1}z_2^{\alpha_2}}{\alpha_1!\alpha_2!}, \qquad \alpha_1,\alpha_2 \in \Z_{\geq 0}, \\
T_{\{1,3\}}^{<N}(F;\, z) &= \sum_{\substack{\alpha_1 < n_1 \\ \alpha_3 < n_3}} f_{\{1,3\},\alpha}(z_2)\frac{z_1^{\alpha_1}z_3^{\alpha_3}}{\alpha_1!\alpha_3!}, \qquad \alpha_1,\alpha_3 \in \Z_{\geq 0}, \\
T_{\{2,3\}}^{<N}(F;\, z) &= \sum_{\substack{\alpha_2 < n_2 \\ \alpha_3 < n_3}} f_{\{2,3\},\alpha}(z_1)\frac{z_2^{\alpha_2}z_3^{\alpha_3}}{\alpha_2!\alpha_3!}, \qquad \alpha_2,\alpha_3 \in \Z_{\geq 0}, \\
T_{\{1,2,3\}}^{<N}(F;\, z) &= \sum_{\substack{\beta_1 < n_1 \\ \beta_2 < n_2 \\ \beta_3 < n_3}} f_{\{1,2,3\},\beta}\frac{z_1^{\beta_1}z_2^{\beta_2}z_3^{\beta_3}}{\beta_1!\beta_2!\beta_3!}, \qquad \beta_1,\beta_2,\beta_3 \in \Z_{\geq 0}, \\
	\operatorname{App}^{<N}(F;\, z) &=
T_{\{1\}}^{<N}(F;\, z) + T_{\{2\}}^{<N}(F;\, z) + T_{\{3\}}^{<N}(F;\, z) \\
& - T_{\{1,2\}}^{<N}(F;\, z) - T_{\{1,3\}}^{<N}(F;\, z)  - T_{\{2,3\}}^{<N}(F;\, z) \\
& + T_{\{1,2,3\}}^{<N}(F;\, z).
\end{aligned}
$$
We say that $f$ is multi-asymptotically developable to $F$ along
$\chi$ on $S=S(W,\epsilon)$ if and only if for any cone
$S'=S(W',\epsilon')$ properly contained in $S$ and for any $N=(n_1,n_2,n_3)
\in \mathbb{Z}^3_{\ge 0}$, there exists a constant $C_{S',N}$ such that
$$
\left|
f(z) - \operatorname{App}^{<N}(F;\, z)
\right|
\le C_{S',N}|z_1|^{n_1}|z_2|^{n_2}|z_3|^{n_3}
 \qquad (z \in S').
$$
The family $F$ is consistent if
\begin{itemize}
\item $f_{\{1\},k}(z_2,z_3)$ is strongly asymptotically developable to
    $$
    \begin{aligned}
    &\left\{\{f_{\{1,2\},(k,\alpha_2)}(z_3)\}_{\alpha_2 \in \Z_{\geq 0}},\{f_{\{1,3\},(k,\alpha_3)}(z_2)\}_{\alpha_3 \in \Z_{\geq 0}}, \right.\\
    &\qquad\qquad\qquad\qquad\left.
    \{f_{\{1,2,3\},(k,\beta_2,\beta_3)}\}_{(\beta_2,\beta_3) \in \Z^2_{\geq 0}}\right\}
    \end{aligned}
    $$
    on $S_{\{1\}}$ for each $k \in \Z_{\geq 0}$,
\item $f_{\{2\},k}(z_1,z_3)$ is strongly asymptotically developable to
    $$
    \begin{aligned}
    &\left\{\{f_{\{1,2\},(\alpha_1,k)}(z_3)\}_{\alpha_1 \in \Z_{\geq 0}},\{f_{\{2,3\},(k,\alpha_3)}(z_1)\}_{\alpha_3 \in \Z_{\geq 0}},\right.\\
    &\qquad\qquad\qquad\qquad\left.\{f_{\{1,2,3\},(\beta_1,k,\beta_3)}\}_{(\beta_1,\beta_3) \in \Z^2_{\geq 0}}\right\}
    \end{aligned}
    $$
    on $S_{\{2\}}$ for each $k \in \Z_{\geq 0}$,
\item $f_{\{3\},k}(z_1,z_2)$ is strongly asymptotically developable to
    $$
    \begin{aligned}
    &\left\{\{f_{\{1,3\},(\alpha_1,k)}(z_2)\}_{\alpha_1 \in \Z_{\geq 0}},\{f_{\{2,3\},(\alpha_2,k)}(z_1)\}_{\alpha_2 \in \Z_{\geq 0}},\right.\\
    &\left.\qquad\qquad\qquad\qquad\{f_{\{1,2,3\},(\beta_1,\beta_2,k)}\}_{(\beta_1,\beta_2) \in \Z^2_{\geq 0}}\right\}
    \end{aligned}
    $$
    on $S_{\{3\}}$ for each $k \in \Z_{\geq 0}$,
\item $f_{\{1,2\},(\alpha_1,\alpha_2)}(z_3)$ is strongly asymptotically developable to
    $$
    \{f_{\{1,2,3\},(\alpha_1,\alpha_2,\beta_3)}\}_{\beta_3 \in \Z_{\geq 0}}
    $$
    on $S_{\{1,2\}}$ for each $(\alpha_1,\alpha_2) \in \Z^2_{\geq 0}$,
\item $f_{\{1,3\},(\alpha_1,\alpha_3)}(z_2)$ is strongly asymptotically developable to
    $$
    \{f_{\{1,2,3\},(\alpha_1,\beta_2,\alpha_3)}\}_{\beta_2 \in \Z_{\geq 0}}
    $$
    on $S_{\{1,3\}}$ for each $(\alpha_1,\alpha_3) \in \Z^2_{\geq 0}$,
\item $f_{\{2,3\},(\alpha_2,\alpha_3)}(z_1)$ is strongly asymptotically developable to
    $$
    \{f_{\{1,2,3\},(\beta_1,\alpha_2,\alpha_3)}\}_{\beta_1 \in \Z_{\geq 0}}
    $$
    on $S_{\{2,3\}}$ for each $(\alpha_2,\alpha_3) \in \Z^2_{\geq 0}$.
\end{itemize}
These are the asymptotics of \cite{Ma84} in $\CC^3$.
%%%%%%%%%%%%%%%%%%%%%%%%%%%%%%%%%%%%%%%%%%%%%
% TAKEUCHI C3
%%%%%%%%%%%%%%%%%%%%%%%%%%%%%%%%%%%%%%%%%%%%%

\

\noindent
$\checkmark\,${\bf{Takeuchi asymptotics.}} In this case
$$
A_\chi=
\left(
\begin{array}{ccc}
1 & 1 & 1 \\
0 & 1 & 1 \\
0 & 0 & 1
\end{array}
\right).
$$
We have $Z_1=\{0\}$, $Z_2=\{z_2=z_3=0\}$, $Z_3=\{z_3=0\}$.  Then we can define a  normal deformation $\widetilde{X}=\CC^3 \times \R^3$ with the map
$p(z;\,t)=(t_1z_1, t_1t_2z_2, t_1t_2t_3z_3)$. By Proposition \ref{prop:bundle-description} the zero section $S_\chi$ of $\widetilde{X}$ has a vector bundle structure
$$
S_\chi \simeq T_{Z_1}Z_2 \underset{X}{\times} T_{Z_2}Z_3 \underset{X}{\times} T_{Z_3}X.
$$
The rational monomials $\imin{\varphi_j}$, $j=1,2,3$, are $\imin{\varphi_1}(\tau)=\tau_1$, $\imin{\varphi_2}(\tau)=\displaystyle{\tau_2 \over \tau_1}$, $\imin{\varphi_3}(\tau)=\displaystyle{\tau_3 \over \tau_2}$. Let $\xi=(\xi_1,\xi_2,\xi_3)$, $\xi_i \neq 0$, $i=1,2,3$. A cofinal family of $C(\xi)$ is $C(\xi,F^0)$ consisting of elements $S(W,\epsilon)$, $W=W_1 \times W_2 \times W_3$, defined as follows
$$
S(W,\epsilon)=
\begin{aligned}
 \left\{z \in X;
\begin{aligned}
& z_k \in W_k \quad (k=1,2,3), \\
& |z_1| < \epsilon, \\
& |z_2| < \epsilon |z_1|, \\
& |z_3| < \epsilon |z_2| \\
\end{aligned}
\right\},
\end{aligned}
$$
where $\epsilon > 0$ and each $W_k$ is an $\mathbb{R}^+$-conic open subset in $\CC$ containing the point $\xi_k$. \\

We now construct the asymptotics.
%We have $K_{\{j\}}=K_j$, $j=1,2$, $K_{\{1,2\}}=K_1 \cup K_2$ and $Z_J=\{z^{(k)}=0,\,k \in K_J\}$ for $\emptyset \neq J \subseteq \{1,2\}$. For any subset $K \subset \{1,2\}$, we denote by $z^{(K)}$ the set of the coordinates $z^{(k)}$'s with $k \in K$ and by $z^{(K)}_{\mathrm{C}}$ the set of the coordinates which do not belong to $z^{(K)}$, i.e., $z^{(k)}$'s with $k \in \{0,1,\dots,m\} \setminus K$. Hence, the coordinates of $Z_J$ are given by $z^{(K_J)}_{\mathrm{C}}$. The $\pi_J$ denotes the canonical projection from $X$ to $Z_J$ defined by $z \to z^{(K_J)}_{\mathrm{C}}$. Given $S(W,\epsilon)$, set $S_J=\pi_J(S(W,\epsilon))$.
We have $S_{\{1\}}=S_{\{1,2\}}=S_{\{1,3\}}=S_{\{1,2,3\}}=\{0\}$ and
$$
S_{\{2\}}=S_{\{2,3\}}=
\begin{aligned}
 \left\{z \in Z_2;
\begin{aligned}
&z_1 \in W_1, \\
& |z_1| < \epsilon
\end{aligned}
\right\},
\end{aligned}
$$
$$
S_{\{3\}}=
\begin{aligned}
 \left\{z \in Z_3;
\begin{aligned}
&z_j \in W_j \quad (j=1,2), \\
&|z_1| < \epsilon, \\
& |z_2| < \epsilon |z_2|
\end{aligned}
\right\}.
\end{aligned}
$$

A total family of coefficients of multi-asymptotic expansion is given by
\begin{eqnarray*}
F & = & \{F_{\{1\}},F_{\{2\}},F_{\{3\}},F_{\{1,2\}},F_{\{1,3\}},F_{\{2,3\}},F_{\{1,2,3\}}\} \\
& = & \left\{\{f_{\{1\},\beta}\}_{\beta \in \Z_{\geq 0}^3},
\{f_{\{2\},\alpha}(z_1)\}_{\alpha \in
\Z_{\geq 0}^2},\{f_{\{3\},k}(z_1,z_2)\}_{k \in
\Z_{\geq 0}}, \right. \\
&& \{f_{\{1,2\},\beta}\}_{\beta \in \Z_{\geq 0}^3},\{f_{\{1,3\},\beta}\}_{\beta \in \Z_{\geq 0}^3},\{f_{\{2,3\},\alpha}(z_1)\}_{\alpha \in \Z_{\geq 0}^2}, \\
&& \left.\{f_{\{1,2,3\},\beta}\}_{\beta \in \Z_{\geq 0}^3}\right\},
\end{eqnarray*}
where $f_{\{3\},k}(z_1,z_2)$ (resp. $f_{\{2\},\alpha}(z_1)$, resp. $f_{\{2,3\},\alpha}(z_1)$) is holomorphic in $S_{\{3\}}$ (resp. $S_{\{2\}}$, resp $S_{\{2,3\}}$) and $f_{\{1\},\beta},f_{\{1,2\},\beta},f_{\{1,3\},\beta},f_{\{1,2,3\},\beta} \in \CC$. An asymptotic expansion $\operatorname{App}^{<N}(F;\,z)$, $N=(n_1,n_2,n_3) \in \Z_{\geq 0}^3$ is given by
$$
\begin{aligned}
	T_{\{1\}}^{<N}(F;\, z) &= \sum_{\beta_1 + \beta_2 + \beta_3 < n_1} f_{\{1\},\beta}\frac{z_1^{\beta_1}z_2^{\beta_2}z_3^{\beta_3}}{\beta_1!\beta_2!\beta_3!}, \qquad \beta_1,\beta_2,\beta_3 \in \Z_{\geq 0}, \\
	T_{\{2\}}^{<N}(F;\, z) &= \sum_{\alpha_2+\alpha_3 < n_2} f_{\{2\},\alpha}(z_1)\frac{z_2^{\alpha_2}z_3^{\alpha_3}}{\alpha_2!\alpha_3!}, \qquad \alpha_2,\alpha_3 \in \Z_{\geq 0}, \\	
	T_{\{3\}}^{<N}(F;\, z) &= \sum_{k < n_3} f_{\{3\},k}(z_1,z_2)\frac{z_3^k}{k!}, \qquad k \in \Z_{\geq 0}, \\	
T_{\{1,2\}}^{<N}(F;\, z) &= \sum_{\substack{\beta_1+\beta_2+\beta_3 < n_1 \\ \beta_2+\beta_3 < n_2}} f_{\{1,2\},\beta}\frac{z_1^{\beta_1}z_2^{\beta_2}z_3^{\beta_3}}{\beta_1!\beta_2!\beta_3!}, \qquad \beta_1,\beta_2,\beta_3 \in \Z_{\geq 0}, \\
T_{\{1,3\}}^{<N}(F;\, z) &= \sum_{\substack{\beta_1 + \beta_2 + \beta_3 < n_1 \\ \beta_3 < n_3}} f_{\{1,3\},\beta}\frac{z_1^{\beta_1}z_2^{\beta_2}z_3^{\beta_3}}{\beta_1!\beta_2!\beta_3!}, \qquad \beta_1,\beta_2,\beta_3 \in \Z_{\geq 0}, \\
T_{\{2,3\}}^{<N}(F;\, z) &= \sum_{\substack{\alpha_2 + \alpha_3 < n_2 \\ \alpha_3 < n_3}} f_{\{2,3\},\alpha}(z_1)\frac{z_2^{\alpha_2}z_3^{\alpha_3}}{\alpha_2!\alpha_3!}, \qquad \alpha_2,\alpha_3 \in \Z_{\geq 0}, \\
T_{\{1,2,3\}}^{<N}(F;\, z) &= \sum_{\substack{\beta_1 + \beta_2 + \beta_3 < n_1 \\ \beta_2 + \beta_3 < n_2 \\ \beta_3 < n_3}} f_{\{1,2,3\},\beta}\frac{z_1^{\beta_1}z_2^{\beta_2}z_3^{\beta_3}}{\beta_1!\beta_2!\beta_3!}, \qquad \beta_1,\beta_2,\beta_3 \in \Z_{\geq 0}, \\
	\operatorname{App}^{<N}(F;\, z) &=
T_{\{1\}}^{<N}(F;\, z) + T_{\{2\}}^{<N}(F;\, z) + T_{\{3\}}^{<N}(F;\, z) \\
& - T_{\{1,2\}}^{<N}(F;\, z) - T_{\{1,3\}}^{<N}(F;\, z)  - T_{\{2,3\}}^{<N}(F;\, z) \\
& + T_{\{1,2,3\}}^{<N}(F;\, z).
\end{aligned}
$$
We say that $f$ is multi-asymptotically developable to $F$ along
$\chi$ on $S=S(W,\epsilon)$ if and only if for any cone
$S'=S(W',\epsilon')$ properly contained in $S$ and for any $N=(n_1,n_2,n_3)
\in \mathbb{Z}^3_{\ge 0}$, there exists a constant $C_{S',N}$ such that
$$
\left|
f(z) - \operatorname{App}^{<N}(F;\, z)
\right|
\le C_{S',N}|z_1|^{n_1-n_2}|z_2|^{n_2-n_3}|z_3|^{n_3}
 \qquad (z \in S').
$$
The family $F$ is consistent if
\begin{itemize}
\item $f_{\{1\},\beta}=f_{\{1,2\},\beta}=f_{\{1,3\},\beta}=f_{\{1,2,3\},\beta}$ for each $\beta \in \Z^3_{\geq 0}$,
\item $f_{\{2\},\alpha}=f_{\{2,3\},\alpha}$ for each $\alpha \in \Z^2_{\geq 0}$,
\item $f_{\{3\},k}(z_1,z_2)$ is strongly asymptotically developable to
    $$
    \begin{aligned}
    &\left\{\{f_{\{1,3\},(\beta_1,\beta_2,k)}\}_{(\beta_1,\beta_2) \in \Z^2_{\geq 0}},\{f_{\{2,3\},(\alpha_2,k)}(z_1)\}_{\alpha_2 \in \Z_{\geq 0}},\right.\\
    &\qquad\qquad\qquad\qquad\left.\{f_{\{1,2,3\},(\beta_1,\beta_2,k)}\}_{(\beta_1,\beta_2) \in \Z^2_{\geq 0}}\right\}
    \end{aligned}
    $$
    on $S_{\{3\}}$ for each $k \in \Z_{\geq 0}$,
\item $f_{\{2,3\},(\alpha_2,\alpha_3)}(z_1)$ is strongly asymptotically developable to
    $$
    \{f_{\{1,2,3\},(\beta_1,\alpha_2,\alpha_3)}\}_{\beta_1 \in \Z_{\geq 0}}
    $$
    on $S_{\{2,3\}}$ for each $(\alpha_2,\alpha_3) \in \Z^2_{\geq 0}$.
\end{itemize}
These are the asymptotics of \cite{HP13} generalizing the bispecialization of \cite{ST94} in $\CC^3$.

%%%%%%%%%%%%%%%%%%%%%%%%%%%%%%%%%%%%%%
% CLEAN INTERSECTION C3
%%%%%%%%%%%%%%%%%%%%%%%%%%%%%%%%%%%%%%

\

\noindent
$\checkmark\,${\bf{Clean intersection (3 lines) asymptotics.}} In this case
$$
A_\chi=
\left(
\begin{array}{ccc}
1 & 1 & 0 \\
0 & 1 & 1 \\
1 & 0 & 1
\end{array}
\right).
$$
We have $Z_1=\{z_1=z_2=0\}$, $Z_2=\{z_2=z_3=0\}$, $Z_3=\{z_1=z_3=0\}$.
Then we can define a  normal deformation $\widetilde{X}=\CC^3 \times \R^3$ with the map $p(z;\,t)=(t_1t_3z_1, t_1t_2z_2, t_2t_3z_3)$. Set $Z=Z_1 \cap Z_2 \cap Z_3$. By Proposition \ref{prop:bundle-description} the zero section $S_\chi$ of $\widetilde{X}$ has a vector bundle structure
$$
S_\chi \simeq T_{Z}Z_1 \underset{X}{\times} T_{Z}Z_2 \underset{X}{\times} T_{Z}Z_3.
$$
The rational monomials $\imin{\varphi_j}$, $j=1,2,3$, are $\imin{\varphi_1}(\tau)=\displaystyle{\left({\tau_1\tau_2 \over \tau_3}\right)^{1 \over 2}}$, $\imin{\varphi_2}(\tau)=\displaystyle{\left({\tau_2\tau_3 \over \tau_1}\right)^{1 \over 2}}$, $\imin{\varphi_3}(\tau)=\displaystyle{\left({\tau_1\tau_3 \over \tau_2}\right)^{1 \over 2}}$. Let $\xi=(\xi_1,\xi_2,\xi_3)$, $\xi_i \neq 0$, $i=1,2,3$. A cofinal family of $C(\xi)$ is $C(\xi,F^0)$ consisting of elements $S(W,\epsilon)$, $W=W_1 \times W_2 \times W_3$, defined as follows
$$
S(W,\epsilon)=
\begin{aligned}
 \left\{z \in X;
\begin{aligned}
& z_k \in W_k \quad (k=1,2,3), \\
& |z_1||z_2| < \epsilon |z_3|, \\
& |z_2||z_3| < \epsilon |z_1|, \\
& |z_1||z_3| < \epsilon |z_2|\\
\end{aligned}
\right\},
\end{aligned}
$$
where $\epsilon > 0$ and each $W_k$ is an $\mathbb{R}^+$-conic open subset in $\CC$ containing the point $\xi_k$. \\

We now construct the asymptotics.
%We have $K_{\{j\}}=K_j$, $j=1,2$, $K_{\{1,2\}}=K_1 \cup K_2$ and $Z_J=\{z^{(k)}=0,\,k \in K_J\}$ for $\emptyset \neq J \subseteq \{1,2\}$. For any subset $K \subset \{1,2\}$, we denote by $z^{(K)}$ the set of the coordinates $z^{(k)}$'s with $k \in K$ and by $z^{(K)}_{\mathrm{C}}$ the set of the coordinates which do not belong to $z^{(K)}$, i.e., $z^{(k)}$'s with $k \in \{0,1,\dots,m\} \setminus K$. Hence, the coordinates of $Z_J$ are given by $z^{(K_J)}_{\mathrm{C}}$. The $\pi_J$ denotes the canonical projection from $X$ to $Z_J$ defined by $z \to z^{(K_J)}_{\mathrm{C}}$. Given $S(W,\epsilon)$, set $S_J=\pi_J(S(W,\epsilon))$.
We have $S_{\{1,2\}}=S_{\{1,3\}}=S_{\{2,3\}}=S_{\{1,2,3\}}=\{0\}$ and
$$
S_{\{1\}}=
\begin{aligned}
 \left\{z \in Z_1;
\begin{aligned}
&z_3 \in W_3, \\
& |z_3| < \epsilon
\end{aligned}
\right\},
\end{aligned}
$$
$$
S_{\{2\}}=
\begin{aligned}
 \left\{z \in Z_2;
\begin{aligned}
&z_1 \in W_1, \\
& |z_1| < \epsilon
\end{aligned}
\right\},
\end{aligned}
$$
$$
S_{\{3\}}=
\begin{aligned}
 \left\{z \in Z_3;
\begin{aligned}
&z_2 \in W_2, \\
& |z_2| < \epsilon
\end{aligned}
\right\}.
\end{aligned}
$$

A total family of coefficients of multi-asymptotic expansion is given by
\begin{eqnarray*}
F & = & \{F_{\{1\}},F_{\{2\}},F_{\{3\}},F_{\{1,2\}},F_{\{1,3\}},F_{\{2,3\}},F_{\{1,2,3\}}\} \\
& = & \left\{\{f_{\{1\},\alpha}(z_3)\}_{\alpha \in \Z_{\geq 0}^2},
\{f_{\{2\},\alpha}(z_1)\}_{\alpha \in
\Z_{\geq 0}^2},\{f_{\{3\},\alpha}(z_2)\}_{\alpha \in
\Z_{\geq 0}^2}, \right. \\
&& \{f_{\{1,2\},\beta}\}_{\beta \in \Z_{\geq 0}^3},\{f_{\{1,3\},\beta}\}_{\beta \in \Z_{\geq 0}^3},\{f_{\{2,3\},\beta}\}_{\beta \in \Z_{\geq 0}^3}, \\
&& \left.\{f_{\{1,2,3\},\beta}\}_{\beta \in \Z_{\geq 0}^3}\right\},
\end{eqnarray*}
where $f_{\{1\},\alpha}(z_3)$ (resp. $f_{\{2\},\alpha}(z_1)$, resp. $f_{\{3\},\alpha}(z_2)$) is holomorphic in $S_{\{1\}}$ (resp. $S_{\{2\}}$, resp $S_{\{3\}}$) and $f_{\{1,2\},\beta},f_{\{1,3\},\beta},f_{\{2,3\},\beta},f_{\{1,2,3\},\beta} \in \CC$. An asymptotic expansion $\operatorname{App}^{<N}(F;\,z)$, $N=(n_1,n_2,n_3) \in \Z_{\geq 0}^3$ is given by
$$
\begin{aligned}
	T_{\{1\}}^{<N}(F;\, z) &= \sum_{\alpha_1 + \alpha_2 < n_1} f_{\{1\},\alpha}(z_3)\frac{z_1^{\alpha_1}z_2^{\alpha_2}}{\alpha_1!\alpha_2!}, \qquad \alpha_1,\alpha_2 \in \Z_{\geq 0}, \\
	T_{\{2\}}^{<N}(F;\, z) &= \sum_{\alpha_2+\alpha_3 < n_2} f_{\{2\},\alpha}(z_1)\frac{z_2^{\alpha_2}z_3^{\alpha_3}}{\alpha_2!\alpha_3!}, \qquad \alpha_2,\alpha_3 \in \Z_{\geq 0}, \\	
	T_{\{3\}}^{<N}(F;\, z) &= \sum_{\alpha_1 + \alpha_3 < n_3} f_{\{3\},\alpha}(z_2)\frac{z_1^{\alpha_1}z_3^{\alpha_3}}{\alpha_1!\alpha_3!}, \qquad \alpha_1,\alpha_3 \in \Z_{\geq 0}, \\	
T_{\{1,2\}}^{<N}(F;\, z) &= \sum_{\substack{\beta_1+\beta_2 < n_1 \\ \beta_2+\beta_3 < n_2}} f_{\{1,2\},\beta}\frac{z_1^{\beta_1}z_2^{\beta_2}z_3^{\beta_3}}{\beta_1!\beta_2!\beta_3!}, \qquad \beta_1,\beta_2,\beta_3 \in \Z_{\geq 0}, \\
T_{\{1,3\}}^{<N}(F;\, z) &= \sum_{\substack{\beta_1 + \beta_2 < n_1 \\ \beta_1 + \beta_3 < n_3}} f_{\{1,3\},\beta}\frac{z_1^{\beta_1}z_2^{\beta_2}z_3^{\beta_3}}{\beta_1!\beta_2!\beta_3!}, \qquad \beta_1,\beta_2,\beta_3 \in \Z_{\geq 0}, \\
T_{\{2,3\}}^{<N}(F;\, z) &= \sum_{\substack{\beta_2 + \beta_3 < n_2 \\ \beta_1 + \beta_3 < n_3}} f_{\{2,3\},\beta}\frac{z_1^{\beta_1}z_2^{\beta_2}z_3^{\beta_3}}{\beta_1!\beta_2!\beta_3!}, \qquad \beta_1,\beta_2,\beta_3 \in \Z_{\geq 0}, \\
T_{\{1,2,3\}}^{<N}(F;\, z) &= \sum_{\substack{\beta_1 + \beta_2 < n_1 \\ \beta_2 + \beta_3 < n_2 \\ \beta_1 + \beta_3 < n_3}} f_{\{1,2,3\},\beta}\frac{z_1^{\beta_1}z_2^{\beta_2}z_3^{\beta_3}}{\beta_1!\beta_2!\beta_3!}, \qquad \beta_1,\beta_2,\beta_3 \in \Z_{\geq 0}, \\
	\operatorname{App}^{<N}(F;\, z) &=
T_{\{1\}}^{<N}(F;\, z) + T_{\{2\}}^{<N}(F;\, z) + T_{\{3\}}^{<N}(F;\, z) \\
& - T_{\{1,2\}}^{<N}(F;\, z) - T_{\{1,3\}}^{<N}(F;\, z)  - T_{\{2,3\}}^{<N}(F;\, z) \\
& + T_{\{1,2,3\}}^{<N}(F;\, z).
\end{aligned}
$$
We say that $f$ is multi-asymptotically developable to $F$ along
$\chi$ on $S=S(W,\epsilon)$ if and only if for any cone
$S'=S(W',\epsilon')$ properly contained in $S$ and for any $N=(n_1,n_2,n_3)
\in \mathbb{Z}^3_{\ge 0}$, there exists a constant $C_{S',N}$ such that
$$
\left|
f(z) - \operatorname{App}^{<N}(F;\, z)
\right|
\le C_{S',N}\left(|z_1|^{n_1+n_3-n_2}|z_2|^{n_1+n_2-n_3}|z_3|^{n_2+n_3-n_1}\right)^{1 \over 2}
$$
on $z \in S'$.
The family $F$ is consistent if
\begin{itemize}
\item $f_{\{1\},\beta}=f_{\{2\},\beta}=f_{\{3\},\beta}=f_{\{1,2,3\},\beta}$ for each $\beta \in \Z^3_{\geq 0}$,
\item $f_{\{1\},(\alpha_1,\alpha_2)}(z_3)$ is strongly asymptotically developable to
    $$
    \{f_{\{1,2,3\},(\alpha_1,\alpha_2,\beta_3)}\}_{\beta_3 \in \Z_{\geq 0}}
    $$
    on $S_{\{1\}}$ for each $(\alpha_1,\alpha_2) \in \Z^2_{\geq 0}$,
\item $f_{\{1,2\},(\alpha_2,\alpha_3)}(z_1)$ is strongly asymptotically developable to
    $$
    \{f_{\{1,2,3\},(\beta_1,\alpha_2,\alpha_3)}\}_{\beta_1 \in \Z_{\geq 0}}
    $$
    on $S_{\{2\}}$ for each $(\alpha_2,\alpha_3) \in \Z^2_{\geq 0}$,
\item $f_{\{1,3\},(\alpha_1,\alpha_3)}(z_2)$ is strongly asymptotically developable to
    $$
    \{f_{\{1,2,3\},(\alpha_1,\beta_2,\alpha_3)}\}_{\beta_2 \in \Z_{\geq 0}}
    $$
    on $S_{\{3\}}$ for each $(\alpha_1,\alpha_3) \in \Z^2_{\geq 0}$.
\end{itemize}
%%%%%%%%%%%%%%%%%%%%%%%%%%%%%%%%%%%%%%%%%%%%%%%%
% MIXED C3
%%%%%%%%%%%%%%%%%%%%%%%%%%%%%%%%%%%%%%%%%%%%%%%%

\

\noindent
$\checkmark\,${\bf{Mixed asymptotics (Majima-Takeuchi).}} In this case
$$
A_\chi=
\left(
\begin{array}{ccc}
1 & 1 & 1 \\
0 & 1 & 0 \\
0 & 0 & 1
\end{array}
\right).
$$
We have $Z_1=\{0\}$, $Z_2=\{z_2=0\}$, $Z_3=\{z_3=0\}$,  Then we can define a  normal deformation $\widetilde{X}=\CC^3 \times \R^3$ with the map $p(z;\,t)=(t_1z_1, t_1t_2z_2, t_1t_3z_3)$. By Proposition \ref{prop:bundle-description} the zero section $S_\chi$ of $\widetilde{X}$ has a vector bundle structure
$$
S_\chi \simeq T_{Z_1}(Z_2 \cap Z_3) \underset{X}{\times} T_{Z_2 \cap Z_3}Z_2 \underset{X}{\times} T_{Z_2 \cap Z_3}Z_3.
$$
The rational monomials $\imin{\varphi_j}$, $j=1,2,3$, are $\imin{\varphi_1}(\tau)=\tau_1$, $\imin{\varphi_2}(\tau)=\displaystyle{\tau_2 \over \tau_1}$, $\imin{\varphi_3}(\tau)=\displaystyle{\tau_3 \over \tau_1}$. Let $\xi=(\xi_1,\xi_2,\xi_3)$, $\xi_i \neq 0$, $i=1,2,3$. A cofinal family of $C(\xi)$ is $C(\xi,F^0)$ consisting of elements $S(W,\epsilon)$, $W=W_1 \times W_2 \times W_3$, defined as follows
$$
S(W,\epsilon)=
\begin{aligned}
 \left\{z \in X;
\begin{aligned}
& z_k \in W_k \quad (k=1,2,3), \\
& |z_1| < \epsilon, \\
& |z_2| < \epsilon |z_1|, \\
& |z_3| < \epsilon |z_1|\\
\end{aligned}
\right\},
\end{aligned}
$$
where $\epsilon > 0$ and each $W_k$ is an $\mathbb{R}^+$-conic open subset in $\CC$ containing the point $\xi_k$. \\

We now construct the asymptotics.
%We have $K_{\{j\}}=K_j$, $j=1,2$, $K_{\{1,2\}}=K_1 \cup K_2$ and $Z_J=\{z^{(k)}=0,\,k \in K_J\}$ for $\emptyset \neq J \subseteq \{1,2\}$. For any subset $K \subset \{1,2\}$, we denote by $z^{(K)}$ the set of the coordinates $z^{(k)}$'s with $k \in K$ and by $z^{(K)}_{\mathrm{C}}$ the set of the coordinates which do not belong to $z^{(K)}$, i.e., $z^{(k)}$'s with $k \in \{0,1,\dots,m\} \setminus K$. Hence, the coordinates of $Z_J$ are given by $z^{(K_J)}_{\mathrm{C}}$. The $\pi_J$ denotes the canonical projection from $X$ to $Z_J$ defined by $z \to z^{(K_J)}_{\mathrm{C}}$. Given $S(W,\epsilon)$, set $S_J=\pi_J(S(W,\epsilon))$.
We have $S_{\{1\}}=S_{\{1,2\}}=S_{\{1,3\}}=S_{\{1,2,3\}}=\{0\}$ and
$$
S_{\{2\}}=
\begin{aligned}
 \left\{z \in Z_2;
\begin{aligned}
&z_j \in W_j \quad (j=1,3), \\
& |z_1| < \epsilon, \\
& |z_1| < \epsilon |z_3|
\end{aligned}
\right\},
\end{aligned}
$$
$$
S_{\{3\}}=
\begin{aligned}
 \left\{z \in Z_3;
\begin{aligned}
&z_j \in W_j \quad (j=1,2), \\
& |z_1| < \epsilon, \\
& |z_1| < \epsilon |z_2|
\end{aligned}
\right\},
\end{aligned}
$$
$$
S_{\{2,3\}}=
\begin{aligned}
 \left\{z \in Z_2 \cap Z_3;
\begin{aligned}
&z_1 \in W_1, \\
&|z_1| < \epsilon
\end{aligned}
\right\}.
\end{aligned}
$$

A total family of coefficients of multi-asymptotic expansion is given by
\begin{eqnarray*}
F & = & \{F_{\{1\}},F_{\{2\}},F_{\{3\}},F_{\{1,2\}},F_{\{1,3\}},F_{\{2,3\}},F_{\{1,2,3\}}\} \\
& = & \left\{\{f_{\{1\},\beta}\}_{\beta \in \Z_{\geq 0}^3},
\{f_{\{2\},k}(z_1,z_3)\}_{k \in
\Z_{\geq 0}},\{f_{\{3\},k}(z_1,z_2)\}_{k \in
\Z_{\geq 0}}, \right. \\
&& \{f_{\{1,2\},\beta}\}_{\beta \in \Z_{\geq 0}^3},\{f_{\{1,3\},\beta}\}_{\beta \in \Z_{\geq 0}^3},\{f_{\{2,3\},\alpha}(z_1)\}_{\alpha \in \Z_{\geq 0}^2}, \\
&& \left.\{f_{\{1,2,3\},\beta}\}_{\beta \in \Z_{\geq 0}^3}\right\},
\end{eqnarray*}
where $f_{\{2\},k}(z_1,z_3)$ (resp. $f_{\{3\},k}(z_1,z_2)$, resp. $f_{\{2,3\},\alpha}(z_1)$) is holomorphic in $S_{\{2\}}$ (resp. $S_{\{3\}}$, resp $S_{\{2,3\}}$) and $f_{\{1\},\beta},f_{\{1,2\},\beta},f_{\{1,3\},\beta},f_{\{1,2,3\},\beta} \in \CC$. An asymptotic expansion $\operatorname{App}^{<N}(F;\,z)$, $N=(n_1,n_2,n_3) \in \Z_{\geq 0}^3$ is given by
$$
\begin{aligned}
	T_{\{1\}}^{<N}(F;\, z) &= \sum_{\beta_1 + \beta_2 + \beta_3 < n_1} f_{\{1\},\beta}\frac{z_1^{\beta_1}z_2^{\beta_2}z_3^{\beta_3}}{\beta_1!\beta_2!\beta_3!}, \qquad \beta_1,\beta_2,\beta_3 \in \Z_{\geq 0}, \\
	T_{\{2\}}^{<N}(F;\, z) &= \sum_{k < n_2} f_{\{2\},k}(z_1,z_3)\frac{z_2^{k}}{k!}, \qquad k \in \Z_{\geq 0}, \\	
	T_{\{3\}}^{<N}(F;\, z) &= \sum_{k < n_3} f_{\{3\},k}(z_1,z_2)\frac{z_3^k}{k!}, \qquad k \in \Z_{\geq 0}, \\	
T_{\{1,2\}}^{<N}(F;\, z) &= \sum_{\substack{\beta_1+\beta_2+\beta_3 < n_1 \\ \beta_2 < n_2}} f_{\{1,2\},\beta}\frac{z_1^{\beta_1}z_2^{\beta_2}z_3^{\beta_3}}{\beta_1!\beta_2!\beta_3!}, \qquad \beta_1,\beta_2,\beta_3 \in \Z_{\geq 0}, \\
T_{\{1,3\}}^{<N}(F;\, z) &= \sum_{\substack{\beta_1 + \beta_2 + \beta_3 < n_1 \\ \beta_3 < n_3}} f_{\{1,3\},\beta}\frac{z_1^{\beta_1}z_2^{\beta_2}z_3^{\beta_3}}{\beta_1!\beta_2!\beta_3!}, \qquad \beta_1,\beta_2,\beta_3 \in \Z_{\geq 0}, \\
T_{\{2,3\}}^{<N}(F;\, z) &= \sum_{\substack{\alpha_2 < n_2 \\ \alpha_3 < n_3}} f_{\{2,3\},\alpha}(z_1)\frac{z_2^{\alpha_2}z_3^{\alpha_3}}{\alpha_2!\alpha_3!}, \qquad \alpha_2,\alpha_3 \in \Z_{\geq 0}, \\
T_{\{1,2,3\}}^{<N}(F;\, z) &= \sum_{\substack{\beta_1 + \beta_2 + \beta_3 < n_1 \\ \beta_2 < n_2 \\ \beta_3 < n_3}} f_{\{1,2,3\},\beta}\frac{z_1^{\beta_1}z_2^{\beta_2}z_3^{\beta_3}}{\beta_1!\beta_2!\beta_3!}, \qquad \beta_1,\beta_2,\beta_3 \in \Z_{\geq 0}, \\
	\operatorname{App}^{<N}(F;\, z) &=
T_{\{1\}}^{<N}(F;\, z) + T_{\{2\}}^{<N}(F;\, z) + T_{\{3\}}^{<N}(F;\, z) \\
& - T_{\{1,2\}}^{<N}(F;\, z) - T_{\{1,3\}}^{<N}(F;\, z)  - T_{\{2,3\}}^{<N}(F;\, z) \\
& + T_{\{1,2,3\}}^{<N}(F;\, z).
\end{aligned}
$$
We say that $f$ is multi-asymptotically developable to $F$ along
$\chi$ on $S=S(W,\epsilon)$ if and only if for any cone
$S'=S(W',\epsilon')$ properly contained in $S$ and for any $N=(n_1,n_2,n_3)
\in \mathbb{Z}^3_{\ge 0}$, there exists a constant $C_{S',N}$ such that
$$
\left|
f(z) - \operatorname{App}^{<N}(F;\, z)
\right|
\le C_{S',N}|z_1|^{n_1-n_2-n_3}|z_2|^{n_2}|z_3|^{n_3}
 \qquad (z \in S').
$$
The family $F$ is consistent if
\begin{itemize}
\item $f_{\{1\},\beta}=f_{\{1,2\},\beta}=f_{\{1,3\},\beta}=f_{\{1,2,3\},\beta}$ for each $\beta \in \Z^3_{\geq 0}$,
\item $f_{\{2\},k}(z_1,z_3)$ is strongly asymptotically developable to
    $$
    \begin{aligned}
    &\left\{\{f_{\{1,2\},(\beta_1,k,\beta_3)}\}_{(\beta_1,\beta_2) \in \Z^2_{\geq 0}},\{f_{\{2,3\},(k,\alpha_3)}(z_1)\}_{\alpha_3 \in \Z_{\geq 0}},\right.\\
    &\left.\qquad\qquad\qquad\qquad\{f_{\{1,2,3\},(\beta_1,k,\beta_3)}\}_{(\beta_1,\beta_3) \in \Z^2_{\geq 0}}\right\}
    \end{aligned}
    $$
    on $S_{\{2\}}$ for each $k \in \Z_{\geq 0}$,
\item $f_{\{3\},k}(z_1,z_2)$ is strongly asymptotically developable to
    $$
    \begin{aligned}
    &\left\{\{f_{\{1,3\},(\beta_1,\beta_2,k)}\}_{(\beta_1,\beta_2) \in \Z^2_{\geq 0}},\{f_{\{2,3\},(\alpha_2,k)}(z_1)\}_{\alpha_2 \in \Z_{\geq 0}},\right.\\
    &\left.\qquad\qquad\qquad\qquad\{f_{\{1,2,3\},(\beta_1,\beta_2,k)}\}_{(\beta_1,\beta_2) \in \Z^2_{\geq 0}}\right\}
    \end{aligned}
    $$
    on $S_{\{3\}}$ for each $k \in \Z_{\geq 0}$,
\item $f_{\{2,3\},(\alpha_2,\alpha_3)}(z_1)$ is strongly asymptotically developable to
    $$
    \{f_{\{1,2,3\},(\beta_1,\alpha_2,\alpha_3)}\}_{\beta_1 \in \Z_{\geq 0}}
    $$
    on $S_{\{2,3\}}$ for each $(\alpha_2,\alpha_3) \in \Z^2_{\geq 0}$.
\end{itemize}
These are the asymptotics of \cite{HP13} (mixed case) in $\CC^3$.
%%%%%%%%%%%%%%%%%%%%%%%%%%%%%%%%%%%%%%%%%%%%
% MIXED C3 NO2
%%%%%%%%%%%%%%%%%%%%%%%%%%%%%%%%%%%%%%%%%%%%

\

\noindent
$\checkmark\,${\bf{Mixed asymptotics (2 lines clean intersection-Takeuchi).}} In this case
$$
A_\chi=
\left(
\begin{array}{ccc}
1 & 1 & 0 \\
0 & 1 & 1 \\
0 & 0 & 1
\end{array}
\right).
$$
We have $Z_1=\{z_1=z_2=0\}$, $Z_2=\{z_2=z_3=0\}$, $Z_3=\{z_3=0\}$.
Then we can define a  normal deformation $\widetilde{X}=\CC^3 \times \R^3$ with the map $p(z;\,t)=(t_1z_1, t_1t_2z_2, t_2t_3z_3)$. Set $Z=Z_1 \cap Z_2$. By Proposition \ref{prop:bundle-description} the zero section $S_\chi$ of $\widetilde{X}$ has a vector bundle structure
$$
S_\chi \simeq T_ZZ_1 \underset{X}{\times} T_ZZ_2 \underset{X}{\times} \dfrac{TX \underset{X}{\times} Z}
{TZ_1 \underset{X}{\times} Z + TZ_2 \underset{X}{\times} Z}.
$$
The rational monomials $\imin{\varphi_j}$, $j=1,2,3$, are $\imin{\varphi_1}(\tau)=\tau_1$, $\imin{\varphi_2}(\tau)=\displaystyle{\tau_2 \over \tau_1}$, $\imin{\varphi_3}(\tau)=\displaystyle{\tau_1\tau_3 \over \tau_2}$. Let $\xi=(\xi_1,\xi_2,\xi_3)$, $\xi_i \neq 0$, $i=1,2,3$. A cofinal family of $C(\xi)$ is $C(\xi,F^0)$ consisting of elements $S(W,\epsilon)$, $W=W_1 \times W_2 \times W_3$, defined as follows:
$$
S(W,\epsilon)=
\begin{aligned}
 \left\{z \in X;
\begin{aligned}
& z_k \in W_k \quad (k=1,2,3), \\
& |z_1| < \epsilon, \\
& |z_2| < \epsilon |z_1|, \\
& |z_1||z_3| < \epsilon |z_2|\\
\end{aligned}
\right\},
\end{aligned}
$$
where $\epsilon > 0$ and each $W_k$ is an $\mathbb{R}^+$-conic open subset in $\CC$ containing the point $\xi_k$. \\

We now construct the asymptotics.
%We have $K_{\{j\}}=K_j$, $j=1,2$, $K_{\{1,2\}}=K_1 \cup K_2$ and $Z_J=\{z^{(k)}=0,\,k \in K_J\}$ for $\emptyset \neq J \subseteq \{1,2\}$. For any subset $K \subset \{1,2\}$, we denote by $z^{(K)}$ the set of the coordinates $z^{(k)}$'s with $k \in K$ and by $z^{(K)}_{\mathrm{C}}$ the set of the coordinates which do not belong to $z^{(K)}$, i.e., $z^{(k)}$'s with $k \in \{0,1,\dots,m\} \setminus K$. Hence, the coordinates of $Z_J$ are given by $z^{(K_J)}_{\mathrm{C}}$. The $\pi_J$ denotes the canonical projection from $X$ to $Z_J$ defined by $z \to z^{(K_J)}_{\mathrm{C}}$. Given $S(W,\epsilon)$, set $S_J=\pi_J(S(W,\epsilon))$.
We have $S_{\{1,2\}}=S_{\{1,3\}}=S_{\{1,2,3\}}=\{0\}$ and
$$
S_{\{1\}}=
\begin{aligned}
 \left\{z \in Z_1;
\begin{aligned}
&z_3 \in W_3, \\
& |z_3| < \epsilon
\end{aligned}
\right\},
\end{aligned}
$$
$$
S_{\{2\}}=S_{\{2,3\}}=
\begin{aligned}
 \left\{z \in Z_2;
\begin{aligned}
&z_1 \in W_1, \\
& |z_1| < \epsilon
\end{aligned}
\right\},
\end{aligned}
$$
$$
S_{\{3\}}=
\begin{aligned}
 \left\{z \in Z_3;
\begin{aligned}
&z_j \in W_j \quad (j=1,2), \\
&|z_1| < \epsilon, \\
& |z_2| < \epsilon |z_1|
\end{aligned}
\right\}.
\end{aligned}
$$

A total family of coefficients of multi-asymptotic expansion is given by
\begin{eqnarray*}
F & = & \{F_{\{1\}},F_{\{2\}},F_{\{3\}},F_{\{1,2\}},F_{\{1,3\}},F_{\{2,3\}},F_{\{1,2,3\}}\} \\
& = & \left\{\{f_{\{1\},\alpha}(z_3)\}_{\alpha \in \Z_{\geq 0}^2},
\{f_{\{2\},\alpha}(z_1)\}_{\alpha \in
\Z_{\geq 0}^2},\{f_{\{3\},k}(z_1,z_2)\}_{k \in
\Z_{\geq 0}}, \right. \\
&& \{f_{\{1,2\},\beta}\}_{\beta \in \Z_{\geq 0}^3},\{f_{\{1,3\},\beta}\}_{\beta \in \Z_{\geq 0}^3},\{f_{\{2,3\},\alpha}(z_1)\}_{\alpha \in \Z_{\geq 0}^2}, \\
&& \left.\{f_{\{1,2,3\},\beta}\}_{\beta \in \Z_{\geq 0}^3}\right\},
\end{eqnarray*}
where $f_{\{1\},\alpha}(z_3)$ (resp. $f_{\{2\},\alpha}(z_1)$, resp. $f_{\{3\},k}(z_1,z_2)$, resp. $f_{\{2,3\},\alpha}(z_1)$) is holomorphic in $S_{\{1\}}$ (resp. $S_{\{2\}}$, resp. $S_{\{3\}}$, resp. $S_{\{2,3\}}$) and $f_{\{1\},\beta}$, $f_{\{1,2\},\beta}$, $f_{\{1,3\},\beta}$, $f_{\{1,2,3\},\beta} \in \CC$. An asymptotic expansion $\operatorname{App}^{<N}(F;\,z)$, $N=(n_1,n_2,n_3) \in \Z_{\geq 0}^3$ is given by
$$
\begin{aligned}
	T_{\{1\}}^{<N}(F;\, z) &= \sum_{\alpha_1 + \alpha_2 < n_1} f_{\{1\},\alpha}(z_3)\frac{z_1^{\alpha_1}z_2^{\alpha_2}}{\alpha_1!\alpha_2!}, \qquad \alpha_1,\alpha_2 \in \Z_{\geq 0}, \\
	T_{\{2\}}^{<N}(F;\, z) &= \sum_{\alpha_2 + \alpha_3 < n_2} f_{\{2\},\alpha}(z_1)\frac{z_2^{\alpha_2}z_3^{\alpha_3}}{\alpha_2!\alpha_3!}, \qquad \alpha_2,\alpha_3 \in \Z_{\geq 0}, \\	
	T_{\{3\}}^{<N}(F;\, z) &= \sum_{k < n_3} f_{\{3\},k}(z_1,z_2)\frac{z_3^k}{k!}, \qquad k \in \Z_{\geq 0}, \\	
T_{\{1,2\}}^{<N}(F;\, z) &= \sum_{\substack{\beta_1+\beta_2 < n_1 \\ \beta_2+\beta_3 < n_2}} f_{\{1,2\},\beta}\frac{z_1^{\beta_1}z_2^{\beta_2}z_3^{\beta_3}}{\beta_1!\beta_2!\beta_3!}, \qquad \beta_1,\beta_2,\beta_3 \in \Z_{\geq 0}, \\
T_{\{1,3\}}^{<N}(F;\, z) &= \sum_{\substack{\beta_1 + \beta_2 < n_1 \\ \beta_3 < n_3}} f_{\{1,3\},\beta}\frac{z_1^{\beta_1}z_2^{\beta_2}z_3^{\beta_3}}{\beta_1!\beta_2!\beta_3!}, \qquad \beta_1,\beta_2,\beta_3 \in \Z_{\geq 0}, \\
T_{\{2,3\}}^{<N}(F;\, z) &= \sum_{\substack{\alpha_2 + \alpha_3 < n_2 \\ \alpha_3 < n_3}} f_{\{2,3\},\alpha}(z_1)\frac{z_2^{\alpha_2}z_3^{\alpha_3}}{\alpha_2!\alpha_3!}, \qquad \alpha_2,\alpha_3 \in \Z_{\geq 0}, \\
T_{\{1,2,3\}}^{<N}(F;\, z) &= \sum_{\substack{\beta_1 + \beta_2 < n_1 \\ \beta_2 + \beta_3 < n_2 \\ \beta_3 < n_3}} f_{\{1,2,3\},\beta}\frac{z_1^{\beta_1}z_2^{\beta_2}z_3^{\beta_3}}{\beta_1!\beta_2!\beta_3!}, \qquad \beta_1,\beta_2,\beta_3 \in \Z_{\geq 0}, \\
	\operatorname{App}^{<N}(F;\, z) &=
T_{\{1\}}^{<N}(F;\, z) + T_{\{2\}}^{<N}(F;\, z) + T_{\{3\}}^{<N}(F;\, z) \\
& - T_{\{1,2\}}^{<N}(F;\, z) - T_{\{1,3\}}^{<N}(F;\, z)  - T_{\{2,3\}}^{<N}(F;\, z) \\
& + T_{\{1,2,3\}}^{<N}(F;\, z).
\end{aligned}
$$
We say that $f$ is multi-asymptotically developable to $F$ along
$\chi$ on $S=S(W,\epsilon)$ if and only if for any cone
$S'=S(W',\epsilon')$ properly contained in $S$ and for any $N=(n_1,n_2,n_3)
\in \mathbb{Z}^3_{\ge 0}$, there exists a constant $C_{S',N}$ such that
$$
\left|
f(z) - \operatorname{App}^{<N}(F;\, z)
\right|
\le C_{S',N}|z_1|^{n_1-n_2+n_3}|z_2|^{n_2-n_3}|z_3|^{n_3}
 \qquad (z \in S').
$$
The family $F$ is consistent if
\begin{itemize}
\item $f_{\{1,2\},\beta}=f_{\{1,3\},\beta}=f_{\{1,2,3\},\beta}$ for each $\beta \in \Z^3_{\geq 0}$,
\item $f_{\{2\},\alpha}=f_{\{2,3\},\alpha}$ for each $\alpha \in \Z^2_{\geq 0}$,
\item $f_{\{1\},(\alpha_1,\alpha_2)}(z_3)$ is strongly asymptotically developable to
    $$
    \{f_{\{1,2,3\},(\alpha_1,\alpha_2,\beta_3)}\}_{\beta_3 \in \Z_{\geq 0}} $$
    on $S_{\{1\}}$ for each $(\alpha_1,\alpha_2) \in \Z^2_{\geq 0}$,
\item $f_{\{2\},(\alpha_2,\alpha_3)}(z_1)$ is strongly asymptotically developable to
    $$
    \{f_{\{1,2,3\},(\beta_1,\alpha_2,\alpha_3)}\}_{\beta_1 \in \Z_{\geq 0}}
    $$
    on $S_{\{2\}}$ for each $(\alpha_2,\alpha_3) \in \Z_{\geq 0}$,
\item $f_{\{3\},k}(z_1,z_2)$ is strongly asymptotically developable to
    $$
    \begin{aligned}
    &\left\{\{f_{\{1,3\},(\beta_1,\beta_2,k)}\}_{(\beta_1,\beta_2) \in \Z^2_{\geq 0}},\{f_{\{2,3\},(\alpha_2,k)}(z_1)\}_{\alpha_2 \in \Z_{\geq 0}},\right.\\
    &\left.\qquad\qquad\qquad\qquad\{f_{\{1,2,3\},(\beta_1,\beta_2,k)}\}_{(\beta_1,\beta_2) \in \Z^2_{\geq 0}}\right\}
    \end{aligned}
    $$
    on $S_{\{3\}}$ for each $k \in \Z_{\geq 0}$.
\end{itemize}

%\end{itemize}

%%%%%%%%%%%%%%%%%%%%%%%%%%%%%%%%%%%%%%%%%%%%
% END OF C3
%%%%%%%%%%%%%%%%%%%%%%%%%%%%%%%%%%%%%%%%%%%%

\subsection{Classification of multi-asymptotics for the case of two submanifolds}

Let us now consider $X=\CC^n$ with variables $z=(z^{(1)},\dots, z^{(m)})$, where $z^{(k)}$, $k=1,\dots,m$, denote the coordinates blocks $(z_{i_{k,1}},\dots,z_{i_{k,n_k}})$ (we forget $z^{(0)}$ to lighten notations). We first consider the normal deformation (and the associated multi-specialization and asymptotics) in the two manifold case,
i.e., $\chi=\{Z_1,Z_2\}$ with $Z_j=\{z^{(k)}=0,\,k \in K_j\}$, $K_j \subseteq \{1,\dots,m\}$, $j=1,2$. Remark that $Z_1$ and $Z_2$ may coincide. Then, we give all possible
non-degenerate multi-asymptotics for the case of two manifolds
with two or three coordinates blocks.

The action associated with $\chi$ is now defined by
$$
\mu_j(z,\lambda)=(\lambda^{a_{j1}}z^{(1)},\dots,\lambda^{a_{jm}}z^{(m)}) \qquad (j=1,2)
$$
with $a_{ji}$ non-negative rationals, $a_{ji} \neq 0$
if $i \in K_j$, $a_{ji}=0$ otherwise.
%Up to multiply by $\sigma_A$, the less common multiple of the denominators of $a_{ji}$, $j=1,2$, $i=1,\dots,m$, we may assume $a_{ji}$ non-negative integer.
The associated matrix $A_\chi$ is
$$
\left(
\begin{array}{ccc}
a_{11} & \cdots & a_{1m} \\
a_{21} & \cdots & a_{2m}
\end{array}
\right)
$$
and we assume for simplicity that all the columns are non zero. We also assume that the action is non-degenerate and that the matrix $A'_\chi$  consisting of the first two columns of $A_\chi$ is invertible. Then we can define a  normal deformation $\widetilde{X}=\CC^n \times \R^2$ with the map $p:\widetilde{X} \to X$ defined by
$$
p(z;\,t)=(z^{(0)},\,\varphi_1(t)z^{(1)}, \dots, \varphi_m(t)z^{(m)})
$$
with
\begin{equation}
\varphi_k(t)=t_1^{a_{1k}}t_2^{a_{2k}} \qquad (k = 1,2,\dots,m).
\end{equation}
Following the notations of Section \ref{section:Multi-cones for a general case} we are going to study the family of multicones $C(\xi)$ associated to a point $\xi=(\xi_1,\dots,\xi_n) \in S_\chi$. The rational monomials $\imin{\varphi_j}$, $j=1,2$ are determined by $A'_\chi{}^{-1}$, where
$$
A'_\chi=
\left(
\begin{array}{cc}
a_{11} & a_{12} \\
a_{21} & a_{22}
\end{array}
\right) \ \ \ \ \mbox{and} \ \ \ \
A'_\chi{}^{-1}={1 \over a_{11}a_{22}-a_{12}a_{21}}
\left(
\begin{array}{cc}
a_{22} & -a_{12} \\
-a_{21} & a_{11}
\end{array}
\right).
$$
Setting $d=\displaystyle{1 \over a_{11}a_{22}-a_{12}a_{21}}$, we have $\imin{\varphi_1}(\tau)=(\tau_1^{a_{22}}\tau_2^{-a_{21}})^{d}$ and $\imin{\varphi_2}(\tau)=(\tau_1^{-a_{12}}\tau_2^{a_{11}})^{d}$. For $k \geq 3$ we have
$$
\psi_k(\tau)={\tau_k \over (\tau_1^{a_{1k}a_{22}-a_{2k}a_{12}}\tau_2^{a_{2k}a_{11}-a_{1k}a_{21}})^d}.
$$
Here we may assume $d > 0$ by exchanging the first and the second columns of $A_\chi$ if necessary.
Suppose that $\xi$ is outside of the set of fixed points.
By definition (up to take a permutation of coordinates blocks), we may assume that $\xi^{(1)},\xi^{(2)} \neq 0$.
Further, for simplicity, we also assume $\xi$ to be normalized, that is,
$|\xi^{(1)}| = |\xi^{(2)}| = 1$.
A cofinal family of $C(\xi)$ is $C(\xi,F^0)$ consisting of elements $S(W,\epsilon)$, $W=W_1 \times \cdots \times W_m$, defined as follows
$$
S(W,\epsilon)=
\begin{aligned}
 \left\{z \in X;
\begin{aligned}
&z^{(k)} \in W_k \quad (k=1,\dots,m), \\
& |z^{(i)}|^{a_{jj}} < \epsilon^{1 \over d}|z^{(j)}|^{a_{ji}} \quad (i,j \in \{1,2\},\,\ i \neq j), \\
&|\xi^{(k)}|-\epsilon < {|z^{(k)}| \over (|z^{(1)}|^{a_{1k}a_{22}-a_{2k}a_{12}}|z^{(2)}|^{a_{2k}a_{11}-a_{1k}a_{21}})^d}  \\
&\qquad\qquad\qquad\qquad< |\xi^{(k)}|+\epsilon \quad (k \geq 3)
\end{aligned}
\right\},
\end{aligned}
$$
where $\epsilon > 0$ and each $W_k$ is an $\mathbb{R}^+$-conic open
subset in $\CC^{n_k}$ containing the point $\xi^{(k)}$. \\

We now construct the asymptotics. We recall the notations of Section \ref{section:Multi-asymptotic expansions}. Set $K_{\{j\}}=K_j$, $j=1,2$, $K_{\{1,2\}}=K_1 \cup K_2$ and $Z_J=\{z^{(k)}=0,\,k \in K_J\}$ for $\emptyset \neq J \subseteq \{1,2\}$. For any subset $J \subset \{1,2\}$,
we denote by $z^{(J)}$ the set of the coordinates $z^{(k)}$'s with $k \in K_J$
and by $z^{(J)}_{\mathrm{C}}$ the set of the coordinates
which do not belong to $z^{(J)}$, i.e.,
$z^{(k)}$'s with $k \in \{0,1,\dots,m\} \setminus K_J$.
Hence, the coordinates of $Z_J$ are given by $z^{(J)}_{\mathrm{C}}$.
The $\pi_J$ denotes the canonical projection from $X$ to $Z_J$ defined
by $z \to z^{(J)}_{\mathrm{C}}$.
\begin{oss}
The notations $z^{(J)}$ and $z_{\mathrm{C}}^{(J)}$ are
denoted by $z^{(K_J)}$ and $z^{(\compl K_J)}$
in Section \ref{section:Multi-asymptotic expansions}, respectively.
To make the notations light, we use
$z^{(J)}$ and $z_{\mathrm{C}}^{(J)}$ instead of
$z^{(K_J)}$ and $z^{(\compl K_J)}$ through the section.
Similary we use $\mathbb{Z}_{\ge 0}^{(J)}$
instead of $\mathbb{Z}_{\ge 0}^{(K_J)}$ hereafter.
\end{oss}

Given $S(W,\epsilon)$, set $S_J=\pi_J(S(W,\epsilon))$. A total family of coefficients of multi-asymptotic expansion is given by
\begin{eqnarray*}
F & = & \{F_{\{1\}},F_{\{2\}},F_{\{1,2\}}\} \\
& = & \left\{\{f_{\{1\},\alpha}(z^{(\{1\})}_{\mathrm{C}})\}_{\alpha \in \Z_{\geq 0}^{(\{1\})}},\,\,
\{f_{\{2\},\alpha}(z^{(\{2\})}_{\mathrm{C}})\}_{\alpha \in \Z_{\geq 0}^{(\{2\})}},\right. \\
& &\qquad\qquad\qquad
\left.
\{f_{\{1,2\},\alpha}(z^{(\{1,2\})}_{\mathrm{C}})\}_{\alpha \in \Z_{\geq 0}^{(\{1,2\})}}
\right\},
\end{eqnarray*}
where $f_{J,\alpha}$ is holomorphic on $S_J$, $\alpha \in \Z_{\geq 0}^{(J)}$, $\emptyset \neq J \subseteq \{1,2\}$. Here
$$
\Z_{\geq 0}^{(J)}=\left\{(\alpha^{(1)},\dots,\alpha^{(m)}) \in \Z^n_{\geq 0};\, \alpha^{(k)}=0,\, k \notin K_J\right\}.
$$
Let $\sigma_A$ be the less common multiple of the denominators of $a_{ji}$, $j=1,2$, $i=1,\dots,m$. We are now ready to define an asymptotic expansion $\operatorname{App}^{<N}(F;\,z)$, $N=(n_1,n_2) \in \Z_{\geq 0}^2$. It is given by
$$
\begin{aligned}
	T_{\{1\}}^{<N}(F;\, z) &= \sum_{a_{11}|\alpha^{(1)}|+\cdots+a_{1m}|\alpha^{(m)}| < n_1/\sigma_A} f_{\{1\},\alpha}(z^{(\{1\})}_{\mathrm{C}})\frac{z^\alpha}{\alpha!}, \qquad \alpha \in \Z_{\geq 0}^{(\{1\})}, \\
	T_{\{2\}}^{<N}(F;\, z) &= \sum_{a_{21}|\alpha^{(1)}|+\cdots+a_{2m}|\alpha^{(m)}| < n_2/\sigma_A} f_{\{2\},\alpha}(z^{(\{2\})}_{\mathrm{C}})\frac{z^\alpha}{\alpha!}, \qquad \alpha \in \Z_{\geq 0}^{(\{2\})}, \\
	T_{\{1,2\}}^{<N}(F;\, z) &= \sum_{\substack{ a_{11}|\alpha^{(1)}|+\cdots+a_{1m}|\alpha^{(m)}| < n_1/\sigma_A \\ a_{21}|\alpha^{(1)}|+\cdots+a_{2m}|\alpha^{(m)}| < n_2/\sigma_A}} f_{\{1,2\},\alpha}(z^{(\{1,2\})}_{\mathrm{C}})\frac{z^\alpha}{\alpha!}, \qquad \alpha \in \Z_{\geq 0}^{(\{1,2\})}, \\
	\operatorname{App}^{<N}(F;\, z) &=
T_{\{1\}}^{<N}(F;\, z) + T_{\{2\}}^{<N}(F;\, z) - T_{\{1,2\}}^{<N}(F;\, z).
\end{aligned}
$$
We say that $f$ is multi-asymptotically developable to $F=\{F_{\{1\}},F_{\{2\}},F_{\{1,2\}}\}$ along
$\chi$ on $S=S(W,\epsilon)$ if and only if for any cone
$S'=S(W',\epsilon')$ properly contained in $S$ (i.e. $\overline{W'_k} \setminus \{0\} \subset W_k$, $k=1,\dots,m$, and $\epsilon'<\epsilon$, cf. Definition \ref{definition:proper sub-cone}) and for any $N=(n_1,n_2)
\in \mathbb{Z}^2_{\ge 0}$, there exists a constant $C_{S',N}$ such that
$$
\left|
f(z) - \operatorname{App}^{<N}(F;\, z)
\right|
\le C_{S',N}(|z^{(1)}|^{a_{22}n_1-a_{12}n_2}|z^{(2)}|^{a_{11}n_2-a_{21}n_1})^{d/\sigma_A}
 \qquad (z \in S').
$$

We are now ready to compute the fibers of the specialization of Whitney holomorphic functions at $\xi$. Let $\OW_X$ be the subanalytic sheaf of Whitney holomorphic functions on $X$. Then
$$
H^0(\nu_\chi\OW_X)_\xi \simeq \lind S H^0(S;\OW_X) \simeq \lind S \lpro {S'}H^0(S';\OW_X),
$$
where $S$ ranges through the family $\{S(W,\epsilon)\}_{W,\epsilon}$,
%$\epsilon>0$, $W=W_1,\dots,W_m$, with $W_k$ $\mathbb{R}^+$-conic open subset in $\CC^{n_k}$ containing $\xi^{(k)}$.
$S'$ ranges through the family of multicones properly contained in $S$. \\

We are now going to classify multi-specializations and associated asymptotics in the two manifolds case for $m \leq 3$. We perform the classification by
\begin{itemize}
\item the number $m$ of columns of $A_\chi$,
\item the number of non-zero entries of $A_\chi$ (denoted by $N$)
assuming that the action is non-degenerate (in that case we are reduced to the one manifold case)
\end{itemize}
If $m=1$ the action is degenerate, so we consider $m=2,3$. Moreover, the fact that row (resp. column) permutation and row (resp. column) multiplication by a positive rational number do not change the multi-normal deformation and the associated asymptotics (up to permutation of variable blocks) , we may classify $A_\chi$ up to these operations. In particular, we assume $a_{11}=a_{22}=1$ and $1-a_{12}a_{21} > 0$. \\

Let us consider the case $m=2$. We have 2 coordinates blocks $z^{(1)}=\{z_1,\dots,z_{m_1}\}$ and $z^{(2)}=\{z_{m_1+1},\dots,z_{m_1+m_2}\}$ ($m_1+m_2=n$). In this case there are 5 possibilities.
%\begin{itemize}

%%%%%%%%%%%%%%%%%%%%%%%%%%%%%%%%%%
% N=0,1
%%%%%%%%%%%%%%%%%%%%%%%%%%%%%%%%%%

\

\noindent
$\checkmark\,$ {\bf{Case $m=2$ and $N=0,1$.}} The action is degenerate.

%%%%%%%%%%%%%%%%%%%%%%%%%%%%%%%%%%
% N=2
%%%%%%%%%%%%%%%%%%%%%%%%%%%%%%%%%%

\

\noindent
$\checkmark\,$ {\bf{Case $m=2,N=2$.}} In this case
$$
A_\chi=
\left(
\begin{array}{cc}
1 & 0 \\
0 & 1
\end{array}
\right).
$$
We have $Z_j=\{z^{(j)}=0\}$, $j=1,2$. Then we can define a  normal deformation $\widetilde{X}=\CC^n \times \R^2$ with the map $p(z;\,t)=(t_1z^{(1)}, t_2z^{(2)})$.
The rational monomials $\imin{\varphi_j}$, $j=1,2$ are $\imin{\varphi_1}(\tau)=\tau_1$ and $\imin{\varphi_2}(\tau)=\tau_2$.  A cofinal family of $C(\xi)$ is $C(\xi,F^0)$ consisting of elements $S(W,\epsilon)$, $W=W_1 \times W_2$, defined as follows
$$
S(W,\epsilon)=
\begin{aligned}
 \left\{z \in X;
\begin{aligned}
&z^{(k)} \in W_k \quad (k=1,2), \\
& |z^{(i)}| < \epsilon \quad (i = 1,2) \\
\end{aligned}
\right\},
\end{aligned}
$$
where $\epsilon > 0$ and each $W_k$ is an $\mathbb{R}^+$-conic open subset in $\CC^{m_k}$ containing the point $\xi^{(k)}$. \\

We now construct the asymptotics.
%We have $K_{\{j\}}=K_j$, $j=1,2$, $K_{\{1,2\}}=K_1 \cup K_2$ and $Z_J=\{z^{(k)}=0,\,k \in K_J\}$ for $\emptyset \neq J \subseteq \{1,2\}$. For any subset $K \subset \{1,2\}$, we denote by $z^{(K)}$ the set of the coordinates $z^{(k)}$'s with $k \in K$ and by $z^{(K)}_{\mathrm{C}}$ the set of the coordinates which do not belong to $z^{(K)}$, i.e., $z^{(k)}$'s with $k \in \{0,1,\dots,m\} \setminus K$. Hence, the coordinates of $Z_J$ are given by $z^{(K_J)}_{\mathrm{C}}$. The $\pi_J$ denotes the canonical projection from $X$ to $Z_J$ defined by $z \to z^{(K_J)}_{\mathrm{C}}$. Given $S(W,\epsilon)$, set $S_J=\pi_J(S(W,\epsilon))$.
We have $S_{\{1,2\}}=\{\mathrm{pt}\}$ and, for $i \neq j \in \{1,2\}$
$$
S_{\{i\}}=
\begin{aligned}
 \left\{z \in Z_i;
\begin{aligned}
&z^{(j)} \in W_j, \\
& |z^{(j)}| < \epsilon
\end{aligned}
\right\}.
\end{aligned}
$$
A total family of coefficients of multi-asymptotic expansion is given by
\begin{eqnarray*}
F & = & \{F_{\{1\}},F_{\{2\}},F_{\{1,2\}}\} \\
& = & \left\{\{f_{\{1\},\alpha}(z^{(2)})\}_{\alpha \in \Z_{\geq 0}^{m_1} \times \{0\}^{m_2}},
\{f_{\{2\},\alpha}(z^{(1)})\}_{\alpha \in \{0\}^{m_1} \times \Z_{\geq 0}^{m_2}},\right. \\
&& \qquad\qquad\qquad\left.\{f_{\{1,2\},\alpha}\}_{\alpha \in \Z_{\geq 0}^n}\right\},
\end{eqnarray*}
where $f_{\{1\},\alpha}$ (resp. $f_{\{2\},\alpha}$) is holomorphic on $S_{\{2\}}$ (resp. $S_{\{1\}}$) and $f_{\{1,2\},\alpha} \in \CC$. An asymptotic expansion $\operatorname{App}^{<N}(F;\,z)$, $N=(n_1,n_2) \in \Z_{\geq 0}^2$ is given by
$$
\begin{aligned}
	T_{\{1\}}^{<N}(F;\, z) &= \sum_{|\alpha^{(1)}| < n_1} f_{\{1\},\alpha}(z^{(2)})\frac{z^\alpha}{\alpha!}, \qquad \alpha \in \Z_{\geq 0}^{m_1} \times \{0\}^{m_2}, \\
	T_{\{2\}}^{<N}(F;\, z) &= \sum_{|\alpha^{(2)}| < n_2} f_{\{2\},\alpha}(z^{(1)})\frac{z^\alpha}{\alpha!}, \qquad \alpha \in \{0\}^{m_1} \times \Z_{\geq 0}^{m_2}, \\	
	T_{\{1,2\}}^{<N}(F;\, z) &= \sum_{|\alpha^{(1)}| < n_1 ,\, |\alpha^{(2)}| < n_2} f_{\{1,2\},\alpha}\frac{z^\alpha}{\alpha!}, \qquad \alpha \in \Z_{\geq 0}^{n}, \\
	\operatorname{App}^{<N}(F;\, z) &=
T_{\{1\}}^{<N}(F;\, z) + T_{\{2\}}^{<N}(F;\, z) - T_{\{1,2\}}^{<N}(F;\, z).
\end{aligned}
$$
We say that $f$ is multi-asymptotically developable to $F=\{F_{\{1\}},F_{\{2\}},F_{\{1,2\}}\}$ along
$\chi$ on $S=S(W,\epsilon)$ if and only if for any cone
$S'=S(W',\epsilon')$ properly contained in $S$ and for any $N=(n_1,n_2)
\in \mathbb{Z}^2_{\ge 0}$, there exists a constant $C_{S',N}$ such that
$$
\left|
f(z) - \operatorname{App}^{<N}(F;\, z)
\right|
\le C_{S',N}|z^{(1)}|^{n_1}|z^{(2)}|^{n_2}
 \qquad (z \in S').
$$
The family $F$ is consistent if
\begin{itemize}
\item $f_{\{1\},(\alpha^{(1)},0)}(z^{(2)})$ is strongly asymptotically developable to $$\{f_{\{1,2\},(\alpha^{(1)},\alpha^{(2)})}\}_{\alpha^{(2)} \in \Z_{\geq 0}^{m_2}}$$ on $S_{\{1\}}$ for each $\alpha^{(1)} \in \Z_{\geq 0}^{m_1}$. Namely, for any cone
$S'_{\{1\}}$ properly contained in $S_{\{1\}}$ and for any $n_2
\in \mathbb{Z}_{\ge 0}$, there exists a constant $C_{S'_{\{1\}},n_2}$ such that
$$
\left|
f_{\{1\},(\alpha^{(1)},0)}(z^{(2)}) - \sum_{|\alpha^{(2)}|<n_2}f_{\{1,2\},(\alpha^{(1)},\alpha^{(2)})}{(z^{(2)})^{\alpha^{(2)}} \over \alpha^{(2)}!}
\right|
\le C_{S'_{\{1\}},n_2}|z^{(2)}|^{n_2},
$$
\item $f_{\{2\},(0,\alpha^{(2)})}(z^{(1)})$ is strongly asymptotically developable to $$\{f_{\{1,2\},(\alpha^{(1)},\alpha^{(2)})}\}_{\alpha^{(1)} \in \Z_{\geq 0}^{m_1}}$$ on $S_{\{2\}}$ for each $\alpha^{(2)} \in \Z_{\geq 0}^{m_2}$. Namely, for any cone
$S'_{\{2\}}$ properly contained in $S_{\{2\}}$ and for any $n_1
\in \mathbb{Z}_{\ge 0}$, there exists a constant $C_{S'_{\{2\}},n_1}$ such that
$$
\left|
f_{\{2\},(0,\alpha^{(2)})}(z^{(1)}) - \sum_{|\alpha^{(1)}|<n_1}f_{\{1,2\},(\alpha^{(1)},\alpha^{(2)})}{(z^{(1)})^{\alpha^{(1)}} \over \alpha^{(1)}!}
\right|
\le C_{S'_{\{2\}},n_1}|z^{(1)}|^{n_1}.
$$
\end{itemize}
%%%%%%%%%%%%%%%%%%%%%%%%%%%%%%%
% N=3
%%%%%%%%%%%%%%%%%%%%%%%%%%%%%%%%

\

\noindent
$\checkmark\,$ {\bf{Case $m=2,N=3$.}} In this case
$$
A_\chi=
\left(
\begin{array}{cc}
1 & b \\
0 & 1
\end{array}
\right),
$$
with $b \in \mathbb{Q}_{>0}$. We have $Z_1= \{0\}$, $Z_2 = \{z^{(2)}=0\}$. Then we can define a  normal deformation $\widetilde{X}=\CC^n \times \R^2$ with the map $p(z;\,t)=(t_1z^{(1)}, t_1^bt_2z^{(2)})$.
The rational monomials $\imin{\varphi_j}$, $j=1,2$, are $\imin{\varphi_1}(\tau)=\tau_1$ and $\imin{\varphi_2}(\tau)=\displaystyle{\tau_2 \over \tau_1^b}$.  A cofinal family of $C(\xi)$ is $C(\xi,F^0)$ consisting of elements $S(W,\epsilon)$, $W=W_1 \times W_2$, defined as follows
$$
S(W,\epsilon)=
\begin{aligned}
 \left\{z \in X;
\begin{aligned}
&z^{(k)} \in W_k \quad (k=1,2), \\
& |z^{(1)}| < \epsilon, \\
& |z^{(2)}| < \epsilon |z^{(1)}|^b
\end{aligned}
\right\},
\end{aligned}
$$
where $\epsilon > 0$ and each $W_k$ is an $\mathbb{R}^+$-conic open subset in $\CC^{m_k}$ containing the point $\xi^{(k)}$. \\

We now construct the asymptotics.
%We have $K_{\{j\}}=K_j$, $j=1,2$, $K_{\{1,2\}}=K_1 \cup K_2$ and $Z_J=\{z^{(k)}=0,\,k \in K_J\}$ for $\emptyset \neq J \subseteq \{1,2\}$. For any subset $K \subset \{1,2\}$, we denote by $z^{(K)}$ the set of the coordinates $z^{(k)}$'s with $k \in K$ and by $z^{(K)}_{\mathrm{C}}$ the set of the coordinates which do not belong to $z^{(K)}$, i.e., $z^{(k)}$'s with $k \in \{0,1,\dots,m\} \setminus K$. Hence, the coordinates of $Z_J$ are given by $z^{(K_J)}_{\mathrm{C}}$. The $\pi_J$ denotes the canonical projection from $X$ to $Z_J$ defined by $z \to z^{(K_J)}_{\mathrm{C}}$. Given $S(W,\epsilon)$, set $S_J=\pi_J(S(W,\epsilon))$.
We have $S_{\{1\}}=S_{\{1,2\}}=\{\mathrm{pt}\}$ and
$$
S_{\{2\}}=
\begin{aligned}
 \left\{z \in Z_2;
\begin{aligned}
&z^{(1)} \in W_1, \\
& |z^{(1)}| < \epsilon
\end{aligned}
\right\}.
\end{aligned}
$$
A total family of coefficients of multi-asymptotic expansion is given by
\begin{eqnarray*}
F & = & \{F_{\{1\}},F_{\{2\}},F_{\{1,2\}}\} \\
& = & \left\{\{f_{\{1\},\alpha}\}_{\alpha \in \Z_{\geq 0}^n},
\{f_{\{2\},\alpha}(z^{(1)})\}_{\alpha \in \{0\}^{m_1} \times \Z_{\geq 0}^{m_2}},
\{f_{\{1,2\},\alpha}\}_{\alpha \in \Z_{\geq 0}^n}\right\},
\end{eqnarray*}
where $f_{\{2\},\alpha}$ is holomorphic on $S_{\{1\}}$ and $f_{\{1\},\alpha},f_{\{1,2\},\alpha} \in \CC$. Let $\sigma_A$ be the denominator of $b$. An asymptotic expansion $\operatorname{App}^{<N}(F;\,z)$, $N=(n_1,n_2) \in \Z_{\geq 0}^2$ is given by
$$
\begin{aligned}
	T_{\{1\}}^{<N}(F;\, z) &= \sum_{|\alpha^{(1)}|+b|\alpha^{(2)}| < n_1/\sigma_A} f_{\{1\},\alpha}\frac{z^\alpha}{\alpha!}, \qquad \alpha \in \Z_{\geq 0}^{n}, \\
	T_{\{2\}}^{<N}(F;\, z) &= \sum_{|\alpha^{(2)}| < n_2/\sigma_A} f_{\{2\},\alpha}(z^{(1)})\frac{z^\alpha}{\alpha!}, \qquad \alpha \in \{0\}^{m_1} \times \Z_{\geq 0}^{m_2}, \\	
	T_{\{1,2\}}^{<N}(F;\, z) &= \sum_{\substack{|\alpha^{(1)}|+b|\alpha^{(2)}| < n_1/\sigma_A \\ |\alpha^{(2)}| < n_2/\sigma_A}} f_{\{1,2\},\alpha}\frac{z^\alpha}{\alpha!}, \qquad \alpha \in \Z_{\geq 0}^{n}, \\
	\operatorname{App}^{<N}(F;\, z) &=
T_{\{1\}}^{<N}(F;\, z) + T_{\{2\}}^{<N}(F;\, z) - T_{\{1,2\}}^{<N}(F;\, z).
\end{aligned}
$$
We say that $f$ is multi-asymptotically developable to $F=\{F_{\{1\}},F_{\{2\}},F_{\{1,2\}}\}$ along
$\chi$ on $S=S(W,\epsilon)$ if and only if for any cone
$S'=S(W',\epsilon')$ properly contained in $S$ and for any $N=(n_1,n_2)
\in \mathbb{Z}^2_{\ge 0}$, there exists a constant $C_{S',N}$ such that
$$
\left|
f(z) - \operatorname{App}^{<N}(F;\, z)
\right|
\le C_{S',N}(|z^{(1)}|^{n_1-bn_2}|z^{(2)}|^{n_2})^{1/\sigma_A}
 \qquad (z \in S').
$$
The family $F$ is consistent if
\begin{itemize}
\item $f_{\{1\},\alpha}=f_{\{1,2\},\alpha}$ for each $\alpha \in \Z_{\geq 0}^{n}$,
\item $f_{\{2\},(0,\alpha^{(2)})}(z^{(1)})$ is strongly asymptotically developable to $$\{f_{\{1,2\},(\alpha^{(1)},\alpha^{(2)})}\}_{\alpha^{(1)} \in \Z_{\geq 0}^{m_1}}$$ on $S_{\{2\}}$ for each $\alpha^{(2)} \in \Z_{\geq 0}^{m_2}$. Namely, for any cone
$S'_{\{2\}}$ properly contained in $S_{\{2\}}$ and for any $n_1
\in \mathbb{Z}_{\ge 0}$, there exists a constant $C_{S'_{\{2\}},n_1}$ such that
$$
\left|
f_{\{2\},(0,\alpha^{(2)})}(z^{(1)}) - \sum_{|\alpha^{(1)}|<n_1}f_{\{1,2\},(\alpha^{(1)},\alpha^{(2)})}{(z^{(1)})^{\alpha^{(1)}} \over \alpha^{(1)}!}
\right|
\le C_{S'_{\{2\}},n_1}|z^{(1)}|^{n_1}.
$$
\end{itemize}
%%%%%%%%%%%%%%%%%%%%%%%%%%%%%%%%
% N=4
%%%%%%%%%%%%%%%%%%%%%%%%%%%%%%%%

\

\noindent
$\checkmark\,$ {\bf{Case $m=2,N=4$.}} In this case
$$
A_\chi=
\left(
\begin{array}{cc}
1 & b \\
c & 1
\end{array}
\right),
$$
with $b,c \in \mathbb{Q}_{>0}$ and $1-bc>0$. We have $Z_1=Z_2 = \{0\}$. Then we can define a  normal deformation $\widetilde{X}=\CC^n \times \R^2$ with the map $p(z;\,t)=(t_1t_2^cz^{(1)}, t_1^bt_2z^{(2)})$. Set $d=\displaystyle{1 \over 1-bc}$.
The rational monomials $\imin{\varphi_j}$, $j=1,2$ are $\imin{\varphi_1}(\tau)=\left(\displaystyle{\tau_1 \over \tau_2^c}\right)^d$ and $\imin{\varphi_2}(\tau)=\displaystyle\left({\tau_2 \over \tau_1^b}\right)^d$.  A cofinal family of $C(\xi)$ is $C(\xi,F^0)$ consisting of elements $S(W,\epsilon)$, $W=W_1 \times W_2$, defined as follows
$$
S(W,\epsilon)=
\begin{aligned}
 \left\{z \in X;
\begin{aligned}
&z^{(k)} \in W_k \quad (k=1,2), \\
& |z^{(1)}| < \epsilon^{1 \over d} |z^{(2)}|^c, \\
& |z^{(2)}| < \epsilon^{1 \over d} |z^{(1)}|^b
\end{aligned}
\right\},
\end{aligned}
$$
where $\epsilon > 0$ and each $W_k$ is an $\mathbb{R}^+$-conic open subset in $\CC^{m_k}$ containing the point $\xi^{(k)}$. \\

We now construct the asymptotics.
%We have $K_{\{j\}}=K_j$, $j=1,2$, $K_{\{1,2\}}=K_1 \cup K_2$ and $Z_J=\{z^{(k)}=0,\,k \in K_J\}$ for $\emptyset \neq J \subseteq \{1,2\}$. For any subset $K \subset \{1,2\}$, we denote by $z^{(K)}$ the set of the coordinates $z^{(k)}$'s with $k \in K$ and by $z^{(K)}_{\mathrm{C}}$ the set of the coordinates which do not belong to $z^{(K)}$, i.e., $z^{(k)}$'s with $k \in \{0,1,\dots,m\} \setminus K$. Hence, the coordinates of $Z_J$ are given by $z^{(K_J)}_{\mathrm{C}}$. The $\pi_J$ denotes the canonical projection from $X$ to $Z_J$ defined by $z \to z^{(K_J)}_{\mathrm{C}}$. Given $S(W,\epsilon)$, set $S_J=\pi_J(S(W,\epsilon))$.
We have $S_{\{1\}}=S_{\{2\}}=S_{\{1,2\}}=\{0\}$.
%and
%$$
%S_{\{1\}}=
%\begin{aligned}
% \left\{z \in X;
%\begin{aligned}
%&z^{(1)} \in W_1, \\
%& |z^{(1)}| < \epsilon
%\end{aligned}
%\right\},
%\end{aligned}
%$$
A total family of coefficients of multi-asymptotic expansion is given by
\begin{eqnarray*}
F & = & \{F_{\{1\}},F_{\{2\}},F_{\{1,2\}}\} \\
& = & \left\{\{f_{\{1\},\alpha}\}_{\alpha \in \Z_{\geq 0}^n},
\{f_{\{2\},\alpha}\}_{\alpha \in \Z_{\geq 0}^n},
\{f_{\{1,2\},\alpha}\}_{\alpha \in \Z_{\geq 0}^n}\right\},
\end{eqnarray*}
where $f_{\{1\},\alpha},f_{\{2\},\alpha},f_{\{1,2\},\alpha} \in \CC$. Let $\sigma_A$ be the less common multiple of the denominators of $b,c$. An asymptotic expansion $\operatorname{App}^{<N}(F;\,z)$, $N=(n_1,n_2) \in \Z_{\geq 0}^2$ is given by
$$
\begin{aligned}
	T_{\{1\}}^{<N}(F;\, z) &= \sum_{|\alpha^{(1)}|+b|\alpha^{(2)}| < n_1/\sigma_A} f_{\{1\},\alpha}\frac{z^\alpha}{\alpha!}, \qquad \alpha \in \Z_{\geq 0}^{n}, \\
	T_{\{2\}}^{<N}(F;\, z) &= \sum_{c|\alpha^{(1)}|+|\alpha^{(2)}| < n_2/\sigma_A} f_{\{2\},\alpha}\frac{z^\alpha}{\alpha!}, \qquad \alpha \in \Z_{\geq 0}^{n}, \\	 
	T_{\{1,2\}}^{<N}(F;\, z) &= \sum_{\substack{|\alpha^{(1)}|+b|\alpha^{(2)}| < n_1/\sigma_A \\ c|\alpha^{(1)}|+|\alpha^{(2)}| < n_2/\sigma_A}} f_{\{1,2\},\alpha}\frac{z^\alpha}{\alpha!}, \qquad \alpha \in \Z_{\geq 0}^{n}, \\
	\operatorname{App}^{<N}(F;\, z) &=
T_{\{1\}}^{<N}(F;\, z) + T_{\{2\}}^{<N}(F;\, z) - T_{\{1,2\}}^{<N}(F;\, z).
\end{aligned}
$$
We say that $f$ is multi-asymptotically developable to $F=\{F_{\{1\}},F_{\{2\}},F_{\{1,2\}}\}$ along
$\chi$ on $S=S(W,\epsilon)$ if and only if for any cone
$S'=S(W',\epsilon')$ properly contained in $S$ and for any $N=(n_1,n_2)
\in \mathbb{Z}^2_{\ge 0}$, there exists a constant $C_{S',N}$ such that
$$
\left|
f(z) - \operatorname{App}^{<N}(F;\, z)
\right|
\le C_{S',N}(|z^{(1)}|^{n_1-bn_2}|z^{(2)}|^{n_2-cn_1})^{d/\sigma_A}
 \qquad (z \in S').
$$
The family $F$ is consistent if
\begin{itemize}
\item $f_{\{1\},\alpha}=f_{\{2\},\alpha}=f_{\{1,2\},\alpha}$ for each $\alpha \in \Z^n_{\geq 0}$.
\end{itemize}
%\end{itemize}

\

Let us consider the case $m=3$. We have 3 coordinates blocks
\begin{eqnarray*}
z^{(1)} & = & \{z_1,\dots,z_{m_1}\}, \\
z^{(2)} & = & \{z_{m_1+1},\dots,z_{m_1+m_2}\}, \\
z^{(3)} & = & \{z_{m_1+m_2+1},\dots,z_{m_1+m_2+m_3}\},
\end{eqnarray*}
with $m_1+m_2+m_3=n$. Up to row (resp. column) permutation and row (resp. column) multiplication by a positive rational number there are 5 interesting cases.
%\begin{itemize}
%%%%%%%%%%%%%%%%%%%
% N=3
%%%%%%%%%%%%%%%%%%%

\

\noindent
$\checkmark\,$ {\bf{Case $m=3,N=3$.}} In this case
$$
A_\chi=
\left(
\begin{array}{ccc}
1 & 0 & e \\
0 & 1 & 0
\end{array}
\right)
$$
with $1 \neq e \in \mathbb{Q}_{>0}$. We have $Z_1=\{z^{(1)}=z^{(3)}=0\}$ and $Z_2 = \{z^{(2)}=0\}$. Then we can define a  normal deformation $\widetilde{X}=\CC^n \times \R^2$ with the map $p(z;\,t)=(t_1z^{(1)}, t_2z^{(2)},t_1^ez^{(3)})$.
The rational monomials $\imin{\varphi_j}$, $j=1,2$ are $\imin{\varphi_1}(\tau)=\tau_1$, $\imin{\varphi_2}(\tau)=\tau_2$ and $\psi_3=\displaystyle{\tau_3 \over \tau_1^e}$.  A cofinal family of $C(\xi)$ is $C(\xi,F^0)$ consisting of elements $S(W,\epsilon)$, $W=W_1 \times W_2 \times W_3$, defined as follows
$$
S(W,\epsilon)=
\begin{aligned}
 \left\{z \in X;
\begin{aligned}
&z^{(k)} \in W_k \quad (k=1,2,3), \\
& |z^{(j)}| < \epsilon \quad (j=1,2), \\
& |z^{(3)}| < \epsilon |z^{(1)}|^e
\end{aligned}
\right\},
\end{aligned}
$$
where $\epsilon > 0$ and each $W_k$ is an $\mathbb{R}^+$-conic open subset in $\CC^{m_k}$ containing the point $\xi^{(k)}$. \\

We now construct the asymptotics.
%We have $K_{\{j\}}=K_j$, $j=1,2$, $K_{\{1,2\}}=K_1 \cup K_2$ and $Z_J=\{z^{(k)}=0,\,k \in K_J\}$ for $\emptyset \neq J \subseteq \{1,2\}$. For any subset $K \subset \{1,2\}$, we denote by $z^{(K)}$ the set of the coordinates $z^{(k)}$'s with $k \in K$ and by $z^{(K)}_{\mathrm{C}}$ the set of the coordinates which do not belong to $z^{(K)}$, i.e., $z^{(k)}$'s with $k \in \{0,1,\dots,m\} \setminus K$. Hence, the coordinates of $Z_J$ are given by $z^{(K_J)}_{\mathrm{C}}$. The $\pi_J$ denotes the canonical projection from $X$ to $Z_J$ defined by $z \to z^{(K_J)}_{\mathrm{C}}$. Given $S(W,\epsilon)$, set $S_J=\pi_J(S(W,\epsilon))$.
We have $S_{\{1,2\}}=\{0\}$ and
$$
S_{\{1\}}=
\begin{aligned}
 \left\{z \in Z_1;
\begin{aligned}
&z^{(2)} \in W_2, \\
& |z^{(2)}| < \epsilon
\end{aligned}
\right\}.
\end{aligned}
$$
$$
S_{\{2\}}=
\begin{aligned}
 \left\{z \in Z_2;
\begin{aligned}
&z^{(j)} \in W_j \quad (j=1,3), \\
& |z^{(1)}| < \epsilon, \\
& |z^{(3)}| < \epsilon |z^{(1)}|^e
\end{aligned}
\right\}.
\end{aligned}
$$
A total family of coefficients of multi-asymptotic expansion is given by
\begin{eqnarray*}
F & = & \{F_{\{1\}},F_{\{2\}},F_{\{1,2\}}\} \\
& = & \left\{\{f_{\{1\},\alpha}(z^{(2)})\}_{\alpha \in \Z_{\geq 0}^{m_1} \times \{0\}^{m_2} \times \Z_{\geq 0}^{m_3}},
\right. \\
&&\qquad \left.
\{f_{\{2\},\alpha}(z^{(1)},z^{(3)})\}_{\alpha \in \{0\}^{m_1} \times \Z_{\geq 0}^{m_2} \times \{0\}^{m_3}},\,\,
\{f_{\{1,2\},\alpha}\}_{\alpha \in \Z_{\geq 0}^n}\right\},
\end{eqnarray*}
where $f_{\{j\},\alpha}$ is holomorphic on $S_{\{j\}}$, $j=1,2$ and $f_{\{1,2\},\alpha} \in \CC$. Let $\sigma_A$ be the denominator of $e$. An asymptotic expansion $\operatorname{App}^{<N}(F;\,z)$, $N=(n_1,n_2) \in \Z_{\geq 0}^2$, is given by
$$
\begin{aligned}
	T_{\{1\}}^{<N}(F;\, z) &= \sum_{|\alpha^{(1)}|+e|\alpha^{(3)}| < n_1/\sigma_A} f_{\{1\},\alpha}(z^{(2)})\frac{z^\alpha}{\alpha!}, \quad \alpha \in \Z_{\geq 0}^{m_1} \times \{0\}^{m_2} \times \Z_{\geq 0}^{m_3}, \\
	T_{\{2\}}^{<N}(F;\, z) &= \sum_{|\alpha^{(2)}| < n_2/\sigma_A} f_{\{2\},\alpha}(z^{(1)},z^{(3)})\frac{z^\alpha}{\alpha!}, \quad \alpha \in \{0\}^{m_1} \times \Z_{\geq 0}^{m_2} \times \{0\}^{m_3}, \\	
	T_{\{1,2\}}^{<N}(F;\, z) &= \sum_{\substack{|\alpha^{(1)}|+e|\alpha^{(3)}| < n_1/\sigma_A \\ |\alpha^{(2)}| < n_2/\sigma_A}} f_{\{1,2\},\alpha}\frac{z^\alpha}{\alpha!}, \quad \alpha \in \Z_{\geq 0}^{n}, \\
	\operatorname{App}^{<N}(F;\, z) &=
T_{\{1\}}^{<N}(F;\, z) + T_{\{2\}}^{<N}(F;\, z) - T_{\{1,2\}}^{<N}(F;\, z).
\end{aligned}
$$
We say that $f$ is multi-asymptotically developable to $F=\{F_{\{1\}},F_{\{2\}},F_{\{1,2\}}\}$ along
$\chi$ on $S=S(W,\epsilon)$ if and only if for any cone
$S'=S(W',\epsilon')$ properly contained in $S$ and for any $N=(n_1,n_2)
\in \mathbb{Z}^2_{\ge 0}$, there exists a constant $C_{S',N}$ such that
$$
\left|
f(z) - \operatorname{App}^{<N}(F;\, z)
\right|
\le C_{S',N}(|z^{(1)}|^{n_1}|z^{(2)}|^{n_2})^{1/\sigma_A}
 \qquad (z \in S').
$$
The family $F$ is consistent if
\begin{itemize}
\item $f_{\{1\},(\alpha^{(1)},0,\alpha^{(3)})}(z^{(2)})$ is strongly asymptotically developable to $$\{f_{\{1,2\},(\alpha^{(1)},\alpha^{(2)},\alpha^{(3)})}\}_{\alpha^{(2)} \in \Z_{\geq 0}^{m_2}}$$ on $S_{\{1\}}$ for each $\alpha^{(i)} \in \Z_{\geq 0}^{m_i}$, $i=1,3$. Namely, for any cone
$S'_{\{1\}}$ properly contained in $S_{\{1\}}$ and for any $n_2
\in \mathbb{Z}_{\ge 0}$, there exists a constant $C_{S'_{\{1\}},n_2}$ such that
$$
\begin{aligned}
&\left|
f_{\{1\},(\alpha^{(1)},0,\alpha^{(3)})}(z^{(2)}) - \sum_{|\alpha^{(2)}|<n_2}f_{\{1,2\},(\alpha^{(1)},\alpha^{(2)},\alpha^{(3)})}{(z^{(2)})^{\alpha^{(2)}} \over \alpha^{(2)}!}
\right| \\
&\qquad\qquad\qquad\qquad\qquad\qquad\le C_{S'_{\{1\}},n_2}|z^{(2)}|^{n_2},
\end{aligned}
$$
\item $f_{\{2\},(0,\alpha^{(2)},0)}(z^{(1)},z^{(3)})$ is strongly asymptotically developable to $$\{f_{\{1,2\},(\alpha^{(1)},\alpha^{(2)},\alpha^{(3)})}\}_{\alpha^{(1)} \in \Z_{\geq 0}^{m_1} , \alpha^{(3)} \in \Z_{\geq 0}^{m_3}}$$ on $S_{\{2\}}$ for each $\alpha^{(2)} \in \Z_{\geq 0}^{m_2}$. Namely, for any cone
$S'_{\{2\}}$ properly contained in $S_{\{2\}}$ and for any $n_1
\in \mathbb{Z}_{\ge 0}$, there exists a constant $C_{S'_{\{2\}},n_1}$ such that
$$
\begin{aligned}
& \big|
f_{\{2\},(0,\alpha^{(2)},0)}(z^{(1)},z^{(3)}) \\
&\quad - \sum_{|\alpha^{(1)}|+e|\alpha^{(3)}|<n_1/\sigma_A}f_{\{1,2\},(\alpha^{(1)},\alpha^{(2)},\alpha^{(3)})}{(z^{(1)})^{\alpha^{(1)}}(z^{(3)})^{\alpha^{(3)}} \over \alpha^{(1)}!\alpha^{(3)}!}
\big| \\
& \qquad\qquad\qquad\qquad\le C_{S'_{\{2\}},n_1}|z^{(1)}|^{n_1/\sigma_A}.
\end{aligned}
$$
\end{itemize}
%%%%%%%%%%%%%%%%%%%%%%%%%%%%
% N=4
%%%%%%%%%%%%%%%%%%%%%%%%%%%%

\

\noindent
$\checkmark\,$ {\bf{Case $m=3,N=4$.}} In this case we have two possibilities
\begin{eqnarray*}
\mathrm{(a)} && A_\chi=
\left(
\begin{array}{ccc}
1 & b & e \\
0 & 1 & 0
\end{array}
\right)
\quad  b,e \in \mathbb{Q}_{>0},\, e \neq 1,
\\
\mathrm{(b)} && A_\chi=
\left(
\begin{array}{ccc}
1 & b & 0 \\
0 & 1 & e
\end{array}
\right)
\quad  b,e \in \mathbb{Q}_{>0}
\end{eqnarray*}

%%%%%%%%%%%%%%%%%%%%%%%%%%%%
% N=4 (a)
%%%%%%%%%%%%%%%%%%%%%%%%%%%%

\

\noindent
$\checkmark\,$ {\bf{Case $m=3,N=4$ (a).}}
%\begin{itemize}
%\item[(a)]
We have $Z_1=\{0\}$ and $Z_2 = \{z^{(2)}=0\}$. Then we can define a  normal deformation $\widetilde{X}=\CC^n \times \R^2$ with the map $p(z;\,t)=(t_1z^{(1)}, t_1^bt_2z^{(2)},t_1^ez^{(3)})$.
The rational monomials $\imin{\varphi_j}$, $j=1,2$ are $\imin{\varphi_1}(\tau)=\tau_1$, $\imin{\varphi_2}(\tau)=\displaystyle{\tau_2 \over \tau_1^b}$ and $\psi_3=\displaystyle{\tau_3 \over \tau_1^e}$.  A cofinal family of $C(\xi)$ is $C(\xi,F^0)$ consisting of elements $S(W,\epsilon)$, $W=W_1 \times W_2 \times W_3$, defined as follows
$$
S(W,\epsilon)=
\begin{aligned}
 \left\{z \in X;
\begin{aligned}
&z^{(k)} \in W_k \quad (k=1,2,3), \\
& |z^{(1)}| < \epsilon, \\
& |z^{(2)}| < \epsilon |z^{(1)}|^b, \\
& |z^{(3)}| < \epsilon |z^{(1)}|^e
\end{aligned}
\right\},
\end{aligned}
$$
where $\epsilon > 0$ and each $W_k$ is an $\mathbb{R}^+$-conic open subset in $\CC^{m_k}$ containing the point $\xi^{(k)}$. \\

We now construct the asymptotics.
%We have $K_{\{j\}}=K_j$, $j=1,2$, $K_{\{1,2\}}=K_1 \cup K_2$ and $Z_J=\{z^{(k)}=0,\,k \in K_J\}$ for $\emptyset \neq J \subseteq \{1,2\}$. For any subset $K \subset \{1,2\}$, we denote by $z^{(K)}$ the set of the coordinates $z^{(k)}$'s with $k \in K$ and by $z^{(K)}_{\mathrm{C}}$ the set of the coordinates which do not belong to $z^{(K)}$, i.e., $z^{(k)}$'s with $k \in \{0,1,\dots,m\} \setminus K$. Hence, the coordinates of $Z_J$ are given by $z^{(K_J)}_{\mathrm{C}}$. The $\pi_J$ denotes the canonical projection from $X$ to $Z_J$ defined by $z \to z^{(K_J)}_{\mathrm{C}}$. Given $S(W,\epsilon)$, set $S_J=\pi_J(S(W,\epsilon))$.
We have $S_{\{1\}}=S_{\{1,2\}}=\{0\}$ and
$$
S_{\{2\}}=
\begin{aligned}
 \left\{z \in Z_2;
\begin{aligned}
&z^{(j)} \in W_j \quad (j=1,3), \\
& |z^{(1)}| < \epsilon, \\
& |z^{(3)}| < \epsilon |z^{(1)}|^e
\end{aligned}
\right\}.
\end{aligned}
$$
A total family of coefficients of multi-asymptotic expansion is given by
\begin{eqnarray*}
F & = & \{F_{\{1\}},F_{\{2\}},F_{\{1,2\}}\} \\
& = & \left\{\{f_{\{1\},\alpha}\}_{\alpha \in \Z_{\geq 0}^{n}},
\{f_{\{2\},\alpha}(z^{(1)},z^{(3)})\}_{\alpha \in \{0\}^{m_1} \times \Z_{\geq 0}^{m_2} \times \{0\}^{m_3}},\right. \\
&&\qquad\qquad\qquad\qquad\qquad\qquad\qquad \left.\{f_{\{1,2\},\alpha}\}_{\alpha \in \Z_{\geq 0}^n}\right\},
\end{eqnarray*}
where $f_{\{2\},\alpha}$ is holomorphic on $S_{\{2\}}$ and $f_{\{1\},\alpha},f_{\{1,2\},\alpha} \in \CC$. Let $\sigma_A$ be the less common multiple of the denominators of $b,e$.  An asymptotic expansion $\operatorname{App}^{<N}(F;\,z)$, $N=(n_1,n_2) \in \Z_{\geq 0}^2$ is given by
$$
\begin{aligned}
	T_{\{1\}}^{<N}(F;\, z) &= \sum_{|\alpha^{(1)}|+b|\alpha^{(2)}|+e|\alpha^{(3)}| < n_1/\sigma_A} f_{\{1\},\alpha}\frac{z^\alpha}{\alpha!}, \qquad \alpha \in \Z_{\geq 0}^{n}, \\
	T_{\{2\}}^{<N}(F;\, z) &= \sum_{|\alpha^{(2)}| < n_2/\sigma_A} f_{\{2\},\alpha}(z^{(1)},z^{(3)})\frac{z^\alpha}{\alpha!}, \qquad \alpha \in \{0\}^{m_1} \times \Z_{\geq 0}^{m_2} \times \{0\}^{m_3}, \\		
T_{\{1,2\}}^{<N}(F;\, z) &= \sum_{\substack{|\alpha^{(1)}|+b|\alpha^{(2)}|+e|\alpha^{(3)}| < n_1/\sigma_A \\ |\alpha^{(2)}| < n_2/\sigma_A}} f_{\{1,2\},\alpha}\frac{z^\alpha}{\alpha!}, \qquad \alpha \in \Z_{\geq 0}^{n}, \\
	\operatorname{App}^{<N}(F;\, z) &=
T_{\{1\}}^{<N}(F;\, z) + T_{\{2\}}^{<N}(F;\, z) - T_{\{1,2\}}^{<N}(F;\, z).
\end{aligned}
$$
We say that $f$ is multi-asymptotically developable to $F=\{F_{\{1\}},F_{\{2\}},F_{\{1,2\}}\}$ along
$\chi$ on $S=S(W,\epsilon)$ if and only if for any cone
$S'=S(W',\epsilon')$ properly contained in $S$ and for any $N=(n_1,n_2)
\in \mathbb{Z}^2_{\ge 0}$, there exists a constant $C_{S',N}$ such that
$$
\left|
f(z) - \operatorname{App}^{<N}(F;\, z)
\right|
\le C_{S',N}(|z^{(1)}|^{n_1-bn_2}|z^{(2)}|^{n_2})^{1/\sigma_A}
 \qquad (z \in S').
$$
The consistent family in this case is described as follows:
\begin{itemize}
\item[$-$] $f_{\{1\},\alpha}=f_{\{1,2\},\alpha}$ for each $\alpha \in \Z_{\geq 0}^{n}$,
\item[$-$] $f_{\{2\},(0,\alpha^{(2)},0)}(z^{(1)},z^{(3)})$ is strongly asymptotically developable to $$\{f_{\{1,2\},(\alpha^{(1)},\alpha^{(2)},\alpha^{(3)})}\}_{\alpha^{(1)} \in \Z_{\geq 0}^{m_1} , \alpha^{(3)} \in \Z_{\geq 0}^{m_3}}$$ on $S_{\{2\}}$ for each $\alpha^{(2)} \in \Z_{\geq 0}^{m_2}$. Namely, for any cone
$S'_{\{2\}}$ properly contained in $S_{\{2\}}$ and for any $n_1
\in \mathbb{Z}_{\ge 0}$, there exists a constant $C_{S'_{\{2\}},n_1}$ such that
$$
\begin{aligned}
&\big|
f_{\{2\},(0,\alpha^{(2)},0)}(z^{(1)},z^{(3)})  \\
&\quad - \sum_{|\alpha^{(1)}|+e|\alpha^{(3)}|<n_1/\sigma_A}f_{\{1,2\},(\alpha^{(1)},\alpha^{(2)},\alpha^{(3)})}{(z^{(1)})^{\alpha^{(1)}}(z^{(3)})^{\alpha^{(3)}} \over \alpha^{(1)}!\alpha^{(3)}!}
\big| \\
&\qquad\qquad\qquad\qquad
\le C_{S'_{\{2\}},n_1}|z^{(1)}|^{n_1/\sigma_A}.
\end{aligned}
$$
\end{itemize}

%%%%%%%%%%%%%%%%%%%%%%%%%%%%%%
% N=4 (b)
%%%%%%%%%%%%%%%%%%%%%%%%%%%%%%

\

\noindent
$\checkmark\,${\bf{Case $m=3,N=4$ (b).}}
We have $Z_1=\{z^{(1)}=z^{(2)}=0\}$ and $Z_2 = \{z^{(2)}=z^{(3)}=0\}$. Then we can define a  normal deformation $\widetilde{X}=\CC^n \times \R^2$ with the map $p(z;\,t)=(t_1z^{(1)}, t_1^bt_2z^{(2)},t_2^e z^{(3)})$.
The rational monomials $\imin{\varphi_j}$, $j=1,2$ are $\imin{\varphi_1}(\tau)=\tau_1$, $\imin{\varphi_2}(\tau)=\displaystyle{\tau_2 \over \tau_1^b}$ and $\psi_3=\displaystyle{\tau_3\tau_1^{be} \over \tau_2^e}$.  A cofinal family of $C(\xi)$ is $C(\xi,F^0)$ consisting of elements $S(W,\epsilon)$, $W=W_1 \times W_2 \times W_3$, defined as follows
$$
S(W,\epsilon)=
\begin{aligned}
 \left\{z \in X;
\begin{aligned}
&z^{(k)} \in W_k \quad (k=1,2,3), \\
& |z^{(1)}| < \epsilon, \\
& |z^{(2)}| < \epsilon |z^{(1)}|^b, \\
& |z^{(3)}||z^{(1)}|^{be} < \epsilon |z^{(2)}|^e
\end{aligned}
\right\},
\end{aligned}
$$
where $\epsilon > 0$ and each $W_k$ is an $\mathbb{R}^+$-conic open subset in $\CC^{m_k}$ containing the point $\xi^{(k)}$. \\

We now construct the asymptotics.
%We have $K_{\{j\}}=K_j$, $j=1,2$, $K_{\{1,2\}}=K_1 \cup K_2$ and $Z_J=\{z^{(k)}=0,\,k \in K_J\}$ for $\emptyset \neq J \subseteq \{1,2\}$. For any subset $K \subset \{1,2\}$, we denote by $z^{(K)}$ the set of the coordinates $z^{(k)}$'s with $k \in K$ and by $z^{(K)}_{\mathrm{C}}$ the set of the coordinates which do not belong to $z^{(K)}$, i.e., $z^{(k)}$'s with $k \in \{0,1,\dots,m\} \setminus K$. Hence, the coordinates of $Z_J$ are given by $z^{(K_J)}_{\mathrm{C}}$. The $\pi_J$ denotes the canonical projection from $X$ to $Z_J$ defined by $z \to z^{(K_J)}_{\mathrm{C}}$. Given $S(W,\epsilon)$, set $S_J=\pi_J(S(W,\epsilon))$.
We have $S_{\{1,2\}}=\{0\}$ and
$$
S_{\{1\}}=
\begin{aligned}
 \left\{z \in Z_1;
\begin{aligned}
&z^{(3)} \in W_3, \\
& |z^{(3)}| < \epsilon
\end{aligned}
\right\},
\end{aligned}
$$
$$
S_{\{2\}}=
\begin{aligned}
 \left\{z \in Z_2;
\begin{aligned}
&z^{(1)} \in W_1, \\
& |z^{(1)}| < \epsilon \\
\end{aligned}
\right\}.
\end{aligned}
$$
 A total family of coefficients of multi-asymptotic expansion is given by
\begin{eqnarray*}
F & = & \{F_{\{1\}},F_{\{2\}},F_{\{1,2\}}\} \\
& = & \left\{\{f_{\{1\},\alpha}(z^{(3)})\}_{\alpha \in \Z_{\geq 0}^{m_1} \times \Z_{\geq 0}^{m_2} \times \{0\}^{m_3}}, \right.\\
& & \qquad\left.\{f_{\{2\},\alpha}(z^{(1)})\}_{\alpha \in \{0\}^{m_1} \times \Z_{\geq 0}^{m_2} \times \Z_{\geq 0}^{m_3}},
\{f_{\{1,2\},\alpha}\}_{\alpha \in \Z_{\geq 0}^n}\right\},
\end{eqnarray*}
where $f_{\{j\},\alpha}$ is holomorphic on $S_{\{j\}}$, $j=1,2$ and $f_{\{1,2\},\alpha} \in \CC$. Let $\sigma_A$ be the less common multiple of the denominators of $b,e$. An asymptotic expansion $\operatorname{App}^{<N}(F;\,z)$, $N=(n_1,n_2) \in \Z_{\geq 0}^2$ is given by
$$
\begin{aligned}
	T_{\{1\}}^{<N}(F;\, z) &= \sum_{|\alpha^{(1)}|+b|\alpha^{(2)}| < n_1/\sigma_A} f_{\{1\},\alpha}(z^{(3)})\frac{z^\alpha}{\alpha!}, \qquad \alpha \in \Z_{\geq 0}^{m_1} \times \Z_{\geq 0}^{m_2} \times \{0\}^{m_3}, \\
	T_{\{2\}}^{<N}(F;\, z) &= \sum_{|\alpha^{(2)}|+e|\alpha^{(3)}| < n_2/\sigma_A} f_{\{2\},\alpha}(z^{(1)})\frac{z^\alpha}{\alpha!}, \qquad \alpha \in \{0\}^{m_1} \times \Z_{\geq 0}^{m_2} \times \Z_{\geq 0}^{m_3}, \\		
T_{\{1,2\}}^{<N}(F;\, z) &= \sum_{\substack{|\alpha^{(1)}|+b|\alpha^{(2)}| < n_1/\sigma_A \\ |\alpha^{(2)}|+e|\alpha^{(3)}| < n_2/\sigma_A}} f_{\{1,2\},\alpha}\frac{z^\alpha}{\alpha!}, \qquad \alpha \in \Z_{\geq 0}^{n}, \\
	\operatorname{App}^{<N}(F;\, z) &=
T_{\{1\}}^{<N}(F;\, z) + T_{\{2\}}^{<N}(F;\, z) - T_{\{1,2\}}^{<N}(F;\, z).
\end{aligned}
$$
We say that $f$ is multi-asymptotically developable to $F=\{F_{\{1\}},F_{\{2\}},F_{\{1,2\}}\}$ along
$\chi$ on $S=S(W,\epsilon)$ if and only if for any cone
$S'=S(W',\epsilon')$ properly contained in $S$ and for any $N=(n_1,n_2)
\in \mathbb{Z}^2_{\ge 0}$, there exists a constant $C_{S',N}$ such that
$$
\left|
f(z) - \operatorname{App}^{<N}(F;\, z)
\right|
\le C_{S',N}(|z^{(1)}|^{n_1-bn_2}|z^{(2)}|^{n_2})^{1/\sigma_A}
 \qquad (z \in S').
$$
The family $F$ is consistent if
\begin{itemize}
\item $f_{\{1\},(\alpha^{(1)},\alpha^{(2)},0)}(z^{(3)})$ is strongly asymptotically developable to $$\{f_{\{1,2\},(\alpha^{(1)},\alpha^{(2)},\alpha^{(3)})}\}_{\alpha^{(3)} \in \Z_{\geq 0}^{m_3}}$$ on $S_{\{1\}}$ for each $\alpha^{(i)} \in \Z_{\geq 0}^{m_i}$, $i=1,2$. Namely, for any cone
$S'_{\{1\}}$ properly contained in $S_{\{1\}}$ and for any $n_3
\in \mathbb{Z}_{\ge 0}$, there exists a constant $C_{S'_{\{1\}},n_3}$ such that
$$
\begin{aligned}
&\big|
f_{\{1\},(\alpha^{(1)},\alpha^{(2)},0)}(z^{(3)})  \\
&\quad- \sum_{|\alpha^{(3)}|<n_3}f_{\{1,2\},(\alpha^{(1)},\alpha^{(2)},\alpha^{(3)})}{(z^{(3)})^{\alpha^{(3)}} \over \alpha^{(3)}!}
\big|
\le C_{S'_{\{1\}},n_3}|z^{(3)}|^{n_3}.
\end{aligned}
$$
\item $f_{\{2\},(0,\alpha^{(2)},\alpha^{(3)})}(z^{(1)})$ is strongly asymptotically developable to $$\{f_{\{1,2\},(\alpha^{(1)},\alpha^{(2)},\alpha^{(3)})}\}_{\alpha^{(1)} \in \Z_{\geq 0}^{m_1}}$$ on $S_{\{2\}}$ for each $\alpha^{(i)} \in \Z_{\geq 0}^{m_i}$, $i=2,3$. Namely, for any cone
$S'_{\{2\}}$ properly contained in $S_{\{2\}}$ and for any $n_1
\in \mathbb{Z}_{\ge 0}$, there exists a constant $C_{S'_{\{2\}},n_1}$ such that
$$
\begin{aligned}
&\big|
f_{\{2\},(0,\alpha^{(2)},\alpha^{(3)})}(z^{(1)})  \\
&\qquad - \sum_{|\alpha^{(1)}|<n_1}f_{\{1,2\},(\alpha^{(1)},\alpha^{(2)},\alpha^{(3)})}{(z^{(1)})^{\alpha^{(1)}} \over \alpha^{(1)}!}
\big|
\le C_{S'_{\{2\}},n_1}|z^{(1)}|^{n_1}.
\end{aligned}
$$
\end{itemize}
%\end{itemize}
%%%%%%%%%%%%%%%%%%%%%%%%%%%%
% N=5
%%%%%%%%%%%%%%%%%%%%%%%%%%%%

\

\noindent
$\checkmark\,$ {\bf{Case $m=3,N=5$.}} In this case
$$
A_\chi=
\left(
\begin{array}{ccc}
1 & b & e \\
0 & 1 & f
\end{array}
\right)
$$
with $b,e,f \in \mathbb{Q}_{>0}$, $bf \neq e$. We have $Z_1=\{0\}$ and $Z_2 = \{z^{(2)}=z^{(3)}=0\}$. Then we can define a  normal deformation $\widetilde{X}=\CC^n \times \R^2$ with the map $p(z;\,t)=(t_1z^{(1)}, t_1^bt_2z^{(2)},t_1^et_2^fz^{(3)})$.
The rational monomials $\imin{\varphi_j}$, $j=1,2$ are $\imin{\varphi_1}(\tau)=\tau_1$, $\imin{\varphi_2}(\tau)=\displaystyle{\tau_2 \over \tau_1^b}$ and $\psi_3=\displaystyle{\tau_3 \over \tau_1^{e-bf}\tau_2^f}$. A cofinal family of $C(\xi)$ is $C(\xi,F^0)$ consisting of elements $S(W,\epsilon)$, $W=W_1 \times W_2 \times W_3$, defined as follows
$$
S(W,\epsilon)=
\begin{aligned}
 \left\{z \in X;
\begin{aligned}
&z^{(k)} \in W_k \quad (k=1,2,3), \\
& |z^{(1)}| < \epsilon, \\
& |z^{(2)}| < \epsilon |z^{(1)}|^b, \\
& |z^{(3)}| < \epsilon |z^{(1)}|^{e-bf}|z^{(2)}|^f
\end{aligned}
\right\},
\end{aligned}
$$
where $\epsilon > 0$ and each $W_k$ is an $\mathbb{R}^+$-conic open subset in $\CC^{m_k}$ containing the point $\xi^{(k)}$. \\

We now construct the asymptotics.
%We have $K_{\{j\}}=K_j$, $j=1,2$, $K_{\{1,2\}}=K_1 \cup K_2$ and $Z_J=\{z^{(k)}=0,\,k \in K_J\}$ for $\emptyset \neq J \subseteq \{1,2\}$. For any subset $K \subset \{1,2\}$, we denote by $z^{(K)}$ the set of the coordinates $z^{(k)}$'s with $k \in K$ and by $z^{(K)}_{\mathrm{C}}$ the set of the coordinates which do not belong to $z^{(K)}$, i.e., $z^{(k)}$'s with $k \in \{0,1,\dots,m\} \setminus K$. Hence, the coordinates of $Z_J$ are given by $z^{(K_J)}_{\mathrm{C}}$. The $\pi_J$ denotes the canonical projection from $X$ to $Z_J$ defined by $z \to z^{(K_J)}_{\mathrm{C}}$. Given $S(W,\epsilon)$, set $S_J=\pi_J(S(W,\epsilon))$.
We have $S_{\{1\}}=S_{\{1,2\}}=\{\mathrm{pt}\}$ and
$$
S_{\{2\}}=
\begin{aligned}
 \left\{z \in Z_2;
\begin{aligned}
&z^{(1)} \in W_1, \\
& |z^{(1)}| < \epsilon \\
\end{aligned}
\right\}.
\end{aligned}
$$
A total family of coefficients of multi-asymptotic expansion is given by
\begin{eqnarray*}
F & = & \{F_{\{1\}},F_{\{2\}},F_{\{1,2\}}\} \\
& = & \left\{\{f_{\{1\},\alpha}\}_{\alpha \in \Z_{\geq 0}^{n}},
\{f_{\{2\},\alpha}(z^{(1)})\}_{\alpha \in \{0\}^{m_1} \times \Z_{\geq 0}^{m_2} \times \Z_{\geq 0}^{m_3}},
\{f_{\{1,2\},\alpha}\}_{\alpha \in \Z_{\geq 0}^n}\right\},
\end{eqnarray*}
where $f_{\{2\},\alpha}$ is holomorphic on $S_{\{2\}}$ and $f_{\{1\},\alpha},f_{\{1,2\},\alpha} \in \CC$.  Let $\sigma_A$ be the less common multiple of the denominators of $b,e,f$. An asymptotic expansion $\operatorname{App}^{<N}(F;\,z)$, $N=(n_1,n_2) \in \Z_{\geq 0}^2$ is given by
$$
\begin{aligned}
	T_{\{1\}}^{<N}(F;\, z) &= \sum_{|\alpha^{(1)}|+b|\alpha^{(2)}|+e|\alpha^{(3)}| < n_1/\sigma_A} f_{\{1\},\alpha}\frac{z^\alpha}{\alpha!}, \qquad \alpha \in \Z_{\geq 0}^{n}, \\
	T_{\{2\}}^{<N}(F;\, z) &= \sum_{|\alpha^{(2)}|+f|\alpha^{(3)}| < n_2/\sigma_A} f_{\{2\},\alpha}(z^{(1)})\frac{z^\alpha}{\alpha!}, \qquad \alpha \in \{0\}^{m_1} \times \Z_{\geq 0}^{m_2} \times \Z_{\geq 0}^{m_3}, \\		
T_{\{1,2\}}^{<N}(F;\, z) &= \sum_{\substack{|\alpha^{(1)}|+b|\alpha^{(2)}|+e|\alpha^{(3)}| < n_1/\sigma_A \\ |\alpha^{(2)}|+f|\alpha^{(3)}| < n_2/\sigma_A}} f_{\{1,2\},\alpha}\frac{z^\alpha}{\alpha!}, \qquad \alpha \in \Z_{\geq 0}^{n}, \\
	\operatorname{App}^{<N}(F;\, z) &=
T_{\{1\}}^{<N}(F;\, z) + T_{\{2\}}^{<N}(F;\, z) - T_{\{1,2\}}^{<N}(F;\, z).
\end{aligned}
$$
We say that $f$ is multi-asymptotically developable to $F=\{F_{\{1\}},F_{\{2\}},F_{\{1,2\}}\}$ along
$\chi$ on $S=S(W,\epsilon)$ if and only if for any cone
$S'=S(W',\epsilon')$ properly contained in $S$ and for any $N=(n_1,n_2)
\in \mathbb{Z}^2_{\ge 0}$, there exists a constant $C_{S',N}$ such that
$$
\left|
f(z) - \operatorname{App}^{<N}(F;\, z)
\right|
\le C_{S',N}(|z^{(1)}|^{n_1-bn_2}|z^{(2)}|^{n_2})^{1/\sigma_A}
 \qquad (z \in S').
$$
The family $F$ is consistent if
\begin{itemize}
\item$f_{\{1\},\alpha}=f_{\{1,2\},\alpha}$ for each $\alpha \in \Z_{\geq 0}^{n}$,
\item $f_{\{2\},(0,\alpha^{(2)},\alpha^{(3)})}(z^{(1)})$ is strongly asymptotically developable to $$\{f_{\{1,2\},(\alpha^{(1)},\alpha^{(2)},\alpha^{(3)})}\}_{\alpha^{(1)} \in \Z_{\geq 0}^{m_1}}$$ on $S_{\{2\}}$ for each $\alpha^{(i)} \in \Z_{\geq 0}^{m_i}$, $i=2,3$. Namely, for any cone
$S'_{\{2\}}$ properly contained in $S_{\{2\}}$ and for any $n_1
\in \mathbb{Z}_{\ge 0}$, there exists a constant $C_{S'_{\{2\}},n_1}$ such that
$$
\left|
f_{\{2\},(0,\alpha^{(2)},\alpha^{(3)})}(z^{(1)}) - \sum_{|\alpha^{(1)}|<n_1}f_{\{1,2\},(\alpha^{(1)},\alpha^{(2)},\alpha^{(3)})}{(z^{(1)})^{\alpha^{(1)}} \over \alpha^{(1)}!}
\right|
\le C_{S'_{\{2\}},n_1}|z^{(1)}|^{n_1}.
$$
\end{itemize}
%%%%%%%%%%%%%%%%%%%%%%%%%%%%
% N=6
%%%%%%%%%%%%%%%%%%%%%%%%%%%%

\

\noindent
$\checkmark\,${\bf{Case $m=3,N=6$.}} In this case
$$
A_\chi=
\left(
\begin{array}{ccc}
1 & b & e \\
c & 1 & r
\end{array}
\right)
$$
with $b,c,e,r \in \mathbb{Q}_{>0}$, $bc \neq 1$, $br \neq e$, $ce \neq r$ and $1-bc>0$. We have $Z_1=Z_2 = \{0\}$. Then we can define a  normal deformation $\widetilde{X}=\CC^n \times \R^2$ with the map $p(z;\,t)=(t_1t_2^cz^{(1)}, t_1^bt_2z^{(2)},t_1^et_2^rz^{(3)})$. Set $d=\displaystyle{1 \over 1-bc}$.
The rational monomials $\imin{\varphi_j}$, $j=1,2$ are $\imin{\varphi_1}(\tau)=\displaystyle{\left({\tau_1 \over \tau_2^c}\right)^d}$, $\imin{\varphi_2}(\tau)=\displaystyle{\left(\tau_2 \over \tau_1^b\right)^d}$ and $\psi_3=\displaystyle{\tau_3 \over (\tau_1^{e-br}\tau_2^{r-ce})^d}$. A cofinal family of $C(\xi)$ is $C(\xi,F^0)$ consisting of elements $S(W,\epsilon)$, $W=W_1 \times W_2 \times W_3$, defined as follows
$$
S(W,\epsilon)=
\begin{aligned}
 \left\{z \in X;
\begin{aligned}
&z^{(k)} \in W_k \quad (k=1,2,3), \\
& |z^{(1)}| < \epsilon^{1 \over d} |z^{(2)}|^c, \\
& |z^{(2)}| < \epsilon^{1 \over d} |z^{(1)}|^b, \\
& |z^{(3)}| < \epsilon (|z^{(1)}|^{e-br}|z^{(2)}|^{r-ce})^d
\end{aligned}
\right\},
\end{aligned}
$$
where $\epsilon > 0$ and each $W_k$ is an $\mathbb{R}^+$-conic open subset in $\CC^{m_k}$ containing the point $\xi^{(k)}$. \\

We now construct the asymptotics.
%We have $K_{\{j\}}=K_j$, $j=1,2$, $K_{\{1,2\}}=K_1 \cup K_2$ and $Z_J=\{z^{(k)}=0,\,k \in K_J\}$ for $\emptyset \neq J \subseteq \{1,2\}$. For any subset $K \subset \{1,2\}$, we denote by $z^{(K)}$ the set of the coordinates $z^{(k)}$'s with $k \in K$ and by $z^{(K)}_{\mathrm{C}}$ the set of the coordinates which do not belong to $z^{(K)}$, i.e., $z^{(k)}$'s with $k \in \{0,1,\dots,m\} \setminus K$. Hence, the coordinates of $Z_J$ are given by $z^{(K_J)}_{\mathrm{C}}$. The $\pi_J$ denotes the canonical projection from $X$ to $Z_J$ defined by $z \to z^{(K_J)}_{\mathrm{C}}$. Given $S(W,\epsilon)$, set $S_J=\pi_J(S(W,\epsilon))$.
We have $S_{\{1\}}=S_{\{2\}}=S_{\{1,2\}}=\{0\}$ and
%$$
%S_{\{2\}}=
%\begin{aligned}
% \left\{z \in Z_2;
%\begin{aligned}
%&z^{(1)} \in W_j, \\
%& |z^{(1)}| < \epsilon \\
%\end{aligned}
%\right\}.
%\end{aligned}
%$$
a total family of coefficients of multi-asymptotic expansion is given by
\begin{eqnarray*}
F & = & \{F_{\{1\}},F_{\{2\}},F_{\{1,2\}}\} \\
& = & \left\{\{f_{\{1\},\alpha}\}_{\alpha \in \Z_{\geq 0}^{n}},
\{f_{\{2\},\alpha}\}_{\alpha \in
\Z_{\geq 0}^{n}},
\{f_{\{1,2\},\alpha}\}_{\alpha \in \Z_{\geq 0}^n}\right\},
\end{eqnarray*}
where $f_{\{1\},\alpha},f_{\{2\},\alpha},f_{\{1,2\},\alpha} \in \CC$.  Let $\sigma_A$ be the less common multiple of the denominators of $b,c,e,r$. An asymptotic expansion $\operatorname{App}^{<N}(F;\,z)$, $N=(n_1,n_2) \in \Z_{\geq 0}^2$ is given by
$$
\begin{aligned}
	T_{\{1\}}^{<N}(F;\, z) &= \sum_{|\alpha^{(1)}|+b|\alpha^{(2)}|+e|\alpha^{(3)}| < n_1/\sigma_A} f_{\{1\},\alpha}\frac{z^\alpha}{\alpha!}, \qquad \alpha \in \Z_{\geq 0}^{n}, \\
	T_{\{2\}}^{<N}(F;\, z) &= \sum_{c|\alpha^{(1)}|+|\alpha^{(2)}|+r|\alpha^{(3)}| < n_2/\sigma_A} f_{\{2\},\alpha}(z^{(1)})\frac{z^\alpha}{\alpha!}, \qquad \alpha \in \Z_{\geq 0}^{n}, \\		
T_{\{1,2\}}^{<N}(F;\, z) &= \sum_{\substack{|\alpha^{(1)}|+b|\alpha^{(2)}|+e|\alpha^{(3)}| < n_1/\sigma_A \\ c|\alpha^{(1)}|+|\alpha^{(2)}|+r|\alpha^{(3)}| < n_2/\sigma_A}} f_{\{1,2\},\alpha}\frac{z^\alpha}{\alpha!}, \qquad \alpha \in \Z_{\geq 0}^{n}, \\
	\operatorname{App}^{<N}(F;\, z) &=
T_{\{1\}}^{<N}(F;\, z) + T_{\{2\}}^{<N}(F;\, z) - T_{\{1,2\}}^{<N}(F;\, z).
\end{aligned}
$$
We say that $f$ is multi-asymptotically developable to $F=\{F_{\{1\}},F_{\{2\}},F_{\{1,2\}}\}$ along
$\chi$ on $S=S(W,\epsilon)$ if and only if for any cone
$S'=S(W',\epsilon')$ properly contained in $S$ and for any $N=(n_1,n_2)
\in \mathbb{Z}^2_{\ge 0}$, there exists a constant $C_{S',N}$ such that
$$
\left|
f(z) - \operatorname{App}^{<N}(F;\, z)
\right|
\le C_{S',N}(|z^{(1)}|^{n_1-bn_2}|z^{(2)}|^{n_2-cn_1})^{d/\sigma_A}
 \qquad (z \in S').
$$
The family $F$ is consistent if
\begin{itemize}
\item$f_{\{1\},\alpha}=f_{\{2\},\alpha}=f_{\{1,2\},\alpha}$ for each $\alpha \in \Z_{\geq 0}^{n}$.
\end{itemize}
%\end{itemize}

\

In the other cases one of the following occurs:
\begin{itemize}
\item the action is degenerate,
\item we reduce to $m=2$.
\end{itemize}

\appendix
\begin{section}{Appendix}

In the Appendix we define conic objects and conic sheaves associated to multiple actions of $\RP$. The proofs of the statements of Appendix \ref{A} and \ref{B} can be found in \cite{HP13}. In Appendix \ref{C} we prove a decomposition theorem for open subanalytic subsets which is crucial for the computation of the sections of the multi-specialization.

%\end{section}
\begin{subsection}{Conic sheaves}\label{A}

Let $k$ be a field. Let $X$ be a real analytic manifold endowed
with a subanalytic action $\mu$ of $\RP$. In other words we have a
subanalytic map
$$\mu: X \times \RP \to X,$$
which satisfies, for each $t_1,t_2 \in \RP$:
$$
  \begin{cases}
    \mu(x,t_1t_2)=\mu(\mu(x,t_1),t_2), \\
    \mu(x,1)=x.
  \end{cases}
$$
Note that $\mu$ is open, in fact let $U \in \op(X)$ and $(t_1,t_2)
\in \op(\RP)$. Then $\mu(U,(t_1,t_2))=\bigcup_{t \in
(t_1,t_2)}\mu(U,t)$, and $\mu(\cdot,t):X \to X$ is a homeomorphism
(with inverse $\mu(\cdot,\imin t)$). We have a diagram
$$\xymatrix{X \ar[r]^{\hspace{-5mm}j} & X \times \RP \ar@ <2pt>
[r]^{\hspace{0.3cm}\mu} \ar@ <-2pt> [r]_{\hspace{0.3cm}p}& X,}$$
where $j(x)=(x,1)$ and $p$ denotes the projection. We have $\mu
\circ j=p \circ j=\id$.

\begin{df} (i) Let $S$ be a subset of $X$. We set $\RP S=\mu(S,\RP).$
If $U \in \op(X)$, then $\RP U \in \op(X)$ since $\mu$ is open.

(ii) Let $S$ be a subset of $X$. We say that $S$ is conic if
$S=\RP S$. In other words, $S$ is invariant by the action of
$\mu$.

(iii) An  orbit of $\mu$ is the set $\RP x$ with $x \in X$.
\end{df}

We assume that the orbits of $\mu$ are contractible. For each $x \in X$ there are two possibilities: either $\RP x=x$ or $\RP x \simeq \R$.

\begin{df} We say that a subset $S$ of $X$ is $\RP$-connected if
$S \cap \RP x$ is connected for each $x \in S$.
\end{df}

\begin{lem} \label{lem: RP-connectedness} (i) Let $S_1,S_2 \subset X$ and suppose that $S_2$ is conic. Then $\RP(S_1 \cap S_2)=\RP S_1 \cap S_2$. (ii) If $S_1$ and $S_2$ are $\RP$-connected, then $S_1 \cap S_2$ is $\RP$-connected.
\end{lem}

Let $X,Y$ topological spaces endowed with an action ($\mu_{X}$ and $\mu_{Y}$ respectively) of $\RP$.

\begin{df} A continuous function $f:X \to Y$ is said to be conic if for each $x \in X$, $a \in \RP$ we have
$f(\mu_{X}(x,a))=\mu_{Y}(f(x),a)$.
\end{df}

\begin{lem} \label{lem: conic maps} Let $f:X \to Y$ be a conic map. (i) Suppose that $S \subset Y$ is $\RP$-connected (resp. conic). Then $\imin f(S)$ is $\RP$-connected (resp. conic). (ii) Suppose that $Z \subset X$ is conic. Then $f(Z)$ is conic.
\end{lem}

Let $X$ be a real analytic manifold endowed with a subanalytic action of $\RP$. Denote by $X_{sa}$ the associated subanalytic site.

\begin{df} A sheaf of $k$-modules $F$ on $X_{sa}$ is conic if the restriction morphism
$\Gamma(\RP U;F) \to \Gamma(U;F)$ is an isomorphism for each
$\RP$-connected $U \in \op^c(X_{sa})$ with $\RP U \in
\op(X_{sa})$.
\begin{itemize}
\item[(i)]We denote by $\mod_{\RP}(k_{X_{sa}})$ the subcategory of
$\mod(k_{X_{sa}})$ consisting of conic sheaves.
\item[(ii)] We denote by $D^b_{\RP}(k_{X_{sa}})$,
%where $\sharp=+,b$
the subcategory of $D^b(k_{X_{sa}})$ consisting of objects $F$
such that $H^j(F)$ belongs to $\mod_{\RP}(k_{X_{sa}})$ for all $j
\in \Z$.
\end{itemize}
\end{df}

Assume the hypothesis below:
\begin{equation}\label{hypsa}
  \begin{cases}
 \text{(i) every $U \in \op^c(X_{sa})$ has a finite covering consisting }\\
 \text{\ \ of $\RP$-connected subanalytic open subsets,}\\
 \text{(ii) for any $U \in \op^c(X_{sa})$ we have $\RP U \in \op(X_{sa})$,}\\
 \text{(iii) for any $x \in X$\ \ the set $\RP x$ is contractible,}\\
 \text{(iv) there exists a covering $\{V_n\}_{n \in \N}$ of $X_{sa}$ such that}\\
 \text{\ \ $V_n$ is $\RP$-connected and $V_n \subset\subset V_{n+1}$ for each $n$}.
  \end{cases}
\end{equation}

The following result was proven in \cite{Pr11}.

\begin{prop}\label{RPU} Assume \eqref{hypsa}. Let $U \in \op(X_{sa})$ be $\RP$-connected and such that $\RP U \in \op(X_{sa})$. Let $F \in D^b_{\RP}(k_{X_{sa}})$. Then
$$
R\Gamma(\RP U;F) \iso R\Gamma(U;F).
$$
\end{prop}

\end{subsection}

\begin{subsection}{Multi-conic sheaves}\label{B}

Let $X$ be a topological space
with $\ell$ actions $\{\mu_i\}_{i=1}^\ell$ of $\RP$ such that $\mu_i(\mu_j(x,t_j),t_i)=\mu_j(\mu_i(x,t_i),t_j)$. We have a map
\begin{eqnarray*}
\mu: X \times (\RP)^\ell & \to & X \\
(x,(t_1,\dots,t_\ell)) & \mapsto & \mu_1(\cdots\mu_\ell(x,t_\ell),\dots,t_1).
\end{eqnarray*}

\begin{df} (i) Let $S$ be a subset of $X$. We set $\RP_i S=\mu_i(S,\RP).$
If $U \in \op(X)$, then $\RP_i U \in \op(X)$ since $\mu_i$ is open for each $i=1,\ldots,\ell$.

(ii)  Let $S$ be a subset of $X$. Let $J = \{i_1,\ldots,i_k\} \subseteq \{1,\ldots,\ell\}$. We set $$\RP_JS=\RP_{i_1}\cdots\RP_{i_k}S=\mu_{i_1}(\cdots\mu_{i_k}(S,\RP),\dots,\RP), \ \ i_1,\ldots i_k \in J.$$
We set $(\RP)^\ell S=\RP_{\{1,\ldots,\ell\}}S=\mu(S,(\RP)^\ell).$
If $U \in \op(X)$, then $\RP_J U \in \op(X)$ since $\mu_i$ is open for each $i \in \{1,\ldots,\ell\}$.

(iii) Let $S$ be a subset of $X$. We say that $S$ is $(\RP)^\ell$-conic if
$S=(\RP)^\ell S$. In other words, $S$ is invariant by the action of
$\mu_i$, $i=1,\dots,\ell$.

\end{df}

\begin{df} (i)  We say that a subset $S$ of $X$ is $\RP_i$-connected if
$S \cap \RP_i x$ is connected for each $x \in S$.

(ii) We say that a subset $S$ of $X$ is $(\RP)^\ell$-connected if
there exists a permutation $\sigma:\{1,\ldots,\ell\}\to\{1,\ldots,\ell\}$ such that
\begin{equation}\label{RPl connected}
 \begin{cases}
 \text{$S$ is $\RP_{\sigma(1)}$-connected,}\\
 \text{$\RP_{\sigma(1)}S$ is $\RP_{\sigma(2)}$-connected,}\\
 \text{\ \ $\vdots$}\\
 \text{$\RP_{\sigma(1)}\cdots\RP_{\sigma(\ell-1)}S$ is $\RP_{\sigma(\ell)}$-connected}.
  \end{cases}
\end{equation}
\end{df}

The following results follow from the case $\ell=1$.

\begin{lem} \label{lem: l-connectedness} (i) Let $S_1,S_2 \subset X$ and suppose that $S_2$ is $(\RP)^\ell$-conic. Then $(\RP)^\ell(S_1 \cap S_2)=(\RP)^\ell S_1 \cap S_2$. (ii) If moreover $S_1$ is $(\RP)^\ell$-connected then $S_1 \cap S_2$ is $(\RP)^\ell$-connected.
\end{lem}

\begin{oss} In (ii) of Lemma \ref{lem: l-connectedness} we have to assume that $S_2$ is $(\RP)^\ell$-conic. Indeed it is not true that the intersection of two $(\RP)^\ell$-connected is $(\RP)^\ell$-connected in general.
\end{oss}

Let $X,Y$ topological spaces endowed with $\ell$ actions $\{\mu_{Xi}\}_{i=1}^\ell$, $\{\mu_{Yi}\}_{i=1}^\ell$ of $\RP$.

\begin{df} A continuous function $f:X \to Y$ is said to be $(\RP)^\ell$-conic if for each $x \in X$, $a \in \RP$ we have
$f(\mu_{Xi}(x,a))=\mu_{Yi}(f(x),a)$, $i=1,\ldots,\ell$.
\end{df}

\begin{lem} \label{lem: l-conic maps} Let $f:X \to Y$ be a $(\RP)^\ell$-conic map. (i) Suppose that $S \subset Y$ is $(\RP)^\ell$-connected (resp. $(\RP)^\ell$-conic). Then $\imin f(S)$ is $(\RP)^\ell$-connected (resp. $(\RP)^\ell$-conic). (ii) Suppose that $Z \subset X$ is $(\RP)^\ell$-conic. Then $f(Z)$ is $(\RP)^\ell$-conic.
\end{lem}

Let $X$ be a real analytic manifold and denote by $X_{sa}$ the associated subanalytic site. Assume that $X$ is endowed with $\ell$ subanalytic $\RP$-actions $\mu_1,\dots,\mu_\ell$ commuting with each other.

\begin{df}\label{def:conic sheaf} A sheaf of $k$-modules $F$ on $X_{sa}$ is $(\RP)^\ell$-conic if it is conic with respect to each $\mu_i$.
%This implies that the restriction morphism $\Gamma((\RP)^\ell U;F) \to \Gamma(U;F)$ is an isomorphism for each $(\RP)^\ell$-connected $U \in \op^c(X_{sa})$ with $(\RP U)^\ell \in \op(X_{sa})$.
\begin{itemize}
\item[(i)]We denote by $\mod_{(\RP)^\ell}(k_{X_{sa}})$ the subcategory of
$\mod(k_{X_{sa}})$ consisting of $(\RP)^\ell$-conic sheaves.
\item[(ii)] We denote by $D^b_{(\RP)^\ell}(k_{X_{sa}})$,
%where $\sharp=+,b$
the subcategory of $D^b(k_{X_{sa}})$ consisting of objects $F$
such that $H^j(F)$ belongs to $\mod_{(\RP)^\ell}(k_{X_{sa}})$ for all $j
\in \Z$.
\end{itemize}
\end{df}

Let us assume the following hypothesis
\begin{equation}\label{hypsam}
  \begin{cases}
 \text{(i) the pair $(X,\mu_i)$ satisfies \eqref{hypsa} for each $i=1,\dots,\ell$, }\\
 %\text{(i) every $U \in \op^c(X_{sa})$ has a finite covering consisting }\\
 %\text{\ \ of $\RP_i$-connected subanalytic open subsets, $i=1,\dots,\ell$,}\\
 \text{(ii) every $U \in \op^c(X_{sa})$ has a finite covering consisting }\\
 \text{\ \ of $(\RP)^\ell$-connected subanalytic open subsets,}\\
 \text{(iii) we have $\RP_J U \in \op(X_{sa})$ for any $U \in \op^c(X_{sa})$}\\
 \text{\ \ and any $J \subset \{1,\ldots,\ell\}$.}\\
 %\text{(iv) $\RP_i x$ is contractible for any $x \in X$ and any $i=1,\ldots,\ell$,}\\
 %\text{(v) there exists a covering $\{V_n\}_{n \in \N}$ of $X_{sa}$ such that}\\
 %\text{\ \ $V_n$ is $(\RP)^\ell$-connected and $V_n \subset\subset V_{n+1}$ for each $n$}.
  \end{cases}
\end{equation}
In this situation the orbits of $\mu_i$, $i=1,\dots,\ell$ are either $\RP x \simeq \R$ or $\RP x = x$.\\

\begin{prop}\label{prop: RPUl} Assume \eqref{hypsam}. Let $U \in \op(X_{sa})$ be $(\RP)^\ell$-connected. Let $F \in D^b_{(\RP)^\ell}(k_{X_{sa}})$. Then
$$R\Gamma((\RP)^\ell U;F) \iso R\Gamma(U;F).$$
\end{prop}

If $X$ satisfies \eqref{hypsam} (i)-(iii), then it follows from Proposition \ref{prop: RPUl} that for $(\RP)^\ell$-conic subanalytic sheaves it is enough to study the cohomology of the sections on $(\RP)^\ell$-conic open subsets.\\

\end{subsection}

\subsection{A decomposition theorem for subanalytic open sets}\label{C}

Let $X=\R^n$ with coordinates $x=(x_1,\dots,x_n)$ endowed with the actions $\mu_j:X \times \RP \to X$, $j=1,\dots,\ell$ defined by $\mu_j(x,\lambda)=\mu_j(x_1,\dots,x_n,\lambda)=(\lambda^{\alpha_{j1}}x_1,\dots,\lambda^{\alpha_{jn}},x_n)$, with $\lambda \geq 0$ and $\alpha_{ji} \in \R$, $\alpha_{ji} \neq 0$ if $i \in I_j \subseteq \{1,\dots,n\}$. We can represent the actions $\mu_j$, $j=1,\dots,\ell$ in a matrix form as follows: let $x=(x_1,\dots,x_n)$, $\lambda=(\lambda_1,\dots,\lambda_\ell)$, $A=(\alpha_{ji}) \in M_{\ell,n}(\R)$. Then
\begin{equation}\label{matrix}
\mu(x,\lambda)=\mu_1(\cdots(\mu_\ell(x,\lambda_\ell)\cdots),\lambda_1)=x e^{\log\lambda A}
\end{equation}
where $\log\lambda=(\log\lambda_1,\dots,\log\lambda_\ell)$ and, given $y=(y_1,\dots,y_n)$, $e^y=\delta_{ij}e^{y_i} \in M_{n,n}(\R)$ ($\delta_{ij}$ is the Kronecker's delta) is a diagonal matrix. So the $i$-th diagonal entry of $e^{\log\lambda A}$ is $\prod_{j=1}^\ell\lambda^{\alpha_{ji}}$.

We consider the equivalence classes $\hat{I}_k$ defined as follows: $k_1,k_2 \in \hat{I}_k$ iff $\pi_{k_1}(\mu_j(x,\lambda))=\pi_{k_2}(\mu_j(x,\lambda))$ for any $j \in 1,\dots,\ell$ and any $\lambda \geq 0$. Here $\pi_{k_i}$ denotes the projection on the $k_i$-th coordinate, $i=1,2$. That means $\alpha_{jk_1}=\alpha_{jk_2}$ for each $j=1,\dots,\ell$.

So we may assume that $\{1,\dots,n\}=\bigsqcup_{k=1}^N \hat{I}_k$. Let $m_k=\sharp \hat{I}_k$. Then $\sum_{k=1}^Nm_k=n$. Consider the morphism of manifolds
\begin{eqnarray}
\varphi : \widetilde{X}'=\times_{k=1}^N\mathbb{S}^{m_k-1} \times \R^N  & \to & \R^{n} \label{varphi} \\
(\vartheta_k,r_k)_{k=1}^N & \mapsto & (r_k \iota_k(\vartheta_k))_{k=1}^N, \nonumber
\end{eqnarray}
where $\iota_k:\mathbb{S}^{m_k-1} \hookrightarrow \R^{m_k}$ denotes the
embedding. Endow $\times_{k=1}^N\mathbb{S}^{m_k-1} \times \R^N$ with the actions $\widetilde{\mu}'_j=\mu_j$, $j=1,\dots,\ell-1$, $\widetilde{\mu}'_j((\vartheta_k,r_k)_{k=1}^N,\lambda)=(\vartheta_k,\lambda^{\alpha_{jk}}r_k)_{k=1}^N$. Then $\varphi$ is a $(\RP)^\ell$-conic map.\\

Let $V$ be an $(\RP)^\ell$-conic subanalytic subset in $X$.
Let $S_k=\{x_i=0;\,i \in \hat{I}_k\}$, and let $\pi_j:X \to S_k$ be the projection. We introduce the following conditions Va.~and Vb.~of $V$ for each $k=1,\dots,N$.
\begin{enumerate}
\item[{\bf{Va.}}] $V$ does not intersect $S_k$.
\item[{\bf{Vb.}}] $\pi_j(V) \subset \pi_j(V \cap S_k)$.
\end{enumerate}

\begin{prop}
Let $V$ be a $(\RP)^\ell$-conic subanalytic subset in $X$ and
let $W$ be an open subanalytic neighborhood of $V$.
If $V$ satisfies the condition either Va.~or Vb.~for each $j$,
then there exists a subanalytic subset $W'$ in $X$
satisfying the following conditions.
\begin{enumerate}
\item  $W'$ is an open neighborhood of $V$ and contained in $W$.
\item  $W'$ is $(\RP)^\ell$-connected.
\item  $\RP_1\dots\RP_k W'$
is also subanalytic in $X$ for any $1 \le k \le \ell$.
\end{enumerate}
\end{prop}
\begin{proof}

Thanks to the morphism of manifolds \eqref{varphi} we may reduce to $X=[0,+\infty]^N$ setting $\widetilde{\mu}_j'=\mu_j$ for short and omitting the variables $\vartheta_j$ (they are irrelevant being fixed by each $\mu_j$, $j=1,\dots,\ell$). We argue by induction on the number of actions. We prove the assertion in several steps.\\

(a) Suppose that there exists $i \in I_\ell$ such that $V$ satisfies Va. on $i$. Up to a permutation of coordinates, we may assume $i=N$.
%Moreover we can also assume that there exists $j \in \{1,\dots,\ell\}$ such that $\mu_j$ acts on the $N$-th coordinate.
First, on $\{x_N \neq 0\}$ consider the homeomorphism $\psi$ making the orbits of $\mu_\ell$ orthogonal to $\{x_N=1\}$. Namely, if $\mu_\ell(x_1,\dots,x_{N-1},x_N,\lambda)=(\lambda^{\alpha_{\ell 1}}x_1,\dots,\lambda^{\alpha_{\ell N-1}}x_{N-1},\lambda^{\alpha_{\ell N}}x_N)$ we set $\psi(x_1,\dots,x_{N-1},x_N)=(x_1x_N^{-\alpha_{\ell 1}/\alpha_{\ell N}},\dots,x_{N-1}x_N^{-\alpha_{\ell N-1}/\alpha_{\ell N}},x_N^{1/\alpha_{\ell N}})$. Then $\psi \circ \mu_\ell \circ (\imin\psi \times \id_{\RP})$ acts only on the variable $x_N$. This corresponds to the change of coordinates sending the matrix (representing the action as in \eqref{matrix})
$$
\left(
\begin{array}{cc}
A & B \\
C & d
\end{array}
\right) \ \ \ \ \mbox{to} \ \ \ \
\left(
\begin{array}{cc}
A' & B' \\
0 & 1
\end{array}
\right)
$$
where $A$ (resp. $A'$) is a $(N-1) \times (N-1)$ matrix, $B$ (resp. $B'$) is a $(N-1) \times 1$ matrix, $C$ (resp. $0$) is a $1 \times (N-1)$ matrix and $d \neq 0$.

Replace $V$ with $\psi(V)$ and $\mu_j$ with $\psi \circ \mu_j \circ (\imin\psi \times \id_{\RP})$, $j=1,\dots,\ell$. One can check easily that $V \cap \{x_N=1\}$ (remark that $\psi(\{x_N=1\})=\{x_N=1\}$) is conic with respect to the actions $\widetilde{\mu}_j$ defined as $\widetilde{\mu}_j(x_1,\dots,x_{N-1},1,\lambda)=(\lambda^{\alpha_{j1}}x_1,\dots,\lambda^{\alpha_{jN-1}}x_{N-1},1)$. By the induction hypothesis $W \cap \{x_N=1\}$ contains a subanalytic open neighborhood $W_1$ of $V \cap \{x_N=1\}$ in $\{x_N=1\}$ which is $(\widetilde{\RP})^{\ell-1}$-connected (here $\widetilde{\RP}$ means with respect to $\widetilde{\mu}_j$, $j=1,\dots,\ell-1$).

Let $\pi:(x_1,\dots,x_{N-1},x_N) \mapsto (x_1,\dots,x_{N-1})$ be the projection. Then $W \cap \imin\pi(\pi(p))$ is a disjoint union of intervals. Let us consider the interval $(m(p),M(p))$ containing $x_N=1$. Set
$$
\widetilde{W}=\{(x_1,\dots,x_{N-1},x_N);\;(x_1,\dots,x_{N-1}) \in \pi(W_1),\; m(p)<x_N<M(p)\}.
$$
By construction $\widetilde{W}$ is subanalytic and $\RP_\ell$-connected. Let us prove that it is $(\RP)^\ell$-connected. Remark that $\RP_j\RP_\ell S = \RP_\ell\widetilde{\RP_j} S$ for each $j=1,\dots,\ell-1$ and $S \subseteq \{x_N=1\}$ (here $\widetilde{\RP_j}$ means conic with respect to $\widetilde{\mu}_j$). So by the induction hypothesis we are reduced to prove that $\RP_\ell S$ is $\RP_j$-connected if $S$ is $\widetilde{\RP_j}$-connected, $j=1,\dots,\ell-1$.

Let $x=(x_1,\dots,x_{N-1},x_N), \mu_j(x,a) \in \RP_\ell S$. We shall prove that $\mu_j(x,b) \in \RP_\ell S$ if $b \in [1,a]$ (we assume without loss of generality that $a>1$). We have $(x_1,\dots,x_{N-1},1)=\mu_\ell(x,\imin{x_N}),\widetilde{\mu}_j(x_1,\dots,x_{N-1},1,a)=\mu_\ell(\mu_j(x,a),\imin{(a^{\alpha_{jN}}x_N)}) \in S$. Since $S$ is $\widetilde{\RP_j}$-connected we have $\widetilde{\mu}_j(x_1,\dots,x_{N-1},1,b) \in S$ and hence $\mu_j(x,b)=\mu_\ell(\widetilde{\mu}_j(x_1,\dots,x_{N-1},1,b),b^{\alpha_{jN}}x_N) \in \RP_\ell S$.\\

%(a1) Suppose that $\mu_\ell$ does not act on the $N$-th variable. On $\{x_N=1\}$ $V \cap \{x_N=1\}$ is conic with respect to the actions $\mu_j$, $j \in J_1 \subseteq \{1,\dots,\ell\}$ such that $\mu_j$ does not act on the $N$-th coordinate. In particular it is $\RP_\ell$-conic. By the induction hypothesis $W \cap \{x_N=1\}$ contains a subanalytic open neighborhood $W_1$ of $V \cap \{x_N=1\}$ in $\{x_N=1\}$ which is $\RP_j$-conic, $j \in J_1 \setminus \{\ell\}$ and $\RP_\ell$-connected. The set $\RP_1\cdots\RP_{\ell-1}W_1=\RP_{j_1}\cdots\RP_{j_p}W_1$, $j_1,\dots,j_p \in \{1,\dots,N\} \setminus J_1$ is a neighborhood of $V$. Let us prove that is is $\RP_\ell$-connected.

(b) Suppose that $V$ satisfies Vb. for each $i \in I_\ell$. Up to shrink $W$, we may assume that $W$ satisfies Vb. for each $i \in I_\ell$ as well. Let $J^+=\{i \in I_\ell,\, \alpha_{\ell i}>0\}$, $J^-=\{i \in I_\ell,\, \alpha_{\ell i}<0\}$.\\

We first shrink $W$ on $X \setminus (\{x_i=0,\,i\in J^+\} \sqcup \{x_i=0,\,i \in J^-\})$. First of all we consider the homeomorphism such that $x_i \mapsto x_i^{1/|\alpha_{\ell i}|}$ if $i \in J^- \sqcup J^+$ and $x_i \mapsto x_i$ otherwise. Then we may assume that $\alpha_{\ell i}=\pm 1$ when $\alpha_{\ell i} \neq 0$. Let $x \in X$ and set $|x|_+=\left(\sum_{i \in J^+}x_i^2\right)^{1/2}, |x|_-=\left(\sum_{i \in J^-}x_i^2\right)^{1/2}$. Set $S=\{x \in X;\,|x|_+=|x|_-\}$. One can check easily that with this definition the intersection of $S$ with an orbit of $\mu_\ell$ is a point (namely, given $x_0 \in X$ we have $\RP_\ell x \cap S = \{\mu_\ell(x_0,\lambda)\}$ with $\lambda=\sqrt{|x_0|_-/|x_0|_+}$). Up to take a permutation of coordinates, we may assume that $N \in J^+$. Let us first assume that $x_N \neq 0$. As in (a), we may choose an homeomorphism $\psi$ such that $\psi \circ \mu_\ell \circ (\imin\psi \times \id_{\RP})$ acts only on the variable $x_N$.
%We first shrink $W$ on $X \setminus (\{x_i=0,\,i\in J^+\} \sqcup \{x_i=0,\,i \in J^-\})$. Let $m^+=\operatorname{lcm}(\alpha_i,\,i\in J^+), m^-=\operatorname{lcm}(\alpha_i,\,i\in J^-)$. Set $S=\{\sum_{i \in J^+}x_i^{m^+m^-/\alpha_{i\ell}}=\sum_{i \in J^-}x_i^{-m^+m^+/\alpha_{i\ell}}\}$. One can check easily that with this definition the intersection of $S$ with an orbit of $\mu_\ell$ is a point. Up to take a permutation of coordinates, we may assume that $N \in J^+$. As in (a), we may choose an homeomorphism $\psi$ such that $\psi \circ \mu_\ell \circ (\imin\psi \times \id_{\RP})$ acts only on the variable $x_N$.
%This corresponds to the change of coordinates sending the matrix
%$$
%\left[
%\begin{array}{cc}
%A & B \\
%C & d
%\end{array}
%\right] \ \ \ \ \mbox{to} \ \ \ \
%\left[
%\begin{array}{cc}
%A' & B' \\
%0 & 1
%\end{array}
%\right]
%$$
%where $A$ (resp. $A'$) is a $(N-1) \times (N-1)$ matrix, $B$ (resp. $B'$) is a $(N-1) \times 1$ matrix, $C$ (resp. $0$) is a $1 \times (N-1)$ matrix and $d \neq 0$.

Replace $V$ with $\psi(V)$ and $\mu_j$ with $\psi \circ \mu_j \circ (\imin\psi \times \id_{\RP})$, $j=1,\dots,\ell$. One can check easily that $V \cap S$ is conic with respect to the actions $\widetilde{\mu}_j$ defined as $\widetilde{\mu}_j(x_1,\dots,x_{N-1},x_N,\lambda)=\eta(\lambda^{\alpha_{j1}}x_1,\dots,\lambda^{\alpha_{jN-1}}x_{N-1})$, where $\eta:[0,+\infty]^{N-1} \iso S$. By the induction hypothesis $W \cap S$ contains a subanalytic open neighborhood $W_1$ of $V \cap S$ in $S$ which is $(\widetilde{\RP})^{\ell-1}$-connected (here $\widetilde{\RP}$ means conic with respect to $\widetilde{\mu}_j$, $j=1,\dots,\ell-1$).

Let $\pi:(x_1,\dots,x_{N-1},x_N) \mapsto (x_1,\dots,x_{N-1})$ be the projection. Then $W \cap \imin\pi(\pi(p))$ is a disjoint union of intervals. Let us consider the interval $(m(p),M(p))$ intersecting $S$. Set
$$
\widetilde{W}_1'=\{(x_1,\dots,x_{N-1},x_N);\;(x_1,\dots,x_{N-1}) \in \pi(W_1),\; m(p)<x_N<M(p)\}.
$$
By construction $\widetilde{W}_1'$ is subanalytic and $\RP_\ell$-connected. To show that it is $(\RP)^\ell$-connected one can easily adapt the proof of (a).

Now let us consider $W \cap \{x_N=0\}$. Up to shrink $W$, we may assume that $\widetilde{W}_1' \cup (W \cap \{x_N=0\})$ is open. By induction on the dimension of $J^+$ we may construct $\widetilde{W}_1^0 \subset W \cap \{x_N=0\}$ with the required properties. By construction $\widetilde{W}_1=\widetilde{W}_1' \cup \widetilde{W}_1^0$ is subanalytic and $\RP_\ell$-connected.\\

Let us now shrink $W \cap \{x_i=0,\, i \in J^+\}$. The same method will apply to  $W \cap \{x_i=0,\, i \in J^-\}$. The set $V \cap \{x_i=0,\, i \in J^+\} \cap \{x_i=0,\, i \in J^-\}$ is $\RP_j$-conic, $j=1,\dots,\ell-1$ (the action $\mu_\ell$ is trivial there). By the induction hypothesis we may assume that $W  \cap \{x_i=0,\, i \in J^+\} \cap \{x_i=0,\, i \in J^-\}$ is $(\RP)^{\ell-1}$-connected.
%First of all we consider the homeomorphism such that $x_i \mapsto x_i^{1/|\alpha_{i\ell}|}$ if $i \in J^-$ and $x_i \mapsto x_i$ otherwise. Then we may assume that the orbits of $\mu_\ell$ are half-lines passing through $\{x_i=0,\, i \in J^-\}$.
The intersection of the orbits of $\mu_\ell$ and $W$ are homeomorphic to a disjoint union of intervals. Let us choose the ones whose boundary intersects $\{x_i=0,\, i \in J^-\}$. Let us prove that this provides a set satisfying the hypothesis.

Assume that $N \in J^-$ and $x_N \neq 0$.  As in (a), we may choose an homeomorphism $\psi$ such that $\psi \circ \mu_\ell \circ (\imin\psi \times \id_{\RP})$ acts only on the variable $x_N$. Replace $V$ with $\psi(V)$ and $\mu_j$ with $\psi \circ \mu_j \circ (\imin\psi \times \id_{\RP})$, $j=1,\dots,\ell$. Let $\pi:(x_1,\dots,x_{N-1},x_N) \mapsto (x_1,\dots,x_{N-1})$ be the projection.  One can check easily that $\pi(V)$ is conic with respect to the actions $\widetilde{\mu}_j$ defined as $\widetilde{\mu}_j(x_1,\dots,x_{N-1},\lambda)=(\lambda^{\alpha_{j1}}x_1,\dots,\lambda^{\alpha_{jN-1}}x_{N-1})$ and that $\pi(W)$ is $(\widetilde{\RP})^{\ell-1}$-connected (here $\widetilde{\RP}$ means with respect to $\widetilde{\mu}_j$, $j=1,\dots,\ell-1$). This is because $\imin\psi(W)  \cap \{x_i=0,\, i \in J^+\} \cap \{x_i=0,\, i \in J^-\}$ is $(\RP)^{\ell-1}$-connected and $\imin\psi(W)$ satisfies Vb.

Let $p \in W$. Then $W \cap \imin\pi(\pi(p))$ is a disjoint union of intervals. Let us consider the interval $(m(p),M(p))$ with $m(p)=0$. Set
$$
\widetilde{W}_0^+=\{(x_1,\dots,x_{N-1},x_N);\;(x_1,\dots,x_{N-1}) \in \pi(W),\; 0<x_N<M(p)\}.
$$
By construction $\widetilde{W}_0^+$ is subanalytic and $\RP_\ell$-connected. To show that it is $(\RP)^\ell$-connected one can easily adapt the proof of (a). On $x_N=0$ the proof is similar.\\

To end the proof, let us show that $\widetilde{W}=\widetilde{W}_1 \cup \widetilde{W}_0^+\cup \widetilde{W}_0^-$ is a neighborhood of $\widetilde{W}_0^+\cup \widetilde{W}_0^+$. We argue by contradiction. We suppose that there exists $y \in \widetilde{W}_0^+\cup \widetilde{W}_0^+$ and that for each $\varepsilon>0$ there is $y_\varepsilon \in \{x \in X;\;|x-y|<\varepsilon,\;|x|_+,|x|_- > 0\}$ with $y_\varepsilon \notin \widetilde{W}$. Suppose that $y \in \widetilde{W}_0^-$ (if $y \in \widetilde{W}_0^+$ the proof is similar). Taking $\varepsilon$ small enough, we may assume that $\RP_\ell y_\varepsilon\cap S=\mu_\ell(y_\varepsilon,\lambda_\varepsilon)$ with $\lambda_\varepsilon\in\RP$. By construction of $\widetilde{W}$ there must be $z_\varepsilon \in \mu_\ell(y_\varepsilon,[\lambda_\varepsilon,1])$ (or $z_\varepsilon \in \mu(y_\varepsilon,[1,\lambda_\varepsilon])$) with $z_\varepsilon \notin W$.

%We may also assume that for each $\varepsilon$ small enough, at least countably many $z_\varepsilon$ belong to $K \subset W$ with $K$ compact. The set $K$ is constructed as follows: take a compact neighborhood $B=\{x \in X;\;|x-y|\leq\delta\}$ of $y$ contained in $W$ and set $K=\RP_\ell B \cap \{|x|\leq |y|+\delta\} \cap \{\sum_{i \in J^+}x_i^{m^+m^-/\alpha_{i\ell}} \leq \sum_{i \in J^-}x_i^{-m^+m^+/\alpha_{i\ell}}\}$. Being a neighborhood of $\widetilde{W}_0^+\cup \widetilde{W}_0^+$, $W$ must contain such a $K$ for $\delta$ small enough.

We may also assume that for  $\varepsilon$ small enough, $z_\varepsilon$ belongs to $K \subset W$ with $K$ compact. The set $K$ is constructed as follows: \\
- if $y \notin \widetilde{W}_0^+$, take a compact neighborhood $B=\{x \in X;\;|x-y|\leq\delta\}$ of $y$ contained in $W$ and set $K=\RP_\ell B \cap \{|x|_-\leq |y|_-+\delta\} \cap \{|x|_+ \leq |x|_-\}$,\\
- if $y \in \widetilde{W}_0^- \cap \widetilde{W}_0^+$, take a compact neighborhood $B=\{x \in X;\;|x-y|\leq\delta\}$ of $y$ contained in $W$ and set $K=\RP_\ell B \cap \{|x|_-\leq \delta\} \cap \{|x|_+ \leq \delta\}$.\\
 Being a neighborhood of $\widetilde{W}_0^+\cup \widetilde{W}_0^+$, $W$ must contain such a $K$ for $\delta$ small enough.

Then we can extract a sequence $z_n$, with $z_n \in \mu_\ell(y_{\varepsilon_n},[\lambda_{\varepsilon_n},1])$ (or $z_n \in \mu_\ell(y_{\varepsilon_n},[1,\lambda_{\varepsilon_n}])$)which has a limit $z$ in $K$. This leads to a contradiction since $z \in W$ and $W$ is open.\\

\end{proof}

\subsection{A proof for Propostion \ref{prop:canonical-from}}\label{D}

In this appedix, we give the proof for Propostion \ref{prop:canonical-from}.
Clearly ii.~ of the propostion implies i.~of the propostion. 
We will show the converse implication. We may assume $q=0$ and
$(T^*X)_q = \mathbb{R}^n$. Furthermore,
by a linear coordinates transformation, we also assume that each $V_k$ in the definition
 is the vector subspace spanned by the vectors $dx_i$ $(i \in E_k)$
for some $E_k \subset \{1,\dots,n\}$.
We will prove the claim by the induction with respect to $n$, i.e., the dimension of $X$.
If $n=0$, then the claim clearly holds. Suppose that the claim were true for $0,1,\dots,n-1$ and we will show the
claim for $n > 0$. We need the following lemma.
\begin{lem}{\label{lem:extend_real_analtyic}}
Let $\chi = \{M_1,\dots,M_\ell\}$ be a family of closed submanifolds in $X$,
and let $g_j$ $(j=1,\dots,\ell)$ be a real analytic function on $M_j$ satisfying
$g_{j_1} = g_{j_2}$ on $M_{j_1} \cap M_{j_2}$ for any $1 \le j_1, j_2 \le \ell$.
Suppose Proposition \ref{prop:canonical-from} were true and
$\chi$ is simultaneously linearizable at $q \in M$. Then
there exists a real analytic function $g(x)$ in an open neighborhood of $q$ such that
$g = g_j$ on $M_j$. Furthermore, if $d_{M_j}g_j(q) = 0$ holds for $j=1,\dots,\ell$,
then we can choose such a $g(x)$ with $dg(q) = 0$.
\end{lem}
\begin{proof}
We may assume $q=0$ and $X= \mathbb{R}^n$. We choose a system of coordinates $(x_1, \dots, x_n)$
satisfying that each $M_j$ is give by $\{x \in \mathbb{R}^n;\,x_i = 0\,\, (i \in I_j)\}$
for some $I_j \subset \{1,\dots,n\}$. Let $\pi_j$ be the orthogonal projection from
$\mathbb{R}^n$ to $M_j$. For a subset $\alpha = \{j_1,\dots,j_p\} \subset \{1,\dots,\ell\}$,
we set $M_\alpha = M_{j_1} \cap \dots \cap M_{j_p}$ and
define the function $g_\alpha$ in $M_\alpha$ by $g_j|_{M_\alpha}$
for some $j \in \alpha$.
Note that, by the gluing conditions,  $g_\alpha$ is independent of the choice of $j \in \alpha$.
Let $\pi_{\alpha}$ denote the orthogonal projection from $\mathbb{R}^n$ to
$M_\alpha$. Then
$$
g(x) = \sum_{\alpha \ne \emptyset} (-1)^{|\alpha|+1} g_\alpha(\pi_\alpha(x)),
$$
gives a desired real analytic function
where $\alpha$ ranges through non-empty subsets of $\{1,\dots,\ell\}$ and $|\alpha|$ denotes
the number of elements of $\alpha$.
\end{proof}
Set $X_k := \{x \in \mathbb{R}^n;\, x_i = 0\,\,(i \in E_k)\}$. We denote by
$\pi_k: \mathbb{R}^n \to X_k$ be the orthogonal projection.
Define $\mathcal{J}_k := \{j \in \{1,\dots,\ell\};\, V_k \subset (T^*_{M_j}X)_q\}$.
Note that, in a neighborhood of $0$,  $\pi_k$ gives an isomorphism between
$M_j$ and $\pi_k(M_j) \subset X_k$ for any $j \in \mathcal{J}_k$.

Now let us consider a $1 \le k \le m$ for which $\mathcal{J}_k$ is non-empty.
Suppose that $k=1$ is such a $k$.
For each $j \in \mathcal{J}_1$ and $i \in E_1$, because of $dx_i \in (T^*_{M_j}X)_0$,
we can find a real analytic function $g_{j,i}$
in an open neighborhood of $0$ with $g_{j,i} = 0$ on $M_j$ and $dg_{j,i}(0) = dx_i$.
Then, by the implicit function theorem, we may assume that
each $g_{j,i}$ has a form $x_i - g'_{j,i}(x')$ where
$x'$ denotes the coordinates $x_i$ with $i \in \{1,\dots,n\}\setminus E_1$
(i.e., the coordinates of $X_1$) and $g'_{j,i}(x')$ is a real analytic function in $X_1$
with $g'_{j,i}(0) = 0$ and $d_{X_1}g'_{j,i}(0) = 0$.

Let us consider the family $\chi'=\{\pi_1(M_j)\}_{j \in \mathcal{J}_1}$ of closed submanifolds in
$X_1$. Clearly $(T^*X_1)_0$ contains the vector subspace $V_2 \oplus \dots \oplus V_m$
such that each $(T^*_{\pi_1(M_j)}X)$ is a part of this direct sum for $j \in \mathcal{J}_1$.
\begin{lem}
For any $j_1$ and $j_2$ in $\mathcal{J}_1$,  the subset
$\pi_1(M_{j_1}) \cap \pi_1(M_{j_2})$ is a closed submanifold
in $X_1$ and we have
 $$
 (T^*_{\pi_1(M_{j_1})}X_1)_0 + (T^*_{\pi_1(M_{j_2})}X_1)_0 =
 (T^*_{\pi_1(M_{j_1}) \cap \pi_1(M_{j_2})}X_1)_0.
 $$
\end{lem}
\begin{proof}
We get
\begin{equation}{\label{eq:dim-cross-1}}
\kappa := \dim\, (M_{j_1} \cap M_{j_2}) =
\dim\, X - \dim\, ((T^*_{M_{j_1}}X)_0 + (T^*_{M_{j_2}}X)_0)
\end{equation}
and
\begin{equation}{\label{eq:dim-cross-2}}
\begin{aligned}
&\dim\, X - \dim\, ((T^*_{M_{j_1}}X)_0 + (T^*_{M_{j_2}}X)_0) \\
&\quad =
\dim\, X_1 - \dim\,
((T^*_{\pi_1(M_{j_1})}X_1)_0 + (T^*_{\pi_1(M_{j_2})}X_1)_0)  \\
&\qquad\ge \dim\, (\pi_1(M_{j_1}) \cap \pi_1(M_{j_2})).
\end{aligned}
\end{equation}
Here the last dimension denotes that of an analytic germ at $0$.
As a result, we obtain
$\kappa \ge \dim\, (\pi_1(M_{j_1}) \cap \pi_1(M_{j_2}))$.
Let $L$ be an irreducible component of the germ of
$\pi_1(M_{j_1}) \cap \pi_1(M_{j_2})$. Note that $\dim L \le \kappa$ also holds.
Then, since $M_{j_1} \cap M_{j_2}$ is defined by
$$
\left\{(\{g'_{j_1,i}(x')\}_{i \in E_1}, \,x') \in X;\,
\begin{aligned}
&g'_{j_1,i} (x') = g'_{j_2,i}(x') \,\, (i \in E_1) \\
&x' \in \pi_1(M_{j_1}) \cap \pi_1(M_{j_2})
\end{aligned}
\right\},
$$
if either $\dim\, L < \kappa$ or $g'_{j_1,i} \ne g'_{j_2,i}$ on $L$ for some $i \in E_1$ holds,
then we can find a point $q'$ such that the dimension of the germ of $M_{j_1} \cap M_{j_2}$
at $q'$ is strictly less than $\kappa$, which contradicts the fact that
$M_{j_1} \cap M_{j_2}$ is smooth and its dimension is $\kappa$.
Hence we can conclude that every irreducible components of
$\pi_1(M_{j_1}) \cap \pi_1(M_{j_2})$ has  dimension $\kappa$ and we also have
$g'_{j_1,i} = g'_{j_2,i}$ on $\pi_1(M_{j_1}) \cap \pi_1(M_{j_2})$ for any $i \in E_1$.
Therefore we have
$$
M_{j_1} \cap M_{j_2}
=
\left\{(\{g'_{j_1,i}(x')\}_{i \in E_1}, \,x') \in X;\,
x' \in \pi_1(M_{j_1}) \cap \pi_1(M_{j_2})
\right\},
$$
from which $\pi_1(M_{j_1} \cap M_{j_2}) = \pi_1(M_{j_1}) \cap \pi_1(M_{j_2})$ follows.
Clearly $\pi_1(M_{j_1} \cap M_{j_2})$ is a closed submanifold and its dimension
is still $\kappa$. Therefore
$\pi_1(M_{j_1}) \cap \pi_1(M_{j_2})$ is a closed submanifold in $X_1$.
Furthermore,  by(\ref{eq:dim-cross-1}) and (\ref{eq:dim-cross-2}),
we get
$$
\begin{aligned}
&\dim\, ((T^*_{\pi_1(M_{j_1})}X_1)_0 + (T^*_{\pi_1(M_{j_2})}X_1)_0)  \\
\qquad &=\dim\, X_1 - \kappa =
\dim\, (T^*_{\pi_1(M_{j_1}) \cap \pi_1(M_{j_2})}X_1)_0.
\end{aligned}
$$
Hence we have obtained the last equation in the lemma.
\end{proof}
By these observations,
the family $\chi'=\{\pi_1(M_j)\}_{j \in \mathcal{J}_1}$ of closed submanifolds in
$X_1$ is simultaneously linearizable. Hence,  by the induction hypothesis,
the claim ii.~of the proposition for $\chi'$ holds.
It follows from the proof of the above lemma that, for any $j_1$, $j_2$ in $\mathcal{J}_1$
and $i \in E_1$,
we have $g'_{j_1,i} = g'_{j_2,i}$ on $\pi_1(M_{j_1}) \cap \pi_1(M_{j_2})$.
Hence, for each $i \in E_1$,
by applying Lemma {\ref{lem:extend_real_analtyic}}
to functions $g'_{j,i}$ on $\pi_1(M_j)$ ($j \in \mathcal{J}_1$),
there exists a real analytic function $g'_i(x')$ in $X_1$ such
that $g'_i = g'_{j,i}$ on $\pi_1(M_j)$ for $j \in \mathcal{J}_1$ and $d_{X_1} g'_i(0) = 0$.
Thus we have obtained the real analytic function $f_i(x) := x_i - g'_i(x')$
($i \in E_1$) such that $df_i(0) = dx_i$ and $f_i(x) = 0$ on $M_j$ ($j \in \mathcal{J}_1$).
By repeating this procedure for $E_2$, $\dots$, $E_m$ respectively, we can obtain
$n$-real analytic functions $f_i$ $(i=1,\dots,n$). Here, if
$\mathcal{J}_k$ is empty, then we simply set $f_i(x) := x_i$ ($i \in E_k$). It follows from the construction
that these $f_i$'s form a system of local coordinates of $X$ near $0$ and each
$M_j$ is defined by equations $f_i = 0$ $(i \in \underset{V_k \subset (T^*_{M_j}X)_q}{\bigcup} E_k )$.
Hence we have shown the claim for $n$ and this completes the proof.
\end{section}

\addcontentsline{toc}{section}{

\noindent
\parbox[t]{.48\textwidth}{
Naofumi HONDA \\
Department of Mathematics, \\
Faculty of Science, \\
Hokkaido University, \\
060-0810 Sapporo, Japan. } \hfill
\parbox[t]{.48\textwidth}{
Luca PRELLI\\
Dipartimento di Matematica \\
Universit\`{a} degli Studi di Padova,\\
Via Trieste 63,\\
35121 Padova, Italia. \\
lprelli@math.unipd.it }
\end{document}